\documentclass[11pt,reqno]{article}
\usepackage{ifthen}
\newboolean{shortver}
\setboolean{shortver}{false}
\usepackage{a4wide}
\usepackage{amsmath,amsthm,amssymb}
\usepackage{extarrows}
\usepackage{units}
\usepackage{url}
\usepackage{calc, graphicx} 
\usepackage{bbm}
\usepackage{scalerel}
\usepackage{dsfont}
\usepackage{lmodern} 
\usepackage{nicematrix}
\usepackage{arxiv}
\usepackage{color}
\usepackage[utf8]{inputenc} 
\usepackage[T1]{fontenc}    
\usepackage[hidelinks]{hyperref}       
\usepackage{url}            
\usepackage{booktabs}       
\usepackage{amsfonts}       
\usepackage{nicefrac}       
  \usepackage{pgfplots}
  \pgfplotsset{compat=newest}
  \usepackage{rotating}
\usepackage{tikz}

\usetikzlibrary{matrix, decorations.pathreplacing}
  \usetikzlibrary{plotmarks}
  \usetikzlibrary{arrows.meta}
  \usepgfplotslibrary{patchplots}
  \usepackage{grffile}
  \usepackage{amsmath}
  \usepackage{nicematrix}
  \usepackage{import}
\usepackage{xifthen}
\usepackage{pdfpages}
\usepackage{transparent}
\usepackage{tikz}
\usepackage{extarrows}
\usepackage{caption}
\DeclareCaptionLabelFormat{cont}{#1~#2\alph{ContinuedFloat}}
\usepackage{microtype}      
\usepackage{lipsum}		
\usepackage{graphicx}
\usepackage{float}
\usepackage{multitoc}
\usepackage{enumitem}
\usepackage{subfigure}
\usepackage{doi}
\usepackage{amsfonts}
\usepackage{mathrsfs}
\usepackage{bm}
 \usepackage{nicematrix}
\usepackage{latexsym}
\usepackage{arydshln}
\usepackage[toc,page]{appendix}
\usepackage{cleveref}
\usepackage{thmtools}
\usepackage{thm-restate}

\usepackage{hyperref}

\usepackage{cleveref}
\usepackage{stmaryrd}
\Crefformat{equation}{(#2#1#3)}

\newtheorem{theorem}{Theorem}[section]

\newtheorem{corollary}[theorem]{Corollary}
\newtheorem{proposition}[theorem]{Proposition}
\newtheorem{lemma}[theorem]{Lemma}

\newtheorem{claim}{Claim}[section]
\newtheorem{example}{Example}

\newtheorem{assumption}[theorem]{Assumption}

\newtheorem{notation}[theorem]{Notation}
\numberwithin{equation}{section}
\theoremstyle{definition}
\newtheorem{definitiony}[theorem]{Definition}
\newenvironment{definition}
  {\pushQED{\qed}\definitiony}
  {\popQED\enddefinitiony}

\theoremstyle{remark}
\newtheorem{remarkx}[theorem]{Remark}
\newenvironment{remark}
  {\pushQED{\qed}\remarkx}
  {\popQED\endremarkx}
\newcommand\slv{\ifthenelse{\boolean{shortver}}{You can find the proof in arXiv. }{We defer this proof to \Cref{app_pre_results} of the appendix. }}
\newcommand\slvv{\ifthenelse{\boolean{shortver}}{You can find the proofs in arXiv. }{We defer these proofs to \Cref{app_pre_results} of the appendix. }}
\newcommand\rou[1]{\lfloor{#1}\rfloor}
\newcommand\bc[1]{\left({#1}\right)}
\newcommand\cbc[1]{\left\{{#1}\right\}}
\newcommand\brk[1]{\left\lbrack{#1}\right\rbrack}
\newcommand\abc[1]{\llbracket{#1}\rrbracket}
\newcommand\aabc[1]{\langle{#1}\rangle}
\newcommand\abs[1]{\left|{#1}\right|}
\newcommand\dabs[1]{\left\|{#1}\right\|}
\newcommand\supp[1]{{\rm supp}({#1})}
\newcommand\rk[1]{\rank_{\mathbb{F}}\left({#1}\right)}
\newcommand\trk[1]{\rank_{\mathbb{F}}}
\newcommand\bin[1]{{\rm Bin}\left({#1}\right)}

\newcommand*{\dif}{\mathop{}\!\mathrm{d}}

\newcommand\syn{{\rm Sym}_n(\FF^*)}

\newcommand\teo[1]{\ensuremath{\mathds{1}}{\left\{{#1}\right\}}}

\newcommand\ranges{s\in[\varepsilon,1-\sigma(-\ln\xi)-\varepsilon]}
\newcommand\rangee{\varepsilon\in(0,1/2-\sigma(-\ln\xi)/2)}
\DeclareRobustCommand{\VAN}[3]{#2}

\newcommand\cm[1]{${\mbox{\rm CM}}_{#1}$}
\newcommand\tcm[2]{${\mbox{\rm CM}}_{#1}({#2})$}
\newcommand\scm[2]{${\mbox{\rm CM}}_{#1,\abc{#2}}$}
\newcommand\mcm[2]{${\mbox{\rm CM}}_{#1,\abc{#2}}^{m}$}
\newcommand\mcma[1]{${\mbox{\rm CM}}_{#1}^{m}$}
\newcommand\fmcm[2]{{\mbox{\rm CM}}_{#1,\abc{#2}}^{m}}
\newcommand\fmcma[1]{{\mbox{\rm CM}}_{#1}^{m}}
\newcommand\fcm[1]{{\mbox{\rm CM}}_{#1}}
\newcommand\ftcm[2]{{\mbox{\rm CM}}_{#1}({#2})}
\newcommand\fscm[2]{{\mbox{\rm CM}}_{#1,\abc{#2}}}

\newcommand\pth{[\bm{\theta}]}

\newcommand\verst[1]{\mathcal{V}_{\abc{#1}}}
\newcommand\edgesl[1]{\bm{S}_{\abc{#1}}}
\newcommand\edgeli[1]{\bm{L}_{\abc{#1}}}

\newcommand\verslk[2]{\bm{V}_{#1,\abc{#2}}}

\newcommand\bverslk[2]{\bar{\bm{V}}_{#1,\abc{#2}}}
\newcommand\bcabc[1]{\bm{C}_{\abc{#1}}}
\newcommand\vmm{\bm{m}}

\newcommand\vw{\bm w}
\newcommand\vx{\bm x}
\newcommand\vh{\bm h}

\newcommand\vd{\bm d}

\newcommand\vu{\bm u}
\newcommand\vv{\bm v}
\newcommand\vy{\bm y}
\newcommand\vz{\bm z}

\newcommand\fO{\mathfrak O}
\newcommand\fM{\mathfrak M}
\newcommand\vze{\bm \zeta}

\newcommand\bvet{\bar{\bm \eta}}
\newcommand\vta{\bm \tau}
\newcommand\vnu{\bm{\nu}}

\newcommand\bvta{\bar{\bm \tau}}

\newcommand\vga{\bm{\gamma}}
\newcommand\val{\bm{\alpha}}
\newcommand\vbe{\bm{\beta}}

\newcommand\vth{\bm{\theta}}
\newcommand\hval{\bm{\alpha}^T}

\newcommand\oone{\bar{o}_{\PP}(1)}
\newcommand\ooone[1]{\bar{o}_{\PP}({#1})}

\newcommand\fB{ {\mathfrak F}^c_s}
\newcommand\fBT{{\mathfrak F}{\bc{\rm tr}}^c_s}
\newcommand\ffB{{\mathfrak F}_s}
\newcommand\ffBT{ {\mathfrak F}{\bc{\rm tr}}_s}

\newcommand\cX{\mathcal{X}}
\newcommand\cY{\mathcal{Y}}
\newcommand\cZ{\mathcal{Z}}
\newcommand\cU{\mathcal{U}}
\newcommand\cV{\mathcal{V}}
\newcommand\cW{\mathcal{W}}
\newcommand\ex[1]{{\rm Exp}({#1})}
\newcommand\THETA{\mathbf\Theta}

\newcommand\eul{\mathrm{e}}
\newcommand\eps{\varepsilon}

\newcommand\ZZ{\mathbb{Z}}
\newcommand\FF{\mathbb{F}}
\newcommand\RR{\mathbb{R}}

\newcommand\NN{\mathbb{N}}
\newcommand\PP{\mathbb{P}}

\newcommand\Erw{\mathbb{E}}

\newcommand{\ind}{\ensuremath{\mathds{1}}}

\DeclareMathOperator{\rank}{rk}
\DeclareMathOperator{\PR}{PR}

\usepackage{todonotes}
\usepackage{color}

\DeclareSymbolFont{extraup}{U}{zavm}{m}{n}
\DeclareMathSymbol{\varheart}{\mathalpha}{extraup}{86}
\DeclareMathSymbol{\vardiamond}{\mathalpha}{extraup}{87}
\usepackage{indentfirst}
\usepackage{authblk}

\title{The asymptotic rank of adjacency matrices of weighted configuration models over arbitrary fields \\ \large From graph exploration to fixed-point equations}
\renewcommand{\shorttitle}{The asymptotic rank of adjacency matrices of weighted configuration models over arbitrary fields}

\author[1]{Remco van der Hofstad
\thanks{Email: r.w.v.d.hofstad@tue.nl}}
\author[1]{Noela Müller
\thanks{Email: n.s.muller@tue.nl}}
\author[1]{Haodong Zhu
\thanks{Email: h.zhu1@tue.nl}}
\affil[1]{Department of Mathematics and Computer Science, Eindhoven University of Technology, P.O. Box 513, Eindhoven, 5600 MB, The Netherlands}
\begin{document}
\maketitle
\begin{abstract}
We study the asymptotic rank of adjacency matrices of a large class of edge-weighted configuration models. Here, the weight of a (multi-)edge can be any fixed non-zero element from an arbitrary field, as long as it is independent of the (multi-)graph. Our main result demonstrates that the asymptotic behavior of the normalized rank of the adjacency matrix neither depends on the fixed edge-weights, nor on which field they are chosen from. Our approach relies on a novel adaptation of the component exploration method of \cite{janson2009new}, which enables the application of combinatorial techniques from \cite{coja2022rank, HofMul25}.

\end{abstract}
\section{Introduction}

\subsection{Background and motivation} \label{sec_intro}
Over the past decades, the study of random graphs has developed into a prominent and flourishing field of research 
both within and outside of mathematics, with applications in many different disciplines,  
such as epidemiology \cite{lewis2011network}, information technology  
\cite{barabasi2013network,broder2000graph}, or the social sciences \cite{BOCCALETTI20141,borgatti2009network}. In order to represent the behavior of real-world networks, researchers have explored a great variety of \textit{sparse} random graph models. Among these, the configuration model stands out as one of the most widely used models to generate networks with any specified degree distribution \cite{bollobas1980,frieze2016introduction,van2017random,van2022random,mr1995,mr1998,newman2018networks,riordan2012phase}.  
Adding to its significance, the configuration model generalizes two of the most prominent random graph models: random $d$-regular graphs and Erd\H{o}s-Rényi random graphs.

A convenient way to represent the connectivity structure of a random graph model such as the configuration model is via its (bi-)adjacency matrix. Adjacency matrices are a subject of extensive interest in their own right \cite{bordenave2011rank,coja2022rank,DemGla24,glasgow2023exact,huang2021invertibility, vu2021recent}, and again, also many practically relevant questions about random graphs can be tackled via the properties of the associated adjacency matrices: In the community
detection problem, we calculate the eigenvalues of the adjacency matrix  to recover the community structure  
\cite{decelle2011asymptotic,mcsherry2001spectral}, while the (sub)determinants of the (bi-)adjacency matrix provide a way to solve the total matching problem \cite{FerLuc2024}.

What can be said about the asymptotic properties of the (bi-)adjacency matrix of a given random graph model? Despite their local variability, many global characteristics of the most widely used sparse models converge, such as the size of the largest connected component or the distance between two uniformly chosen vertices \cite{van2022random}. In this sense, also the rank of the (bi-)adjacency matrices of different sparse random graph models, as one of their most important properties, is conjectured or proven to converge under suitable normalization 
\cite{bauer2001exactly,bordenave2011rank,coja2022rank,Bypassing2013}.
 
 To gain some intuition on the rank question for adjacency matrices of sparse random graphs, let us focus on sparse Erd\H{o}s-Rényi random graphs for the time being, where every vertex has a bounded number of neighbors on average. 
In this case, the rank of the adjacency matrix can be upper-bounded quite accurately by 
a peeling process that had originally been invented by Karp and Sipser \cite{karp1981maximum} as an algorithm to find near-optimal matchings
\cite{bauer2001exactly,bordenave2011rank,glasgow2023exact}. This peeling process iteratively removes vertices of degree $1$ along with their unique neighbors. Thus, at the point when the process terminates, a (possibly empty) \textit{Karp-Sipser core} of minimum degree $2$ remains, along with a collection of isolated vertices. By estimating the size of the Karp-Sipser core, Karp and Sipser obtained an asymptotic formula for the matching number of sparse Erd\H{o}s-Rényi random graphs \cite{karp1981maximum}, which later turned out to also yield a sharp upper bound for the rank, as first observed by Bauer and Golinelli \cite{bauer2001exactly}. This upper bound is mostly combinatorial in its nature and does not make use of the $0/1$-entries of the adjacency matrix at all.

It was not until ten years later that Bordenave, Lelarge and Salez found a way to obtain a matching lower bound on the asymptotic rank of sparse Erd\H{o}s-Rényi random graphs \cite{bordenave2011rank}.  Indeed, under an extra condition on the degree distribution, Bordenave, Lelarge and Salez not only derive an asymptotic rank formula for sparse Erd\H{o}s-Rényi random graphs, but also for a much more general class of random graph models. The degree condition has been later removed by Bordenave in \cite{bordenave2016spectrum}. The proof of the lower bound in \cite{bordenave2011rank} is based on the analysis of the spectral measure of the adjacency operator associated to the random graph in question.

As in \cite{bordenave2011rank}, most often, (bi-)adjacency matrices are studied, and their properties are harnessed by regarding them as matrices over the real numbers: For the asymptotic rank of sparse random matrices over finite fields, less is known. 
While this point of view is natural and sufficient for many settings, it leaves out important applications. For example, matrices that naturally arise in problems related to information theory or random constraint satisfaction problems, such as the $k$-XORSAT problem \cite{gwynne2014sat,soos2019bird}, can be treated as adjacency matrices of a suitable random graph (e.g. a  bipartite configuration model \cite{coja2022sparse,Bypassing2013}). In these and related settings, it is more natural to consider the arising (adjacency) matrices over a finite field. While the combinatorial upper bound via matchings is independent of field-specific considerations, the proof of the lower bound in \cite{bordenave2011rank} heavily depends on the field.

Addressing this challenge, Coja-Oghlan, Ergür, Gao, Hetterich and Rolvien \cite{coja2022rank} devised a clean combinatorial argument to derive an asymptotic rank formula for bi-adjacency matrices of \textit{bipartite} configuration models. Their physics-inspired approach can deal with the bi-adjacency matrices of a broad range of bipartite graphs, regardless of the values of the nonzero entries and even the underlying field.

While \cite{coja2022rank} exclusively treats non-symmetric matrices, building upon \cite{coja2022rank}, we  \cite{HofMul25} subsequently extended the result of \cite{bordenave2011rank} to (symmetric) adjacency matrices of sparse Erd\H{o}s-Rényi random graphs, if arbitrary nonzero weights from a fixed field are put on the edges of the graph. The results in \cite{coja2022rank, HofMul25} illustrate that for weighted adjacency matrices of sparse Erd\H{o}s-Rényi random graphs, the precise values of the non-zero entries, or which field they are coming from, do not influence the asymptotic rank. 

However, already the study of the rank of the adjacency matrices of sparse random $d$-regular graphs is more difficult than that of Erd\H{o}s-R\'enyi random graphs, since the edges in these graph are no longer independent. Nevertheless, in \cite{huang2021invertibility}, Huang proved that their adjacency matrices have full rank over $\RR$ with high probability for any fixed $d\geq 3$ by studying the adjacency matrix over finite fields. Other works  \cite{coja2022rank,cook2017singularity,litvak2017adjacency} have focused on the rank of the bi-adjacency matrix of $d$-regular \textit{bipartite} graphs.

These and other results illustrate the huge interest in the rank of sparse random matrices. Most of the mentioned results either concern the rank of asymmetric matrices \cite{coja2022rank,cooper2019rank} or exclusively work over the field $\RR$ \cite{bordenave2011rank,huang2021invertibility}, and thereby take a different route than the clean combinatorial upper bound hinted at by the analysis of the Karp-Sipser core. Given the result in \cite{HofMul25} for sparse Erd\H{o}s-Rényi random graphs, one may wonder whether the asymptotic rank of the adjacency matrix considered over different fields also remains unchanged for a model as complex as the configuration model. What is more, does this still hold true if nonzero entries from an arbitrary field are added to the edges? 
In this article, we prove that under a mild condition on the degree sequence, the asymptotic rank of the adjacency matrix of a weighted configuration model does not rely on the values of the nonzero entries in the adjacency matrix and the field we consider. 
This general result is made possible by a novel implementation of the component exploration method of \cite{janson2009new}, 
specifically tailored to enable the application of the combinatorial techniques developed in  \cite{coja2022rank, HofMul25}. It also extends the earlier work of Bordenave, Lelarge, and Salez \cite{bordenave2011rank}, who employed heavy functional analytic machinery to establish a corresponding result for the real rank of unweighted configuration models.

\section{Model definition}\label{sec_model_def}
\subsection{The configuration model}\label{sec_conmodel_def}
For each positive integer $n$, let $[n]=\cbc{1,\ldots,n}$ denote the vertex set, and $\vd:=(\vd_i)_{i\in[n]}:=(\vd(n,i))_{i\in[n]}$ be a sequence of non-negative, integer-valued random variables or specified constants, which we will refer to as the degree sequence of the graph that is to be constructed. For reasons that will become apparent soon, we additionally assume that the sum of the $\vd_i$ is even. The configuration model \cm{n} is a random \textit{multi-graph} graph on the vertex set $[n]$ such that vertex $i\in[n]$ has degree $\vd_i$. 
To construct a graph with degree sequence $\vd$, we start from a graph on vertex set $[n]$ where $\vd_i$ ``half-edges'', or stubs, are attached to each vertex $i$. To form the edges of the graph, subsequently, half-edges are matched up in pairs to form full edges. More specifically, the half-edges are matched up uniformly at random, with each possible pairing having the same probability. This procedure yields a valid graph as the number of half-edges is even. 
The result of the pairing procedure generally yields a multigraph, i.e., a graph that potentially contains loops and multiedges.

Let $\mathcal N_k:=\mathcal N_k(n):=\cbc{i\in [n]:\vd_i=k}$ be the set of vertices of degree $k$  in the configuration model and let $n_k=\abs{\mathcal N_k}$ denote its size. In the study of sparse random graphs, it is often beneficial to take the perspective of a uniformly chosen vertex: Even though the whole graph may contain cycles, in many prominent random graph models, the finite neighborhood around most vertices does not contain cycles. This gives rise to a local tree structure that can be exploited. To ensure that a uniformly chosen vertex, along with its neighborhood, is well behaved in the configuration model, we impose the following regularity assumptions on $n_k$:
\begin{assumption}[Regularity conditions for the degree sequence]\label{assumption_weaker}
For a degree sequence $\vd=(\vd_i)_{i\in[n]}$ such that $\sum_{i =1}^n \vd_i$ is even for all $n$, assume that there exists a sequence of deterministic nonnegative numbers $(p_k)_{k\geq 0}$ such that $p_0<1$ and the following hold:
\begin{enumerate}[label=(\alph*)]
    \item \label{it_weaker_a_1} 
For all $k \geq 0$: $\lim_{n\to\infty}\Erw\abs{\frac{ n_k }{n}-p_k}= 0$.

    \item  $\lim_{n\to\infty}\Erw\brk{\sum_{k\geq 0}k \frac{ n_k }{n}}= \sum_{k\geq0} k p_k<\infty$.
\end{enumerate}
\end{assumption}
\begin{remark}
\Cref{assumption_weaker} is standard in the study of configuration models, see for example \cite[Condition 1.7 (a)\&(b)]{van2022random}. It also implies that the sequence $(p_k)_{k\geq 0}$ can be regarded as a probability distribution (see \Cref{result_p1}). The assumption $p_0<1$ is made to exclude the trivial case where almost all vertices are isolated.  
\end{remark}

 \Cref{assumption_weaker}\ref{it_weaker_a_1} guarantees that the probability that a uniformly chosen vertex $\vu$ in \cm{n} has degree $k$ is approximately equal to $p_k$. This yields a precise description of the immediate neighborhood ($1$-neighborhood) of a typical vertex. 
 Furthermore, \Cref{assumption_weaker} also allows to look further away than merely distance $1$: One can show that for any fixed $R$, if we examine the $R$-neighborhood of a typical vertex in \cm{n}, as $n$ tends to infinity, it will closely resemble the $R$-neighborhood of the root in a unimodular
branching-process, where the offspring distribution of the root is given by $(p_k)_{k\geq 0}$, while the offspring distribution of the other vertices is given by the size-biased distribution $((k+1) p_{k+1}/\sum_{m\geq 1}m p_m)_{k\geq 0}$ (see \cite[Theorem 4.1 and Corollary 4.5]{van2022random}). More specifically, the just described unimodular
branching-process is the limit of \cm{n} in the sense of \emph{local weak convergence}.

For our purposes, this local limit captures the most important characteristics of \cm{n}, and we define the probability generating function (p.g.f.) $\psi$ of the offspring distribution of the root of the above branching process by setting
\begin{align}\label{def_psi}
    \psi(\alpha):=\sum_{k\geq 0} p_k\alpha^k,\quad \alpha\in[0,1],
\end{align}
and the  probability generating function $\hat{\psi}$ of the offspring distribution of non-root vertices by setting
\begin{align}\label{def_psi_zero}
    \hat{\psi}(\alpha):=\sum_{k\geq 0}  kp_k\alpha^{k-1}/\sum_{k\geq 0} kp_k=\psi'(\alpha)/\psi'(1), \quad \alpha\in[0,1].
\end{align}
Moreover, for $\alpha \in [0,1]$ and  $\phi:[0,1]\mapsto[0,1]$ differentiable, we set
\begin{equation}\label{def_rd}
   R_\phi(\alpha):=2-\phi(1-\phi'(\alpha)/\phi'(1))-\phi(\alpha)-\phi'(\alpha)(1-\alpha).
    \end{equation}
Finally, throughout the article, we will assume that the p.g.f. $\psi$ from \Cref{def_psi} satisfies the following assumption:
\begin{assumption}\label{tech_assumption}
$(p_k)_{k\geq 0}$ from \Cref{assumption_weaker} is such that the second derivative of $\psi$ is log-concave on $(0,1)$. 
\end{assumption}
\Cref{tech_assumption} also appears in \cite[Theorem 13]{bordenave2011rank} and corresponds to the case that there is an almost perfect matching on the Karp-Sipser core of the graph \cite{bordenave2011rank,bordenave2013matchings}. See \Cref{sec-tech_assumption} for a discussion of this condition.

\subsection{Weighted adjacency matrix and main result}\label{sec_intro_adj}
We next define a general weighted version of \cm{n} along with its adjacency matrix. Given an arbitrary field $\FF$, let $\FF^\ast:=\FF \setminus\{0\}$ be its multiplicative group. The elements of $\FF^\ast$ will serve as the edge-weights of our model. 
More specifically, fix a symmetric matrix $J_n\in \syn$\footnote{Here, $\syn$ denotes the set of all symmetric $n \times n$ matrices with entries from $\FF^\ast$.}. For each distinct $i,j\in[n]$, we give weight $J_n(i,j)$ to all edges between vertices $i$ and $j$, if there are any present, and 
define the weighted adjacency matrix of \cm{n} by ignoring multi-edges and self-loops by setting 
\begin{equation}\label{ec1}
  \bm{A}_{n}(i,j):=\begin{cases}
      \teo{{\mbox{there is at least one edge between $i$ and $j$ in \cm{n}}} }J_n(i,j),\quad &i\neq j; \\
      0,&i=j.
  \end{cases}
\end{equation}

To emphasize that we consider the rank of $\bm{A}_n$ over different fields, let $\rank_{\FF}(\bm{A}_n)$ specifically denote the rank of $\bm{A}_n$ over $\FF$. 
The main result of this paper is the following theorem, which states that  under Assumptions \ref{assumption_weaker} and \ref{tech_assumption}, the normalized ranks of $(\bm{A}_n)_n$ over $\FF$ converge in probability to the constant $\min_{\alpha\in [0,1]}R_\psi(\alpha)$, regardless of the field $\FF$ and the symmetric matrices $(J_n)_n$:
\begin{theorem}\label{coro_main_basic}
Assume that the degree sequence $\vd$ satisfies \Cref{assumption_weaker} with a probability  distribution $(p_k)_{k\geq 0}$  satisfying \Cref{tech_assumption}. Then, for any field $\FF$ and any $(J_n)_n \subseteq \syn$, 
\begin{align*}
\frac{1}{n}\rank_{\FF}\bc{\bm{A}_{n}} \stackrel{\mathbb{P}}{\longrightarrow} \min_{\alpha\in [0,1]}R_\psi(\alpha),\quad  n\to\infty.
\end{align*}
\end{theorem}
\begin{remark}
    Note that we assume \( p_0 < 1 \) in \Cref{assumption_weaker}. When \( p_0 = 1 \),
\begin{equation*}
    \frac{1}{n} \rank_{\mathbb{F}}(\bm{A}_n) \xrightarrow{\mathbb{P}} 0,
\end{equation*}
since only a negligible proportion of the rows and columns of \( \bm{A}_n \) are nonzero.

\end{remark}
We finish this section by stating some examples in which the asymptotic rank formula of \Cref{coro_main_basic} holds. Despite the restrictions that \Cref{tech_assumption} imposes, \Cref{coro_main_basic} can be applied to a large group of random graphs of interest, including random $r$-regular and also Erd\H{o}s-Rényi random graphs: 
\begin{example}
The normalized rank of the adjacency matrix of an edge-weighted random $r$-regular graph, where $r\geq 2$, converges to $1$ in probability (see also \cite{huang2021invertibility}).
\end{example}

\begin{example}[\cite{bordenave2011rank,HofMul25}]
For an edge-weighted Erd\H{o}s-Rényi random graph with average degree $\lambda$, 
\begin{align*}
    \frac{1}{n}\rank_{\FF}\bc{\bm{A}_{n}}  \stackrel{\mathbb{P}}{\longrightarrow}  \min_{\alpha\in[0,1]}\bc{2-\eul^{-\lambda\eul^{\lambda(\alpha-1)}}-\eul^{\lambda(\alpha-1)}(\lambda+1-\lambda\alpha)},\quad  n\to\infty.
\end{align*}
Indeed, if $(\vd_i)_{i\in[n]}$ is the degree sequence of an Erd\H{o}s-Rényi random graph on the vertex set $[n]$ with average degree $\lambda$, then using \cite[Theorem 7.18]{van2017random}, \Cref{coro_main_basic} implies the convergence of the normalized rank of the associated weighted adjacency
matrices within that model. Consequently, \Cref{coro_main_basic} extends the asymptotic rank formula of \cite{HofMul25}. 
\end{example}

\begin{example}
For an edge-weighted configuration model where all but a sublinear number of degrees are $1, 2$ or $3$ (i.e., $p_1+p_2+p_3=1$), 
\begin{align*}
    \frac{1}{n}\rank_{\FF}\bc{\bm{A}_{n}} \stackrel{\mathbb{P}}{\longrightarrow} \min_{\alpha\in[0,1]}R_\psi(\alpha),\quad  n\to\infty.
\end{align*}
This example is used in \cite{budzinski2022critical} to study the behavior of the critical Karp-Sipser core. In the special case where $p_3=0$, the limit can be explicitly evaluated as
\begin{align*}
    \frac{1}{n}\rank_{\FF}\bc{\bm{A}_{n}} \stackrel{\mathbb{P}}{\longrightarrow}  \frac{(2-p_1)^2}{4-3p_1},\quad  n\to\infty.
\end{align*}
\end{example}

\section{Proof overview}\label{sec_prov}
We next provide a high-level overview of the proof of \Cref{coro_main_basic}. We actually prove the following stronger result,
which reveals that the convergence in  \Cref{coro_main_basic} is uniform with respect to the choice of the weight matrices:
\begin{theorem}\label{coro_main}
Assume that the degree sequence $\vd$ satisfies \Cref{assumption_weaker} with a probability  distribution $(p_k)_{k\geq 0}$  satisfying \Cref{tech_assumption} or $\sum_{k\geq 0}k(k-2)p_k\leq 0$. Then, for any field $\FF$, $\rank_{\FF}\bc{\bm{A}_n}/n$ converges in probability to $\min_{\alpha\in [0,1]}R_\psi(\alpha)$ uniformly in $\bc{J_n}_{n\geq 1}$ in the sense that, for any $\varepsilon>0$,
\begin{equation}\label{e01a}
\lim_{n\to\infty}\sup_{J_n\in \syn}\PP\bc{\abs{\frac{1}{n}\rank_{\FF}\bc{\bm{A}_{n}}-\min_{\alpha\in [0,1]}R_\psi(\alpha)}\geq \varepsilon}=0.
\end{equation}
\end{theorem}

The organization of the rest of this section is as follows: In \Cref{sec_local_str} we explain informally how the local structure of \cm{n} determines our proof approach. In \Cref{sec_LR,sec_permat} we introduce general concepts from the study of random constraint satisfaction problems which we will employ in the rank computation. Finally, in \Cref{sec_cm_exp}, we describe the graph exploration process from \cite{janson2009new}, and explain how it can be combined with the machinery and general ideas from \Cref{sec_LR,sec_permat}. We end  with a discussion in \Cref{sec_discussion}.

\subsection{From local structure to rank}\label{sec_local_str}
\begin{figure}[ht]
\begin{minipage}[r]{0.3\linewidth}
\end{minipage}\hfill
\begin{minipage}[r]{0.3\linewidth}
\includegraphics[width=1\linewidth]{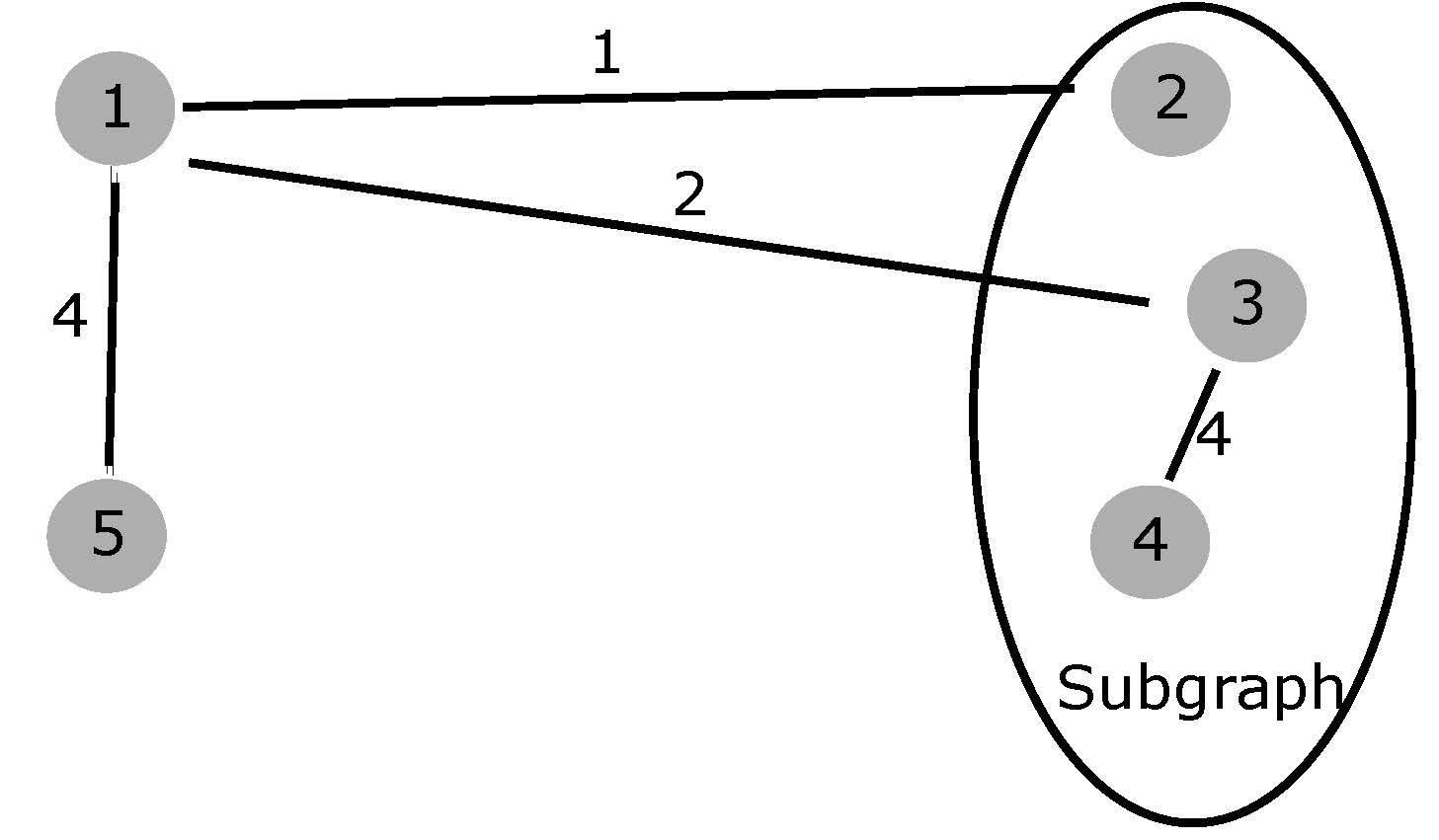}
\end{minipage}\hfill
\begin{minipage}[l]{0.5\linewidth}
\begin{align*}
    \begin{bNiceArray}{ccw{c}{1cm}cc}[margin]
0 & 1 & 2 & 0 & 4\\
 1& \Block[draw]{3-3}{\mbox{Adjacency}\\ \mbox{matrix of}\\ \mbox{the subgraph}} & & & 0 \\
2 & & & & 0 \\
0 & & & & 0 \\
4 & 0 & 0& 0 & 0
\end{bNiceArray}
\end{align*}
\end{minipage}
\caption{An example of a weighted graph (left) and its adjacency matrix (right). Let us look at the effect of the removal of the leaf $5$ along with its unique neighbor $1$ on the rank of the adjacency matrix, where we first remove the neighbor of the leaf. In the adjacency matrix, the first row (corresponding to vertex $1$) is linearly independent of the other rows in the matrix, as it is the only row with a non-zero entry in position $5$ (the column corresponding to the leaf). 
  The same argument applies to the removal of the first column, and is clearly independent of the edge weights. Subsequently removing the row and column corresponding to vertex $5$ does not change the rank any further, as both have been turned into zero-vectors by then.}
\label{f1}
\end{figure}

As a warm-up, assume that we aim to compute the rank of the weighted adjacency matrix of a finite tree graph.  
In this case, a clean combinatorial decomposition argument that is solely based on the Karp-Sipser algorithm, as described in \Cref{sec_intro}, is already sufficient to show that 
the rank of the weighted adjacency matrix is independent of the edge-weights (always assuming that they are non-zero): 
Indeed, if the original graph is a tree, the successive application of the Karp-Sipser leaf-removal algorithm will terminate either when the graph is empty, or when all remaining vertices are isolated. At this point, the adjacency matrix of the remaining graph is either empty or a zero-matrix. In either case, its rank is $0$, and therefore independent of the original field or edge weights. 
Moreover, 
the removal of each degree-1 vertex (a \textit{leaf})  and its unique neighbor along the way reduces the rank of the adjacency matrix by $2$ \cite{bauer2001exactly}. The corresponding argument applies both to graphs with and without edge weights, 
as is illustrated in \Cref{f1}.
As the removal procedure, the rank reduction and the final graph, are all  insensitive to the edge weights, the rank of the weighted adjacency matrix of the original tree graph must be too.

One might hope that the previous discussion already lays out a viable combinatorial proof strategy, as  \Cref{assumption_weaker} holds, at least \textit{locally}, \cm{n} looks like a tree: If we sample a vertex $\bm{u}$ uniformly at random, its local neighborhood is approximated by the unimodular branching process specified in \Cref{sec_conmodel_def}. In the cases where this branching process dies out almost surely, and therefore is a \textit{finite} tree with probability one, the Karp-Sipser heuristic is indeed sufficient to determine the asymptotic rank, as 
 a large proportion of the vertices in \cm{n} are part of \textit{finite} tree components as $n$ approaches infinity, 
and the overall rank is the sum of the ranks of the adjacency matrices of the different connected components. 

The degree sequences for which the neighborhoods of most vertices are (or are not) finite trees can be compactly described in terms of the  
sequence $\bc{p_k}_{k\geq 0}$ from Assumption \ref{assumption_weaker} (see e.g. \cite{van2022random,mr1995,mr1998}):
\begin{enumerate}
    \item \label{case_subc_less1} \textbf{Subcritical degree sequence}:
      $\sum_{k\geq 0} k(k-2)p_k\leq 0$ and $p_0+p_2< 1$: The associated  branching process is a finite tree with probability $1$.
\item \label{case_subc_equal1}\textbf{Critical degree sequence}: $p_0+p_2= 1$: The associated 
branching process is a single vertex or an infinite line.

    \item \textbf{Supercritical degree sequence}: $\sum_{k\geq 0} k(k-2)p_k> 0$: The associated
branching process is infinite with positive probability.
\end{enumerate}

Inspired by this local behavior, we divided the proof of \Cref{coro_main} into three parts, depending on the sign of $\sum_{k\geq 0}k(k-2)p_k$ and the value of $p_2$ as follows:
\begin{proposition}\label{coro_main_prop_1}
Assume that the degree sequence $\vd$ satisfies \Cref{assumption_weaker} with a probability  distribution $(p_k)_{k\geq 0}$  satisfying $\sum_{k\geq 0}k(k-2)p_k\leq 0$ and $p_2\neq 1$. Then, for any field $\FF$, $\rank_{\FF}\bc{\bm{A}_n}/n$ converges in probability to $\min_{\alpha\in [0,1]}R_\psi(\alpha)$ uniformly in $\bc{J_n}_{n\geq 1}$ in the sense that, for any $\varepsilon>0$,
\begin{equation*}
\lim_{n\to\infty}\sup_{J_n\in \syn}\PP\bc{\abs{\frac{1}{n}\rank_{\FF}\bc{\bm{A}_{n}}-\min_{\alpha\in [0,1]}R_\psi(\alpha)}\geq \varepsilon}=0.
\end{equation*}
\end{proposition}
\begin{proposition}\label{coro_main_prop_2}
Assume that the degree sequence $\vd$ satisfies \Cref{assumption_weaker} with a probability  distribution $(p_k)_{k\geq 0}$  satisfying  $p_2= 1$, i.e., $\psi(\alpha)=\alpha$. Then, for any field $\FF$, $\rank_{\FF}\bc{\bm{A}_n}/n$ converges in probability to $\min_{\alpha\in [0,1]}R_\psi(\alpha)=1$ uniformly in $\bc{J_n}_{n\geq 1}$ in the sense that, for any $\varepsilon>0$,
\begin{equation*}
\lim_{n\to\infty}\sup_{J_n\in \syn}\PP\bc{\abs{\frac{1}{n}\rank_{\FF}\bc{\bm{A}_{n}}-1}\geq \varepsilon}=0.
\end{equation*}
\end{proposition}
\begin{proposition}\label{coro_main_prop_3}
Assume that the degree sequence $\vd$ satisfies \Cref{assumption_weaker} with a probability  distribution $(p_k)_{k\geq 0}$  satisfying  \Cref{tech_assumption} and $\sum_{k\geq 0}k(k-2)p_k> 0$. Then, for any field $\FF$, $\rank_{\FF}\bc{\bm{A}_n}/n$ converges in probability to $\min_{\alpha\in [0,1]}R_\psi(\alpha)$ uniformly in $\bc{J_n}_{n\geq 1}$ in the sense that, for any $\varepsilon>0$,
\begin{equation*}
\lim_{n\to\infty}\sup_{J_n\in \syn}\PP\bc{\abs{\frac{1}{n}\rank_{\FF}\bc{\bm{A}_{n}}-\min_{\alpha\in [0,1]}R_\psi(\alpha)}\geq \varepsilon}=0.
\end{equation*}
\end{proposition}
\begin{proof}[Proof of \Cref{coro_main} subject to \Cref{coro_main_prop_1,coro_main_prop_2,coro_main_prop_3}]
    \Cref{coro_main} follows directly from the combination of \Cref{coro_main_prop_1,coro_main_prop_2,coro_main_prop_3}.
\end{proof}

Note that the subcritical degree sequences are those for which the typical local neighborhood is described by a finite tree.
Moreover, critical degree sequences can be essentially reduced to the same setting, as we can remove a small proportion of the edges uniformly at random such that the remaining
degree distribution becomes subcritical, while the rank of the adjacency matrix changes only slightly. 
Therefore,  \cref{e01a} for (sub-)critical configuration models is essentially a corollary of the asymptotic real rank formula for adjacency matrices of unweighted configuration models from \cite{bordenave2011rank}. 
For the sake of completeness, we carry out the detailed derivations of \Cref{coro_main_prop_1,coro_main_prop_2} 
in \Cref{sec_proof_main}.

However, for supercritical degree sequences, the previous, relatively simple, approach is deemed to fail. While a tight almost sure upper bound on the normalized rank of $\bm{A}_n$ can still be obtained in terms of the matching number (through an adaptation of the proof in \cite{Bypassing2013} from the bi-adjacency matrix to the adjacency matrix), obtaining a matching lower bound on the rank poses a bigger problem. Correspondingly, the main contribution of this article is the derivation of a tight lower bound on the asymptotic rank for the supercritical configuration model. 
Our proof of \Cref{coro_main_prop_3} contains two main ingredients: The first is a statistical-physics inspired combinatorial perspective on the problem, which had been adapted to the study of random, non-symmetric matrices in \cite{coja2022rank};  
the second is a \textit{graph exploration} process that has been used to study the largest connected component of configuration models in \cite{janson2009new}. An involved, but ultimately successful, combination of core concepts from \cite{coja2022rank} with the graph exploration enables the derivation of a matching lower bound.

\subsection{Rank and random constraint satisfaction} \label{sec_LR}
We begin with a description of the first component. 
 Any homogeneous linear system $Ax = 0$, $A \in \FF^{m \times n}$, can be naturally regarded as a random constraint satisfaction problem (rCSP), where the vertices $[n]$ constitute the variables, while the $m$ rows of $A$ take over the role of the (linear) constraints. 
To be completely explicit about the connection, assume that $\FF = \FF_q$ for a finite field $\FF_q$ with $q$ elements, even though this is not assumed in the rest of the article.
In this case, if the rank of $A$ over $\FF_q$ is equal to $r$, the number of solutions of the linear system $Ax=0$  is $q^{n-r}$, and determining the rank of a matrix over $\FF_q$ becomes equivalent to calculating the number of solutions of the associated rCSP. Therefore, ideas that have been developed for rCSPs can be adapted to the study of random matrices \cite{coja2022rank}, while on the other hand, specific random matrices may shed light on possible phenomena in rCSPs \cite{coja2022sparse}. 

\paragraph{The cavity approach.} Particularly successful contemporary approaches in the mathematical study of rCSPs include statistical-physics inspired proofs that are based on the celebrated cavity method of Mézard, Parisi and Virasoro (see e.g. \cite{mezard2009information}). Simply speaking, when employed in the context of random matrices, the cavity method would prescribe a comparison of two linear systems whose size differs by one. While, due to technical obstructions, we do not follow this prescription blindly in the present article, we will employ it from an alternative perspective in the sense that we will successively shrink $\bm{A}_n$. On a high level, to determine the rank of $\bm{A}_n$, we will iteratively remove rows and columns from $\bm{A}_n$, and then sum the rank decreases to  obtain the rank of the original matrix. Since $\bm{A}_n$ is a symmetric matrix, and it is natural to maintain a certain level of self-similarity, and therefore symmetry, at each step, we will simultaneously eliminate rows with their symmetric columns.

While we will not compare two copies $\bm{A}_{n}$ and $\bm{A}_{n+1}$ of the original model, but instead employ a graph decomposition approach, we will still build upon two main ingredients for a successful implementation of the cavity method: A sufficiently good understanding of the marginal distributions of the uniform distribution over solutions, and  an appropriate correlation decay property. In this section, we focus on the marginal structure, while we will show how to use decorrelation properties in our proof without actually proving decorrelation in the original model in the following section.

\paragraph{Marginals and frozen variables.} The assumption of linear constraints, which is present in rCSPs associated to random matrices, allows for key simplifications in comparison to rCSPs whose constraints do not have algebraic structure. Returning to the study of solutions to the homogeneous linear system $Ax=0$ over $\FF_q$, it is straightforward to see (see e.g. \cite[Lemma 2.3]{ayre2020satisfiability}) that each coordinate $x_i$ has one of the following two types:
\begin{enumerate}
    \item The coordinate is zero in all solutions;
    \item For any $s \in \FF_q$, there is exactly the same number of solutions in which the coordinate assumes the value $s$. 
\end{enumerate}
This observation encodes that the uniform measure over the solutions to $Ax=0$ has reasonably simple coordinate marginal distributions.
Moreover, in the language of rCSPs, coordinates of the first type correspond to so-called \textit{frozen variables}. These are more generally defined as variables that take the same value within one solution cluster 
\cite{budzynski2019biased,mezard2009information}.

While the actual idea is to remove rows (and columns), for the proofs, we will set matrix rows and columns to zero instead of removing them. In this sense, in the following, for any matrix $A\in \mathbb{F}^{m\times n}$, we write $A\abc{i_1,\ldots,i_j;k_1,\ldots,k_\ell}$ for the matrix that is obtained from $A$ by replacing all rows $i_1,\ldots,i_j$ and columns $k_1,\ldots,k_\ell$ of $A$ by zero rows and columns. Even though we consider the operation of simultaneously removing row and column $i$, for reasons that will become apparent in the next section, we do not assume that $A$ is symmetric in the following.

\paragraph{Types of variables and rank.} In our formalisation of the decomposition procedure, the structure of the solution space, and particularly (some, even though not complete) knowledge of the proportion of frozen variables, will play a key role. 
\begin{remark}\label{re_unitv-frozen}
    Let $e_n(i)\in\FF^{1\times n}$ be the $i$th unit vector, i.e., $e_n(i)$ has a $1$ in the $i$th coordinate and zeros everywhere else. Then, a more explicit connection between frozen variables and the rank of $A$ is as follows: 

\begin{enumerate}
    \item If $e_n(i)$ belongs to the row space of $A$, then $e_n(i)\cdot x=0$ for every solution of $Ax = 0$. Hence, coordinate $i$ is frozen. Conversely, if the $i$th coordinate of all solutions is $0$, then $e_n(i)$ is in the row space of $A$. Thus, $i$ is frozen in $A$ if and only if $e_n(i)$ belongs to the row space of $A$.

    \item On the other hand, certainly, $e_n(i)$ is not in the row space of $A\abc{;i}$. Thus, if $e_n(i)$ belongs to the row space of $A$, the latter strictly contains the row space of $A\abc{;i}$. Hence, $i$ is frozen in $A$ if and only if $\rk{A}-\rk{A\abc{;i}}=1$.
\end{enumerate}
This observation is valid for all fields $\FF$, and the distinction of coordinate marginals can be related to the rank change upon removal of a row and its corresponding column via the telescoping sum
\begin{align}\label{eq-re-unitv-frozen}
    \rank_{\FF}\bc{A}-\rank_{\FF}\bc{A\abc{i;i}}=\rank_{\FF}\bc{A^T}-\rank_{\FF}\bc{A^T\abc{;i}}+\rank_{\FF}\bc{A\abc{i;}}-\rank_{\FF}\bc{A\abc{i;i}}.
\end{align}
\end{remark}
While knowledge of the marginal type (frozen / not frozen) of coordinate $i$ in $A^T$ is sufficient to evaluate $\rank_{\FF}\bc{A^T}-\rank_{\FF}\bc{A^T\abc{;i}}$, the first summand in \eqref{eq-re-unitv-frozen}, for the second one, we apparently need to consider its type with respect to the matrix $A\abc{i;}$. Instead of doing so directly, we introduce the following more fine-grained partition of frozen variables in $A$  that also takes the marginal type of $i$ in $A\abc{i;}$ into account: 

\begin{definition}[Frozen status {\cite[Definitions 2.3 and 2.12]{HofMul25}}]\label{dsc}
For any matrix $A\in \mathbb{F}^{m\times n}$ and $i \in [m\wedge n]$, we say that
\begin{enumerate}[label=(\roman*)]
  \item $i$ is \textbf{frozen} in $A$, or $i$ is a frozen variable in $A$, if the unit row vector $e_n(i)$ is in the row space of $A$ and we denote by $\mathcal{F}(A)$ the set of all frozen variables in $A$;
  \item $i$ is \textbf{ frailly frozen} in $A$ if $i\in\mathcal{F}(A)\backslash \mathcal{F}\bc{A\abc{i;}}$;
  \item $i$ is \textbf{firmly frozen} in $A$ if $i\in \mathcal{F}\bc{A\abc{i;}}$.
  \end{enumerate}
We refer to whether $i$ is not frozen, frailly frozen, or firmly frozen in matrix $A$ as its \textit{frozen status}.
\end{definition}
The theoretically possible fourth case where $i\in \mathcal{F}\bc{A\abc{i;}} \backslash \mathcal{F}(A)$ is void, as $i$ being frozen in $A\abc{i;}$ actually implies that $i$ is frozen in $A$. Indeed, if $e_n(i)$ is in the row space of $A\abc{i;}$, then it also lies in the row space of $A$. 
As a consequence, any coordinate $i$ is either not frozen, frailly frozen, or firmly frozen in $A$.

Therefore, to compute the rank difference $\rank_{\FF}\bc{A}-\rank_{\FF}\bc{A\abc{i;i}}$, it is sufficient to determine the more fine-grained frozen status of $i$ in $A$ and $A^T$.  Since there are three different frozen statuses for either matrix, it seems that there are $9$ different categories to be considered. However, by \cite[Proposition 4.5]{HofMul25}, $i$ is frailly frozen in $A$ if and only if $i$ is frailly frozen in $A^T$. Therefore, we can partition the set of coordinates into five disjoint sets as follows:    
\begin{definition}(Typecasting of variables {\cite[Definition 2.13]{HofMul25}})\label{dxyzuv}
For any matrix $A\in \mathbb{F}^{m\times n}$, we partition the set $[m\wedge n]$ into
\begin{enumerate}[label=(\roman*)]
    \item the set $\cX(A)$ of frailly frozen variables;
  \item the set $\cY(A)$ of variables that are firmly frozen in both $A$ and $A^T$;
  \item the set $\cZ(A)$ of variables that are neither frozen in $A$ or $A^T$;
  \item the set $\cU(A)$ of variables that are not (firmly) frozen in $A$ and firmly frozen in $A^T$;
  \item the set $\cV(A)$ of variables that are firmly frozen in $A$ and not (firmly) frozen in $A^T$.
\end{enumerate}
For each $i \in [m \wedge n]$ and $W\in\cbc{X,Y,Z,U,V}$, we refer to the type of $i$ to be $W$ if $i\in\cW(A)$.
\end{definition}
Finally, with \Cref{dxyzuv}, we have the following explicit relation between the rank decrease and the variable types in $A$:
\begin{lemma}[{\cite[
Lemma 4.7]{HofMul25}}]\label{lr}
For any $A \in \FF^{m \times n}$ and $ i\in [m\land n]$,
\begin{enumerate}[label=(\roman*)]
  \item \makebox[14em][l]{$i \in \cY(A)$}$ \Longleftrightarrow  \qquad\rank_{\FF}(A)-\rank_{\FF}(A\abc{i;i})=2$;
  \item \makebox[14em][l]{$i \in \cX(A) \cup \cU(A) \cup \cV(A) $}$ \Longleftrightarrow \qquad \rank_{\FF}(A)-\rank_{\FF}(A\abc{i;i})=1$;
  \item \makebox[14em][l]{$i \in \cZ(A)$}$  \Longleftrightarrow  \qquad\rank_{\FF}(A)-\rank_{\FF}(A\abc{i;i})=0$.
\end{enumerate}
As a consequence, the rank difference between $A$ and $A\abc{i;i}$ satisfies
\begin{equation}\label{lr1}
\begin{aligned}
\rank_{\FF}(A)-\rank_{\FF}(A\abc{i;i})&=\ind\{i \in \cX(A)\}+2\cdot \ind\{i \in \cY(A)\}+\ind\{i \in \cU(A)\}+\ind\{i \in \cV(A)\}.
\end{aligned}
\end{equation}
\end{lemma}

\paragraph{Relating the type of a vertex to the types of its neighbors.} \Cref{lr} leaves us with the task of determining the type distribution in $\bm A_n$, which at first sight does not appear to be any easier than tracing the original rank change. Indeed, we will not aim to determine the type distribution completely, but rather to characterize it via a set of fixed-point equations by relating the variable types in $A\abc{i,i}$ to the variables types in $A$. 
For example, if  ${\rm supp}(b)={\rm supp}(b^T)=\left\{i\in [n]\colon b_i\neq 0\right\}$ for $b\in\FF^{1\times n}$, then firmly frozen variables in $A$ can be (almost) characterized as follows:
\begin{align}
\mbox{$i$ is firmly frozen in $A$}\label{im_fnf}\quad &\xLongleftrightarrow{1}  \quad
     \mbox{ $e_n(i)$ is in the row space of $A\abc{i;}$}\\ &\xLongleftrightarrow{2}  \quad \mbox{$A\abc{i;}(,i)$, the $i$th column of $A\abc{i;}$, is not in the column space of $A\abc{i;i}$}\nonumber\\
     &\xLongrightarrow{\ \ 3\ \ \ }  \quad \mbox{$\exists j\in\supp{A\abc{i;}(,i)}$ such that $e_n(j)^T$ is not in the column space of $A\abc{i;i}$}\nonumber\\
     &\xLongleftrightarrow{4}  \quad \mbox{$\exists j\in\supp{A\abc{i;}(,i)}$ such that $j$ is not frozen in $A\abc{i;i}^T$.}\nonumber
\end{align}
Note that the implications $\xLongleftrightarrow{1}$ and $\xLongleftrightarrow{4}$  follow directly from \Cref{dsc}.  
We now explain how to derive $\xLongrightarrow{2}$:  
If $e_n(i)$ lies in the row space of $A\abc{i;}$, then if we append $e_n(i)$ to the bottom of $A\abc{i:}$, the rank of the resulting matrix will remain unchanged. However, since the $i$th column of the augmented matrix cannot be expressed as a linear combination of the other columns in that matrix, this preservation of rank holds true only if the $i$th column of $A$ cannot be expressed as a linear combination of the other columns in $A\abc{i;}$ as well.  
The implication $\xLongleftarrow{2}$ follows from the reverse of the argument for $\xLongrightarrow{2}$.  
The implication $\xLongrightarrow{3}$ follows from the fact that $A\abc{i;}(,i)$ can be expressed as a linear combination of $\{e_n(j)^T\}_{j \in \supp{A\abc{i;}(,i)}}$.

If the converse of the third implication in \Cref{im_fnf} (i.e., $\xLongleftarrow{3}$) was true, we could infer whether $i$ is firmly frozen in $A$ based on whether all $j\in\supp{A\abc{i;}(,i)}$ are frozen in $A\abc{i;i}^T$. Moreover, if we could establish similar properties for non-frozen and frailly frozen variables, and if the type distributions in $A$ and $A\abc{i;i}$ were similar, we would be in a good position to derive \textit{fixed-point equations} for the type distribution of $A$. 

Sadly, the converse of the third implication of \Cref{im_fnf} is \textbf{not} generally true. However, if the uniform measure over solutions does not exhibit excessive long-range correlations between subsets of variables of bounded size,
it is likely to be true. In the next section, we will explain how to turn $A$ into a matrix with such a convenient solution space structure by an appropriate perturbation.

Under the assumption of full equivalence of all statements in \Cref{im_fnf}, if we successively turn the rows and columns in $\bm{A}_n$ into zero rows and columns, 
\Cref{lr1} provides a route towards a successful  rank estimation of $\bm{A}_n$. The first remaining challenge lies in determining an appropriate order in which to replace the rows and columns by zero rows and columns, that simultaneously allows inference of the types of the associated coordinates. 
The second remaining challenge lies in the derivation of meaningful fixed-point equations for the type distributions.

\subsection{Perturbation matrices}\label{sec_permat}
In this section, we introduce a matrix perturbation that will enable us to establish the validity of the converse of the third implication of \Cref{im_fnf}. We start this section with the following concept of linear relations from \cite{coja2022rank}:
\begin{definition}[Linear relations {\cite[Definition 2.1]{coja2022rank}}]\label{d2}
  Let $A \in \FF^{m \times n}$.
\begin{enumerate}[label=(\roman*)]
  \item A set $\emptyset \not= I\subseteq [n]$ is a \textbf{relation} of $A$ if there exists a row vector $y\in \mathbb{F}^{1\times m}$ such that $\emptyset\neq{\rm supp}(y A) \subseteq I$. If furthermore ${\rm supp}(y A) = I$, then we call $y$ a \textbf{representation} of $I$ in $A$.
  \item A relation $I \subseteq [n]$ is a \textbf{proper relation} of $A$ if $I\backslash \mathcal{F}(A)$ is a relation of $A$. We denote by $\PR(A)$ the set of proper relations of $A$.\label{it_def_proper_rela}
  \item\label{it_deltalfree} For $\delta>0,\ell\geq 2$, we say that $A$ is \textbf{$(\delta,\ell)$-free} if there are no more than $\delta n^\ell$ proper relations $I \subseteq [n]$ of size $\abs{I}=\ell$.
\end{enumerate}
\end{definition} 
For our purposes, matrices with few proper relations have the following desirable properties: Assume that specifically, the matrix $A\abc{i;i}^T$ is $(\delta,\ell)$-free, where $\ell$ is the size of the support of $A\abc{i;}(,i)$, the $i$th column of $A\abc{i;}$. If $\supp{A\abc{i;}(,i)}$ is chosen approximately uniformly over all subsets of $[n]$ of size $\ell$, then \ref{it_deltalfree} yields that $\supp{A\abc{i;}(,i)}$ is a proper relation of $A\abc{i,i}^T$ with probability tending to $0$ as $\delta$ tends to $0$. We may thus assume that it either forms no relation and therefore contains no frozen variables, or that \textit{all} of its elements are frozen in $A\abc{i;i}^T$. In the first case, clearly both the assumption and the consequence of $\xLongleftarrow{3}$ in \Cref{im_fnf} are true. In the second case for $j\in\supp{A\abc{i;}(,i)}$, \emph{all} $e_n(j)^T$ are in the column space of $A\abc{i;i}$, so that by an appropriate linear combination, both the assumption and the consequence of $\xLongleftarrow{3}$ in \Cref{im_fnf} are false. Morally, we thus have the desired equivalence in the converse of the third implication in \Cref{im_fnf}.

While it is far from clear, or even true, that the matrices we are working with are $(\delta,\ell)$-free for appropriate choices of $\delta$ and $\ell$, \cite{coja2022rank,HofMul25} introduce a way to make \emph{any} matrix and its transpose  $(\delta,\ell)$-free without a major change in its rank. The main idea is to attach  perturbation matrices as follows:
\begin{definition}[Perturbation matrices and canonical perturbation]\label{d34}
Fix $P>0$ and let $\bm{\theta} :=(\bm{\theta}_r,\bm{\theta}_c)$ where $\bm{\theta}_r$ and $\bm{\theta}_c$ are independent and uniformly distributed on $[P]$. 
\begin{enumerate}
    \item Let $\THETA_r[\bm{\theta}_r, n]\in \FF^{\bm{\theta}_r\times n}$ be a matrix such that each \emph{row} has exactly one $1$ in a  uniformly chosen position among the $n$ possibilities. We assume that this choice is independent of everything else. All other entries of the matrix are $0$. 
\item Let $\THETA_c[n, \bm{\theta}_c]\in \FF^{n\times\bm{\theta}_c}$ be a matrix such that each \emph{column} has exactly one $1$ 
in a  uniformly chosen position among the $n$ possibilities. We assume that this choice is independent of everything else. All other entries of the matrix are $0$.
\end{enumerate}
For $A \in \FF^{n \times n}$, we write
$$A[\bm{\theta}]=\begin{pmatrix}A&\THETA_c[n, \bm{\theta}_c]\\\THETA_r[\bm{\theta}_r, n ] & 0_{\bm{\theta}_r\times \bm{\theta}_c} \end{pmatrix}.$$
We further call $A[\bm{\theta}]$ the perturbation of $A$ and abbreviate $\THETA:=(\THETA_r[\theta_r, n],\THETA_c[n, \theta_c])$.
\end{definition}
Indeed, the following proposition shows that after attachment of the perturbation matrices, both the perturbation of $A$ and its transpose will become $(\delta,\ell)$-free with high probability as the dimension tends to infinity:
\begin{proposition}[Perturbation eliminates most short proper relations  {\cite[Proposition 2.10]{HofMul25}}]\label{p1}
Fix $\delta>0, L\in\NN_{\geq 2}$. Then with $A[\vth]$ as in \Cref{d34},
    \begin{equation}\label{e1_rep}
    \lim_{P\to\infty}\lim_{n\to\infty}\sup_{A\in \FF^{n \times n}}\mathbb{P}\left( \text{\rm $A[\bm{\theta}]$  or $A[\bm{\theta}]^T$ is not $(\delta,\ell)$-free for some $2\leq \ell\leq L$}\right) = 0.
    \end{equation}
\end{proposition}
Hence, in our proofs, we will work with the perturbation of $\bm{A}_n$ rather than $\bm A_n$ itself to assume $(\delta,\ell)$-freeness. We can do so, as conveniently, the rank difference between $\bm{A}_n$ and $\bm{A}_n[\bm{\theta}]$ is bounded by $\vth_r+\vth_c\leq 2P$ and therefore asymptotically negligible. Consequently, the lower bound of \Cref{coro_main} follows if
\begin{align}\label{ineq_rk_l}
\liminf_{P\to\infty}\liminf_{n\to\infty}\Erw\brk{\frac{1}{n}\rk{\bm{A}_n\pth}}
    \geq\min_{\alpha\in [0,1]}R_\psi(\alpha), \mbox{ uniformly in $\cbc{J_n}_{n\geq 1}$.}
\end{align}

\subsection{The choice of the next vertex: Graph exploration}\label{sec_cm_exp}
Now that we have the decomposition idea in place, along with an expression of the rank decrease in terms of the different variable types \cref{lr1} and a morally true characterization of the type of the removed vertex in terms of the types of its neighbors as in \cref{im_fnf}, the `only' thing that needs specification is the choice of the vertex that is to be removed at each step.  
Here, the first idea that might spring to mind is to remove vertices uniformly at random, and indeed, this is what was done in \cite{HofMul25}.  
However, in the present case, this turns out to not align well with the properties we require, 
 as we aim 
to select the next vertex $i$ such that the probability for $i$ to have any type in $\bm{A}_n[\vth]$
is approximately the same as the probability for a uniformly chosen neighbor of $i$ in $\bm{A}_n\abc{i;i}[\vth]$ (see \cref{im_fnf}).   

We can get an idea as to what the crucial point is by the following consideration: As the degree is the distinguishing feature of vertices in the configuration model, it seems natural to suspect that the key point is to choose the next vertex $\bm{i}$ in such a way that its \emph{degree distribution} in \cm{n} is the same as the degree distribution of its neighbors $\bm{j}$ in the subgraph of \cm{n} excluding vertex $\bm{i}$. Indeed, if such a property holds, the type distributions of $\bm{i}$ and $\bm{j}$ should be close. This similarity will then allow us to derive fixed-point equations on the type distribution of $\bm i$, which will eventually tell us what this type distribution looks like. From this, we could trace the rank decrease when setting the $\bm i$th row and column to zero using \Cref{lr1}. 

However, as consideration of the special case of a random $r$-regular graph reveals, we cannot expect that such a procedure always exists, since the degree of $\bm i$ is always equal to $r$ for $r$-regular graphs, while the degrees of its neighbors excluding $\bm i$ will be $r-1$ in the subgraph with high probability. More generally, as we are interested in randomly chosen vertices and their neighborhoods, the local limit of \cm{n} might give us a hint on how to solve our problem. In this local limit, which is the branching process described in \Cref{sec_conmodel_def}, the offspring distributions of the root and the other vertices in the tree differ, unless the former distribution follows a Poisson distribution, as in the Erd\H{o}s-R\'enyi random graph case \cite{HofMul25}. Indeed, \emph{in all other cases}, 
if we select a uniform vertex $\bm i$ from the set $[n]$ and a uniform neighbor $\bm j$ of $\bm i$, then the number of neighbors of $\bm i$ 
and the number of neighbors of $\bm j$ other than $\bm i$ 
do not asymptotically follow the same distribution. Therefore, one cannot expect that the type distributions of $\bm i$ and $\bm j$ in the corresponding perturbed adjacency matrices are similar, contrary to the Erd\H{o}s-R\'enyi random graph case.

However,  
all \emph{non-root} vertices have the same offspring distribution. So we somehow want to remove the vertices in the branching process other than the root, which leads to the following removal algorithm:
\begin{enumerate}
    \item Choose a vertex $v$ according to some law (for example, uniformly at random or with probability proportional to its degree), and remove it from the graph.

    \item Choose a vertex from the remaining vertices among the neighbors of the previously removed vertices in \cm{n} and remove it. Repeat this step until there is no such vertex.

    \item Repeat from the first step until all vertices have been removed. 
\end{enumerate}
From the limiting branching process we see that each vertex that has been removed \emph{by the second step} has approximately the same degree distribution as its children in the corresponding subgraph, so we may be able to determine its type distribution as outlined previously.
So, if \emph{almost all} the vertices in the removal algorithm are removed in the second step, we might consider our task done, as we should be able to determine the type distribution of almost all the removed vertices. Using equation \eqref{lr1}, we could compute their contributions to the rank decrease when successively replacing the corresponding rows and columns by zero rows and columns. 

Sadly, if we actually remove \emph{all} vertices in this way, the vertices removed by the first step do form a non-negligible proportion. On the other hand, if we stop the removal procedure early, 
the proportion of vertices that have been removed in the first step until this point might still be negligible if the graph has a giant component containing a positive proportion of the graph. Indeed, as it turns out, the largest proportion we can choose to ensure this property is equal to the proportion of vertices in the largest connected component of \cm{n}, which we refer to as the \textit{giant component} \cite{van2022random}.

\paragraph{Graph exploration in continuous time \cite{janson2009new}.}
We formalize the previous considerations in terms of the graph exploration from \cite{janson2009new}, which had originally been introduced to study the size of the giant component of the configuration model. We will use it to perform an educated decomposition of the graph, which translates to successively turning rows and columns of the adjacency matrix to all-zero rows and columns. The exploration process, which is continuous-time, is defined as follows:

At any time $t \geq 0$, a vertex can either be \textbf{sleeping} or \textbf{awake}, while a half-edge can be \textbf{sleeping}, \textbf{active} or \textbf{dead}. Sleeping or active half-edges are also called \textbf{living}. We interpret awake vertices as the vertices that have been explored and dead half-edges as paired half-edges. In contrast, active half-edges are connected to explored vertices, but have not yet been paired.

To introduce the pairing order, we assign each half-edge $h$ an i.i.d.\ (independent and identically distributed) random lifetime $\bm{E}_h$, where $\bm{E}_h$ 
is an exponential random variable with mean $1$. Each half-edge dies  spontaneously when the time exceeds its lifetime. At time $t=0$, all the half-edges and vertices are sleeping and we explore the graph following  three steps, which are the same as those in \cite[Section 4]{janson2009new}:
\begin{enumerate}[label=\textbf{Step \arabic*}]
\item\label{exploration_step1} When there is no \textbf{active} half-edge (as in the beginning), we instantaneously choose a half-edge  uniformly at random among all \textbf{sleeping} half-edges.  We awaken the vertex it belongs to and activate all its half-edges.
If there is no sleeping half-edge left, the process stops; the remaining sleeping vertices are all
isolated and we have explored all other  components.

\item \label{exploration_step2} Pick an \textbf{active} half-edge uniformly at random among 
all \textbf{active} half-edges and change its status to dead.
\item \label{exploration_step3}Wait until the next half-edge dies because of the time exceeding its lifetime. This half-edge is
joined (on paired) to the one killed in the previous \ref{exploration_step2} to form an edge of the
graph. If the vertex it belongs to is \textbf{sleeping}, we change the status of this vertex
to awake and all its remaining adjacent half-edges to active. Repeat from \ref{exploration_step1}.
\end{enumerate}
While it might not be completely obvious that this procedure is equivalent to the definition of the configuration model from \Cref{sec_model_def}, the equivalence is a consequence of the fact that the order in which the half-edges are paired does not matter (see for example \cite[Section 7.2, Lemma 7.6]{van2017random}). 

Further, comparing the graph exploration and the previous removal algorithm, one can see that they essentially describe the same thing:  The next awakened vertex in the graph exploration is just the next removed vertex in the above removal algorithm.  In other words, the graph exploration meets our needs.

Returning to the terminology of the graph exploration, when exploring the giant component, excluding the first vertex that is awakened, each subsequent vertex $v$ is awakened according to \ref{exploration_step3}, which corresponds to the second step in the removal algorithm. This allows us to derive fixed-point equations to estimate the distribution of the type of vertex $v$ in the perturbed adjacency matrix when we zero out the elements in the corresponding rows and columns of the previously awakened vertices. As a result, we obtain a good estimate of the rank decrease. 

As mentioned earlier, this approach is not applicable once we have finished exploring the giant component, as the local limit of a vertex chosen uniformly from the vertices outside the giant component falls into the subcritical case discussed in \Cref{sec_local_str} \cite{van2022random}. Consequently, the proportion of vertices awakened in \ref{exploration_step1}, corresponding to the first step in the removal algorithm, is no longer negligible, and we cannot use our previous approach.

Fortunately, after the giant component has been explored, the degree distribution in the subgraph induced by the sleeping vertices falls into the subcritical case \cite{van2022random}, which we can handle as we have shown in \Cref{sec_local_str}. In this way, we eventually obtain the expected rank of $\bm{A}_n$.

\subsection{Discussion and open problems}\label{sec_discussion}
\subsubsection{Other types of configuration models}
Several variants of the configuration model exist (see \cite{van2017random}), out of which we work with the so-called \emph{erased configuration model}. In this short section, we explain how \Cref{coro_main_basic} applies to other variants as well. The reason for this is that a slight modification of the graph or the adjacency matrix does not change the rank of the adjacency matrix by too much: 
\begin{remark}\label{remark_basic_remove_add} 
    For any matrix $A$, replacing a zero row or column by an arbitrary row or column will increase the rank by $0$ or $1$. Assume that graph $G_2$ is obtained from graph $G_1$ by removing  $m_1$ arbitrary edges $e_1,\ldots,e_{m_1}$ and adding $m_2$ arbitrary edges $e_1',\ldots,e_{m_2}'$. Denote by $x_i$ one of the endpoints of $e_i$ and by $y_j$ one of the endpoints of $e_j'$, and let $I:=\bc{\cup_{i=1}^{m_1}\cbc{x_i}}\cup\bc{\cup_{j=1}^{m_2}\cbc{y_j}}$. If we set all elements in the rows and columns with coordinates in $I$ to zero in the adjacency matrices of $G_1$ and $G_2$, the resulting matrices become identical. As we have set at most $2(m_1+m_2)$ rows and columns to zero, the rank difference between the two adjacency matrices is upper bounded by $2(m_1+m_2)$.
\end{remark}
\textbf{Configuration model with self-loops:} By \Cref{remark_basic_remove_add}, instead of consider the erased \cm{n}, we can take the self-loops in \cm{n} into account. If we define 
    \[\bm{A}_{n}(i,i):=\teo{{\mbox{there is at least a self-loop on $i$ in \cm{n}}} }J_n(i,i)\]
    rather than $0$ in \cref{ec1}, then \cref{e01a} still holds. Indeed, the rank difference between the matrices under the old and the new definitions is upper bounded by the number of self-loops in the graph since we can make all the elements zero in the rows with coordinates in $\bm{SL}:=\cbc{i\in[n]:\mbox{there is at least a self-loop on $i$ in \cm{n}}}$ and the two adjacency matrices become the same. By \cite[Proposition 7.11]{van2017random}, the size of $\bm{SL}$ is negligible compared to $n$. Then \Cref{remark_basic_remove_add} yields that the rank difference between the two adjacency matrices is negligible as well, so their asymptotic ranks are the same. An analogous result applies to multiple edges. Indeed, \cite[Proposition 7.11]{van2017random} also shows that the number of such edges is negligible compared to $n$. Hence, regardless of the values of the corresponding entries, the asymptotic rank remains unchanged.

\textbf{Configuration model conditioned on being simple:} Suppose that we condition \cm{n} on being a simple graph, i.e., not having any self-loops and multiple edges, and denote by $\bm{A}_n$ the (weighted) adjacency matrix of the resulting model. Under the additional assumption $\sum_{k\geq k} k^2p_k<\infty$, \Cref{coro_main_basic} for the conditioned model follows directly from the unconditional result and \cite[Corollary 7.17]{van2017random}. 

\subsubsection{The log-concavity \Cref{tech_assumption}}\label{sec-tech_assumption}
While \Cref{tech_assumption} is not needed in the derivation of the asymptotic rank formula for subcritical and critical configuration models, as well as in our derivation of a tight upper bound for supercritical configuration models, our lower bound in the supercritical case depends crucially on \Cref{tech_assumption} (see the proof of \Cref{imp_pro} in \Cref{sec_progt}). In particular, it does not appear to be easy to get rid of \Cref{tech_assumption} following our approach.

As mentioned earlier, in the important case that $\FF=\RR$ and $J_n(i,j)\equiv 1$, the normalized rank of $\bm{A}_{n}$ has been derived in \cite{bordenave2011rank}. Interestingly, \cite{bordenave2011rank}  make an  assumption like \Cref{tech_assumption} in their derivation of a tight \emph{upper bound}, while we need it to get the \emph{lower bound}. It is worth noting that the proof of the upper bound in~\cite{bordenave2011rank}  
relies on the Karp--Sipser algorithm, and in particular does not depend on the specific values of the nonzero entries in the graph. By following the arguments in~\cite[Theorem 13, Lemma 14, Proposition 15]{bordenave2011rank}, one can verify that under Assumptions~\ref{assumption_weaker} and~\ref{tech_assumption}, the following inequality holds even for weighted adjacency matrices over an arbitrary field~$\FF$:
\[
\liminf_{n \to \infty} \Erw\brk{\frac{1}{n} \rk{\bm{A}_n}} \leq \min_{\alpha \in [0,1]} R_\psi(\alpha),
\]
and thus obtain an upper bound for our case as well, \emph{provided} that \Cref{tech_assumption} holds.

In the case   $\FF=\RR$ and $J_n(i,j)\equiv 1$, Bordenave \cite{bordenave2016spectrum} later removed the log-concavity assumption that was used in \cite{bordenave2011rank} to derive the upper bound, by using the convergence of the spectral measure under the Kolmogorov-Smirnov distance. In this case, the assumption is thus not necessary for the rank formula to hold true. 
It would be very interesting to see if \Cref{tech_assumption} can be removed in \Cref{coro_main} as well.

\subsubsection{The connection between the Karp-Sipser core and the adjacency matrix for Erd\H{o}s-R\'enyi random graphs}
For (unweighted) Erd\H{o}s-R\'enyi random graphs, the connection between their Karp-Sipser core and adjacency matrix has been made explicit fairly recently. In \cite{DemGla24}, DeMichele, Glasgow and Moreira derive a precise formula for the rank in the near sparse regime when $p=\omega(1/n)$ by showing that with high probability, the corank of the adjacency matrix is exactly \emph{equal} to the number of isolated vertices that remain after the Karp-Sipser leaf-removal procedure. In the more challenging regime where $p=\Theta(1/n)$, among other things, Glasgow, Kwan, Sah and Sawhney \cite{glasgow2023exact}  prove that with high probability, the corank of the adjacency matrix is equal to the number of isolated vertices that remain after the Karp-Sipser leaf-removal procedure, \emph{plus} an extra term that is due to the presence of certain cycles in the Karp-Sipser core.  
These two works reveal the intricate relation between the Karp-Sipser core and the rank of the adjacency matrix of the Erd\H{o}s-R\'enyi random graph. Glasgow, Kwan, Sah and Sawhney \cite{glasgow2023exact} also remark that extending their methods to the general case, where both the field of concern and the nonzero entries in the matrix can be chosen arbitrarily, is not an easy task. It would be very interesting to investigate whether such a precise formula still holds in this general case.

\section{Preliminaries}
\subsection{Simpler model with deterministic degree sequence}\label{sec_sim_model}
\Cref{assumption_weaker} imposes very general regularity conditions on the degree sequence for the configuration model. However, this level of generality makes our proofs quite involved. In this section, we introduce a simpler configuration model, which imposes a much stronger restriction on the degree sequence and simplifies the proof of \Cref{coro_main}. We will generalize the result to the setting described in \Cref{assumption_weaker} in \Cref{sec_proof9.1}.

To handle the case supercritical case when $\sum_{k \geq 0} k(k - 2)p_k > 0$, we fix a number $K \in \NN_{\geq 3}$, which will serve as a uniform bound on the vertex degrees (when $K \leq 2$, the inequality $\sum_{k \geq 0} k(k - 2)p_k > 0$ does not hold.). Apart from \Cref{assumption_weaker}, we further assume that the degree sequence $\vd$ satisfies the following assumption:
\begin{assumption}[Extra assumption on the degree sequence]\label{assumption_stronger}
Let $\vd$ be a degree sequence that satisfies \Cref{assumption_weaker} with limiting distribution $(p_k)_{k \geq 0}$. Additionally, assume that there exists $K\in\NN_{\geq 3}$ such that the following hold:
    \begin{enumerate}
        \item $\vd$ is non-random. In this case, we write $d$, $d_i$ or $d(n,i)$.

        \item For any positive integer $n$ and $i\in[n]$, $d_i\in\cbc{0, \ldots, K}$.

        \item The probability distribution $(p_k)_{k\geq 0}$ satisfies 
 $\min_{0 \leq k \leq K}p_k:= q>0$.
    \end{enumerate}
\end{assumption}

As we discussed in \Cref{sec_local_str}, the main difficulty of the proof comes from the supercritical case where the local structure of a typical vertex might be an infinite tree. With our simpler model, our main work focuses on the proof of the following proposition:
\begin{proposition}\label{t_main}
    Assume that the degree sequence $d$ satisfies Assumption \ref{assumption_stronger} with a probability  distribution $(p_k)_{k\geq 0}$  satisfying \Cref{tech_assumption} and $\sum_{k\geq 0}k(k-2)p_k>0$. For any field $\FF$, $\Erw\brk{\rank_{\FF}\bc{\bm{A}_n}/n}$ is asymptotically lower bounded by  $\min_{\alpha\in [0,1]}R_\psi(\alpha)$ uniformly in $\bc{J_n}_{n\geq 1}$ in the sense that 
\begin{align*}
\liminf_{n\to\infty}\inf_{J_n\in \syn}\Erw\brk{\frac{1}{n}\rank_{\FF}\bc{\bm{A}_{n}}}\geq \min_{\alpha\in [0,1]}R_\psi(\alpha).
\end{align*}
\end{proposition}
Although it may be possible to prove \Cref{t_main} under a condition that does not require  \Cref{assumption_stronger}, like \cite[Condition 2.1]{janson2009new}, a standard (albeit somewhat tedious) approach in \Cref{sec_proof9.1} shows that the result for our simpler model suffices for the target lower bound under Assumptions \ref{assumption_weaker} and \ref{tech_assumption}. Furthermore, we will demonstrate in \Cref{sec_proof_main} that the proof of \Cref{t_main} is a key component in establishing \Cref{coro_main}. Therefore, we will focus on the proof of \Cref{t_main} from now until \Cref{sec_rank_est}.
\subsection{Notation}\label{sec_pre}
This section can be used as a reference for recurring notation that is used throughout the article.

\paragraph{Sets.}
We write $[\ell] = \cbc{1,2,\ldots,\ell}$ and denote the cardinality of a set $B$ by $|B|$. For two sets $B_1$ and $B_2$, we denote their symmetric difference by
$B_1 \Delta B_2$. If $B$ is a set and $\ell \leq |B|$, then we write $\binom{B}{\ell}$ for the collection of $\ell$-subsets of $B$.

\paragraph{Real numbers and fields.} For $a, b \in \RR$, we write $a\vee b = \max\cbc{a,b}$ and $a\wedge b = \min\cbc{a,b}$. $\FF$ is reserved to denote a generic field, and $\FF^\ast = \FF \setminus \{0\}$ its multiplicative group.

\paragraph{Vectors and matrices.} For $A \in \FF^{m \times n}$, we denote its transpose by $A^T$. For a vector $b=\bc{b_1,b_2,\ldots,b_n}\in \mathbb{F}^{1\times n}$, we let ${\rm supp}(b)={\rm supp}(b^T)=\left\{i\in [n]\colon b_i\neq 0\right\}$. We denote the $i$th standard unit vector in $ \mathbb{F}^{1\times n}$ by $e_n(i)$.

For $s=(s_1,s_2,\ldots,s_\ell) \in \RR^{1 \times \ell}$, define $\|s\|_\infty=\sup_{i\in [\ell]}|s_i|$ and $\|s\|_k=(\sum_{i=1}^{\ell}|s_i|^k)^{1/k}$.

For $A \in \FF^{m \times n}$, we denote
\begin{enumerate}[label=(\roman*)]
    \item the $i$th row of $A$ by $A(i,)$ and the $j$th column of $A$ by $A(,j)$.

    \item the matrix obtained by replacing rows $\ell_1,\ell_2,\ldots,\ell_s$ and columns $\ell_1',\ell_2',\ldots,\ell_t'$ in $A$ by zero rows and columns, respectively, by 
$A\abc{\ell_1,\ell_2,\ldots,\ell_s;\ell_1',\ell_2',\ldots,\ell_t'}$. 
\end{enumerate}

\paragraph{Notions of convergence.}
Throughout the article, 
the order in which limits are taken matters significantly. 
For fixed $\varepsilon >0$ and families of real numbers $(a_{n,P,J_n,s})_{n, P\in \ZZ^+,J_n\in \syn,\ranges}$ (where the parameters $\sigma(\cdot)$ and $\xi$ are introduced in \Cref{sec_graphexp}), we write
\begin{enumerate}[label=(\roman*)]
\item  \makebox[7em][l]{\mbox{$a_{n,P,J_n,s}=o_n(1)$}} $\quad \Longleftrightarrow \quad $ 
 $\forall P\geq 1,J_n\in \syn,\ranges:~$ 
$\lim_{n\to\infty}a_{n,P,J_n,s}=0$;
\item  \makebox[7em][l]{\mbox{$a_{n,P,J_n,s}=\bar{o}_n(1)$}} $\quad \Longleftrightarrow \quad $ $\forall P\geq 1: \quad $ $\lim_{n\to\infty}\sup_{J_n\in\syn,\ranges}\abs{a_{n,P,J_n,s}}=0$;
  \item \makebox[7em][l]{\mbox{$a_{n,P,J_n,s}=\bar{o}_{n,P}(1)$}} $ \quad \Longleftrightarrow \quad  \lim_{P\to\infty}\limsup_{n\to\infty}\sup_{J_n\in\syn,\ranges}\abs{a_{n,P,J_n,s}}=0$.
\end{enumerate}
For a family of \textit{uniformly bounded} random variables $(\bm{b}_{n,P,J_n,s})_{n, P\in \ZZ^+,J_n\in \syn, \ranges}$, we write
\begin{enumerate}[label=(\roman*)]
  \item $  \bm{b}_{n,P,J_n,s}=\oone\qquad \Longleftrightarrow \qquad \mathbb{E}\abs{\bm{b}_{n,P,J_n,s}}=\bar{o}_{n,P}(1)$;
  \item $\bm{b}_{n,P,J_n,s}\geq \oone\qquad \Longleftrightarrow \qquad \bc{\bm{b}_{n,P,J_n,s}}^-=\oone$.
  \item $\bm{b}_{n,P,J_n,s}\leq \oone\qquad \Longleftrightarrow \qquad \bc{\bm{b}_{n,P,J_n,s}}^+ = \oone$.
\end{enumerate}
For a family of events $(\mathfrak B_{n,P,J_n,s})_{n, P\in \ZZ^+,J_n\in \syn,\ranges}$, we say that $\mathfrak B_{n,P,J_n,s}$ occurs with high probability if  $\mathbb{P}\bc{\mathfrak B_{n,P,J_n,s}}=1+\bar{o}_{n,P}(1)$.

We extend the above notions of convergence to families of numbers and events that only depend on subsets of the parameters. For example, for a family of real numbers $(c_{n,P})_{n, P\in\ZZ^+}$, by treating it as constant on the unspecified parameters, we write $c_{n,P}=\bar{o}_{n,P}(1)$ whenever $\lim_{P\to\infty}\limsup_{n\to\infty}\abs{c_{n,P}}=0$.

\section{Graph decomposition}\label{sec_graphexp}
In this section, we summarize the main properties of the continuous-time graph exploration presented in \Cref{sec_cm_exp}, and reformulate it as a discrete vertex-removal model. Within the latter, in \Cref{sec_current_degrees}, we calculate the degree distribution of the induced graph of remaining vertices after removal of $\lfloor sn\rfloor$ vertices and in \Cref{sec_type_new}, we investigate how the next vertex that is to be removed is chosen given the induced graph of remaining vertices. 
As explained, we assume throughout that the configuration model is supercritical, based on the simpler model described in \Cref{sec_sim_model}, so that in particular $p_k=0$ for $k>K$, and
\begin{align}\label{con_super_c}
    \sum_{k=0}^K k(k-2)p_k>0.
\end{align}

\subsection{Preliminaries: Living half-edges and sleeping vertices} 
We first summarize the main results from \cite{janson2009new} that we will build upon in the subsequent analysis. 
The degree distribution of the sleeping vertices at a given time $t$ of the graph exploration has been studied by Janson and Luczak in \cite{janson2009new}, which conveniently applies to our setting:
\begin{remark}\label{rem_janson_applies}
Any degree sequence $\vd$ that satisfies \Cref{assumption_stronger} also satisfies Condition 2.1 from \cite{janson2009new}. Thus, all results of \cite{janson2009new} apply to the present setting.
\end{remark}

To build some intuition, consider the early stages of the exploration. In this time regime, there is a comparatively simple and accurate approximation of the number of sleeping vertices based on the following reasoning: 
If, instead of following the graph exploration, vertices were awakened at the very moment when their first adjacent half-edge dies, then at time $t$, any vertex of degree $k$ would be sleeping \textit{independently} with probability exactly equal to $\eul^{-kt}$. In fact, spontaneously dying half-edges are the dominant mechanism for awakening vertices also in the graph exploration, as long as \ref{exploration_step1} is not performed too often: Let
\begin{align}\label{def_sigma_lambda}
    \sigma(t)=\sum_{k=0}^Kp_k \eul^{-kt} \qquad \text{and} \qquad \lambda(t)=\sum_{k=0}^K kp_k \eul^{-kt}.
\end{align}
Then informally, in a time regime where the effect of \ref{exploration_step1} is negligible, we expect the number of sleeping vertices at time $t$ divided by $n$ to converge to $\sigma(t)$, while we expect the number of sleeping half-edges divided by $n$ to converge to $\lambda(t)$, as $n \to \infty$. 

The time regime where the effect of \ref{exploration_step1} is negligible has been determined in \cite{janson2009new}; we provide an alternative formulation in terms of the generating function $\hat{\psi}$ from \eqref{def_psi_zero} here. A comparison of \cref{def_sigma_lambda} with \eqref{def_psi_zero} reveals that for $t \geq 0$, 
\begin{align}\label{eq_psi_lamda}
    \hat{\psi}(\eul^{-t})=\frac{\lambda(t)\eul^t}{\lambda(0)}.
\end{align}
In \cite{janson2009new} it has been shown when the conditions assumed hold that $\hat{\psi}$ has exactly one fixed point in $(0,1)$:
\begin{lemma}[{\cite[Lemma 5.5]{janson2009new}}]\label{lem_psi0}
     For $(p_k)_{0\leq k\leq K}$ satisfying  \Cref{assumption_stronger} and Condition \cref{con_super_c}, the function $\hat{\psi}$ defined in \eqref{def_psi_zero} has exactly one fixed point $\xi$ in the open interval $(0,1)$. Moreover, $\hat{\psi}(\alpha){>}\alpha$ for $\alpha\in(0,\xi)$ and $\hat{\psi}(\alpha){<}\alpha$ for $\alpha\in(\xi,1)$. 
\end{lemma} 
 Indeed, the value $-\ln\xi$ is relevant in the graph exploration as the limiting time when the giant component is fully explored \cite{van2022random}, 
    which is also the time regime relevant for our rank approach.
To provide a solid basis for the intuition from the beginning of this section, we observe that $-\ln\xi$ also marks a time region where \ref{exploration_step1} is not performed w.h.p.:
\begin{lemma}[{\cite[Proof of Theorem 2.3(i)]{janson2009new}}]\label{lem_janson_s}
 For  any degree sequence $d$ satisfying Assumption \ref{assumption_stronger} and any $\varepsilon\in (0,-(\ln \xi)/2)$, 
\begin{align}\label{s1_notperform_t}
    \PP\bc{\mbox{\rm  \ref{exploration_step1} is performed in time interval $[\varepsilon,-\ln \xi-\varepsilon]$}}=o_n(1).
\end{align}
\end{lemma}
Finally, let us return to the numbers of living half-edges and sleeping vertices throughout the initial phase of the exploration. Let $\bm{S}(t)$ and $\bm{L}(t)$ denote the numbers of sleeping and living half-edges at time $t\geq 0$, respectively, and, for $0 \leq k \leq K$, let $\bm{V}_{k}(t)$ be the number of sleeping degree-$k$ vertices at time $t$. Then $\bm{V}(t):=\sum_{k=0}^K\bm{V}_{k}(t)$ gives the total number of sleeping vertices. 
In line with common practice, we define all these random functions to be right-continuous. From \cite{janson2009new}, we infer the following precise descriptions of $\bm{L}(t)$ and $\bm{V}_{k}(t)$, at least up to time $-\ln\xi$:
\begin{lemma}[{\cite[Lemma 5.1]{janson2009new}}]\label{lem_janson_l}
For any degree sequence $d$ satisfying \Cref{assumption_stronger}, 
\begin{align*}
    \sup_{t\geq 0}\abs{n^{-1}\bm{L}(t)-\lambda(0)\eul^{-2t}}\stackrel{\mathbb{P}}{\longrightarrow} 0.
\end{align*}
\end{lemma}
\begin{lemma}[{\cite[Proofs of Lemmas 5.2 and 5.3]{janson2009new}}]\label{lem_janson_l_2}
For any degree sequence $d$ satisfying \Cref{assumption_stronger}, 
\begin{align*}
 \sup_{t\in [0,-\ln \xi]}\abs{n^{-1}\bm{V}_{k}(t)-p_k\eul^{-kt}}\stackrel{\mathbb{P}}{\longrightarrow}0.
\end{align*}
\end{lemma}
\Cref{lem_psi0,lem_janson_s,lem_janson_l,lem_janson_l_2} are proved in \cite{janson2009new}. However, \Cref{lem_janson_s,lem_janson_l_2} are not stated as separate results, but rather proved along the way. For the sake of completeness, in Appendix \ref{sec_app_graphexp}, we therefore explain in more detail how these two lemmas are derived based on \cite{janson2009new}.
\subsection{From graph exploration to vertex removal} 
\Cref{lem_janson_l_2} shows that during the graph exploration, up to the point where the giant has been explored, we have good control over the numbers $\bm{V}_{k}(t)$ of sleeping vertices of degree $k$. The mechanism that enables this approximation is that whether a vertex is awake or not is independent of the randomness at the other vertices, unless the vertex is awakened by \ref{exploration_step1}. And indeed, in the time regime $[0,-\ln\xi]$, according to \Cref{lem_janson_s}, the effect of \ref{exploration_step1} is negligible.  
However, for the implementation of the intended decomposition procedure and the derivation of type fixed-point equations, we need precise knowledge of the unexplored graph after a fixed number of vertices have been declared awake, 
such as the degree distribution among the sleeping vertices. 
More specifically, we will remove vertices from the graph, following the order that is dictated by the awakening of vertices during the graph exploration. As mentioned, instead of removing rows and columns from the adjacency matrix, we will set all their entries to zero. 
Therefore, in the remainder of this section, we transform Lemmas \ref{lem_janson_s}, \ref{lem_janson_l} and \ref{lem_janson_l_2} from an exploration model into a vertex removal model that is indexed by the proportion of vertices awakened so far (in continuous time).

For $s \in [1/n,1]$, let $\vnu_s$ denote the $\rou{sn}$th vertex that is awakened during the graph exploration and 
$$\verst{s}:=[n]\backslash\cbc{\vnu_r:r\in[1/n,s]} $$ 
the set of sleeping vertices after $\rou{sn}$ vertices have been awakened. 
The induced subgraph of \cm{n} on the vertex set $\verst{s}$ will be denoted by \scm{n}{s}. 

Similarly, let \tcm{n}{t} be the induced subgraph of \cm{n} on the set of sleeping vertices at continuous time $t$. 
\begin{notation}[Vertices and degrees in the reduced models] \label{notation_reduced}
   We call the degree of vertex $i \in \verst{s}$ in 
   \scm{n}{s} its \textbf{current} degree and denote it by $\bar{\vd}_{i,s}$. We also use the term `current degree' to refer to the degree of a vertex $i$ in \tcm{n}{t}. The degree $d_i$ of vertex $i \in [n]$ in \cm{n} will be called the \textbf{original} degree. Moreover, 
\begin{itemize}
    \item[(i)] for $0\leq k\leq K$, $\bar{\bm{V}}_k(t)$ denotes the number of vertices of current degree $k$ in \tcm{n}{t}.
    \item[(ii)] $\verslk{k}{s}$ and $\bverslk{k}{s}$ denote the numbers of vertices of original and current degree $k$ in  \scm{n}{s}, respectively. 
    \item[(iii)] $\edgesl{s}:=\sum_{k=0}^K k\verslk{k}{s}$ and $\edgeli{s}$ denote the numbers of sleeping and living half-edges in \cm{n}, respectively, at the time when the $\rou{sn}$th vertex is awakened.
\end{itemize}
\end{notation}
We next define the stopping times
\begin{align*}
    \vta_{n,s}:=\inf\cbc{t\geq 0:\bm{V}(t)\leq n-\lfloor sn\rfloor}.
\end{align*}
Recall that in the graph exploration, vertices are either awakened one-by-one or in pairs. The latter happens if and only if \ref{exploration_step1} is performed instantaneously after \ref{exploration_step3}. Thus, there are two cases: 
If only \ref{exploration_step3} is performed at time $\vta_{n,s}$, at this time, there are exactly $n-\lfloor sn \rfloor$ sleeping vertices. Hence, $\ftcm{n}{\vta_{n,s}}=\fscm{n}{s}$ and for all $0\leq k\leq K$,
\begin{align}\label{def_vtans} \bm{V}(\vta_{n,s})=n-\lfloor n s\rfloor ,\quad \bm{V}_{k}(\vta_{n,s})=\verslk{k}{s} \quad\mbox{and}\quad \bar{\bm{V}}_k(\vta_{n,s})=\bverslk{k}{s}.
\end{align}
On the other hand, if both \ref{exploration_step3} and \ref{exploration_step1} are performed at time $\vta_{n,s}\neq 0$, at this time, there are exactly $n-\lfloor sn \rfloor-1$ sleeping vertices. Hence, $\ftcm{n}{\vta_{n,s}}=\fscm{n}{s+1/n}$ and for all $0\leq k\leq K$,
\begin{align}\label{def_vtans_2} \bm{V}(\vta_{n,s})=n-\lfloor n s\rfloor-1,\quad \bm{V}_{k}(\vta_{n,s})=\verslk{k}{s+1/n} \quad\mbox{and}\quad \bar{\bm{V}}_k(\vta_{n,s})=\bverslk{k}{s+1/n}.
\end{align}

\begin{definition}[{Approximation of $\vta_{n,s}$}]\label{def_t_s}
Consider the parametrised functions
\begin{align}\label{def_H_s}
    H_s:[0, \infty) \to \RR, \quad H_s(t):= \sigma(t)-1+s, \qquad \text{where } s \in [0,1-p_0).
\end{align}
Observe that for $s\in [0,1-p_0)$, $H_s$ is strictly decreasing with $H_s(0) = s \geq 0$ and $\lim_{t \to \infty} H_s(t) = s+p_0-1 < 0$. We denote the unique zero of $H_s$ by $t_s$. 
\end{definition}

Using \Cref{def_t_s}, we are in the position to transform \Cref{lem_janson_l,lem_janson_l_2,lem_janson_s}  into corresponding results for the vertex-removal process:
\begin{lemma}[Concentration in vertex-removal process]\label{l_concentration}
Fix $\varepsilon\in (0,1/2-\sigma(-\ln\xi)/2)$. Then for any degree sequence $d$ satisfying \Cref{assumption_stronger}:
\begin{itemize}
\item[(i)] For any $\varepsilon'\in (0,\varepsilon)$, $\PP\bc{\exists s\in [\varepsilon',1-\sigma(-\ln\xi)-\varepsilon']: \vta_{n,s} \notin \brk{t_{s-\varepsilon'}, t_{s+\varepsilon'}}}
        = o_n(1)$; 
    \item[(ii)] \label{con_re_2}$\sup_{\ranges}\abs{n^{-1}\verslk{k}{s}- \eul^{-k t_s}p_k}\stackrel{\mathbb{P}}{\longrightarrow}0$;
\item[(iii)] $\sup_{\ranges}\abs{n^{-1}\edgeli{s}-\lambda(0)\eul^{-2t_s}}\stackrel{\mathbb{P}}{\longrightarrow}0$; 
    \item[(iv)] $\PP\bc{\mbox{\rm \ref{exploration_step1} is performed in the process of awakening vertices $\varepsilon n, \ldots, (1-\sigma(-\ln\xi)-\varepsilon) n$}} =o_n(1).$
\end{itemize}
\end{lemma}

\begin{proof}[Proof of \Cref{l_concentration}]
\textit{Proof of (i):} By \Cref{lem_janson_l_2},
    \begin{align*}
        \sup_{t\in [0,-\ln \xi]}\abs{n^{-1}\bm{V}(t)-\sigma(t)}\leq \sum_{k=0}^K\sup_{t\in [0,-\ln \xi]}\abs{n^{-1}\bm{V}_{k}(t)-p_k\eul^{-kt}} \stackrel{\mathbb{P}}{\longrightarrow}0. 
    \end{align*}
With $H_s$ and $t_s$ as in \Cref{def_t_s}, we have the relation $s=1-\sigma(t_s)$. Furthermore, as $\sigma(0)=1$, $t_0=0$. Hence,
    \begin{align*}
        \sup_{s\in [0,1-\sigma(-\ln\xi)]}\abs{n^{-1}\bm{V}(t_s)-1+s}=\sup_{t_s\in [0,-\ln \xi]}\abs{n^{-1}\bm{V}(t_s)-\sigma(t_s)} \stackrel{\mathbb{P}}{\longrightarrow}0.
    \end{align*}
Let $a_n:= \Erw[\sup_{t_s\in [0,-\ln \xi]}\abs{n^{-1}\bm{V}(t_s)-\sigma(t_s)}]$. Since $n^{-1}\bm{V}(t)\in [0,1]$ for all $t \geq 0$, by the dominated convergence theorem, $a_n = o_n(1)$ and by 
Markov's inequality, for any $\varepsilon'\in (0,\varepsilon)$,
    \begin{align}\label{eq_newadd1}
        \PP\bc{\sup_{s\in [0,1-\sigma(-\ln\xi)-\varepsilon']}\bc{n^{-1}\bm{V}(t_{s+\varepsilon'})-1+s+\varepsilon'}\geq \varepsilon'/2} \leq \frac{2 a_n}{\varepsilon'},
    \end{align}
and
    \begin{align}\label{eq_newadd2}
        \PP\bc{\inf_{s\in [\varepsilon',1-\sigma(-\ln\xi)]}\bc{n^{-1}\bm{V}(t_{s-\varepsilon'})-1+s-\varepsilon'}\leq -\varepsilon'/2} \leq \frac{2 a_n}{\varepsilon'}.
    \end{align}
Recall from \cref{def_vtans} and \cref{def_vtans_2} that 
    $n^{-1}\bm{V}(\vta_{n,s})=1-n^{-1}\lfloor ns\rfloor$ or $1-n^{-1}(\lfloor ns\rfloor+1)$ almost surely. We then conclude from \Cref{eq_newadd1} and \Cref{eq_newadd2} that, for $n$ large enough,  
    \begin{align}\label{ineq_vtau}
        &\PP\bc{\exists s\in [\varepsilon',1-\sigma(-\ln\xi)-\varepsilon']: n^{-1}\bm{V}(\vta_{n,s}) \notin \bc{n^{-1}\bm{V}(t_{s+\varepsilon'})+\varepsilon'/4, n^{-1}\bm{V}(t_{s-\varepsilon'})-\varepsilon'/4}} \leq \frac{4 a_n}{\varepsilon'}.
    \end{align}
Since $\bm{V}(\cdot)$ is a nonincreasing function, \cref{ineq_vtau} immediately gives that    
\begin{align}\label{ineq_vtau2}
\PP\bc{\exists s\in [\varepsilon,1-\sigma(-\ln\xi)-\varepsilon]: \vta_{n,s} \notin \brk{t_{s-\varepsilon}, t_{s+\varepsilon}}} = o_n(1).
    \end{align}
    
\textit{Proof of (ii):} By \cref{def_vtans} and \cref{def_vtans_2},  $\verslk{k}{s}=\bm{V}_{k}(\vta_{n,s})$ or $\bm{V}_{k}(\vta_{n,s})+1$ almost surely. The combination of \Cref{lem_janson_l_2}, \cref{ineq_vtau2} and the
dominated convergence theorem gives that for any $\varepsilon' \in (0, \varepsilon)$ and $\tilde a_n := 3\max\{4a_n, 1/n\}$,
\begin{align}\label{eq_difntau}
    &\Erw\brk{\sup_{s\in [\varepsilon,1-\sigma(-\ln\xi)-\varepsilon]}\abs{\frac{\verslk{k}{s}}{n}- \eul^{-k t_s}p_k}}\nonumber\\
    \leq& \Erw\brk{\sup_{s\in [\varepsilon,1-\sigma(-\ln\xi)-\varepsilon]}\abs{\frac{\bm{V}_{k}(t_{s-\varepsilon'})}{n}- \eul^{-k t_s}p_k}}+\Erw\brk{\sup_{s\in [\varepsilon,1-\sigma(-\ln\xi)-\varepsilon]}\abs{\frac{\bm{V}_{k}(t_{s+\varepsilon'})}{n}- \eul^{-k t_s}p_k}}\nonumber\\
    &+\PP\bc{\exists\ s\in [\varepsilon',1-\sigma(-\ln\xi)-\varepsilon']:\vta_{n,s} \notin \brk{t_{s-\varepsilon'}, t_{s+\varepsilon'}}
    }+o_n(1)\nonumber\\
    \leq &\bc{\eul^{-k t_{s-\varepsilon'}}p_k- \eul^{-k t_s}p_k}+\bc{\eul^{-k t_s}p_k-\eul^{-k t_{s+\varepsilon'}}p_k}+ \frac{\tilde a_n}{\varepsilon'}+o_n(1).
\end{align}
On the other hand, since $t_s$ is the unique zero of $H_s(t)=\sigma(t)-1+s$ in $[0,1-p_0)$, by the inverse function theorem, $t_s$ is a continuously differentiable function with respect to $s$. Differentiating both sides of $H_s(t_s)=0$ with respect to $s$ gives that
\begin{align}\label{dif_Ts}
    \frac{\dif t_s}{\dif s}=\bc{\sum_{k=0}^K p_k k \eul^{-k t_s}}^{-1}=\lambda(t_s)^{-1}\in [K^{-1}(1-s-p_0)^{-1},(1-s-p_0)^{-1}],
\end{align}
where we use the fact that $K^{-1}\lambda(t) \leq \sigma(t)-p_0\leq \lambda(t)$ (recall \eqref{def_sigma_lambda}).
 The combination of \cref{dif_Ts} and the mean value theorem then gives that
\begin{align}\label{eq_difexpts-ep}
    \eul^{-k t_{s-\varepsilon'}}p_k- \eul^{-k t_s}p_k\leq k \eul^{-k t_{s-\varepsilon'}}p_k(t_s-t_{s-\varepsilon'})\leq k p_k\varepsilon'\bc{\sigma(-\ln\xi)-p_0}^{-1},
\end{align}
and
\begin{align}\label{eq_difexpts+ep}
   \eul^{-k t_s}p_k-\eul^{-k t_{s+\varepsilon'}}p_k \leq k \eul^{-k t_{s}}p_k(t_{s+\varepsilon'}-t_s)\leq k p_k\varepsilon'\bc{\sigma(-\ln\xi)-p_0}^{-1}.
\end{align}
The desired result now follows from the combination of \cref{eq_difntau} to \cref{eq_difexpts+ep} and taking $\varepsilon' = \varepsilon'(n)=\min\{\varepsilon,\tilde a_n^{1/2}\}$ in \cref{eq_difntau}.

\textit{Proof of (iii):} Given \Cref{lem_janson_l}, the proof is analogous to the proof of item (ii).

\textit{Proof of (iv):}
    By \cref{def_vtans} and \cref{def_vtans_2},
    \begin{align*}
        &\PP\bc{\mbox{\rm \ref{exploration_step1} is performed in the process of awakening vertices $\varepsilon n, \ldots, (1-\sigma(-\ln\xi)-\varepsilon) n$}}\\
        &\leq \PP\bc{\mbox{\rm \ref{exploration_step1} is not performed between $\vta_{n,\varepsilon-1/n}$ and $\vta_{n,1-\sigma(-\ln\xi)-\varepsilon+1/n}$}}.
    \end{align*}  
    Then the desired result is a direct consequence of the combination of \cref{s1_notperform_t} and \cref{ineq_vtau2}.
\end{proof}

\subsection{Current degrees in the removal model}\label{sec_current_degrees}
In this section, 
we prove the following lemma, which provides a good estimate of the number of vertices with current degree $k$ at the time when $\lfloor sn\rfloor$ vertices have been awakened in the graph exploration:
\begin{lemma}\label{lem_bvks}
For each integer $0\leq k\leq K$ and $\varepsilon >0$, uniformly in $\ranges$, 
\begin{align*}
    \mathbb{E}\abs{\frac{1}{n}\bverslk{k}{s}-\sum_{m=k}^K \binom{m}{k}\bc{\frac{\lambda(t_s)\eul^{2t_s}}{\lambda(0)}}^k\bc{1-\frac{\lambda(t_s)\eul^{2t_s}}{\lambda(0)}}^{m-k} \eul^{-m t_s} p_m}= \bar{o}_n(1).
 \end{align*}
\end{lemma}
The proof of \Cref{lem_bvks} requires the following lemma on estimating the ratio of two random variables:
\begin{restatable}{lemma}{appdivide}
\label{lem_app_divide}
    Let $\bc{\bm{X}_n}_{n\geq 1}$, $\bc{\bm{Y}_n}_{n\geq 1}$ and $\bc{\bm{Z}_n}_{n\geq 1}$ be three sequences of random variables defined on the same probability space such that $0\leq \bm{X}_n\leq C\bm{Y}_n$ for some constant $C>0$, $\bm{Y_n}>0$ for all $n$ and $\bc{\bm{Z}_n}_{n\geq 1}\subseteq [0,1]$. Let  
     $\bc{\mathfrak{H}_n}_{n\geq 0}$ be a sequence of events
    with $\PP\bc{\mathfrak{H}_n}\geq 1-a_n$ for all $n \geq 1$ and a sequence $(a_n)_{n\geq 1} \subseteq (0,1]$. Finally, assume that there exist $x \geq 0, y>0$ and $\beta >0$ such that
    \begin{align}\label{ass_lem_app_divide}
        \mathbb{E}\brk{\ensuremath{\mathds{1}}_{\mathfrak{H}_n}\abs{\bm{X}_n-x n^{\beta}}}\leq a_n n^{\beta}\quad\mbox{and}\quad \mathbb{E}\brk{\ensuremath{\mathds{1}}_{\mathfrak{H}_n}\abs{\bm{Y}_n-y n^{\beta}}}\leq a_n n^{\beta}.
    \end{align}
    Then for $n$ such that $a_n\leq y^2$,
    \begin{align*}
        \mathbb{E}\abs{\frac{\bm{X}_n}{\bm{Y}_n}\bm{Z}_n-\frac{x}{y}\bm{Z}_n}\leq \mathbb{E}\abs{\frac{\bm{X}_n}{\bm{Y}_n}-\frac{x}{y}}\leq\frac{(x+y) \sqrt{a_n}}{y(y-\sqrt{a_n})}+\bc{C+\frac{x}{y}}(2\sqrt{a_n}+ a_n).
    \end{align*}
\end{restatable}
\Cref{lem_app_divide} can be seen as a generalization of the comparison of drawing with and without replacement: For example, fixing positive integers $u,v,w,k_1$ and $k_2$, \Cref{lem_app_divide} with $\bm{X}_n = \binom{un}{k_1}\binom{vn}{k_2}$, $\bm{Y}_n = \binom{wn}{k_1+k_2}$, $\bm{Z}_n=1$ gives that
\begin{align*}
    \abs{\frac{\binom{un}{k_1}\binom{vn}{k_2}}{\binom{wn}{k_1+k_2}}-\binom{k_1+k_2}{k_1}\frac{u^{k_1}v^{k_2}}{w^{k_1+k_2}}}=o_n(1).
\end{align*}
The proof of \Cref{lem_app_divide} is deferred to Appendix \ref{app_needproof}. 

\begin{proof}[Proof of \Cref{lem_bvks}]
Let $\mathfrak{L}_{i}$ denote the event that there is no self-loop attached to vertex $i$ in \cm{n} for $j\in[n]$. Similarly, for distinct $i,j\in[n]$, let $\mathfrak{L}_{i,j}$ denote the event that there are no edges between vertices $i$ and $j$ in \cm{n}. Thanks to \Cref{assumption_stronger}, 
\begin{align}\label{est_event_L}
\PP\bc{\mathfrak{L}_{i}^{\mathrm{c}}} \leq \frac{d_i(d_i-1)}{\sum_{k=1}^n d_k-1} 
= \bar{o}_n\bc{1}, 
\qquad
\PP\bc{\mathfrak{L}_{i,j}^{\mathrm{c}}}\leq \frac{d_i d_j}{\sum_{k=1}^n d_k-1} = \bar{o}_n(1), 
\end{align}
i.e., $\mathfrak{L}_{j}$ and $\mathfrak{L}_{i,j}$ are w.h.p. events. Note that both estimates are uniform over the choice of $i, j \in [n]$. 
We next use the second-moment method to prove that $\bverslk{k}{s}$ concentrates around its average.

\textit{Analysis of $\mathbb{E}\brk{\bverslk{k}{s}}$:} 
  Recall that $\verst{s}$ denotes the set of sleeping vertices after $\lfloor sn \rfloor$ vertices have been awakened and $\edgesl{s}, \edgeli{s}$ from \Cref{notation_reduced} (iii). 
  By the tower property and \Cref{est_event_L}, uniformly in $\ranges$,
        \begin{align*}
            \mathbb{E}\brk{\bverslk{k}{s}}=& \Erw\brk{\sum_{i\in \verst{s}}\Erw\brk{\mathds{1}\mathfrak{L}_{i}  \cdot \teo{\bar{\vd}_{i,s}=k}\big\vert \verst{s},~\edgesl{s},~\edgeli{s}}} +\bar{o}_n(n).
        \end{align*}
Conditionally on $\mathfrak{L}_{i}$ and the values of $\verst{s}, \edgesl{s}$ and $\edgeli{s}$, for $i \in \verst{s}$, the neighbors of its half-edges are chosen uniformly from $\edgesl{s}-d_i$ sleeping half-edges and $\edgeli{s}-\edgesl{s}$ active half-edges without replacement. As a consequence, 
\begin{align*}
\Erw\brk{\teo{\bar{\vd}_{i,s}=k}|\mathfrak{L}_{i},\verst{s},\edgesl{s},\edgeli{s}}=\Erw\brk{\frac{\binom{\edgesl{s}-d_i}{k}\binom{\edgeli{s}-\edgesl{s}}{d_i-k}}{\binom{\edgeli{s}-d_i}{d_i}}\mid \mathfrak{L}_{i},\verst{s},\edgesl{s},\edgeli{s}}.
\end{align*}
Therefore, another application of the tower property yields that
\begin{align*}
    \mathbb{E}\brk{\bverslk{k}{s}}
            =&\sum_{i\in [n]}\Erw\brk{\teo{i\in \verst{s}}\frac{\binom{\edgesl{s}-d_i}{k}\binom{\edgeli{s}-\edgesl{s}}{d_i-k}}{\binom{\edgeli{s}-d_i}{d_i}}}+\bar{o}_n(n).
\end{align*}
By \Cref{l_concentration} and \Cref{lem_app_divide} in Appendix \ref{app_needproof} on the comparison between sampling with and without replacement, 
\begin{align*}
    \mathbb{E}\abs{\teo{i\in \verst{s}}\frac{\binom{\edgesl{s}-d_i}{k}\binom{\edgeli{s}-\edgesl{s}}{d_i-k}}{\binom{\edgeli{s}-d_i}{d_i}}-\teo{i\in \verst{s}}\binom{d_i}{k}\frac{\lambda(t_s)^k(\lambda(0)\eul^{-2 t_s}-\lambda(t_s))^{d_i-k} }{(\lambda(0)\eul^{-2 t_s})^{d_i}}}=\bar{o}_n(1).
\end{align*}
On the other hand, the number of vertices in $\verst{s}$ with original degree $m$ is $\verslk{m}{s}$, whose expectation has been estimated in \Cref{l_concentration}. Hence, we conclude that 
\begin{align*}
    \mathbb{E}\brk{\bverslk{k}{s}}
            &=\sum_{i\in [n]}\Erw\brk{\teo{i\in \verst{s}}\binom{d_i}{k}\frac{\lambda(t_s)^k(\lambda(0)\eul^{-2 t_s}-\lambda(t_s))^{d_i-k} }{(\lambda(0)\eul^{-2 t_s})^{d_i}}}+\bar{o}_n(n)\\
            =&\sum_{m=k}^K \binom{m}{k}\bc{\frac{\lambda(t_s)\eul^{2t_s}}{\lambda(0)}}^k\bc{1-\frac{\lambda(t_s)\eul^{2t_s}}{\lambda(0)}}^{m-k} \mathbb{E}\brk{\verslk{m}{s}}+\bar{o}_n(n)\\
            =&n\sum_{m=k}^K \binom{m}{k}\bc{\frac{\lambda(t_s)\eul^{2t_s}}{\lambda(0)}}^k\bc{1-\frac{\lambda(t_s)\eul^{2t_s}}{\lambda(0)}}^{m-k} \eul^{-m t_s} p_m+\bar{o}_n(n).
\end{align*}

\textit{Analysis of $\mathbb{E}\brk{\bverslk{k}{s}^2}$:} By \cref{est_event_L}, we can compute the second moment of $\bverslk{k}{s}$ as
\begin{align*}
            \mathbb{E}\brk{\bverslk{k}{s}^2}=&\sum_{i,j\in [n]}\PP\bc{i,j\in \verst{s},\bar{\vd}_{i,s}=\bar{\vd}_{j,s}=k}\\
            =&\sum_{i,j\in [n],i\neq j}\Erw\brk{\Erw\brk{\teo{i,j\in \verst{s},\mathfrak{L}_{i},\mathfrak{L}_{j},\mathfrak{L}_{i,j}}\teo{\bar{\vd}_{i,s}=\bar{\vd}_{j,s}=k}|\edgesl{s},\edgeli{s}}}+\bar{o}_n(n^2).
\end{align*}
As in the first moment computation, using \Cref{l_concentration} and \Cref{lem_app_divide}, we deduce that  
    \begin{align*}
    & \mathbb{E}\brk{\bverslk{k}{s}^2} \\
            =&\sum_{i,j\in [n],i\neq j}\Erw\brk{\Erw\Big[
            \frac{\teo{i,j\in \verst{s}}\binom{\edgesl{s}-d_i-d_j}{k}\binom{\edgesl{s}-d_i-d_j-k}{k}\binom{\edgeli{s}-\edgesl{s}}{d_i-k}\binom{\edgeli{s}-\edgesl{s}-d_i+k}{d_j-k}}{\binom{\edgeli{s}-d_i-d_j}{d_i}\binom{\edgeli{s}-2d_i-d_j}{d_j}} \mid \edgesl{s},\edgeli{s}\Big]} +\bar{o}_n(n^2)\\
            =&n^2\bc{\sum_{m=k}^K \binom{m}{k}\bc{\frac{\lambda(t_s)\eul^{2t_s}}{\lambda(0)}}^k\bc{1-\frac{\lambda(t_s)\eul^{2t_s}}{\lambda(0)}}^{m-k} \eul^{-m t_s} p_m}^2+\bar{o}_n(n^2).
\end{align*} 
As
\begin{align*}
    &\mathbb{E}\abs{\bverslk{k}{s}-n\sum_{m=k}^K \binom{m}{k}\bc{\frac{\lambda(t_s)\eul^{2t_s}}{\lambda(0)}}^k\bc{1-\frac{\lambda(t_s)\eul^{2t_s}}{\lambda(0)}}^{m-k} \eul^{-m t_s} p_m}\\
    \leq& \sqrt{\Erw\brk{\abs{\bverslk{k}{s}-n\sum_{m=k}^K \binom{m}{k}\bc{\frac{\lambda(t_s)\eul^{2t_s}}{\lambda(0)}}^k\bc{1-\frac{\lambda(t_s)\eul^{2t_s}}{\lambda(0)}}^{m-k} \eul^{-m t_s} p_m}^2}}=\bar{o}_n(n),
\end{align*}
the desired result follows from the second-moment method.
\end{proof}
In light of \Cref{lem_bvks}, we define 
\begin{align*}
\psi_t(\alpha):&=\sigma(t)^{-1}\sum_{k=0}^K\alpha^k\sum_{m=k}^K \binom{m}{k}\bc{\frac{\lambda(t)\eul^{2t}}{\lambda(0)}}^k\bc{1-\frac{\lambda(t)\eul^{2t}}{\lambda(0)}}^{m-k} \eul^{-m t} p_m\\
&=\sigma(t)^{-1}\sum_{k=0}^K p_k \eul^{-kt}\bc{1+\frac{\lambda(t)\eul^{2t}}{\lambda(0)}(\alpha-1)}^k.
\end{align*}
By \Cref{lem_bvks}, we can regard $\psi_{t_s}$ as the limiting generating function of the current degree distribution of  \scm{n}{s}. Finally, we set
\begin{align}\label{def_psi_hat}
\hat{\psi}_t(\alpha):=\lambda(t)^{-1}\sum_{k=0}^K k p_k \eul^{-kt}\bc{1+\frac{\lambda(t)\eul^{2t}}{\lambda(0)}(\alpha-1)}^{k-1} = \psi_t'(\alpha)/\psi_t'(1).
\end{align}

\subsection{The choice of the next awakened vertex}\label{sec_type_new}
In this section, we take a closer look at vertex $\vnu_{s+1/n}$ 
and its probabilistic properties. Recall that the current degree of vertex $j\in \verst{s}$ in \scm{n}{s} is denoted by $\bar{\vd}_{j,s}$. The current degree $\bar{\vd}_{j,s}$ of vertex $j\in \verst{s}$ deviates from its original degree $d_j$ if and only if $j$ is adjacent to one of the previously awakened vertices $\{\vnu_{1/n}, \ldots, \vnu_{s}\}$.  According to the graph exploration, the first half-edge of a vertex in $\verst{s}$ to be killed is either a half-edge connected to a vertex outside $\verst{s}$ when $\sum_{i \in \verst{s}} (d_i - \bar{\vd}_{i,s}) > 0$, or an arbitrary half-edge in \scm{n}{s} when $\sum_{i \in \verst{s}} (d_i - \bar{\vd}_{i,s}) = 0$. This leads to the following definition, which characterizes all possible choices of $\vnu_{s + 1/n}$ given \scm{n}{s}:

\begin{definition}[Hypnopompic half-edge] 
    Given \scm{n}{s},  
a half-edge $h$ incident to vertex $j \in \verst{s}$ is called hypnopompic if one of the following holds:
\begin{itemize}
    \item $\sum_{i\in \verst{s}}(d_i-\bar{\vd}_{i,s})=0$;
    \item $d_j-\bar{\vd}_{j,s}>0$ and $h$ is not one of the $\bar{\vd}_{j,s}$ half-edges that connect $j$ to a vertex in $\verst{s}$.
\end{itemize}
\end{definition}
The following lemma shows that given \scm{n}{s} and thus a set of hypnopompic half-edges incident to vertices in $\verst{s}$, each of the latter is equally likely to be the half-edge that is killed in the step of the graph exploration that awakens $\vnu_{s+1/n}$:

\begin{lemma}\label{lem-half-edge-exchangeable}
    Fix $\ranges$. Given $\fscm{n}{s}=G$, 
    let $h_1$ and $h_2$ be two hypnopompic half-edges, and $\vh$ be the first half-edge incident to a vertex in $\verst{s}$ that is declared dead. Then
    \begin{align*}
        \PP\bc{\vh=h_1\mid \fscm{n}{s}=G}=\PP\bc{\vh=h_2\mid \fscm{n}{s}=G}.
    \end{align*}
    Moreover, for any half-edge $h$ incident to a vertex $j \in \verst{s}$ that is not hypnopompic, $\PP\bc{\vh=h\mid \fscm{n}{s}=G} = 0$.
\end{lemma}

\begin{proof}
Observe that the set of non-hypnopompic half-edges $h$ is only non-empty if $\fscm{n}{s}$ is such that $\sum_{j\in \verst{s}}(d_j-\bar{\vd}_{j,s})>0$. In this case, conditioning on $\fscm{n}{s}=G$ entails that in \cm{n}, by definition, the subgraph \scm{n}{s} is connected to the vertices in $\{\vnu_{1/n}, \ldots, \vnu_s\}$ via a hypnopompic half-edge. Therefore, in the graph exploration, such an edge will be the first to be explored and thus killed. This implies the second claim.

For the remainder of the proof, fix $G$ and two hypnopompic half-edges $h_1$ and $h_2$. Let $\textbf{Hist}_s$ denote the (random) history of the vertex-removal graph exploration up to, but excluding the step when $\vnu_{s+1/n}$ is awakened. Thus, each possible realization of $\textbf{Hist}_s$ records the sequence of steps that have been carried out, together with the half-edges that are activated or killed and the vertices that are awakened in each step. Its last component either consists of the step in which $\vnu_s$ is awakened and its adjacent half-edges are activated or a step in which two half-edges between vertices in $\{\vnu_{1/n},\ldots,\vnu_s\}$ are paired and thus killed, but where no new vertex is awakened.

By Bayes' Theorem, for $i \in \{1,2\}$,
\begin{align*}
   \PP\bc{\vh=h_i\mid \fscm{n}{s}=G}&= \sum_{H}\frac{\PP\bc{\fscm{n}{s}=G\mid \textbf{Hist}_s=H,~\vh=h_i}\PP\bc{\vh=h_i\mid \textbf{Hist}_s=H}\PP\bc{\textbf{Hist}_s=H}}{\PP\bc{ \fscm{n}{s}=G}}.
\end{align*}
Thus, to show the claim, it is sufficient to show that for all histories $H$ that are compatible with $\fscm{n}{s}=G$, 
\begin{align}
    \PP\bc{\vh=h_1\mid \textbf{Hist}_s=H} &= \PP\bc{\vh=h_2\mid \textbf{Hist}_s=H} \quad \text{and} \label{eq_bayes_h1}\\
     \PP\bc{\fscm{n}{s}=G\mid \textbf{Hist}_s=H,~\vh=h_1} &= \PP\bc{\fscm{n}{s}=G\mid \textbf{Hist}_s=H,~\vh=h_2}. \label{eq_bayes_h2}
\end{align}
Given a history $H$, it is already determined whether $\vh$, the next half-edge to be killed, will be killed by \ref{exploration_step2} (if all half-edges adjacent to $\{\vnu_{1/n}, \ldots, \vnu_s\}$ have been explored) or \ref{exploration_step3} (otherwise). If $H$ is such that $\bm h$ is killed by \ref{exploration_step2}, then conditionally on $H$, the next step of the graph exploration is to awaken $\vnu_{s+1/n}$ through \ref{exploration_step1} by choosing a sleeping half-edge \textit{uniformly} at random. Thus, as both $h_1$ and $h_2$ are sleeping, they are chosen with the same conditional probability. If $H$ is such that $\bm h$ is killed by \ref{exploration_step3}, then conditionally on $H$, the next step of the graph exploration is to kill a half-edge in \ref{exploration_step3} based on its life-time. Conditionally on $H$, we know that the half-edge with the lowest life-time is adjacent to a vertex in $\fscm{n}{s}$, but not more than that. Again, by exchangeability of the lifetimes, $h_1$ and $h_2$ are chosen with the same conditional probability. This establishes \eqref{eq_bayes_h1}.

Next, also in the consideration of \eqref{eq_bayes_h2}, we distinguish two cases. First, assume that $H$ is such that $\bm h$ is chosen according to \ref{exploration_step2}, which means that in \cm{n}, the subgraph \scm{n}{s} is disconnected from the vertices in $\{\vnu_{1/n}, \ldots, \vnu_s\}$, and that $\bm h$ is part of \scm{n}{s}. Conditioning on $\textbf{Hist}_s=H,~\vh=h_i$ therefore potentially starts the exploration of the components of \scm{n}{s} from different half-edges, but at each future step, the probability to choose the `right' half-edge to pair with the last half-edge killed by \ref{exploration_step2} to complete $\fscm{n}{s}=G$ is uniformly distributed over the remaining half-edges.

 Second, assume that $H$ is such that $\bm h$ is chosen according to \ref{exploration_step3}, which means that in \cm{n}, the subgraph \scm{n}{s} is connected to the vertices in $\{\vnu_{1/n}, \ldots, \vnu_s\}$, and that $\bm h$ is not part of \scm{n}{s}. In this case, after conditioning on $\textbf{Hist}_s=H,~\vh=h_i$, not only unpaired half-edges in $\fscm{n}{s}$ are remaining, but also unpaired half-edges adjacent to vertices in $\{\vnu_{1/n}, \ldots, \vnu_s\}$. However, again, at each future step, the probability to choose a particular second half-edge to pair with the last half-edge killed by \ref{exploration_step2} is uniformly distributed over the remaining half-edges, and the number of possible pairings between hypnopomic half-edges and unpaired half-edges in $\{\vnu_{1/n}, \ldots, \vnu_s\}$ is the same given $\textbf{Hist}_s=H,~\vh=h_1$ and given $\textbf{Hist}_s=H,~\vh=h_2$. This establishes \eqref{eq_bayes_h2}.
 \end{proof}

In particular, \Cref{lem-half-edge-exchangeable} implies that given \scm{n}{s}, the probability that any sleeping vertex in $\verst{s}$ is chosen as the next vertex to be awakened is either proportional to its original degree or to the number of half-edges that connect it to the set of awake vertices: 
\begin{corollary}\label{cor-half-edge-exchangeable}
   Given \scm{n}{s},  the following two cases can occur:
\begin{enumerate}[label=(\roman*)]
    \item \label{it-cor-half-edge-exchangeable-1} $\sum_{i\in \verst{s}}(d_i-\bar{\vd}_{i,s})=0$: In this case, $\vnu_{s+1/n}$ is chosen among the vertices in $\verst{s}$
  with probability proportional to the original degree sequence $(d_i)_{i \in \verst{s}}$.

    \item \label{it-cor-half-edge-exchangeable-2} $\sum_{i\in \verst{s}}(d_i-\bar{\vd}_{i,s})\neq 0$: In this case, $\vnu_{s+1/n}$ is chosen among the vertices in $\verst{s}$ with probability proportional to $(d_i-\bar{\vd}_{i,s})_{i\in\verst{s}}$, i.e., the number of half-edges that connect it to the vertices $ \{\vnu_1, \ldots, \vnu_s\}$.
\end{enumerate}
\end{corollary}

\section{Graph decomposition and rank}

\subsection{The type of the next awakened vertex}
Given the previous results about the choice of $\vnu_{s+1/n}$, we next return to the associated adjacency matrices. Denote by $\bm{A}_{n,s}$ the matrix in which the rows and columns corresponding to $\vnu_{1/n},\vnu_{2/n},\ldots,\vnu_s$ in $\bm{A}_n$ are replaced by zero rows and columns, respectively. Then excluding those zero rows and columns, $\bm{A}_{n,s}$ can be viewed as the adjacency matrix of \scm{n}{s}. By a slight abuse of terminology, from now on, we refer to $\bm{A}_{n,s}$ as the  adjacency matrix of \scm{n}{s}.

According to the proof strategy in \Cref{sec_cm_exp}, our primary goal is to compute the probability that the newly awakened vertex $\vnu_{s+1/n}$ is of a specific type, conditional on
$\bm{A}_{n,s}[\vth]$, the perturbed adjacency matrix of the induced subgraph of sleeping vertices before its awakening. Recall that the matrix $\bm{A}_{n,s}[\bm{\theta}]$ is fully determined by \scm{n}{s} and $\THETA$. 
In this context, \Cref{cor-half-edge-exchangeable} inspires the following definitions, which correspond to the conditional probability that $\vnu_{s+1/n}$ 
has a certain type in $\bm{A}_{n,s}[\bm{\theta}]$:
\begin{definition}[Size-biased type proportions]\label{def_aw}
Given \scm{n}{s} and $\THETA=(\THETA_r[\theta_r, n],\THETA_c[n, \theta_c])$, define
 \begin{equation*}
\begin{aligned}
    \val_s:=&\begin{cases}
        \frac{\sum_{i\in \verst{s}}\teo{ i\in\mathcal{F}\bc{\bm{A}_{n,s}[\bm{\theta}]}}(d_i-\bar{\vd}_{i,s})}{\sum_{i\in \verst{s}}(d_i-\bar{\vd}_{i,s})}, &\sum_{i\in \verst{s}}(d_i-\bar{\vd}_{i,s})\neq 0;\\
         \frac{\sum_{i\in \verst{s}}\teo{ i\in\mathcal{F}\bc{\bm{A}_{n,s}[\bm{\theta}]}}d_i}{\sum_{i\in \verst{s}}d_i},&\sum_{i\in \verst{s}}(d_i-\bar{\vd}_{i,s})=0;
    \end{cases}\\
    \vw_s:=&\begin{cases}
        \frac{\sum_{i\in \verst{s}}\teo{ i\in\mathcal{W}\bc{\bm{A}_{n,s}[\bm{\theta}]}}(d_i-\bar{\vd}_{i,s})}{\sum_{i\in \verst{s}}(d_i-\bar{\vd}_{i,s})}, &\sum_{i\in \verst{s}}(d_i-\bar{\vd}_{i,s})\neq 0;\\
         \frac{\sum_{i\in \verst{s}}\teo{ i\in\mathcal{W}\bc{\bm{A}_{n,s}[\bm{\theta}]}}d_i}{\sum_{i\in \verst{s}}d_i},&\sum_{i\in \verst{s}}(d_i-\bar{\vd}_{i,s})=0,
    \end{cases}
\end{aligned}
 \end{equation*}
for $\vw\in\cbc{\vx,\vy,\vz,\vu,\vv}$. We also abbreviate $\vze_s=(\vx_s,\vy_s,\vz_s,\vu_s,\vv_s)$ and call the components of this vector the size-biased type proportions.     
\end{definition}
 
 In terms of the size-biased type proportions, \Cref{cor-half-edge-exchangeable} reads as follows:
\begin{lemma}\label{lem_con_typp}
For any $\vw\in\cbc{\vx,\vy,\vz,\vu,\vv}$ and $\ranges$,
    \begin{align*}
\PP\bc{\vnu_{s+1/n}\in\mathcal{W}\bc{\bm{A}_{n,s}[\bm{\theta}]} \mid \vze_s}=\PP\bc{\vnu_{s+1/n}\in \mathcal{W}\bc{\bm{A}_{n,s}[\bm{\theta}]} \mid \fscm{n}{s},\THETA}=\vw_s.
\end{align*}
\end{lemma}

\begin{proof}
Given $\fscm{n}{s}$, its vertex set $\verst{s}$ and induced degrees $(\bar{\vd}_{i,s})_{i\in \verst{s}}$ are fixed.
\begin{enumerate}[label=\textbf{Case \arabic*:}]
    \item\label{lem_con_typp_case1} If $\sum_{i\in \verst{s}}(d_i-\bar{\vd}_{i,s})=0$, by \Cref{cor-half-edge-exchangeable} \ref{it-cor-half-edge-exchangeable-1}, for each sleeping vertex $j \in \verst{s}$,
\begin{align}\label{eq_prob_vbej}
\PP\bc{\vnu_{s+1/n}=j \mid \fscm{n}{s}}=\frac{d_j}{\sum_{i\in\verst{s}}d_i}.
\end{align}
Since the perturbation matrices $\THETA$ are independent of the exploration and the configuration model \cm{n}, \cref{eq_prob_vbej} gives that
\begin{align*}
&\PP\bc{\vnu_{s+1/n}\in \mathcal{W}\bc{\bm{A}_{n,s}[\bm{\theta}]} \mid \fscm{n}{s},\THETA}
=\sum_{j\in \mathcal{W}\bc{\bm{A}_{n,s}[\bm{\theta}]} \cap \verst{s}}\PP\bc{\vnu_{s+1/n}=j \mid \fscm{n}{s},\THETA}=\vw_s.\nonumber
\end{align*}

\item\label{lem_con_typp_case2} If $\sum_{i\in \verst{s}}(d_i-\bar{\vd}_{i,s})\neq 0$, by \Cref{cor-half-edge-exchangeable} \ref{it-cor-half-edge-exchangeable-2}, for each sleeping vertex $j \in \verst{s}$,
\begin{align}\label{eq_prob_vbej2}
\PP\bc{\vnu_{s+1/n}=j \mid \fscm{n}{s}}=\frac{d_j-\bar{\vd}_{j,s}}{\sum_{i\in \verst{s}}(d_i-\bar{\vd}_{i,s})}.
\end{align}
As in \textbf{Case $\bm 1$}, it follows that $\PP\bc{\vnu_{s+1/n}\in \mathcal{W}\bc{\bm{A}_{n,s}[\bm{\theta}]} \mid \fscm{n}{s},\THETA}=\vw_s$.
\end{enumerate}
Finally, the missing identity follows from the tower property:
\begin{align*}
\PP\bc{\vnu_{s+1/n}\in\mathcal{W}\bc{\bm{A}_{n,s}[\bm{\theta}]} \mid \vze_s}=\mathbb{E}\brk{\PP\bc{\vnu_{s+1/n}\in \mathcal{W}\bc{\bm{A}_{n,s}[\bm{\theta}]} \mid \fscm{n}{s},\THETA} \mid \vze_s}
= \vw_s.
\end{align*}
\end{proof}

\subsection{Rank and types} \label{sec_rank_types}
We now relate the graph decomposition procedure and the size-biased proportions of types to our main objective, the derivation of a lower bound on the expected rank for supercritical configuration models. Abbreviating $\iota:=1-\sigma(-\ln \xi)-\varepsilon$, we obtain the lower bound 
\begin{align}\label{l_decom}
\Erw\brk{\rk{\bm{A}_n\pth}}\geq&\sum_{j=\lfloor \varepsilon n \rfloor}^{\lceil n\iota\rceil -1 }\Erw\brk{\rk{\bm{A}_{n,j/n}\pth}-\rk{\bm{A}_{n,(j+1)/n}\pth}}+\Erw\brk{\rk{\bm{A}_{n,\lceil n\iota\rceil/n}\pth}}.
\end{align}
Further, by the choice of $\iota$, $\Erw\brk{\rk{\bm{A}_{n,\lceil n\iota\rceil/n }\pth}}$ effectively corresponds to a rank computation for a subcritical configuration model.
Thus, we will be mainly concerned with estimating the rank differences in \Cref{l_decom}.

For this, we define the event
\begin{align}\label{event_pns}
\mathfrak{P}_{n,s}=\cbc{\THETA_r[\bm{\theta}_r,n](,\vnu_{s+1/n})=0_{\bm{\theta}_r\times 1}, \THETA_c[n,\bm{\theta}_c](\vnu_{s+1/n},)=0_{1\times\bm{\theta}_c}}.
\end{align}
For $\THETA$ in $\mathfrak{P}_{n,s}$, $\bm{A}_{n,s}[\vth]\abc{\vnu_{s+1/n};\vnu_{s+1/n}}=\bm{A}_{n,s+1/n}[\vth]$.
Moreover, \Cref{d34}  ensures  that 
\begin{align}\label{p_event_P}
\PP\bc{\mathfrak{P}_{n,s}}=1+\bar{o}_{n,P}(1).
\end{align}
Then by the deterministic rank relation \Cref{lr1} and \Cref{lem_con_typp} on the type of the next awakened vertex, for $\THETA \in \mathfrak{P}_{n,s}$,
\begin{align*}
&\Erw\brk{\rk{\bm{A}_{n,s}[\vth]} -\rk{\bm{A}_{n,s+1/n}[\vth]} \mid \THETA}  \\
&\qquad=\PP\bc{\vnu_{s+1/n} \in \cX(\bm{A}_{n,s+1/n}[\vth]) \mid \THETA}+2\PP\bc{\vnu_{s+1/n} \in \cY(\bm{A}_{n,s+1/n}[\vth]) \mid \THETA}\\
&\qquad\quad+\PP\bc{\vnu_{s+1/n} \in \cU(\bm{A}_{n,s+1/n}[\vth]) \mid \THETA}+\PP\bc{\vnu_{s+1/n} \in \cV(\bm{A}_{n,s+1/n}[\vth]) \mid \THETA}\\
&\qquad=\Erw\brk{\vx_s+2\vy_s+\vu_s+\vv_s \mid \THETA}.
\end{align*}
Combining the last identity with \Cref{p_event_P} gives that
\begin{align}\label{eq_re_rd}
    \Erw\brk{\rk{\bm{A}_{n,s}[\vth]}-\rk{\bm{A}_{n,s+1/n}[\vth]}}=\Erw\brk{\vx_s+2\vy_s+\vu_s+\vv_s}+\bar{o}_{n,P}(1).
\end{align}
\eqref{eq_re_rd} illustrates why it is crucial for our approach to estimate the values of the size-biased proportions of frozen types. 
In the spirit of \cite{HofMul25}, we aim 
to derive fixed-point equations for the size-biased type proportions. These equations will eventually help us to lower bound the asymptotic rank through \cref{eq_re_rd}. 

Our route towards a derivation of the fixed-point equations is the following proposition: 
\begin{proposition}\label{prop_allcomestogether}
 Fix $\rangee$.   Assume that 
    \begin{align}\label{eq_allcomestogether_a1}
        \dabs{\vze_{s}-\vze_{s+1/n}}_1 = \oone 
    \end{align}
    and that there exist functions $Y,Z,U,V$ such that for each $W \in \{Y,Z,U,V\}$, $W=W(\zeta,\hat{\psi}_{t_s})$ is differentiable with respect to $\zeta=(x,y,z,u,v)\in [0,1]^5$ with uniformly bounded gradient, and that
\begin{align}\label{eq_allcomestogether_a2}
   \mathbb{P}\bc{\vnu_{s+1/n}\in\mathcal{W}\bc{\bm{A}_{n,s}[\bm{\theta}]} \mid \vze_{s+1/n}}-W\bc{\vze_{s+1/n},\hat{\psi}_{t_s}}=\oone. 
    \end{align}
    Then, for any $\vw \in \{ \vy, \vz, \vu, \vv\}$,
    \begin{align}\label{eq_allcomestogether_re}
      \vw_s - W(\vze_s,\hat{\psi}_{t_s}) = \oone.  
    \end{align}
    If  only 
    \begin{align}\label{eq_allcomestogether_a2_2}
       \mathbb{P}\bc{\vnu_{s+1/n}\in\mathcal{W}\bc{\bm{A}_{n,s}[\bm{\theta}]} \mid \vze_{s+1/n}}-W\bc{\vze_{s+1/n},\hat{\psi}_{t_s}}\geq \oone, 
    \end{align}
    then
\begin{align}\label{eq_allcomestogether_re2}
      \vw_s - W(\vze_s,\hat{\psi}_{t_s})\geq \oone.  
    \end{align}
\end{proposition}
The proof of \Cref{prop_allcomestogether} is a direct consequence of \cite[Proposition C.1]{HofMul25}. Its derivation from that proposition is given in Appendix \ref{sec_proof_allct}.

\section{Stability of types under graph decomposition}\label{Sec_stability_type}

To derive the fixed-point equations for the type proportions, we first show that in the beginning of the graph decomposition, the latter remain relatively stable under the removal of a single vertex and its incident edges. Fix $\rangee$, and recall that we aim to derive uniform bounds for $s \in [\varepsilon, 1-\sigma(-\ln \xi)-\varepsilon]$, which is captured in our $o_n$, $\bar{o}_n$ and $\bar{o}_{n,P}$ notation. The main result of the current section demonstrates that the proportions of types remain relatively stable under a single step of the graph exploration:
\begin{proposition}[Stability of types]\label{pro_sta_vze}
Fix $\rangee$. For any $\vw \in \{ \vx, \vy, \vz, \vu, \vv, \vze\}$,
\[\mathbb{E}\dabs{\vw_{s}-\vw_{s+1/n}}_1 = \bar{o}_{n,P}(1). \]
\end{proposition}
We prove \Cref{pro_sta_vze} in the remainder of this section.
\subsection{Reduction to average type changes}\label{sec_sop}
In this section, we reduce the proof of \Cref{pro_sta_vze} to the task of bounding the probability that a uniformly chosen vertex changes its type in $\bm{A}_{n,s}[\bm{\theta}]$ as compared to $\bm{A}_{n,s+1/n}[\bm{\theta}]$. 

For $\rangee$, define 
\begin{align}\label{def_c}
   c = c(\varepsilon):=\min_{r\in [\varepsilon,1-\sigma(-\ln\xi)-\varepsilon]}\frac{\lambda(t_r)(\lambda(0)\eul^{-2t_r}-\lambda(t_r))}{2K\lambda(0)\eul^{-2t_r}}, 
\end{align}
where $t_s$ is the unique zero of $H_s(t)=\sigma(t)-1+s$ in $[0,1-p_0)$ from \Cref{def_t_s}. \Cref{l_dcs} in the Appendix shows that $n\lambda(t_r)(\lambda(0)\eul^{-2t_r}-\lambda(t_r))\lambda(0)^{-1}\eul^{2t_r}$ is a good approximation of the number of half-edges joining vertices in $\verst{r}$ to vertices outside of $\fscm{n}{r}$. Therefore, $cn$ functions as a uniform lower bound on the reservoir of half-edges joining the set of sleeping vertices to the set of awake vertices during the current exploration epoch.
The combination of \cref{eq_psi_lamda} and \Cref{lem_psi0} yields that $c>0$. 

In most proofs of this section, we will work on the event
\begin{align}\label{event_gs}
    \mathfrak{G}_{\varepsilon,s} := \cbc{\text{the number of $j\in \verst{s}$ such that $\bar{\vd}_{j,s}<d_j$ is greater than $cn$}}.
\end{align}
This ensures that the number of half-edges of vertices in $\verst{s}$ connected to vertices in $[n] \backslash \verst{s}$ remains of order $n$, thereby guaranteeing that the next vertex to be awakened will be chosen according to \ref{exploration_step3}. 
\Cref{l_dcs} in  Appendix \ref{ap_eventg} shows that $\mathfrak{G}_{\varepsilon,s}$ occurs with high probability: 
\begin{align}\label{ineq_dcs}\PP\bc{\mathfrak{G}^c_{\varepsilon,s}}=\bar{o}_n(1).
\end{align}

\subsubsection{Proof of Proposition~\ref{pro_sta_vze}, Step 1.}
Recall \Cref{def_aw} for the type proportions $\vw \in \{ \vx, \vy, \vz, \vu, \vv\}$. We first relate the difference of the proportions $\vw_{s}$ and $\vw_{s+1/n}$ to the type change of a uniformly chosen vertex.
Set $\bar{\mathfrak G}_{\eps,s}:= \mathfrak{G}_{\eps,s}\cap \mathfrak{G}_{\eps,s+1/n}$. As $\abs{\vw_{s}-\vw_{s+1/n}}$ is bounded above by $1$, by \Cref{ineq_dcs},
\begin{align*} 
&\mathbb{E}\abs{\vw_{s}-\vw_{s+1/n}} = \bar{o}_n(1) \\
+ &
\mathbb{E}\Bigg[\mathds{1} \bar{\mathfrak G}_{\eps,s} \abs{\frac{\sum_{i\in \verst{s}}\teo{ i\in\mathcal{W}\bc{\bm{A}_{n,s}[\bm{\theta}]}}(d_i-\bar{\vd}_{i,s})}{\sum_{i\in \verst{s}}(d_i-\bar{\vd}_{i,s})}-\frac{\sum_{i\in \verst{s+1/n}}\teo{ i\in\mathcal{W}\bc{\bm{A}_{n,s+1/n}[\bm{\theta}]}}(d_i-\bar{\vd}_{i,s+1/n})}{\sum_{i\in \verst{s+1/n}}(d_i-\bar{\vd}_{i,s+1/n})}}\Bigg].\nonumber 
\end{align*}
Moreover,  
\begin{flalign*}
    &\sum_{i\in \verst{s}}\teo{ i\in\mathcal{W}\bc{\bm{A}_{n,s}[\bm{\theta}]}}(d_i-\bar{\vd}_{i,s})-\sum_{i\in \verst{s+1/n}}\teo{ i\in\mathcal{W}\bc{\bm{A}_{n,s}[\bm{\theta}]}}(d_i-\bar{\vd}_{i,s})\\
    &\qquad=\teo{ \vnu_{s+1/n}\in\mathcal{W}\bc{\bm{A}_{n,s}[\bm{\theta}]}}(d_{\vnu_{s+1/n}}-\bar{\vd}_{\vnu_{s+1/n}} ,s)\leq K,
\end{flalign*}
so that
\begin{align*}
&\mathbb{E}\abs{\vw_{s}-\vw_{s+1/n}}\\
=&
\mathbb{E}\Bigg[\mathds{1} \bar{\mathfrak G}_{\eps,s}\Bigg|\sum_{i\in \verst{s+1/n}}\Bigg(\frac{\teo{ i\in\mathcal{W}\bc{\bm{A}_{n,s}[\bm{\theta}]}}(d_i-\bar{\vd}_{i,s})}{\sum_{j\in \verst{s}}(d_j-\bar{\vd}_{j,s})}-\frac{\teo{ i\in\mathcal{W}\bc{\bm{A}_{n,s+1/n}[\bm{\theta}]}}(d_i-\bar{\vd}_{i,s+1/n})}{\sum_{j\in \verst{s+1/n}}(d_j-\bar{\vd}_{j,s+1/n})}\Bigg)\Bigg|\Bigg]+\bar{o}_n(1).\nonumber
\end{align*}
We next use that for any $m \geq 1$ and sequences $(\hat a_i)_{i \in [m]},  (\hat b_i)_{i \in [m]}, (\hat c_i)_{i \in [m]}, (\hat d_i)_{i \in [m]}$ of numbers,
\begin{align}\label{eq_prod_deff}
    \Big|\sum_{i\in [m]}(\hat{a}_i \hat{b}_i-\hat{c}_i \hat{d}_i)\Big| &= \Big|\sum_{i\in [m]}(\hat{a}_i\hat{b}_i-\hat a_i\hat d_i+ \hat d_i \hat a_i - \hat d_i\hat{c}_i)\Big| \\
    &\leq \sup_{i\in [m]}\abs{\hat{a}_i}\sum_{i\in [m]}\abs{\hat{b}_i-\hat{d}_i}+\sup_{i\in [m]}\abs{\hat{d}_i}\sum_{i\in [m]}\abs{\hat{a}_i-\hat{c}_i}.\nonumber
\end{align}
Moreover, on $\bar{\mathfrak{G}}_{\eps,s}$, $\frac{d_i-\bar{\vd}_{i,s+1/n}}{\sum_{j\in \verst{s+1/n}}(d_j-\bar{\vd}_{j,s+1/n})} \leq K/(cn)$, so that by \cref{eq_prod_deff}, again on $ \bar{\mathfrak G}_{\eps,s}$,
\begin{align*}
    & \abs{\sum_{i\in \verst{s+1/n}}\Bigg( \frac{\teo{ i\in\mathcal{W}\bc{\bm{A}_{n,s}[\bm{\theta}]}}(d_i-\bar{\vd}_{i,s})}{\sum_{j\in \verst{s}}(d_j-\bar{\vd}_{j,s})}-\frac{\teo{ i\in\mathcal{W}\bc{\bm{A}_{n,s+1/n}[\bm{\theta}]}}(d_i-\bar{\vd}_{i,s+1/n})}{\sum_{j\in \verst{s+1/n}}(d_j-\bar{\vd}_{j,s+1/n})}\Bigg)}\\
   \leq &\sum_{i\in \verst{s+1/n}}\abs{\frac{d_i-\bar{\vd}_{i,s}}{\sum_{j\in \verst{s}}(d_j-\bar{\vd}_{j,s})}-\frac{d_i-\bar{\vd}_{i,s+1/n}}{\sum_{j\in \verst{s+1/n}}(d_j-\bar{\vd}_{j,s+1/n})}}\\
    &    +\frac{K }{cn}\sum_{i\in \verst{s+1/n}}\teo{ i\in\mathcal{W}\bc{\bm{A}_{n,s}[\bm{\theta}]}\Delta \mathcal{W}\bc{\bm{A}_{n,s+1/n}[\bm{\theta}]}}.
\end{align*}
Furthermore, $\sum_{i\in \verst{s+1/n}}\abs{\bar{\vd}_{i,s}-\bar{\vd}_{i,s+1/n}}\leq K$ since for each $i\in \verst{s+1/n}$, $\bar{\vd}_{i,s}-\bar{\vd}_{i,s+1/n}$ is equal to the number of half-edges $\vnu_{s+1/n}$ connected to vertex $i$. Then, on $\bar{\mathfrak G}_{\eps,s}$, 
\begin{align*}
    &\sum_{i\in \verst{s+1/n}}\abs{\frac{d_i-\bar{\vd}_{i,s+1/n}}{\sum_{j\in \verst{s+1/n}}(d_j-\bar{\vd}_{j,s+1/n})}-\frac{d_i-\bar{\vd}_{i,s}}{\sum_{j\in \verst{s}}(d_j-\bar{\vd}_{j,s})}}\\
    \leq &\sum_{i\in \verst{s+1/n}}(d_i-\bar{\vd}_{i,s+1/n})\abs{\frac{1}{\sum_{j\in \verst{s+1/n}}(d_j-\bar{\vd}_{j,s+1/n})}-\frac{1}{\sum_{j\in \verst{s}}(d_j-\bar{\vd}_{j,s})}}+\frac{\sum_{i\in \verst{s+1/n}}\abs{\bar{\vd}_{i,s}-\bar{\vd}_{i,s+1/n}}}{\sum_{i\in \verst{s}}(d_i-\bar{\vd}_{i,s})}\\
    = &\abs{\frac{\sum_{i\in \verst{s+1/n}}\bc{\bar{\vd}_{i,s+1/n}-\bar{\vd}_{i,s}}+d_{\vnu_{s+1/n}}-\bar{\vd}_{\vnu_{s+1/n},s}}{\sum_{i\in \verst{s}}(d_i-\bar{\vd}_{i,s})}}+\frac{\sum_{i\in \verst{s+1/n}}\abs{\bar{\vd}_{i,s}-\bar{\vd}_{i,s+1/n}}}{\sum_{i\in \verst{s}}(d_i-\bar{\vd}_{i,s})}\leq \frac{3K}{cn}=\bar{o}_n(1),
\end{align*}
where in the last inequality we use $\sum_{i\in \verst{s+1/n}}\abs{\bar{\vd}_{i,s}-\bar{\vd}_{i,s+1/n}}\leq K$. 
Thus we conclude that
\begin{align}\label{eq_l91_x1}
    \mathbb{E}\abs{\vw_{s}-\vw_{s+1/n}}\leq &\frac{K}{cn}\mathbb{E}\Big[\sum_{i\in \verst{s+1/n}}\teo{ i\in\mathcal{W}\bc{\bm{A}_{n,s}[\bm{\theta}]}\Delta \mathcal{W}\bc{\bm{A}_{n,s+1/n}[\bm{\theta}]}}\Big]+\bar{o}_n(1)\nonumber\\
    \leq&\frac{K}{cn}\sum_{i\in [n]}\PP\bc{ i\in\mathcal{W}\bc{\bm{A}_{n,s}[\bm{\theta}]}\Delta \mathcal{W}\bc{\bm{A}_{n,s+1/n}[\bm{\theta}]}}+\bar{o}_n(1).
\end{align}

\subsection{Average type changes}
The bound \cref{eq_l91_x1} illustrates that in order to prove the stability result \Cref{pro_sta_vze}, it is sufficient to obtain a bound on the probability of type change of an average vertex. In this sense, we show that advancing the graph exploration by a single step has a negligible effect on whether a vertex is frozen or not with respect to the associated perturbed matrices $\bm{A}_{n,s}[\bm{\theta}], \bm{A}_{n,s+1/n}[\bm{\theta}]$  \textit{on average}. To be able to deal with all types of proportions $\vw \in \{ \vx, \vy, \vz, \vu, \vv\}$, we also investigate freezing/unfreezing with respect to a modified matrix. 
The key steps in the derivation of the desired stability result are the following two lemmas:
\begin{lemma}[One-step matrix zeroing, original matrix]\label{l7}
Fix $\delta >0$ and $\rangee$. Then
\[\frac{1}{n}\sum_{i\in [n]}\mathbb{P}\left(i\in \mathcal{F}\bc{\bm{A}_{n,s}[\bm{\theta}]}\Delta\mathcal{F}\bc{\bm{A}_{n,s+1/n}[\bm{\theta}]}\right)\leq 4\delta (K+1)K^{2K+2}c^{-K-1}+\bar{o}_{n,P}(1). \]
\end{lemma}
\begin{lemma}[One-step matrix zeroing, row-zeroed matrix]\label{l71}
Fix $\delta >0$ and $\rangee$. Then
\[\frac{1}{n}\sum_{i\in [n]}\mathbb{P}\left(i\in\mathcal{F}\bc{\bm{A}_{n,s}[\bm{\theta}]\abc{i;}}\Delta \mathcal{F}\bc{\bm{A}_{n,s+1/n}[\bm{\theta}]\abc{i;}}\right) \leq 4\delta (K+1)K^{2K+2}c^{-K-1}+\bar{o}_{n,P}(1).\]
\end{lemma}

\begin{proof}[Proof of \Cref{pro_sta_vze}]
  Given the bound \cref{eq_l91_x1} along with \Cref{l7,l71}, the only missing step in the proof of \Cref{pro_sta_vze} is relating the events $\cbc{ i\in\mathcal{W}\bc{\bm{A}_{n,s}[\bm{\theta}]}\Delta \mathcal{W}\bc{\bm{A}_{n,s+1/n}[\bm{\theta}]}}$ to freezing or unfreezing events in various matrices as studied in  \Cref{l7,l71}. However, this analysis is only a minor adaption of the proof of \cite[Proposition 4.11]{HofMul25}. \slv
\end{proof}

The rough proof strategy for \Cref{l7,l71} is similar to the approach that has been used for the weighted Erd\H{o}s-Rényi random graph in  \cite[Lemmas 4.17 and 4.18]{HofMul25}. However, the details differ significantly, as edges in the configuration model are not independent of each other. Therefore, we have to resort to alternative approaches, such as the moment method, to establish concentration. 

\subsection{One-step matrix zeroing, original matrix: Proof of \Cref{l7}}
This proof is an adaptation of the proof of \cite[Lemma 4.17]{van2022random}. First, observe that 
it is sufficient to bound
\begin{align}\label{eq_delta_backslash}
\mathbb{P}\bc{i \in \mathcal{F}\bc{\bm{A}_{n,s+1/n}[\bm{\theta}]}\backslash \mathcal{F}\bc{\bm{A}_{n,s}[\bm{\theta}]}} \quad \text{and} \quad \mathbb{P}\bc{i \in  \mathcal{F}\bc{\bm{A}_{n,s}[\bm{\theta}]}\backslash \mathcal{F}\bc{\bm{A}_{n,s+1/n}[\bm{\theta}]}}.
\end{align}
Our proof strategy is to relate both freezing and unfreezing of $i$ to the existence of a proper relation in either $\bm{A}_{n,s}[\bm{\theta}]$ or $\bm{A}_{n,s+1/n}[\bm{\theta}]$ that involves $i$: 
 On the event $\mathfrak{P}_{n,s}$ defined in \eqref{event_pns},
\begin{align}\label{implication_unfreeze}
  i\in \mathcal{F}\bc{\bm{A}_{n,s+1/n}[\bm{\theta}]}\backslash \mathcal{F}\bc{\bm{A}_{n,s}[\bm{\theta}]} \quad \Longrightarrow \quad 
  \{i,\vnu_{s+1/n}\} \text{ is a proper relation of } \bm{A}_{n,s}[\bm{\theta}],\end{align}
and 
\begin{align}\label{implication_freeze}
    & i\in \mathcal{F}\bc{\bm{A}_{n,s}[\bm{\theta}]}\backslash \mathcal{F}\bc{\bm{A}_{n,s+1/n}[\bm{\theta}]},  \quad  i\notin \{\vnu_{s+1/n}\}\cup\supp{\bm{A}_{n,s}[\bm{\theta}](\vnu_{s+1/n},)}\nonumber \\
     \Longrightarrow  \quad &\cbc{i}\cup \supp{\bm{A}_{n,s}[\bm{\theta}](\vnu_{s+1/n},)} \text{ is a proper relation of } \bm{A}_{n,s+1/n}[\bm{\theta}].
 \end{align}
We defer the proofs of the implications \Cref{implication_unfreeze} and \Cref{implication_freeze} to Appendix \ref{app_pre_results}, as they are analogous to \cite[Lemma 4.17]{HofMul25}.  
Given these deterministic statements, we estimate the two probabilities in \cref{eq_delta_backslash} on the event
\begin{align}\label{event_R}
    \mathfrak{R}_{n,s}=\cbc{\text{{both} $\bm{A}_{n,s+1/n}[\bm{\theta}]$ and $\bm{A}_{n,s+1/n}[\bm{\theta}]^T$ are $(\delta,\ell)$-free for  $2\leq \ell \leq 2K$}}.
\end{align}
The benefit of $\mathfrak{R}_{n,s}$ is that we can \textit{almost} disregard the impact of short proper linear relations of $\bm{A}_{n,s+1/n}[\bm{\theta}]$ and its transpose. Moreover, by \Cref{p1},
\begin{equation}\label{oe}
\mathbb{P}\left({\mathfrak{R}_{n,s}^c}\right)=\bar{o}_{n,P}(1).
\end{equation} 
Given \eqref{implication_unfreeze} - \eqref{oe}, we now proceed with bounding the probabilties in \eqref{eq_delta_backslash} separately.

\begin{enumerate}[label=(\roman*),left=1pt] 
\item \textit{Freezing.} 
We first show that, on the high-probability event $\mathfrak R_{n,s-1/n}$, the probability that $\{i,\vnu_{s+1/n}\}$ forms a proper relation of $\bm{A}_{n,s}[\bm{\theta}]$ is small when averaged over $i\in [n]$. 
Given $\fscm{n}{s} = G \in \mathfrak{G}_{\eps,s}$, it holds true that $\sum_{j\in\verst{s}}(d_j-\bar{d}_{j,s})\geq cn$ and $\vnu_{s+1/n}$ is chosen according to \ref{exploration_step3}. Then, for each vertex $j \in \verst{s}$, 
\begin{align}\label{eq_conp_nus}
    \PP\bc{\vnu_{s+1/n}=j \mid \fscm{n}{s} = G}\leq \frac{d_j-\bar{d}_{j,s}}{\sum_{k\in\verst{s}}\bc{d_k-\bar{d}_{k,s}}}\leq \frac{K}{cn}.
\end{align}
Recall from \Cref{sec_permat} that $\PR(\bm{A}_{n,s}[\bm{\theta}])$ is the set of all proper relations of $\bm{A}_{n,s}[\bm{\theta}]$. 
By \eqref{eq_conp_nus},
\begin{align}\label{eq_ex_nu}
   & \sum_{i\in [n]}\PP\bc{\cbc{i, \vnu_{s+1/n}}\in \PR(\bm{A}_{n,s}[\bm{\theta}]) \mid \fscm{n}{s}=G,\THETA}\nonumber\\
&\qquad=\sum_{j \in \verst{s}} \sum_{j\in\verst{s}} \teo{\cbc{i, j}\in \PR(\bm{A}_{n,s}[\bm{\theta}])}\PP\bc{\vnu_{s+1/n}=j \mid \fscm{n}{s}=G} \nonumber\\
    &\qquad\leq  \frac{K}{cn}\sum_{i,j\in [n]}\teo{\cbc{i, j}\in \PR(\bm{A}_{n,s}[\bm{\theta}])}.
\end{align}
Recall \Cref{{d2}} and \cref{oe}.  Then \Cref{ineq_dcs} and \cref{eq_ex_nu} give that
\begin{align}\label{eq_ivnave}
\frac{1}{n}\sum_{i\in[n]}\mathbb{P}&\bc{\cbc{i, \vnu_{s+1/n}}\in \PR(\bm{A}_{n,s}[\bm{\theta}])}\\
    \leq &\frac{1}{n}\sum_{i\in[n]}\mathbb{E}\brk{\ind{\mathfrak{G}_{\eps,s}}\PP\bc{\cbc{i,\vnu_{s+1/n}}\in \PR(\bm{A}_{n,s}[\bm{\theta}]) \mid \fscm{n}{s},\THETA}}+\bar{o}_n(1)\nonumber\\
    \leq& \frac{K}{c n^2}\sum_{i,j\in [n]}\PP\bc{\cbc{i,j}\in \PR(\bm{A}_{n,s}[\bm{\theta}]),\mathfrak{R}_{n,s-1/n}}+\bar{o}_{n,P}(1).\nonumber
\end{align}

On the other hand,  on $\mathfrak{R}_{n,s-1/n}$, there are at most $\delta (n+P)^2$ proper relation of size $2$. 
Since each proper relation is counted twice in the sum $\sum_{i,j\in [n]}\teo{\cbc{i,j}\in \PR(\bm{A}_{n,s}[\bm{\theta}]),\mathfrak{R}_{n,s-1/n}}$, we conclude that this sum is upper bounded by $2\delta (n+P)^2$. Hence,
\begin{align}\label{proper_relation_counting}
\sum_{i,j\in [n]}\PP\bc{\cbc{i,j}\in \PR(\bm{A}_{n,s}[\bm{\theta}]),\mathfrak{R}_{n,s-1/n}}\leq 2\delta (n+P)^2.
\end{align}

Combining \cref{implication_unfreeze}, \Cref{eq_ivnave} and \Cref{proper_relation_counting} yields that 
\begin{align}\label{eq_imp_unfreeze}
    \frac{1}{n}\sum_{i\in [n]}\PP\bc{i\in \mathcal{F}\bc{\bm{A}_{n,s+1/n}[\bm{\theta}]}\backslash \mathcal{F}\bc{\bm{A}_{n,s}[\bm{\theta}]}}\leq 2\delta K c^{-1}+\bar{o}_{n,P}(1).
\end{align}

  \item 
  \textit{Unfreezing.} 
 Similarly to part (i), we show that the probability that $\{i\}\cup\supp{\bm{A}_{n,s}[\bm{\theta}](\vnu_{s+1/n},)}$ forms a proper relation of $\bm{A}_{n,s+1/n}[\bm{\theta}]$ is small when averaged over $i\in [n]$.  
 To this end, observe that
\begin{align}
    &\frac{1}{n}\sum_{i\in[n]}\PP\bc{\{i\}\cup\supp{\bm{A}_{n,s+1/n}[\bm{\theta}](\vnu_{s+1/n},)}\in \PR(\bm{A}_{n,s+1/n}[\bm{\theta}]),\mathfrak{G}_{\eps,s+1/n}} \nonumber \\
\leq&\frac{1}{n}\sum_{i\in[n]} \sum_{k=1}^K\sum_{1\leq j_1<j_2\ldots<j_k\leq n} \PP\bc{\cbc{i, j_1,\ldots,j_k}\in \PR(\bm{A}_{n,s+1/n}[\bm{\theta}])} \nonumber \\
&\times\PP\bc{\supp{\bm{A}_{n,s}(\vnu_{s+1/n},)}=\cbc{j_1,\ldots,j_k},\mathfrak{G}_{\eps,s+1/n} \mid \cbc{i, j_1,\ldots,j_k}\in \PR(\bm{A}_{n,s+1/n}[\bm{\theta}])}. \label{eq_unfreezing}
\end{align}
We next bound the two probability terms in \eqref{eq_unfreezing} separately, by firstly showing an anticoncentration property of $\supp{\bm{A}_{n,s}(\vnu_{s+1/n},)}$ on $\mathfrak{G}_{\eps,s+1/n}$, and secondly using that there are only few proper relations on $\mathfrak R_{n,s}$.

Let $\bm{k}$ be the number of half-edges of $\vnu_{s+1/n}$ that connect to vertices in $\verst{s}$, which is an upper bound on the number of non-zero entries of row $\vnu_{s+1/n}$ in $\bm A_{n,s}$. We next determine the distribution of their neighbors given  \scm{n}{s+1/n}. 
For this, order these $\bm{k}$ half-edges arbitrarily as half-edges $1,2,\ldots,\bm{k}$. 
We call a half-edge incident to vertex $j\in \verst{s}$ \textit{free after  
$m \in [\bm{k}]$ pairings of $\vnu_{s+1/n}$} 
if it is neither paired with half-edges $1, \ldots m-1$ of $\vnu_{s+1/n}$ nor with any vertex from $\verst{s+1/n}$ in \scm{n}{s+1/n}. These are exactly the half-edges whose neighbors are unknown when conditioned on \scm{n}{s+1/n}, $\vnu_{s+1/n}$ and the neighbors of half-edges $1, \ldots, m-1$.

\begin{claim}\label{claim_extra}
Given \scm{n}{s+1/n},  $\vnu_{s+1/n}$, $\bm{k}$ and the neighbors of half-edges $1, \ldots, m-1$ for $m\in [\bm{k}]$, the neighbor of half-edge $m$ is chosen uniformly from all half-edges that are free after 
$m-1$ pairings of $\vnu_{s+1/n}$.
\end{claim}

\begin{proof}[Proof of \Cref{claim_extra}]
    The proof uses Bayes' theorem and is analogous to the derivation of \Cref{lem-half-edge-exchangeable}.
\end{proof}

Moreover, for $\fscm{n}{s+1/n}$ such that $\mathfrak{G}_{\eps,s+1/n}$ holds, the number of free half-edges at any stage $m$ is at least $ cn-K$. Let $\bm v_1, \ldots, \bm v_{\bm k} \in \verst{s+1/n}$ denote the vertices that half-edges $1, \ldots, \bm k$ are connected to, and let $v$ be a vertex incident to a half edge that is free at stage $m$. Then  
by \Cref{claim_extra},
\begin{align}\label{eq_qlvs1n}
  &\ind{\mathfrak{G}_{\eps,s+1/n}} \PP\bc{\bm v_m = v \mid \fscm{n}{s+1/n},\bm{k}, \bm v_1, \ldots, \bm v_{m-1}} \leq \frac{K}{cn-K}.
\end{align}
Multiplying \eqref{eq_qlvs1n} for all $m\in [\bm{k}]$ gives that for all $v_1, \ldots, v_{\bm k} \in \verst{s+1/n}$,
\begin{align}\label{eq_prob_q}
 \ind{\mathfrak{G}_{\eps,s+1/n}}  \PP\bc{\bm v_1 = v_1, \ldots, \bm v_{\bm k} = v_{\bm k} \mid \fscm{n}{s+1/n},\bm{k}} \leq \bc{\frac{K}{cn-K}}^{\bm{k}}.
\end{align}
Now, fix $k \in [K]$ and $j_1, \ldots, j_k \in [n]$. There are at most $\bm k^{\bm k}$ choices of $(v_1, \ldots, v_{\bm k})$ that give rise to $\supp{\bm{A}_{n,s}(\vnu_{s+1/n},)}$ $=\cbc{j_1,\ldots,j_{k}}$. Summing over such $(v_1, \ldots, v_{\bm k})$ in \eqref{eq_prob_q} and using that $\bm k \leq K$, we arrive at
\begin{align*}
   \ind{\mathfrak{G}_{\eps,s+1/n}}  \PP\bc{\supp{\bm{A}_{n,s}(\vnu_{s+1/n},)}=\cbc{j_1,\ldots,j_{k}} \mid \fscm{n}{s+1/n},\bm{k}}\leq \frac{K^{2K}}{(cn-K)^{\bm{k}}}.
\end{align*}
Since $k\leq \bm{k}\leq K$ and $c<1$, we can use the tower property to remove the conditioning on $\bm k$:
\begin{align}\label{eq_con_halfedge}
 \ind{\mathfrak{G}_{\eps,s+1/n}}  \PP\bc{\supp{\bm{A}_{n,s}(\vnu_{s+1/n},)}=\cbc{j_1,\ldots,j_k} \mid \fscm{n}{s+1/n}}\leq \frac{K^{2K}}{c^K n^k}(1+\bar{o}_n(1)).
\end{align}
 Furthermore, \Cref{eq_con_halfedge}, the tower property and the fact that the support of $\bm{A}_{n,s}(\vnu_{s+1/n},)$ is independent of the perturbation matrices $\THETA$ give that
\begin{align}\label{eq_con_halfedge21}
    &\PP\bc{\supp{\bm{A}_{n,s}(\vnu_{s+1/n},)}=\cbc{j_1,\ldots,j_k},\mathfrak{G}_{\eps,s+1/n} \mid \cbc{i, j_1,\ldots,j_k}\in \PR(\bm{A}_{n,s+1/n}[\bm{\theta}])}\nonumber\\
    =&\Erw[\PP(\supp{\bm{A}_{n,s}(\vnu_{s+1/n},)}=\{j_1,\ldots,j_k\},\mathfrak{G}_{\eps,s+1/n} \mid \fscm{n}{s+1/n},\THETA) \mid \{i,\ldots,j_k\}\in \PR(\bm{A}_{n,s+1/n}[\bm{\theta}])]\nonumber\\
    \leq & \frac{K^{2K}}{c^K n^k}+\bar{o}_n(n^{-k}).
\end{align}
Returning to \eqref{eq_unfreezing}, we have thus obtained the upper bound
\begin{align*}
    &\frac{1}{n}\sum_{i\in[n]}\PP\bc{\{i\}\cup\supp{\bm{A}_{n,s+1/n}[\bm{\theta}](\vnu_{s+1/n},)}\in \PR(\bm{A}_{n,s+1/n}[\bm{\theta}]),\mathfrak{G}_{\eps,s+1/n}}\\
\leq &\sum_{k=1}^K\frac{K^{2K}}{c^K n^{k+1}}\sum_{1\leq j_1<\ldots<j_k\leq n}\sum_{i\in[n]} \PP\bc{\cbc{i, j_1,\ldots,j_k}\in \PR(\bm{A}_{n,s+1/n}[\bm{\theta}]),\mathfrak{R}_{n,s}}+\bar{o}_{n,P}(1).
\end{align*}  
As in the proof of \cref{proper_relation_counting}, the sum $\sum_{1\leq j_1<\ldots<j_k\leq n}\sum_{i\in[n]}  \teo{\cbc{i,j_1,\ldots,j_k}\in \PR(\bm{A}_{n,s+1/n}[\bm{\theta}]),\mathfrak{R}_{n,s}}$ counts each proper relation of size $k$ exactly $k$ times and each proper relation of size $k+1$ exactly $k+1$ times. However, on $\mathfrak R_{n,s}$, $ \bm{A}_{n,s+1/n}[\bm{\theta}]$ is 
$(\delta,\ell)$-free for $2\leq \ell\leq 2K$, so that 
\begin{align}\label{proper_relation_counting2}
&\sum_{1\leq j_1<\ldots<j_k\leq n}\sum_{i\in[n]}  \PP\bc{\cbc{i,j_1,\ldots,j_k}\in \PR(\bm{A}_{n,s+1/n}[\bm{\theta}]),\mathfrak{R}_{n,s}}\\
&\qquad\leq \delta (k (n+P)^{k}+(k+1)(n+P)^{k+1}).\nonumber
\end{align}
Therefore, 
    \begin{align*}
   &  \frac{1}{n}\sum_{i\in[n]}\PP\bc{\{i\}\cup\supp{\bm{A}_{n,s+1/n}[\bm{\theta}](\vnu_{s+1/n},)}\in \PR(\bm{A}_{n,s+1/n}[\bm{\theta}])}
\leq \sum_{k=1}^{K}\frac{\delta K^{2K}(K+1)}{c^K}+\bar{o}_{n,P}(1)\nonumber\\
\leq & \hspace{0.1 cm} \delta K^{2K+1}(K+1)c^{-K}+\bar{o}_{n,P}(1).
\end{align*}
Finally, by the deterministic upper bound $\sum_{i\in[n]}\teo{i\in \{\vnu_{s+1/n}\}\cup\supp{\bm{A}_{n,s}[\bm{\theta}](\vnu_{s+1/n},)}} \leq K+P+1$, 
 \begin{align}\label{eq_ex_supp_i}
      &\frac{1}{n}\sum_{i\in[n]}\PP\bc{i\in \{\vnu_{s+1/n}\}\cup\supp{\bm{A}_{n,s}[\bm{\theta}](\vnu_{s+1/n},)}}
   \leq  \frac{K+P+1}{n}=\bar{o}_n(1). 
  \end{align}

Hence, by \cref{implication_freeze}, 
\begin{equation}\label{el702}
\frac{1}{n}\sum_{i\in[n]}\mathbb{P}\left(i\in \mathcal{F}\bc{\bm{A}_{n,s}[\bm{\theta}]} \setminus \mathcal{F}\bc{\bm{A}_{n,s+1/n}[\bm{\theta}]}\right)\leq \delta K^{2K+1}(K+1)c^{-K}+\bar{o}_{n,P}(1).
\end{equation}
\end{enumerate}

Combining \cref{eq_imp_unfreeze} and \cref{el702}
finishes the proof of \Cref{l7}.\qed

\subsection{One-step matrix zeroing, row-zeroed matrix: Proof of \Cref{l71}}
The proof is an adaption of the proof of \cite[Lemma 4.18]{van2022random}. Again, we relate the `type change' events in question to the existence of proper relations. More precisely, for $i \in [n]$, we will show that on  a reasonably likely event $\mathfrak S_n(i)$,
\begin{align}\label{eq_rel_imp}
    &i\in \mathcal{F}\bc{\bm{A}_{n,s}[\bm{\theta}]\abc{i;}} \Delta \mathcal{F}\bc{\bm{A}_{n,s+1/n}[\bm{\theta}]\abc{i;}}\nonumber \\
    \Longrightarrow \quad & \supp{\bm{A}_{n,s}(,i)} \text{ is a proper relation of } \bm{A}_{n,s}\abc{i;i}[\bm{\theta}]^T \text{ or of } \bm{A}_{n,s+1/n}\abc{i;i}[\bm{\theta}]^T.
\end{align}
The proof now proceeds in four steps: First, we define the events $\mathfrak S_n(i)$. Second, we derive an upper bound on the `average' probability of $\mathfrak S_n(i)^c$ over $i \in [n]$. This requires the bulk of the work. Third, we prove the implication \eqref{eq_rel_imp}. Finally, we show that the probability of the events on the right-hand side of \eqref{eq_rel_imp} is uniformly bounded in $i$.

\textit{Definition of $\mathfrak S_n(i)$.} The event that we work on is an intersection of three events: First, set 
\begin{align*}
    \mathfrak{S}_{n,1}(i)= \mathfrak P_{n,s} \cap \cbc{i\neq \vnu_{s+1/n},\bm{A}_n(\vnu_{s+1/n},i)=0, \THETA_r[\bm{\theta}_r,n](,i)=0_{\bm{\theta}_r \times 1}, \THETA_c[n,\bm{\theta}_c](i,)=0_{1\times \bm{\theta}_c}
   }.
\end{align*}
On $\mathfrak S_{n,1}(i)$, 
the non-zero entries of column (or row) $i$ in all involved matrices are contained in $[n]$, i.e.,
\begin{align}\label{eq_eq_supp}
    \supp{\bm{A}_{n,s}(,i)}=\supp{\bm{A}_{n,s}[\bm{\theta}](,i)}=\supp{\bm{A}_{n,s+1/n}[\bm{\theta}](,i)}=\supp{\bm{A}_{n,s+1/n}(,i)}.
\end{align}
Next, let
\begin{align*}
    \mathfrak{S}_{n,2}(i)=\cbc{\text{{for} all } j\in \supp{\bm{A}_{n,s}(,i)}: j \notin \mathcal{F}\bc{\bm{A}_{n,s}\abc{i;i}[\bm{\theta}]^T} \Delta \mathcal{F}\bc{\bm{A}_{n,s+1/n}\abc{i;i}[\bm{\theta}]^T}}
\end{align*}
be the event that no element of the support of column $i$ has a different type in $\bm{A}_{n,s}\abc{i;i}[\bm{\theta}]^T$ than in $\bm{A}_{n,s+1/n}\abc{i;i}[\bm{\theta}]^T$. 
Finally, let
\begin{align*}
    \mathfrak{S}_{n,3}(i)& =\cbc{\bm{A}_{n,s}\abc{i;i}[\bm{\theta}]^T  (\delta,\ell)\text{-free for } 2\leq \ell\leq K+1} \cap \cbc{ \bm{A}_{n,s+1/n}\abc{i;i}[\bm{\theta}]^T  (\delta,\ell)\text{-free for } 2\leq \ell\leq K+1} \\
    &= \mathfrak{S}_{n,3,s}(i) \cap \mathfrak{S}_{n,3,s+1/n}(i).
\end{align*}
We set $\mathfrak S_n(i) = \mathfrak S_{n,1}(i) \cap \mathfrak S_{n,2}(i) \cap \mathfrak S_{n,3}(i)$. 

\textit{Bounding the average of $\PP(\mathfrak S_n(i))^c$.} We next derive an upper bound on
\begin{align}\label{eq_l71_0}
    \frac{1}{n}\sum_{i\in [n]}\PP(\mathfrak S_n(i)^c) \leq  \frac{1}{n}\sum_{i\in [n]}\PP(\mathfrak S_{n,1}(i)^c) + \frac{1}{n}\sum_{i\in [n]}\PP(\mathfrak S_{n,1}(i) \cap \mathfrak S_{n,2}(i)^c) + \frac{1}{n}\sum_{i\in n}\PP(\mathfrak S_{n,3}(i)^c). 
\end{align}
In doing so, observe that \eqref{p_event_P} and the construction of the perturbation in \Cref{d34} ensure that
\begin{align}
\frac{1}{n}\sum_{i\in[n]}\PP\bc{\mathfrak S_{n,1}(i)^c}\leq \bar{o}_{n,P}(1) + \frac{K}{n}+\frac{1}{n} = \bar{o}_{n,P}(1).
\label{eq_sum_bound_eve1}    
\end{align}
Moreover, \Cref{p1} implies that
\begin{align}\label{upperbound_set3}
\frac{1}{n}\sum_{i\in[n]}\PP\bc{\mathfrak S_{n,3}(i)^c}= \bar{o}_{n,P}(1).
\end{align}
We are thus left with bounding the middle sum on the r.h.s. of \eqref{eq_l71_0}. As $\mathfrak S_{n,1}(i) \cap \mathfrak S_{n,2}(i)^c$ involves a `status change' event with respect to a matrix in which row and column $\vnu_{s+1/n}$ are replaced by a zero row and a zero column, the basic approach that we pursue to bound the last missing sum in \cref{eq_l71_0} is similar to the one of \Cref{l7}. However, in the present case, the vertices for which we investigate a status change are not chosen uniformly among all vertices, but from the support of the $i$th column of $\bm A_{n,s}$, and we will need to address this complication. 
But first, a derivation analogous to the one of implications (\ref{implication_unfreeze}) and (\ref{implication_freeze}) shows that
 \begin{align*}
    \mathfrak S_{n,1}(i)\cap \mathfrak S_{n,2}(i)^c &  \subseteq \mathfrak S_{n,1}(i)\cap \Big(\big\{\ \exists\ j\in\supp{\bm{A}_{n,s}(,i)} :  \cbc{j,\vnu_{n,s+1/n}} \in \PR\bc{\bm{A}_{n,s}\abc{i;i}[\bm{\theta}]^T}\big\}\nonumber\\
     & \cup \big\{ \ \exists\ j\in\supp{\bm{A}_{n,s}(,i) }: \cbc{j}\cup \supp{\bm{A}_{n,s}\abc{i;i}[\bm{\theta}]^T(\vnu_{s+1/n},)} \in\PR\bc{\bm{A}_{n,s+1/n}\abc{i;i}[\bm{\theta}]^T}\big\}\nonumber\\
     & \cup \big\{ \ \exists\ j\in\supp{\bm{A}_{n,s}(,i)} :  j\in   \{\vnu_{s+1/n}\}\cup\supp{\bm{A}_{n,s}\abc{i;i}[\bm{\theta}]^T(\vnu_{s+1/n},)}\big\}\Big).
\end{align*}
Since $\bm{A}_n(i,\vnu_{s+1/n})=\bm{A}_n(\vnu_{s+1/n},i)=0$ as well as $\THETA_c[n,\bm{\theta}_c](i,)=0_{1\times \bm{\theta}_c}$ on $\mathfrak S_{n,1}(i)$, on this event, $\supp{\bm{A}_{n,s}\abc{i;i}[\bm{\theta}]^T(\vnu_{s+1/n},)} = \supp{\bm{A}_{n,s}(\vnu_{s+1/n},)}$, so that we may simplify the last inclusion as stating that
 \begin{align}
     \mathfrak S_{n,1}(i)\cap \mathfrak S_{n,2}(i)^c &  \subseteq \mathfrak S_{n,1}(i)\cap \big( \mathfrak E_{n,1}(i)\cup \mathfrak E_{n,2}(i) \cup \mathfrak E_{n,3}(i) \big), \label{es1cs2i}
\end{align}
where
 \begin{align*}
     \mathfrak E_{n,1}(i) &= \big\{\ \exists\ j\in \supp{\bm{A}_{n,s}(,i)} :  \cbc{j,\vnu_{n,s+1/n}} \in \PR\bc{\bm{A}_{n,s}\abc{i;i}[\bm{\theta}]^T}\big\},\nonumber\\
     \mathfrak E_{n,2}(i) &=   \big\{ \ \exists\ j\in\supp{\bm{A}_{n,s}(,i)}: \cbc{j}\cup \supp{\bm{A}_{n,s}(\vnu_{s+1/n},)} \in\PR\bc{\bm{A}_{n,s+1/n}\abc{i;i}[\bm{\theta}]^T}\big\},\nonumber\\
     \mathfrak E_{n,3}(i) &=   \big\{ \ \exists\ j\in\supp{\bm{A}_{n,s}(,i)}:  j\in   
     \supp{\bm{A}_{n,s}(\vnu_{s+1/n},)} \big\}.
\end{align*}
 As these events involve the matrices $\bm{A}_{n,s}\abc{i;i}, \bm{A}_{n,s+1/n}\abc{i;i}$ that arise from the removal of vertex $i$, we define a corresponding graph as follows: For $i\in [n]$, let \scm{n}{s,\aabc{i}} be the graph that is induced by \cm{n} on the vertex set $\verst{s}\backslash\cbc{i}$. Moreover,  for $j\in \verst{s}\backslash\cbc{i}$, let $\bar{\vd}_{j,s,\aabc{i}}$ be the (current) degree of vertex $j$ in \scm{n}{s,\aabc{i}}. 
We emphasize that via its vertex count, the graph \scm{n}{s,\aabc{i}} contains the information whether or not $i\in \verst{s}$. 
Finally, 
on the event $\mathfrak{G}_{\eps,s}$ from \cref{event_gs}, 
\begin{align}\label{lower_degree_irem}
\sum_{j\in \verst{s}\backslash\cbc{i}}(d_j-\bar{\vd}_{j,s,\aabc{i}})\geq \sum_{j\in \verst{s}\backslash\cbc{i}}(d_j-\bar{\vd}_{j,s}) \geq cn - (d_i-\bar{\vd}_{i,s} )  \geq cn-K.
\end{align}
We next bound the probability of the three events on the right-hand side of \eqref{es1cs2i} separately.
\begin{enumerate}
    \item \label{it_b_e_1} We first establish an upper bound on $\PP\bc{\mathfrak E_{n,1}(i)}$ from \eqref{es1cs2i}.
  Note that by \Cref{p1} and \cref{ineq_dcs},
\begin{align}\label{eq_set2_1}
    \PP\bc{\mathfrak E_{n,1}(i)} \leq & \sum_{j,\ell\in[n]}\PP\bc{j\in\supp{\bm{A}_{n,s}(,i)},\vnu_{n,s+1/n}=\ell,\mathfrak{G}_{\eps,s},\cbc{j,\ell} \in \PR\bc{\bm{A}_{n,s}\abc{i;i}[\bm{\theta}]^T},\mathfrak S_{n,3,s}(i)} \nonumber \\
    & \quad + \bar{o}_{n,P}(1).
\end{align}
By the tower property and independence of the perturbation matrix and the graph structure, for each $j, \ell \in [n]$,
\begin{align}\label{eq_set2_0}
    &\PP\bc{j\in\supp{\bm{A}_{n,s}(,i)},\vnu_{n,s+1/n}=\ell,\mathfrak{G}_{\eps,s},\cbc{j,\ell} \in \PR\bc{\bm{A}_{n,s}\abc{i;i}[\bm{\theta}]^T},\mathfrak S_{n,3,s}(i)}\\
    &\quad=\Erw\Big[\PP\bc{j\in\supp{\bm{A}_{n,s}(,i)},\vnu_{n,s+1/n}=\ell,\mathfrak{G}_{\eps,s} \mid \fscm{n}{s,\aabc{i}}}\nonumber\\
    &\qquad\times\teo{\cbc{j,\ell} \in \PR\bc{\bm{A}_{n,s}\abc{i;i}[\bm{\theta}]^T},\mathfrak S_{n,3,s}(i)}\Big]. \nonumber
\end{align}
By \cref{eq_conp_nus}, for each $j,\ell\in [n]$,
\begin{align}\label{eq_remove_i}
    &\PP\bc{j\in\supp{\bm{A}_{n,s}(,i)},\vnu_{s+1/n}=\ell,\mathfrak{G}_{\eps,s} \mid \fscm{n}{s,\aabc{i}}}\nonumber\\
    &\qquad= \Erw\brk{\teo{j\in\supp{\bm{A}_{n,s}(,i)},\mathfrak{G}_{\eps,s}}\PP\bc{\vnu_{s+1/n}=\ell \mid \fscm{n}{s}} \mid \fscm{n}{s,\aabc{i}}}\nonumber\\
    &\qquad\leq \frac{K}{cn}\PP\bc{j\in \supp{\bm{A}_{n,s}(,i)},\mathfrak{G}_{\eps,s} \mid \fscm{n}{s,\aabc{i}}}.
\end{align}
As in \Cref{lem-half-edge-exchangeable} and \Cref{claim_extra}, it can be shown that conditionally on \scm{n}{s,\aabc{i}}, each vertex $j\in \verst{s}\backslash\cbc{i}$ has a probability proportional to $d_j-\bar{\vd}_{j,s,\aabc{i}}$ of being chosen as a neighbor of a half-edge of  $i$.   
Combined with \eqref{lower_degree_irem}, this yields 
\begin{align}\label{eq_con_halfedge2}
    \PP\bc{j\in\supp{\bm{A}_{n,s}(,i)},\mathfrak{G}_{\eps,s} \mid \fscm{n}{s,\aabc{i}}}\leq d_i\frac{d_j-\bar{\vd}_{j,s,\aabc{i}}}{\sum_{j\in \verst{s}\backslash\cbc{i}}(d_j-\bar{\vd}_{j,s,\aabc{i}})}\leq \frac{K^2}{cn-K}.
\end{align}
Plugging \cref{eq_con_halfedge2} into \cref{eq_remove_i}, we arrive at 
\begin{align}\label{ineq_suppvs_precise}
    \PP\bc{j\in\supp{\bm{A}_{n,s}(,i)},\vnu_{s+1/n}=\ell,\mathfrak{G}_{\eps,s}  \mid \fscm{n}{s,\aabc{i}}}
    \leq \frac{K^3}{(cn-K)^2}.
\end{align} 
Plugging  \eqref{ineq_suppvs_precise} into \eqref{eq_set2_0} yields that
\begin{align}\label{eq_set2_2}
    &\PP\bc{j\in\supp{\bm{A}_{n,s}(,i)},\vnu_{n,s+1/n}=\ell, \mathfrak{G}_{\eps,s},\cbc{j,\ell} \in \PR\bc{\bm{A}_{n,s}\abc{i;i}[\bm{\theta}]^T},\mathfrak S_{n,3,s}(i) }\\
    &\leq \frac{K^3}{(cn-K)^2}\PP\bc{\cbc{j,\ell} \in \PR\bc{\bm{A}_{n,s}\abc{i;i}[\bm{\theta}]^T},\mathfrak S_{n,3,s}(i)}.\nonumber
\end{align}
Moreover, analogously to the proof of \Cref{proper_relation_counting},
\begin{align}\label{proper_relation_counting_2}
    \sum_{j,\ell\in [n]}\PP\bc{\cbc{i, j}\in \PR(\bm{A}_{n,s}\abc{i;i}[\bm{\theta}]^T),\mathfrak S_{n,3,s}(i)}\leq 2\delta (n+P)^2.
\end{align}
Combining \Cref{ineq_dcs}, \Cref{upperbound_set3}, \Cref{eq_set2_1}, \Cref{eq_set2_2} and \Cref{proper_relation_counting_2} gives
\begin{align}\label{eq_it_b_e_1}
    \PP&\bc{\mathfrak E_{n,1}(i)}
    \leq 2\delta K^3c^{-2}+\bar{o}_{n,P}(1).
\end{align}

\item  We next establish an upper bound on $\PP\bc{\mathfrak E_{n,2}(i)}$ from \eqref{es1cs2i}. As in the previous case, by  \Cref{p1},
\begin{align}\label{proper_relation_counting2_0}
    \PP\bc{\mathfrak E_{n,2}(i)} &
  \leq  \sum_{k=0}^K\sum_{1\leq j_1<j_2\ldots<j_k\leq n}\sum_{j\in[n]} \PP\big(j\in\supp{\bm{A}_{n,s+1/n}(,i)},\supp{\bm{A}_{n,s}(,\vnu_{s+1/n})}=\cbc{j_1,\ldots,j_{k}}, \nonumber \\
  & \qquad  \cbc{j,j_1,\ldots,j_k}\in \PR\bc{\bm{A}_{n,s+1/n}\abc{i;i}[\bm{\theta}]^T},\mathfrak{G}_{\eps,s+1/n},\mathfrak S_{n,3, s+1/n}(i)\big) + \bar{o}_{n,P}(1).
\end{align}
Using the same conditioning approach as in \cref{eq_set2_0} and \eqref{eq_remove_i} (first conditioning on $\fscm{n}{s+1/n,\aabc{i}}$, then on $\fscm{n}{s+1/n}$), we may first bound the following conditional probabilities on the support of column $i$ in $\bm{A}_{n,s+1/n}$ and column $\vnu_{s+1/n}$ in $\bm{A}_{n,s}$.
Given $0\leq k\leq K$ and $1\leq j_1<\ldots< j_k\leq n$, 
by \Cref{eq_con_halfedge},
\begin{align}\label{cond_step_0}
  \mathds{1}\mathfrak{G}_{\eps,s+1/n}  \PP\bc{\supp{\bm{A}_{n,s}(,\vnu_{s+1/n})}=\cbc{j_1,\ldots,j_{k}} \mid \fscm{n}{s+1/n}}\leq \frac{K^{2K}}{c^K n^{k}}(1+\bar{o}_n(1)).
\end{align}
Furthermore, for $j\in[n]$, by \Cref{eq_con_halfedge2}, 
\begin{align}\label{eq_con_halfedge2__2}
    \PP\bc{j\in\supp{\bm{A}_{n,s+1/n}(,i)}, \mathfrak{G}_{\eps,s+1/n} \mid \fscm{n}{s+1/n,\aabc{i}}}\leq \frac{K^2}{cn-K}.
\end{align}
Combining \cref{cond_step_0} and \cref{eq_con_halfedge2__2} with the outlined conditioning, we obtain the upper bound
\begin{align}\label{eq_set2_3}
    \frac{K^{2K+2}}{c^{K+1}n^{k+1}}(1+\bar{o}_n(1))\PP\bc{\cbc{j,j_1,\ldots,j_k} \in \PR\bc{\bm{A}_{n,s+1/n}\abc{i;i}[\bm{\theta}]^T},\mathfrak S_{n,3, s+1/n}(i)}
\end{align}
for the single summands in \cref{proper_relation_counting2_0}.
Finally, analogously to the proof of \Cref{proper_relation_counting2},
\begin{align}\label{proper_relation_counting2_2}
  \sum_{1\leq j_1<j_2\ldots<j_k\leq n}&\sum_{j\in[n]} \PP\bc{\cbc{j,j_1,\ldots,j_k} \in \PR\bc{\bm{A}_{n,s+1/n}\abc{i;i}[\bm{\theta}]^T},\mathfrak S_{n,3, s+1/n}(i)} \nonumber \\
  \leq & \delta (k (n+P)^{k}+(k+1)(n+P)^{k+1}).
\end{align}
Combining  \Cref{eq_con_halfedge2__2}, \Cref{eq_set2_3} and \Cref{proper_relation_counting2_2} gives that
\begin{align}
   \PP&\bc{\mathfrak E_{n,2}(i)} \leq \delta (K+1)K^{2K+2}c^{-K-1}+\bar{o}_{n,P}(1).\label{eq_it_b_e_2}
\end{align}

\item Finally, we establish an upper bound on $\PP\bc{\mathfrak S_{n,1}(i) \cap \mathfrak E_{n,3}(i)}$ from \eqref{es1cs2i}. 
Observe that on $\mathfrak G_{n,1}(i)$, $\supp{\bm{A}_{n,s}(,i)} = \supp{\bm{A}_{n,s+1/n}(,i)}$, so that, as in the previous two cases, by \cref{ineq_dcs},
\begin{align*}
 \PP\bc{\mathfrak S_{n,1}(i) \cap \mathfrak E_{n,3}(i)}\leq  \sum_{j\in[n]}\PP\bc{j\in\supp{\bm{A}_{n,s+1/n}(,i)},j \in \supp{\bm{A}_{n,s}(,\vnu_{s+1/n})},\mathfrak{G}_{\eps ,s+1/n}} + \bar{o}_n(1).
\end{align*}
By the tower property,
\begin{align}\label{eq_tower_twosupport}
    &\PP\bc{j\in\supp{\bm{A}_{n,s+1/n}(,i)},j \in \supp{\bm{A}_{n,s}(,\vnu_{s+1/n})},\mathfrak{G}_{\eps,s+1/n}}\\
    &\qquad=\Erw\brk{ \teo{j\in\supp{\bm{A}_{n,s+1/n}(,i)},\mathfrak{G}_{\eps,s+1/n}}\PP\bc{j \in \supp{\bm{A}_{n,s}(\vnu_{s+1/n},)}\mid \fscm{n}{s+1/n}}}.\nonumber
\end{align}
On the other hand, by \Cref{eq_qlvs1n}, 
\begin{align}\label{eq_con_halfedge3}
 \mathds{1} \mathfrak{G}_{\eps,s+1/n}   \PP\bc{j \in \supp{\bm{A}_{n,s}(,\vnu_{s+1/n})}\mid \fscm{n}{s+1/n}} \leq \frac{K^2}{cn-K}\mathds{1} \mathfrak{G}_{\eps,s+1/n}.
\end{align}
Hence, by 
\Cref{eq_tower_twosupport} and \Cref{eq_con_halfedge3},
\begin{align*}
  \PP\bc{\mathfrak S_{n,1}(i) \cap \mathfrak E_{n,3}(i)}\leq \frac{K^2}{cn-K} \Erw\brk{ \sum_{j\in[n]} \teo{j \in \supp{\bm{A}_{n,s}(,\vnu_{s+1/n})} } } + \bar{o}_n(1) \leq \frac{K^2}{cn-K} + \bar{o}_n(1).
\end{align*}
\end{enumerate}

Plugging in  \eqref{eq_sum_bound_eve1}, \eqref{upperbound_set3} and the last three cases into \eqref{eq_l71_0} 
finally gives that
\begin{align}\label{eq_l71_1}
    \frac{1}{n}\sum_{i\in [n]}\PP(\mathfrak S_n(i)^c) \leq 2\delta (K+1)K^{2K+2}c^{-K-1}+\bar{o}_{n,P}(1).
\end{align}
\textit{Proof of implication (\ref{eq_rel_imp}) by contraposition.} Observe that if supp$(\bm{A}_{n,s}(,i))=\emptyset$, given $\mathfrak P_{n,s}$, clearly the premise of \cref{eq_rel_imp} is not satisfied and there is nothing to prove. We thus assume supp$(\bm{A}_{n,s}(,i))\not=\emptyset$. Suppose that $\mathfrak S_n(i)$ 
holds and that supp$(\bm{A}_{n,s}(,i))$ is neither a proper relation of $\bm{A}_{n,s}\abc{i;i}[\bm{\theta}]^T$ nor of  $ \bm{A}_{n,s+1/n}\abc{i;i}[\bm{\theta}]^T$. 
On $\mathfrak S_n(i)$,
\begin{align*}
  \supp{\bm{A}_{n,s}(,i)}=\supp{\bm{A}_{n,s}[\bm{\theta}]\abc{i;}(,i)}=\supp{\bm{A}_{n,s+1/n}[\bm{\theta}]\abc{i;}(,i)}, 
\end{align*}
so that we may regard supp$(\bm{A}_{n,s}(,i))$ as the set of non-zero entries of the $i$th column in any of these matrices, as this is unaffected by the perturbation and replacing rows $i, \bm \nu_{s+1/n}$ by a zero rows and column $\bm \nu_{s+1/n}$ by a zero column. Moreover, again on $\mathfrak S_n(i)$,
\begin{align*}
   \bm{A}_{n,s}\abc{i;i}[\bm{\theta}] = \bm{A}_{n,s}[\bm{\theta}]\abc{i;i}= \bc{\bm{A}_{n,s}[\bm{\theta}]\abc{i;} }\abc{;i},  \bm{A}_{n,s+1/n}\abc{i;i}[\bm{\theta}]=  \bm{A}_{n,s+1/n}[\bm{\theta}]\abc{i;i} = \bc{\bm{A}_{n,s+1/n}[\bm{\theta}]\abc{i;} }\abc{;i},
\end{align*}
so that the assumption on proper relations concerns the $i$th column of $\bm{A}_{n,s}[\bm{\theta}]\abc{i;}$ and $\bm{A}_{n,s+1/n}[\bm{\theta}]\abc{i;}$, respectively.
Since the support may still form a (non-proper) relation, we distinguish four cases:
\begin{itemize}
    \item[\bf{Case 1:}] $\supp{\bm{A}_{n,s}(,i)}$ has no representation in  $ \bm{A}_{n,s}[\bm{\theta}]\abc{i;i}^T$ and no representation  in $ \bm{A}_{n,s+1/n}[\bm{\theta}]\abc{i;i}^T$. \\
    As there is no representation of the non-zero entries of the $i$th column of $\bm{A}_{n,s}[\bm{\theta}]\abc{i;}$ in $\bm{A}_{n,s}[\bm{\theta}]\abc{i;i}^T$, 
    this column cannot lie in the column space of $\bm{A}_{n,s}[\bm{\theta}]\abc{i;i}$. Thus, by the second equivalence in \Cref{im_fnf}, $i$ is frozen in $\bm{A}_{n,s}[\bm{\theta}]\abc{i;}$. Applying the same reasoning to $\bm{A}_{n,s+1/n}[\bm{\theta}]\abc{i;}$ yields $ i\in \mathcal{F}\bc{\bm{A}_{n,s}[\bm{\theta}]\abc{i;}} \cap \mathcal{F}\bc{\bm{A}_{n,s+1/n}[\bm{\theta}]\abc{i;}}$.
    \item[\bf{Case 2:}]   supp$(\bm{A}_{n,s}(,i))$ has a representation both in  $ \bm{A}_{n,s}[\bm{\theta}]\abc{i;i}^T$ and in $ \bm{A}_{n,s+1/n}[\bm{\theta}]\abc{i;i}^T$. \\
   Since we also assume that supp$(\bm{A}_{n,s}(,i))$ is not a proper relation in $ \bm{A}_{n,s}[\bm{\theta}]\abc{i;i}^T$ and in $ \bm{A}_{n,s+1/n}[\bm{\theta}]\abc{i;i}^T$, all variables in $\supp{\bm{A}_{n,s}(,i)}$ are frozen both in $ \bm{A}_{n,s+1/n}[\bm{\theta}]\abc{i;i}^T$ and in $ \bm{A}_{n,s}[\bm{\theta}]\abc{i;i}^T$. In particular, column $i$ of $\bm{A}_{n,s+1/n}[\bm{\theta}]\abc{i;}$ is contained in the column space of $\bm{A}_{n,s+1/n}[\bm{\theta}]\abc{i;i}$ and column $i$ of $\bm{A}_{n,s}[\bm{\theta}]\abc{i;}$ is contained in the column space of $\bm{A}_{n,s}[\bm{\theta}]\abc{i;i}$. Thus, by the second equivalence in \Cref{im_fnf}, $i$ is neither frozen in $\bm{A}_{n,s+1/n}[\bm{\theta}]\abc{i;}$ nor in $\bm{A}_{n,s}[\bm{\theta}]\abc{i;}$.
    \item[\bf{Case 3:}]  supp$(\bm{A}_{n,s}(,i))$ has a representation in  $ \bm{A}_{n,s+1/n}[\bm{\theta}]\abc{i;i}^T$, but none in $ \bm{A}_{n,s}[\bm{\theta}]\abc{i;i}^T$. \\
   As in \textbf{Case 2}, all variables in $\supp{\bm{A}_{n,s}(,i)}$ must be frozen in $ \bm{A}_{n,s+1/n}[\bm{\theta}]\abc{i;i}^T$, but there is a variable that is not frozen in $ \bm{A}_{n,s}[\bm{\theta}]\abc{i;i}^T$.
    However, this cannot happen on $\mathfrak S_{n,2}(i)$.
    \item[\bf{Case 4:}] supp$(\bm{A}_{n,s}(,i))$ has a representation in  $ \bm{A}_{n,s}[\bm{\theta}]\abc{i;i}^T$, but none in $ \bm{A}_{n,s+1/n}[\bm{\theta}]\abc{i;i}^T$. 
    
    By the same reasoning as in Case 3, this cannot happen on $\mathfrak S_{n,2}(i)$.
\end{itemize}
\textit{Excluding proper relations.} By (\ref{eq_rel_imp}), 
\begin{align}
& \PP(i\in \mathcal{F}\bc{\bm{A}_{n,s}[\bm{\theta}]\abc{i;}} \Delta \mathcal{F}\bc{\bm{A}_{n,s+1/n}[\bm{\theta}]\abc{i;}})\label{eq_bound_symdif_rowi}\\
&\qquad\leq 
     \PP\bc{\mathfrak S_n(i), \text{supp}(\bm{A}_{n,s}(,i)) \text{ is a proper relation in } \bm{A}_{n,s}[\bm{\theta}]\abc{i;i}^T \text{ or in } \bm{A}_{n,s+1/n}[\bm{\theta}]\abc{i;i}^T} + \PP\bc{\mathfrak S_n(i)^c}.\nonumber
\end{align}
Employing the same conditioning approach as in \cref{eq_set2_1} and \cref{eq_remove_i} once more, we obtain
\begin{align}\label{final_step_0}
    & \PP\bc{\mathfrak S_n(i), \text{supp}(\bm{A}_{n,s}(,i)) \text{ is a proper relation of } \bm{A}_{n,s}[\bm{\theta}]\abc{i;i}^T}  \\
   &\qquad\leq\sum_{k=2}^K\sum_{1\leq j_1\ldots<j_k\leq n}\Erw\Bigg[ \teo{\cbc{j_1,\ldots,j_k}\in \PR\bc{\bm{A}_{n,s}[\bm{\theta}]\abc{i;i}^T}} \ind{\mathfrak{S}_{n,3,s}(i)} \nonumber\\
   & \qquad\quad \times\PP\bc{\mathfrak G_{\eps, s},  \text{supp}(\bm{A}_{n,s}(,i)) =\cbc{j_1,\ldots,j_k}\mid \fscm{n}{s, \aabc{i}} }\Bigg] + \bar{o}_{n,P}(1).\nonumber
\end{align}
Analogously to \Cref{eq_con_halfedge}, 
\begin{align}\label{final_step_1}
  \mathds{1} \mathfrak G_{\eps, s}  \PP\bc{\supp{\bm{A}_{n,s}(,i)}=\cbc{j_1,\ldots,j_k} \mid \fscm{n}{s}}\leq\frac{K^{2K}}{c^K n^k}(1+\bar{o}_n(1)).
\end{align}
Moreover, by \Cref{d2},
\begin{align}\label{final_step_2}
\sum_{1\leq j_1\ldots<j_k\leq n}\Erw\Bigg[ \teo{\cbc{j_1,\ldots,j_k}\in \PR\bc{\bm{A}_{n,s}[\bm{\theta}]\abc{i;i}^T}} \ind{\mathfrak{S}_{n,3,s}(i)} \Bigg] \leq \delta K^{2K} c^{-K}+\bar{o}_{n,P}(1).
\end{align}
Using the tower property in \eqref{final_step_0} (conditioning on $\fscm{n}{s}$) in combination with \eqref{final_step_1} and \eqref{final_step_2},  we conclude that
\begin{align}
\label{eq_useinapp}
  \PP\bc{\mathfrak S_n(i), \text{supp}(\bm{A}_{n,s}(,i)) \in \PR\bc{\bm{A}_{n,s}[\bm{\theta}]\abc{i;i}^T}}    \leq \delta K^{2K+1} c^{-K}+\bar{o}_{n,P}(1).
\end{align}
Analogously, 
\begin{align*}
    \PP\bc{\mathfrak S_n(i), \supp{\bm{A}_{n,s+1/n}(,i)}\in \PR\bc{\bm{A}_{n,s+1/n}[\bm{\theta}]\abc{i;i}^T}}\leq \delta K^{2K+1} c^{-K}+\bar{o}_{n,P}(1).
\end{align*}
Hence, 
\begin{equation}\label{eq_l71_2}\begin{aligned}
&\PP(\mathfrak S_n(i), \text{supp}(\bm{A}_{n,s}(,i)) 
\in \PR\bc{\bm{A}_{n,s}[\bm{\theta}]\abc{i;i}^T} 
\cup \PR\bc{\bm{A}_{n,s+1/n}[\bm{\theta}]\abc{i;i}^T}) 
\leq 2\delta K^{2K+2} c^{-K}+\bar{o}_{n,P}(1),  
\end{aligned}
\end{equation}
Combining \eqref{eq_bound_symdif_rowi}, (\ref{eq_l71_1}) and (\ref{eq_l71_2}) yields the claim.\qed

\section{{Conditional degree of the last awakened vertex in \scm{n}{s+1/n}}}\label{complex_conexp}
Recall that the vector $\vze_{s+1/n}$ represents the type distribution with respect to the perturbed adjacency matrix $\bm{A}_{n,s+1/n}[\vth]$, which is determined by \scm{n}{s+1/n} and the perturbation $\THETA$. 
In order to apply \Cref{prop_allcomestogether} to derive the type fixed-point equations, the next step is to estimate $\PP(\vnu_{s+1/n}\in\mathcal{W}\bc{\bm{A}_{n,s}[\bm{\theta}]} \mid \vze_{s+1/n})$ from \cref{eq_allcomestogether_a2}. 

The present section prepares this estimate:  
To determine the type distribution of $\vnu_{s+1/n}$ in $\bm{A}_{n,s}[\vth]$ given $\vze_{s+1/n}$, we first compute the current degree distribution of $\vnu_{s+1/n}$ in \scm{n}{s} conditionally on \scm{n}{s+1/n}. 
Let 
\begin{align}\label{def_q_k}
    q_k:=q_k(s)=\lambda(t_s)^{-1}\sum_{\ell=k+1}^K \ell\binom{\ell-1}{k}\bc{\frac{\lambda(t_s)\eul^{2t_s}}{\lambda(0)}}^k\bc{1-\frac{\lambda(t_s)\eul^{2t_s}}{\lambda(0)}}^{\ell-k-1} \eul^{-\ell t_s} p_\ell
\end{align}
be the coefficient of $\alpha^k$ in $\hat{\psi}_{t_s}(\alpha)$ as defined in \cref{def_psi_hat}. Intuitively, $q_k$ approximates the unconditional probability that a uniform vertex in a uniformly chosen edge in \scm{n}{s} has degree $k$. The main result of the present section is the following:
\begin{proposition}\label{lem_vscm}
Uniformly in $\ranges$,
    \begin{align*}
\mathbb{E}\abs{\PP\bc{\bar{\vd}_{\vnu_{s+1/n},s}=k \mid \fscm{n}{s+1/n}}-q_k}=\bar{o}_n(1).
\end{align*}
\end{proposition}

\Cref{lem_vscm} states that  
the distribution of the current degree of $\vnu_{s+1/n}$ in \scm{n}{s}
is approximately independent of \scm{n}{s+1/n}. 
This result reflects the fact that the macroscopic structure of the graph \scm{n}{s+1/n} is constrained by inequalities such as \Cref{lem_janson_l,lem_janson_l_2}, or informally, laws of large numbers, and therefore that \scm{n}{s+1/n} only provides limited additional information about $\bar{\vd}_{\vnu_{s+1/n},s}$.

Our proof strategy towards \Cref{lem_vscm} is as follows: Recall that for $\ranges$, with high probability, $\vnu_{s+1/n}$ is awakened by \ref{exploration_step3}. Therefore, its current degree in \scm{n}{s} is strictly smaller than its original degree. Therefore, by the law of total probability, on the event that $\vnu_{s+1/n}$ is awakened by \ref{exploration_step3},
\begin{align}\label{eq_degree_LTP}
  &  \PP\bc{\bar{\vd}_{\vnu_{s+1/n},s}=k \mid \fscm{n}{s+1/n}} \nonumber \\
    = &  \sum_{\ell=k+1}^K \PP\bc{\bar{\vd}_{\vnu_{s+1/n},s}=k \mid \fscm{n}{s+1/n}, d_{\vnu_{s+1/n}}=\ell} \PP\bc{d_{\vnu_{s+1/n}}=\ell \mid \fscm{n}{s+1/n}}.
\end{align}
The benefit of the decomposition of \eqref{eq_degree_LTP} is as follows:  While the current degree $\bar{\vd}_{\vnu_{s+1/n},s}$ is a priori rather complicated to work with, the original degree $d_{\vnu_{s+1/n}}$ is more closely related to the graph exploration, as sleeping vertices are awakened with probability proportional to their \textit{original degrees} in \ref{exploration_step1} and \ref{exploration_step3}. Thus, coming from the graph exploration, it seems like a more reasonable task to compute $\PP(d_{\vnu_{s+1/n}}=\ell \mid \fscm{n}{s+1/n})$.

On the other hand, the probability that the current degree $\bar{\vd}_{\vnu_{s+1/n},s}$ of $\vnu_{s+1/n}$ in \scm{n}{s} takes on a certain value given its original degree can be computed by a similar counting argument as in \Cref{lem_bvks}.
While we provide details on this counting argument in the proof of \Cref{lem_vscm} in Section \ref{sec_proof_vscm}, the main work in this section goes into the analysis of the conditional distribution of the original degree of $\vnu_{s+1/n}$. 
 Through a careful examination of the graph exploration,  we prove in \Cref{lem_vcm} that $\PP(d_{\vnu_{s+1/n}}=\ell \mid \fscm{n}{s+1/n})$ is close to a constant:
\begin{lemma}\label{lem_vcm}
For any $\ell\in \cbc{0}\cup[K]$, uniformly in $\ranges$,
    \begin{align}\label{eq-lem_vcm-new}
\mathbb{E}\abs{\PP\bc{d_{\vnu_{s+1/n}}=\ell \mid \fscm{n}{s+1/n}}-\ell p_\ell \eul^{-t_s \ell}\lambda(t_s)^{-1}}=\bar{o}_n(1).
\end{align}
\end{lemma}

In the same spirit as \Cref{lem_vscm}, \Cref{lem_vcm} demonstrates that the conditional probability $\PP(d_{\vnu_{s+1/n}}=\ell \mid \fscm{n}{s+1/n})$ is close to the unconditional one. We first prove \Cref{lem_vcm} in Section \ref{sec_proof_vcm} and then \Cref{lem_vscm} in Section \ref{sec_proof_vscm}. The next two sections serve as a preparation towards the proof of \Cref{lem_vcm}.

\begin{remark}
  In terms of the graph exploration, the conditional probability in \Cref{lem_vcm} still is not easy to compute, 
  as it conditions on the future graph \scm{n}{s+1/n}, while seeking information about the past (the last removed vertex $\vnu_{s+1/n}$).  
  One could consider applying Bayes' rule to exchange the event and the condition.  However, the probability $\PP\bc{\fscm{n}{s+1/n}=G \mid d_{\vnu_{s+1/n}}=\ell }$ remains challenging to compute, and we still need to prove that $d_{\vnu_{s+1/n}}$ and $\fscm{n}{s+1/n}$ are approximately independent. So ultimately, this approach would not differ significantly from our current one.
\end{remark}

\subsection{The original and a modified graph exploration}
Throughout this section, we fix $\ranges$ and an integer $r>0$. The proof of \Cref{lem_vcm}, which is the main technical part in this section, is based on a `law-of-large-numbers' argument. 
For a start, it might not be too far-stretched to believe that from the perspective of \scm{n}{s+1/n}, the original degrees of the last $r$ removed vertices all behave similarly in the sense that for $u\in[r]$, 
\begin{align}\label{app_exc_d}
    \PP\bc{d_{\vnu_{s+2/n-u/n}}=\ell \mid \fscm{n}{s+1/n}}\approx \PP\bc{d_{\vnu_{s+1/n}}=\ell \mid \fscm{n}{s+1/n}}.
\end{align}
Thus, one could hope that the desired probability can be computed as an average over these last $r$ removed vertices:
\begin{align}\label{appro_od_con}
    \PP\bc{d_{\vnu_{s+1/n}}=\ell \mid \fscm{n}{s+1/n}}\approx\frac{1}{r}\Erw\brk{\abs{\cbc{u\in [r]:d_{\vnu_{s+2/n-u/n}}=\ell}} \mid \fscm{n}{s+1/n}}.
\end{align}
More specifically, at least informally, as we first let $n\to\infty$ and then $r\to\infty$,  $\frac{1}{r}\abs{\cbc{u\in [r]:d_{\vnu_{s+2/n-u/n}}=\ell}}$ should converge to a constant by `the law of large numbers'. Hence, the right-hand side of \cref{appro_od_con} converges to a constant, and so does the left-hand side.
The proof of \Cref{lem_vcm} is essentially based on this idea. 

In establishing a rigorous version of \cref{appro_od_con}, the definition of a slightly modified graph exploration will be instrumental. For this purpose, we first revisit the graph exploration from \cite{janson2009new}, as explained in \Cref{sec_cm_exp}, and make all involved randomness completely explicit -- that is, the i.i.d.~exponential life-times and the uniform choices of half-edges in \ref{exploration_step1} and \ref{exploration_step2}.  Let $\mathcal{H}$ be the set of half-edges in \cm{n}. We furthermore set 
\begin{itemize}
    \item $\bc{\bm{E}_h}_{h\in\mathcal{H}}$ to be the i.i.d.\ \textit{lifetimes} of the half-edges used in \ref{exploration_step3}, each following an $\ex{1}$ distribution; and independently given these lifetimes,
    \item $(\bm{U}(h,t,v))_{h\in \mathcal{H},t\in \cbc{0}\cup\cbc{\bm{E}_a:a\in\mathcal{H}},v\in\cbc{\text{active},\text{dead}}}$ to be  i.i.d.\ \textit{decision variables} associated to triples consisting of a half-edge $h$, a potential time $t$ at which \ref{exploration_step1} or \ref{exploration_step2} can be called (initially or after performance of \ref{exploration_step3}), and state $v$, each following a uniform distribution on $[0,1]$.
\end{itemize}
 The existence of all these random variables on the same probability space is guaranteed by the Kolmogorov extension theorem. We imagine that the decision random variables $\bm{U}(h,t,v)$ associated to half-edge $h$ can be used to set $h$ to state $v$ at time $t$. Finally, recall that in the \emph{original} graph exploration, as time $t \geq 0$ evolves, the process successively goes through \ref{exploration_step1} to \ref{exploration_step3}, where \ref{exploration_step1} and \ref{exploration_step2} are performed instantaneously, while it takes an exponential time to complete \ref{exploration_step3}. Thus, in terms of the exponential and decision random variables, the original graph exploration can be formulated as follows:
\begin{enumerate}[label=\textbf{Step \arabic*}]
\item\label{exploration_step11} At times $t$ when there is no active half-edge (such as at time $t=0$), we instantaneously choose the  half-edge associated to the \textit{ smallest decision variable} among the subset of decision variables that are indexed by half-edges that are sleeping at time $t$, time $t$ and state $v=$ active. We awaken its adjacent vertex and activate all its half-edges.
The process stops if there is no sleeping half-edge left.

\item \label{exploration_step21}  At times $t$ when there is an active half-edge, pick the one associated to the \textit{smallest decision variable} among all decision variables indexed by half-edges that are active at time $t$, time $t$ and state $v=$ dead. Change the status of this half-edge to dead.
\item \label{exploration_step31}Wait until the next half-edge dies  because of the time exceeding its lifetime. This half-edge is
joined to the one killed in \ref{exploration_step21} to form an edge of the
graph. If the vertex it belongs to is sleeping, we change the status of this vertex
to awake and all its remaining adjacent half-edges to active. Repeat from \ref{exploration_step11}.
\end{enumerate} 
We call $\bm{\eta}:=(\bm{E}_h,\bm{U}(h,t,v))_{h\in \mathcal{H},t\in \cbc{0}\cup\cbc{\bm{E}_a:a\in\mathcal{H}},v\in\cbc{\text{active},\text{dead}}}$ the \textit{decision vector} that almost surely uniquely determines the original graph exploration as described above. Unless stated otherwise, in this section, all variables in the graph exploration, such as $\vnu_s$, are derived from the original graph exploration determined by the decision vector $\bm{\eta}$. When we explicitly want to emphasize this, we also write $\vnu_s(\bm \eta)$.

We next define a modified graph exploration that is coupled to the original one. While it is rather similar to the latter, it will be easier to analyze a certain \textit{switching} procedure in the modified exploration. 
More specifically, given a decision vector $\bm \eta$, the  \textit{modified} graph exploration will proceed exactly as the original one, using \ref{exploration_step1}, \ref{exploration_step22} and \ref{exploration_step3}, where \ref{exploration_step22} is the following modification of \ref{exploration_step21}:
\begin{enumerate}[label=\textbf{Step 2'}]
    \item \label{exploration_step22} If the number of awakened vertices in the graph is at least $\lfloor ns+2-r\rfloor$ and at most $\lfloor ns\rfloor$ at time $t$, and there are active half-edges belonging to the first $\lfloor ns+1-r\rfloor$ awakened vertices, pick the active half-edge associated to the smallest decision variable among all decision variables indexed by half-edges that are active at time $t$ and belong to one of \textit{the first $\lfloor ns+1-r\rfloor$ awake vertices}, time $t$ and state $v=$ dead. Change the status of this half-edge to dead.
    Otherwise, pick the active half-edge associated to the smallest decision variable among all decision variables indexed by half-edges that are active at time $t$, time $t$ and state $v=$ dead and change the status of this half-edge to dead.
\end{enumerate}
 \ref{exploration_step22} proceeds exactly as \ref{exploration_step21}, except when the number of awakened vertices in the graph lies between $\lfloor ns+2-r\rfloor$ and  $\lfloor ns\rfloor$.  In the time-window between the awakening of $\vnu_{s+2/n-r/n}$ and $\vnu_{s+1/n}$ (including the steps where $\vnu_{s+2/n-r/n}$ and $\vnu_{s+1/n}$ are awakened), the aim of \ref{exploration_step22} is to  actively suppress the possibility that half-edges adjacent to vertices $\vnu_{s-(r-2)/n},\ldots, \vnu_{s}$, are paired by other decision mechanisms than \ref{exploration_step3}.
In particular, vertices are awakened in the same order both in the original and in the modified graph explorations, so that we can unambiguously use the notation $\vnu_s$ for both.

On a high level, the modified graph exploration facilitates a proof of \cref{app_exc_d} because it is more stable with respect to small modifications in $\bm \eta$ in the following sense: Fix $\bm \eta$ such that \ref{exploration_step11} is not performed between the awakening of $\vnu_{s+2/n-r/n}$ and $\vnu_{s+1/n}$ in either graph exploration. For $u\in [r]$, denote by $\bm h_u:=h_u(\bm{\eta})$ the first dead half-edge of $\vnu_{s+2/n-u/n}$ with respect to the original graph exploration, so that among the half-edges in $\cbc{\bm h_u:u\in[r]}$, $\bm{h}_r$ dies first and $\bm h_1$ dies last. Now given $u_1,u_2\in [r]$ with $u_1<u_2$, suppose that we modify $\bm\eta$ to $\bm \eta'$ by keeping everything fixed, apart from exchanging the exponential life-times $E_{\bm h_{u_1}}, E_{\bm h_{u_2}}$. This change has the effect that in the original graph exploration with respect to $\bm \eta'$, as \ref{exploration_step1} is not performed, the $\lfloor ns+2-u_2\rfloor$th awakened vertex is $\vnu_{s+2/n-u_1/n}(\bm \eta)$.
Consequently, in the next passes through \ref{exploration_step2}, in comparison to the exploration with respect to $\bm \eta$, other decision variables are available (first only the ones corresponding to edges adjacent to $\vnu_{s+2/n-u_1/n}(\bm \eta)$), and if these are chosen, the seemingly innocent swap might lead to an avalanche of changes, finally resulting in a significantly altered graph. 
On the other hand, in the modified graph exploration, \textit{if} there are always enough active half-edges, \ref{exploration_step22} excludes the possibility to choose from the altered decision variables, so that effectively, the lifetime-exchange only swaps  $\vnu_{s+2/n-u_1/n}$ and $\vnu_{s+2/n-u_2/n}$ along with the neighbors of $\bm h_{u_1}$ and $\bm h_{u_2}$,  but otherwise yields the same graph as the modified graph exploration with the original life-times.

 Denote by $\bm{H}_r:=\bm{H}_r(\bm{\eta})$ the set of all half-edges of $\vnu_{s+2/n-r/n},\ldots,\vnu_{s+1/n}$ according to the original and modified graph explorations. Observe that also the first dead half-edge of $\vnu_{s+2/n-u/n}$, $u \in [r]$, does not depend on the specific exploration, so that we again use identical notation $\bm h_u$ in both cases.

Next, we 
identify an appropriate event $\fO_{s,r}$ which ensures that the original and modified graph exploration are identical: 
Define $\fO_{s,r}$ ($\fM_{s,r}$) to be the set of all decision vectors $\eta$ such that, for the original (modified) graph exploration it holds true that:
\begin{enumerate}
    \item\label{set_a_1} All vertices $\nu_{s+2/n-r/n},\ldots,\nu_{s+1/n}$ are activated by \ref{exploration_step3};

    \item\label{set_a_2} At the time when the awakening step of $\nu_{s+1/n}$ finishes, exactly one of its half-edges has been killed for each vertex $\nu_{s+2/n-u/n} (\mbox{for }u\in[r])$;
\end{enumerate}
By definition, the lifetime of $h_1$,  which is  the first dead
half-edge of $\nu_{s+1/n}$, is the maximum of the lifetimes of all half-edges in $\cbc{h_u:u\in [r]}$. On the other hand, item \ref{set_a_2} in the definition of $\fO_{s,r}$ ($\fM_{s,r}$) implies that the lifetime of $h_1$ is smaller than the minimum of the lifetimes of the half-edges in $H_r\backslash\cbc{h_u:u\in [r]}$.

\begin{claim}\label{claim_1}
    On $\fO_{s,r}$, both the original and the modified explorations perform exactly the same. 
In particular, $\fO_{s,r}\subseteq \fM_{s,r}$. 
\end{claim}
\begin{proof}
    By definition, both explorations perform exactly the same until the $\lfloor ns+2-r \rfloor$th vertex has been awakened. For $\eta\in \fO_{s,r}$, this vertex is awakened by \ref{exploration_step3}. In the time interval until vertex $\nu_{s+1/n}$ is awakened, \ref{exploration_step21} and \ref{exploration_step22} only differ if there are both active half-edges belonging to the first $\lfloor ns+1-r\rfloor$ awakened vertices and active half-edges belonging to the already awakened vertices among $\nu_{s+2/n-r/n}, \ldots, \nu_{s}$, and in \ref{exploration_step21}, one of the latter is killed. However, in that case, the adjacent vertex would  have at least two killed half-edges at the time when the awakening step of $\nu_{s+1/n}$ finishes, which cannot happen on $\fO_{s,r}$. Thus, the original and modified graph explorations agree until vertex $\nu_{s+1/n}$ is awakened. As they also  perform the same afterwards, the claim follows. 
\end{proof}
We will now show that $\fO_{s,r}$ happens with high probability. Thus, with decent probability, the original and modified graph explorations are the same, and we can transfer any exchangeability properties from the modified graph exploration back to the original one.

\begin{claim}\label{claim-app_as}
For any integer $r>0$,
\begin{align}\label{eq_app_as}
\mathbb{P}(\bm{\eta}\in \fO_{s,r})=1+\bar{o}_n(1).
\end{align}
\end{claim}

\begin{proof}
By \Cref{l_concentration} and 
\Cref{ineq_dcs}, 
\begin{align}\label{eq_eventa_1}
   \mathbb{P}\bc{\mbox{$\vnu_{s+2/n-r/n},\ldots,\vnu_{s+1/n}$ are activated by \ref{exploration_step31}},\mathfrak{G}_{\eps,s+2/n}}=1+\bar{o}_n(1).
\end{align}
It thus remains to show that 
\begin{align*}
    \mathbb{P}\bc{\exists\ u\  \in[r]\mbox{ s.t. at least two half-edges of $\vnu_{s+2/n-u/n}$ are killed when $\vnu_{s+1/n}$ is awakened}}=\bar{o}_n(1).
\end{align*}
 For this, we make the following two observations:
\begin{enumerate}
    \item\label{case_step2_kill}  In any iteration of \ref{exploration_step21} after the awakening of $\nu_{s+2/n-r/n}$, an active half-edge is chosen uniformly at random among all currently \emph{active} half-edges. As there is at least one active half-edge each time \ref{exploration_step21} is called, at any such moment between the awakening of $\vnu_{s+2/n-r/n}$ and $\vnu_{s+1/n}$, the number of active half-edges is lower bounded by $\max\cbc{1,\edgeli{s+2/n}-\edgesl{s+2/n-r/n}}$. 
    As a consequence, for each iteration $i$ of \ref{exploration_step21} after the awakening of $\nu_{s+2n-r/n}$, by \Cref{l_concentration,lem_app_divide},
\begin{align}\label{pro_halfedge_kill1}
        \PP&\bc{\exists \mbox{$h \in \bm{H}_r\backslash\cbc{h_u:u\in[r]}: h$ is killed in the $i$th iteration of \ref{exploration_step21} after the awakening of $\vnu_{s+2/n-r/n}$}}\nonumber\\
        \leq &\Erw\brk{\frac{rK}{\max\cbc{1,\edgeli{s+2/n}-\edgesl{s+2/n-r/n}}}}=\bar{o}_n(1).
    \end{align}

    \item In any iteration of \ref{exploration_step31} after the awakening of $\nu_{s+2/n-r/n}$, an active half-edge is chosen uniformly at random among all currently \emph{living} half-edges.  
    Analogously to the previous observation, for each iteration, by \Cref{lem_app_divide},
    \begin{align}\label{pro_halfedge_kill2}
    \PP&\bc{\exists \mbox{$h \in \bm{H}_r\backslash\cbc{h_u:u\in[r]}: h$ is killed in the $i$th iteration of \ref{exploration_step31} after the awakening of $\vnu_{s+2/n-r/n}$}}\nonumber\\
        \leq &\Erw\brk{\frac{rK}{\edgeli{s+2/n}}}=\bar{o}_n(1).
    \end{align}
\end{enumerate}
Now, fix a number $N>0$ and let $\bm B = B(\bm \eta)$ be the number of iterations of \ref{exploration_step21} between the awakening of $\vnu_{s+2/n-r/n}$ and $\vnu_{s+1/n}$. Then summing \cref{pro_halfedge_kill1} and \cref{pro_halfedge_kill2} over $i\in[Nr]$, by a union bound,
\begin{align}\label{ineq_twohalfedge}
    \PP&\bc{\exists\ u\  \in[r],\mbox{ s.t. at least two half-edges of $\vnu_{s+2/n-u/n}$ have been killed when $\vnu_{s+1/n}$ is awakened}}\nonumber\\
    \leq &\PP(\bm B\geq Nr)+\bar{o}_n(1).
\end{align}
On the event that $\vnu_{s+2/n-r/n},\ldots,\vnu_{s+1/n}$ are activated by \ref{exploration_step3} and $\mathfrak{G}_{\eps,s+2/n}$ holds, $\bm B\geq Nr$ implies that between the awakening of $\vnu_{s+2/n-r/n}$ and that of $\vnu_{s+1/n}$, at least $Nr$ half-edges were killed in iterations of \ref{exploration_step21}. On the other hand, the number of living half-edges cannot exceed $Kn$, while on $\mathfrak{G}_{\eps,s+2/n}$, the number of sleeping half-edges in the given time frame is always bounded below by $cn$. Consequently, each time a half-edge is killed in \ref{exploration_step21}, with probability at least $c/K$ it is a sleeping half-edge and a new vertex is awakened. By Chebyshev's inequality,
\begin{align}\label{bound_bnr}
    \PP\Big(\bm B\geq Nr,\mbox{ $\vnu_{s+2/n-r/n},\ldots,\vnu_{s+1/n}$ are activated by \ref{exploration_step3}},\mathfrak{G}_{\eps,s+2/n}\Big)
    \leq &\PP\bc{\bin{Nr,c/K}\leq r}\\
    \leq& Nr/(cN-K)^2,\nonumber
\end{align}
where $\bin{n,p}$ denotes a binomial random variable with $n$ trials and success probability $p$. 
By \eqref{eq_eventa_1}, \cref{ineq_twohalfedge} and \cref{bound_bnr}, taking $N\to\infty$,
\begin{align}\label{eq_eventa_2}
  \PP&\bc{\exists\ u\  \in[r],\mbox{ s.t. at least two half-edges of $\vnu_{s+2/n-u/n}$ are killed when $\vnu_{s+1/n}$ is awakened}}\nonumber\\
    \leq& \limsup_{N\to\infty}\PP(\bm B\geq Nr)+\bar{o}_n(1)=\bar{o}_n(1).
\end{align}
Combining \eqref{eq_eventa_1} and \eqref{eq_eventa_2} 
gives \cref{eq_app_as}.
\end{proof}

\subsection{Switching in the modified graph exploration} \label{sec_exchange}

In this section, we define a switching operation on the $r$ first killed half-edges of $\vnu_{s+2/n-r/n}, \ldots, \vnu_{s+1/n}$ that works well with the modified graph exploration. 
Fix a permutation $\rho$ on $[r]$. Then for each decision vector $\eta$, we define a new decision vector $\varrho(\eta)$ by replacing the lifetime of half-edge $h_{u}$ by $\bm{E}_{h_{\rho(u)}}$, for $u\in[r]$. 
Finally, as we are interested in degrees, let $(w_u)_{u\in[r]}$ be a sequence of integers in $[K]$ and $\fM_{s,(w_u)_{u\in[r]}} \subseteq\fM_{s,r}$ be the subset of all decision vectors such that the original degree of $\vnu_{s+2/n-u/n}$ is equal to $w_u$ for all $u\in[r]$. Let \mcma{n} be the graph obtained through the modified graph exploration procedure. The following claim then clarifies the effect of this operation on the modified graph exploration and the degrees of $\vnu_{s+2/n-r/n}, \ldots, \vnu_{s+1/n}$:
\begin{lemma}\label{claim_2}
For each permutation $\rho$ on $[r]$ and $\eta\in \fM_{s,(w_u)_{u\in[r]}}$,
\begin{enumerate}
    \item  $\varrho(\eta)\in \fM_{s,(w_{\rho(u)})_{u\in[r]}}$;
    \item\label{claim2_it3} $\fmcma{n}(\varrho(\eta))$ agrees with $\fcm{n}(\eta)$, except that for $u \in [r]$, the neighbor of $h_u$ in  $\fmcma{n}(\eta)$ becomes the neighbor of $h_{\rho^{-1}(u)}$ in $\fmcma{n}(\varrho(\eta))$. 
\end{enumerate}
\end{lemma}
As we will show in the following proof, for $\eta\in \fM_{s,(w_u)_{u\in[r]}}$, the same half-edges are killed in the modified graph exploration with respect to $\eta$ and with respect to $\varrho(\eta)$,  except that when $h_u\in H_r$ is killed in the  modified graph exploration with respect to $\eta$, $h_{\rho(u)}\in H_r$ is killed in the  modified graph exploration with respect to $\varrho(\eta)$.
\begin{proof} 
In the following, let `step' refer to \ref{exploration_step11}, \ref{exploration_step22} or \ref{exploration_step31} in the modified graph exploration, where \ref{exploration_step11} is also counted there are active half-edges remaining. We note that the random variables $\bm{E}_h$ and $\bm{U}(h,t,v)$ in the decision set $\bm{\eta}$ are distinct almost surely. 
Therefore, almost surely, the trajectory of the set of \emph{lifetimes} of living half-edges with respect to the number of taken steps is in bijection with the evolution of the set of living half-edges and therefore the pairing order of the whole graph exploration. The main ingredient in the proof of  \Cref{claim_2} is thus the following claim:

\begin{claim}\label{claim_traj_living}
    For $\eta \in \fM_{s,(w_u)_{u\in[r]}}$ such that all its components are distinct, at any step $\ell$, the set of lifetimes of living half-edges is the same for the modified graph explorations determined by $\eta$ and $\varrho(\eta)$.
\end{claim}

\begin{proof}[Proof of \Cref{claim_traj_living}]
   We use induction on the number of steps $\ell \geq 0$ that have already been taken.
   
For $\ell =0$, no step has been taken yet, and the set of lifetimes of living half-edges is $\cbc{E_h:h\in\mathcal{H}}$ in both modified graph explorations.

To advance the induction, assume that the claim is true when at most $(\ell-1)$ steps have been taken, so that the sets of lifetimes of living half-edges after the $(\ell-1)$st step are the same in both explorations. 
If the $\ell$th step is an iteration of \ref{exploration_step11}, the sets of lifetimes of living half-edges do not change from step $\ell-1$ to step $\ell$, since no half-edges are killed in \ref{exploration_step11}. If the $\ell$th step is an iteration of  \ref{exploration_step31}, then in both explorations, the half-edge with the smallest lifetime among all living half-edges is killed and removed from the sets of lifetimes of living half-edges.

We now consider the case that the $\ell$th step is an iteration of \ref{exploration_step22} and emphasize that $h_u$, $H_r$ and $\nu_{t}$ always refer to the modified graph exploration with respect to $\eta$. 
Given the set of lifetimes of living half-edges after the $(\ell-2)$th step, which has been an iteration of \ref{exploration_step31}, the following four sets are determined for both $\eta$ and $\varrho(\eta)$:
\begin{itemize}
    \item The sets of \emph{living} half-edges: This is due to the fact that the lifetimes stand in (known) bijection to the half-edges (even though the bijections may be different for $\eta$ and $\varrho(\eta)$);
    \item The sets of \emph{dead} half-edges as the complements of the sets of living half-edges with respect to $\mathcal{H}$;
    \item The sets of \emph{sleeping} half-edges: This is due to the fact that  after each iteration of \ref{exploration_step31}, the adjacent half-edges of vertices with only living half-edges can be unambiguously identified as sleeping.
    \item The sets of \emph{active} half-edges as the complements of the sets of sleeping half-edges with respect to the set of living half-edges.
\end{itemize}
Given this information, we in particular know the sets of awake vertices and next distinguish the cases whether we are in an epoch between the awakening of $\nu_{s+2/n-r/n}$ and $\nu_{s+1/n}$ or not.
\begin{enumerate}
    \item Suppose that in the modified graph exploration for $\eta$, after step $\ell-2$, there  are at most $\lfloor ns+1-r\rfloor$ or at least $\lfloor ns +1\rfloor$ awakened vertices. Then up to this stage, none or all of the half-edges in $\cbc{h_u:u\in [r]}$ have been killed in the modified graph exploration w.r.t. $\eta$ and thus all or none of their lifetimes are in the set of living half-edges after the $(\ell-2)$th step. Since the bijections between a half-edge and its lifetime agree on $\mathcal{H}\backslash \cbc{h_u:u\in [r]}$ for both modified graph explorations, we conclude from the induction hypothesis that the sets of living half-edges are the same after the $(\ell-2)$th step, and thus that all four sets above are the same. 
    
    Then in the $(\ell-1)$st step, which is an iteration of \ref{exploration_step11}, if the set of active half-edges was non-empty in the $(\ell-2)$th step, no vertex is awakened, in which case it is immediately clear that all four sets above stay the same. On the other hand,  if the set of active half-edges was non-empty in the $(\ell-2)$th step, a sleeping half-edge is activated according to the value of its decision variable, all of which agree in $\eta$ and $\rho(\eta)$. As $\eta\in \fM_{s,r}$, the corresponding awakened vertex is not one of $\nu_{s+2/n-r/n},\ldots,\nu_{s+1/n}$ in the modified graph exploration w.r.t $\eta$, and thus the lifetimes of sleeping and active vertices change in the same way. Again, we conclude that after the $(\ell-1)$st step, all four sets above agree for $\eta$ and $\rho(\eta)$.

As the sets of active half-edges agree after the $(\ell-1)$st step and the choice of the killed half-edge in step $\ell$ is based on identical decision variables, the same half-edge is killed in both modified graph explorations. Moreover, as $\eta\in \fM_{s,r}$, the killed half-edge is not in $\cbc{h_u:u\in[r]}$, so that it has the same lifetime in $\eta$ and $\rho(\eta)$. 
We conclude that the sets of lifetimes of living half-edges are the same after the $\ell$th step.

\item Suppose that in the modified graph exploration for $\eta$, after step $\ell-2$, there  are at least $\lfloor ns+2-r\rfloor$ and at most $\lfloor ns\rfloor$ awakened vertices. 
As execution of \ref{exploration_step11} would entail awakening of a new vertex,  by \cref{set_a_1} in the definition of $\fM_{s,r}$, nothing happens in the $(\ell-1)$st step in the modified graph exploration for $\eta$. Moreover, by the induction hypothesis, after step $\ell-2$, the sets of lifetimes of living, and therefore also those of dead, half-edges agree in the modified explorations w.r.t. $\eta$ and $\rho(\eta)$. From this we conclude that every half-edge not in $\cbc{h_u:u\in [r]}$ has the same status \emph{active}, \emph{sleeping} or \emph{dead} in both explorations after step $\ell-2$. Therefore, if the set of active half-edges w.r.t. $\eta$ is non-empty, so is the one w.r.t. $\rho(\eta)$, and nothing happens in the $(\ell-1)$st step in the modified graph exploration for $\varrho(\eta)$ either.

Moreover, as the set of active half-edges \textbf{not in $\cbc{h_u:u\in [r]}$} in the $\ell$th step in both explorations is the same and non-empty, in step $\ell$, the same active half-edge is chosen according to the value of its decision variable and killed in both explorations. 
As this half-edge is not in $\cbc{h_u:u\in [r]}$, the same value is removed from both sets of lifetimes of living half-edges. 
\end{enumerate}
As a consequence, the sets of lifetimes of living half-edges are the same after the completion of the $\ell$th step. 
\end{proof}

From the proof we see that at each step, every half-edge not in $\cbc{h_u:u\in [r]}$ has the same status \emph{active}, \emph{sleeping} or \emph{dead} in both explorations. Therefore, both explorations perform the same, except when the vertices $\nu_{s+2/n-r/n}, \ldots, \nu_{s+1/n}$ are awakened by \ref{exploration_step31}. At these times, $h_u$ is paired in the exploration w.r.t. $\eta$, while $h_{\rho^{-1} (u)}$ is paired in the exploration w.r.t. $\rho(\eta)$. This permutes the order in which the $\lfloor ns +2-r\rfloor$th, $\ldots, \lfloor ns+1\rfloor$th vertices are awakened and therefore their degrees. 
It is then very easy to see that $\varrho(\eta)\in \fM_{s,(w_{\rho(u)})_{u\in[r]}}$, and $\fmcma{n}(\varrho(\eta))$ agrees with $\fmcma{n}(\eta)$, except that for $u \in [r]$, the neighbor of $h_u$ in  $\fmcma{n}(\eta)$ becomes the neighbor of $h_{\rho(u)}$ in $\fmcma{n}(\varrho(\eta))$. 
\end{proof}

\subsection{Proof of \Cref{lem_vcm}} \label{sec_proof_vcm}
\begin{proof}[Proof of \Cref{lem_vcm}]
In the following, denote by \mcm{n}{s+1/n} the induced graph on the last $n-\lfloor ns+1\rfloor$ sleeping vertices in the graph that results from the modified graph exploration. 
Note that for each permutation $\rho$ of $[r]$, $\rho^{-1}$ is a permutation of $[r]$ as well. Hence, $\varrho(\eta)\in \fM_{s,(w_{\rho(u)})_{u\in[r]}}$ only if $\eta=\varrho^{-1}(\varrho(\eta))\in \fM_{s,(w_u)_{u\in[r]}}$ by \Cref{claim_2}, and thus $\varrho:\fM_{s,(w_u)_{u\in[r]}}\to \fM_{s,(w_{\rho(u)})_{u\in[r]}}$ is a bijection.
Moreover,  \cref{claim2_it3} in \Cref{claim_2} gives that $\fmcm{n}{s+1/n}(\eta) = \fmcm{n}{s+1/n}(\rho(\eta))$.
Hence,
\begin{align*}
    &\eta\in \fM_{s,(w_u)_{u\in[r]}}\cap\cbc{\gamma:\fmcm{n}{s+1/n}(\gamma)=G} \quad
    \Longleftrightarrow \quad\varrho(\eta)\in \fM_{s,(w_{\rho(u)})_{u\in[r]}}\cap\cbc{\gamma:\fmcm{n}{s+1/n}(\gamma)=G}.
\end{align*}
On the other hand, since we merely exchange the lifetimes of the half-edges in $\cbc{h_u}_{u\in[r]}$, the probability density function of both $\bm{\eta}$ and $\varrho(\bm{\eta})$ under the condition $\bm{E}_h = e_h\ (h\in\mathcal{H})$ equals to $\prod_{h\in\mathcal{H}} \eul^{-e_h}$.
Hence,
\begin{align}\label{eq_ex_eta}
    \frac{\dif\PP( \bm{\eta})}{\dif\PP( \varrho(\bm{\eta}))}=1.
\end{align}
As a consequence, for any graph $G$, by \Cref{eq_ex_eta}, 
  \begin{align*}
     \PP\bc{\bm{\eta}\in \fM_{s,(w_u)_{u\in[r]}},\fmcm{n}{s+1/n}(\bm \eta)=G}
     &=\int_{\bm{\eta}\in \fM_{s,(w_u)_{u\in[r]}}\cap\cbc{\vga:\fmcm{n}{s+1/n}(\vga)=G}}\dif\PP( \bm{\eta})\\
     &=\int_{\bm{\eta}\in \fM_{s,(w_u)_{u\in[r]}}\cap\cbc{\vga:\fmcm{n}{s+1/n}(\vga)=G}}\dif\PP( \varrho(\bm{\eta}))\\
&=\int_{\varrho(\bm{\eta})\in \fM_{s,(w_{\rho(u)})_{u\in[r]}}\cap\cbc{\vga:\fmcm{n}{s+1/n}(\vga)=G}}\dif\PP( \varrho(\bm{\eta}))\\
     &=\mathbb{P}\bc{\bm{\eta}\in \fM_{s,(w_{\rho(u)})_{u\in[r]}},\fmcm{n}{s+1/n}(\bm \eta)=G}.
 \end{align*}
Hence,
\begin{align}\label{eq_proaper}
\mathbb{P}\bc{\bm{\eta}\in \fM_{s,(w_u)_{u\in[r]}} \mid \fmcm{n}{s+1/n}}=\mathbb{P}\bc{\bm{\eta}\in \fM_{s,(w_{\rho(u)})_{u\in[r]}} \mid \fmcm{n}{s+1/n}}.
\end{align}
Let $(m_k)_{k\in[K]}$ be a sequence of non-negative integers such that $\sum_{k=1}^Km_k=r$. Let $\hat{\fM}_{s,(m_k)_{k\in[K]}}$ denote the disjoint union of all  $\fM_{s,(w_u)_{u\in[r]}}$  such that the number of $u\in[r]$ with $w_u=k$ is equal to $m_k$ for all $k\in [K]$. There are $\frac{r!}{\prod_{i\in[K]} m_i!}$ such disjoint $\fM_{s,(w_u)_{u\in[r]}}$ that belong to each $\hat{\fM}_{s,(m_k)_{k\in[K]}}$ and  $\mathbb{P}(\bm{\eta}\in \fM_{s,(w_u)_{u\in[r]}} \mid \fmcm{n}{s+1/n})$ is the same for all such $(w_u)_{u\in[r]}$ by \cref{eq_proaper}. As a consequence, 
\begin{align}\label{relation_seta_setb}
    \mathbb{P}\bc{\bm{\eta}\in \hat{\fM}_{s,(m_k)_{k\in[K]}} \mid \fmcm{n}{s+1/n}}=\frac{r!}{\prod_{i\in[K]} m_i!}\mathbb{P}\bc{\bm{\eta}\in \fM_{s,(w_u)_{u\in[r]}} \mid \fmcm{n}{s+1/n}},
\end{align}
where on the right-hand side the number of $u\in[r]$ such that  $w_u=k$ is equal to $m_k$ for all $k\in [K]$.

On the other hand,
\begin{align*}
\mathbb{P}\bc{\bm{\eta}\in \fM_{s,r},d_{\vnu_{s+1/n}}=\ell \mid \fmcm{n}{s+1/n}}=&\sum_{(w_u)_{u\in[r]}\in [K]^r,w_1=\ell}\mathbb{P}\bc{\bm{\eta}\in \fM_{s,(w_u)_{u\in[r]}} \mid \fmcm{n}{s+1/n}}.
\end{align*}
In the summation on the right-hand side, the number of $(w_u)_{u\in[r]}$ such that the number of $u\in[r]$ with $w_u=k$ is equal to $m_k$ for all $k\in [K]$ equals $\frac{(r-1)!}{(m_\ell-1)!\prod_{i\in[K]\backslash \cbc{\ell}} m_i!}$. Hence, by \Cref{relation_seta_setb},
\begin{align}
&\mathbb{P}\bc{\bm{\eta}\in\fM_{s,r},d_{\vnu_{s+1/n}}=\ell \mid \fmcm{n}{s+1/n}} =\sum_{\substack{(m_k)_{k\in[K]}:\\ \sum_{k=1}^Km_k=r}}\mathbb{P}\bc{\bm{\eta}\in \hat{\fM}_{s,(m_k)_{k\in[K]}} \mid \fmcm{n}{s+1/n}}\frac{m_\ell}{r}.\label{eq_cond_prob_sum_hatm}
\end{align}
Let $\vmm_\ell$ and $\vmm_{\ell,m}$ denote the number of vertices in $\{\nu_{s+2/n-r/n}, \ldots, \nu_{s+1/n}\}$ that have 
original degree $\ell$ in the original and modified graph explorations, respectively. Then \eqref{eq_cond_prob_sum_hatm} yields that
\begin{align}\label{eq_core_claim}
\mathbb{P}\bc{\bm{\eta}\in\fM_{s,r},d_{\vnu_{s+1/n}}=\ell \mid \fmcm{n}{s+1/n}}
=&\Erw\brk{\frac{\vmm_{\ell,m}}{r}\teo{\bm{\eta}\in\fM_{s,r}} \mid \fmcm{n}{s+1/n}}.
\end{align}
The l.h.s. of \Cref{eq_core_claim} is already close to that of \Cref{eq-lem_vcm-new}, except that the graph is modified and we are considering $\bm{\eta} \in\fM_{s,r}$. These two issues can be resolved if $\bm{\eta} \in\fO_{s,r}$ w.h.p., by applying \Cref{claim_1}. To prove \Cref{eq-lem_vcm-new}, we also need to express the r.h.s. of \Cref{eq_core_claim} in a suitable form. For this purpose, we formulate the following claim:
\begin{claim}\label{claim-eq_app_gl}
For any $\ell\in \cbc{0}\cup[K]$, uniformly in $\ranges$,
\begin{align}\label{eq_app_gl}
    \Erw\abs{\frac{\vmm_\ell}{r}-\ell p_\ell \eul^{-t_s \ell}\lambda(t_s)^{-1}}=\bar{o}_{r}(1)+\bar{o}_n(1),
\end{align}
where $a_r=\bar{o}_{r}(1)$  means that $(\abs{\bm{a}_r})_{r \in \NN, \ranges}$ is uniformly bounded and $\lim_{r\to\infty} a_r=0$.
\end{claim}
The proof of Claim~\ref{claim-eq_app_gl} is deferred to the end of this section. 
Now, recall that \Cref{claim_1} states that, on $\fO_{s,r}$, both explorations perform exactly the same. Hence,
\begin{align}\label{eq_c_claim1}
\Erw\brk{\frac{\vmm_\ell}{r}\teo{\bm{\eta}\in\fO_{s,r}} \mid \fscm{n}{s+1/n}}=\Erw\brk{\frac{\vmm_{\ell,m}}{r}\teo{\bm{\eta}\in\fO_{s,r}} \mid \fmcm{n}{s+1/n}}.
\end{align}
Since $\fO_{s,r}\subseteq\fM_{s,r}$ by  \Cref{claim_1} and $\vmm_\ell,\vmm_{\ell,m}\leq r$,   \Cref{eq_c_claim1} yields that
\begin{align}\label{eq_exp_gl2}
&\Erw\abs{\Erw\brk{\frac{\vmm_\ell}{r} \mid \fscm{n}{s+1/n}}-\Erw\brk{\frac{\vmm_{\ell,m}}{r}\teo{\bm{\eta}\in\fM_{s,r}} \mid \fmcm{n}{s+1/n}}}\nonumber\\
 &\qquad\leq\Erw\abs{\Erw\brk{\frac{\vmm_\ell}{r}\teo{\bm{\eta}\in\fO_{s,r}} \mid \fscm{n}{s+1/n}}-\Erw\brk{\frac{\vmm_{\ell,m}}{r}\teo{\bm{\eta}\in\fO_{s,r}} \mid \fmcm{n}{s+1/n}}}\nonumber\\
&\qquad\quad+1-\PP\bc{\bm{\eta}\in\fO_{s,r}}+\PP\bc{\bm{\eta}\in\fM_{s,r}}-\PP\bc{\bm{\eta}\in\fO_{s,r}}\nonumber\\
 &\qquad=1+\PP\bc{\bm{\eta}\in\fM_{s,r}}-2\PP\bc{\bm{\eta}\in\fO_{s,r}}\leq 2-2\PP\bc{\bm{\eta}\in\fO_{s,r}}=\bar{o}_n(1).
\end{align}
Further, \cref{eq_app_as}, along with $\fO_{s,r}\subseteq\fM_{s,r}$, gives that $\mathbb{P}(\bm{\eta}\in\fM_{s,r})=1+\bar{o}_n(1)$. Hence,
\begin{align}\label{eq_conp_dnus}
&\Erw\abs{\mathbb{P}\bc{d_{\vnu_{s+1/n}}=\ell \mid \fscm{n}{s+1/n}}-\mathbb{P}\bc{\bm{\eta}\in\fM_{s,r},d_{\vnu_{s+1/n}}=\ell \mid \fscm{n}{s+1/n}}}\leq \mathbb{P}(\bm{\eta}\not\in\fM_{s,r})=\bar{o}_n(1).
\end{align}
Consequently, by \cref{eq_core_claim}, \cref{eq_app_gl}, \cref{eq_exp_gl2} and \cref{eq_conp_dnus},
\begin{align*}
\mathbb{E}\abs{\mathbb{P}\bc{d_{\vnu_{s+1/n}}=\ell \mid \fscm{n}{s+1/n}}-\ell p_\ell \eul^{-t_s \ell}\lambda(t_s)^{-1}}=\bar{o}_n(1)+\bar{o}_{r}(1).
\end{align*}
Since $\mathbb{E}\abs{\mathbb{P}\bc{d_{\vnu_{s+1/n}}=\ell \mid \fscm{n}{s+1/n}}-\ell p_\ell \eul^{-t_s \ell}\lambda(t_s)^{-1}}$ does not depend on $r$, taking the limit $\lim_{r\to \infty}$ on both sides of the above equation, we conclude that the left-hand side is equal to $\bar{o}_n(1)$, as desired.
\end{proof}

Hence, we are left to prove \cref{eq_app_gl}.

\begin{proof}[Proof of \Cref{claim-eq_app_gl}]
We prove \Cref{claim-eq_app_gl} by using a second-moment method. 
For any $b\in [s+2/n-r/n,s+1/n]$,  by  \Cref{l_concentration} and the fact that $s \mapsto t_s$ is continuous, 
\begin{align}\label{eq_app_db}
\mathbb{P}\bc{d_{\vnu_{b}}=\ell}=&\Erw\brk{\mathbb{P}\bc{d_{\vnu_{b}}=\ell \mid  \verslk{0}{b-1/n}, \ldots, \verslk{K}{b-1/n}}}=\Erw\brk{\frac{\ell\verslk{\ell}{b-1/n}}{\sum_{k\in [K]}k \verslk{k}{b-1/n}}}\nonumber\\
=&\ell p_\ell \eul^{-t_s \ell}\lambda(t_s)^{-1}+\bar{o}_n(1),
\end{align}
where in the second equation we use the fact that a sleeping vertex in the graph exploration is awakened with probability proportional to its original degree.

Furthermore, for any $b_1,b_2\in [s+2/n-r/n,s+1/n]$ such that $b_1\leq b_2-1/n$, by \Cref{l_concentration} and the tower property, 
\begin{align*}
&\PP\bc{d_{\vnu_{b_1}}=d_{\vnu_{b_2}}=\ell \mid \verslk{k}{b_1-1/n},k\in[K]}\\
=&\PP\bc{d_{\vnu_{b_1}}=\ell \mid \verslk{k}{b_1-1/n},k\in[K]}\PP\bc{d_{\vnu_{b_2}}=\ell \mid d_{\vnu_{b_1}}=\ell,\verslk{k}{b_1-1/n},k\in[K]}\\
=&\frac{\ell\verslk{\ell}{b_1-1/n}}{\sum_{k\in [K]}k \verslk{k}{b_1-1/n}}\Erw\brk{\mathbb{P}\bc{d_{\vnu_{b_2}}=\ell \mid \verslk{k}{b_i-1/n},i\in [2],k\in[K],d_{\vnu_{b_1}}=\ell} \mid \verslk{k}{b_1-1/n},k\in[K],d_{\vnu_{b_1}}=\ell} \\
=&\mathbb{E}\brk{\frac{\ell\verslk{\ell}{b_1-1/n}}{\sum_{k\in [K]}k \verslk{k}{b_1-1/n}}\frac{\ell\verslk{\ell}{b_2-1/n}}{\sum_{k\in [K]}k \verslk{k}{b_2-1/n}} \mid \verslk{k}{b_1-1/n},k\in[K],d_{\vnu_{b_1}}=\ell}.
\end{align*}
Therefore, by \Cref{l_concentration},
\begin{align}\label{eq_app_db12} \mathbb{P}\bc{d_{\vnu_{b_1}}=d_{\vnu_{b_2}}=\ell}=\Erw\brk{\frac{\ell\verslk{\ell}{b_1-1/n}}{\sum_{k\in [K]}k \verslk{k}{b_1-1/n}}\frac{\ell\verslk{\ell}{b_2-1/n}}{\sum_{k\in [K]}k \verslk{k}{b_2-1/n}}}
=(\ell p_\ell \eul^{-t_s \ell}\lambda(t_s)^{-1})^2+\bar{o}_n(1).
\end{align}
On the other hand, the definition of $\vmm_\ell$ yields that  $\vmm_\ell=\sum_{u\in[r]}\teo{ d_{\vnu_{s+2/n-u/n}}=\ell}$ . Thus,
\begin{align}\label{eq_app_difglr}
    &\bc{\Erw\abs{\frac{\vmm_\ell}{r}-\ell p_\ell \eul^{-t_s \ell}\lambda(t_s)^{-1}}}^2\leq \Erw\abs{\frac{\vmm_\ell}{r}-\ell p_\ell \eul^{-t_s \ell}\lambda(t_s)^{-1}}^2\\
    =&\frac{1}{r^2}\Big(\sum_{u\in[r]}\mathbb{P}\bc{ d_{\vnu_{s+2/n-u/n}}=\ell}+2\sum_{u\in[r]}\sum_{v\in [u-1]}\mathbb{P}\bc{ d_{\vnu_{s+2/n-u/n}}=d_{\vnu_{s+2/n-v/n}}=\ell}\nonumber\\
    &-2r\ell p_\ell \eul^{-t_s \ell}\lambda(t_s)^{-1}\sum_{u\in[r]}\mathbb{P}\bc{ d_{\vnu_{s+2/n-u/n}}=\ell}+\bc{r\ell p_\ell \eul^{-t_s \ell}\lambda(t_s)^{-1}}^2\Big).\nonumber
\end{align}
The equality in \cref{eq_app_gl} then follows directly from the combination of \cref{eq_app_db}, \cref{eq_app_db12} and \cref{eq_app_difglr}.
\end{proof}

\subsection{Proof of \Cref{lem_vscm}} \label{sec_proof_vscm}
Recall the decomposition \cref{eq_degree_LTP}. While \Cref{lem_vcm} yields the desired estimate for the second factor $\PP\bc{  d_{\vnu_{s+1/n}}=\ell \mid \fscm{n}{s+1/n}}$ on the right hand side, we now derive the asymptotics of the first factor $\PP\bc{\bar{\vd}_{\vnu_{s+1/n},s}=k \mid \fscm{n}{s+1/n}, d_{\vnu_{s+1/n}}=\ell}$. 
Let $\bcabc{s+1/n}$ be the number of half-edges of vertices in $\verst{s+1/n}$ that are joined to a half-edge outside of \scm{n}{s+1/n}. Then
    \begin{align*}
        \mathbb{E}\abs{n^{-1}\bcabc{s+1/n}-\lambda(t_s)+\frac{\lambda^2(t_s)\eul^{2t_s}}{\lambda(0)}}=\bar{o}_n(1) \quad\mbox{(for a proof, see  Appendix \ref{ap_eventg}).}
    \end{align*}
Next, recall \Cref{notation_reduced} and the event $\mathfrak G_{\eps, s+1/n}$ from \cref{event_gs}. On the event $\mathfrak G_{\eps, s+1/n}$, for sufficiently large $n$, vertex $\vnu_{s+1/n}$ is awakened via \ref{exploration_step3}. 
At the time when vertex $\vnu_{s+1/n}$ 
is awakened, the number of active half-edges $\edgeli{s+1/n}-\edgesl{s+1/n}$ is lower bounded by $\bcabc{s+1/n}$, as all half-edges of vertices in $\verst{s+1/n}$ that are joined to a half-edge outside of \scm{n}{s+1/n} will connect to active half-edges. Given $\fscm{n}{s+1/n}$ such that $\mathfrak G_{\eps, s+1/n}$ holds, $\edgeli{s+1/n}, \edgesl{s+1/n}$ and $d_{\vnu_{s+1/n}}=\ell$, the neighbors of the $\bcabc{s+1/n}$ half-edges are chosen uniformly at random among the active half-edges.  
As $\vnu_{s+1/n}$ has $\ell-1$ active half-edges when it is awakened, for $k\leq \ell-1$,
\begin{align*}
    &\PP\bc{\bar{\vd}_{\vnu_{s+1/n},s}=k \mid \fscm{n}{s+1/n},d_{\vnu_{s+1/n}}=\ell,\edgeli{s+1/n},\edgesl{s+1/n} }\\
    =&\frac{\binom{\edgeli{s+1/n}-\edgesl{s+1/n}-\ell+1}{\bcabc{s+1/n}-k}\binom{\ell-1}{k}}{\binom{\edgeli{s+1/n}-\edgesl{s+1/n}}{\bcabc{s+1/n}}}
    =\frac{\binom{\ell-1}{k}\prod_{i=0}^{k-1}(\bcabc{s+1/n}-i)\prod_{i'=0}^{\ell-2-k}{(\edgeli{s+1/n}-\edgesl{s+1/n}-\bcabc{s+1/n}-i')}}{\prod_{i''=0}^{\ell-2}(\edgeli{s+1/n}-\edgesl{s+1/n}-i'')}.
\end{align*}
Define $r_{k,\ell}:=\binom{\ell-1}{k}\bc{\frac{\lambda(t_s)\eul^{2t_s}}{\lambda(0)}}^k\bc{1-\frac{\lambda(t_s)\eul^{2t_s}}{\lambda(0)}}^{\ell-k-1}$. By \Cref{l_concentration} and \Cref{lem_app_divide}, 
\begin{align*}
\Erw \brk{\mathds{1}\mathfrak{G}_{\varepsilon,s+1/n}\bigg\lvert \binom{\ell-1}{k}\frac{\prod_{i=0}^{k-1}(\bcabc{s+1/n}-i)\prod_{i'=0}^{\ell-2-k}{(\edgeli{s+1/n}-\edgesl{s+1/n}-\bcabc{s+1/n}-i')}}{\prod_{i''=0}^{\ell-2}(\edgeli{s+1/n}-\edgesl{s+1/n}-i'')} 
-r_{k,\ell} \bigg\rvert} =\bar{o}_n(1).
\end{align*}
Then the tower property yields that 
\begin{align*}
&\mathbb{E}\brk{\ind{\mathfrak G_{\eps, s+1/n}} \abs{\PP\bc{\bar{\vd}_{\vnu_{s+1/n},s}=k
\mid \fscm{n}{s+1/n},d_{\vnu_{s+1/n}}=\ell}-r_{k,\ell}}}=\bar{o}_n(1).
\end{align*}
Thus, 
\begin{align}\label{eq_con_cur_degree}
& \Erw\abs{\PP\bc{\bar{\vd}_{\vnu_{s+1/n},s}=k \mid \fscm{n}{s+1/n},d_{\vnu_{s+1/n}}=\ell}-r_{k,\ell}}\\
=&\mathbb{E}\brk{\ind{\mathfrak G_{\eps, s+1/n}} \abs{\PP\bc{\bar{\vd}_{\vnu_{s+1/n},s}=k,  
 \mid \fscm{n}{s+1/n},d_{\vnu_{s+1/n}}=\ell}-r_{k,\ell}}}+\bar{o}_n(1) =\bar{o}_n(1).\nonumber
\end{align}
Note that $q_k=\sum_{\ell=k+1}^K \ell\eul^{-\ell t_s}p_\ell\lambda(t_s)^{-1}r_{k,\ell}$. The desired result then follows directly from the combination of \cref{eq_prod_deff}, \cref{eq_degree_LTP}, \Cref{lem_vcm}, \cref{eq_con_cur_degree},  and the tower property.

\section{Fixed-point equations}\label{sec_fixeqs}
Building upon the stability results from Section \ref{Sec_stability_type} and the conditional degree analysis from \Cref{complex_conexp}, we may now derive the type fixed-point equations. In this sense, the main result of the present section is the following:
\begin{proposition}[Type fixed-point equations]\label{al11}
Fix $\rangee$. Then uniformly in $\ranges$,
\begin{align}
   & \vy_s=1-\hat{\psi}_{t_s}\bc{\vx_s+\vy_s+\vu_s}-\hat{\psi}_{t_s}\bc{\vx_s+\vy_s+\vv_s}+\hat{\psi}_{t_s}\bc{\vx_s+\vy_s}+\oone;\label{eq_fpy}\\
  &  \vu_s=\hat{\psi}_{t_s}\bc{\vx_s+\vy_s+\vu_s}-\hat{\psi}_{t_s}\bc{\vx_s+\vy_s}+\oone;\label{eq_fpu}\\
 &   \vv_s=\hat{\psi}_{t_s}\bc{\vx_s+\vy_s+\vv_s}-\hat{\psi}_{t_s}\bc{\vx_s+\vy_s}+\oone;\label{eq_fpv}\\
  &  \vz_s\geq \hat{\psi}_{t_s}\bc{\vy_s}+\oone. \label{eq_fpz}
\end{align}
\end{proposition}

Throughout the section, to treat all four equations in a unified way, 
we use the following suggestive notation for 
the functions on the right-hand sides of (\ref{eq_fpy}) to (\ref{eq_fpz}):

\begin{definition}[Type functions]\label{def_typ_fun}
Let $\mathcal G$ denote the set of non-decreasing functions $g:[0,1]\to [0,1]$ and $\Delta^4$ be the four-dimensional standard simplex. We then define the following three functions $Y,U,V:\Delta^4\times \mathcal{G} \to [0,1]$ by setting
\begin{enumerate}[label=(\roman*)]
  \item $Y\bc{\zeta,g}=1-g(x+y+u)-g(x+y+v)+g(x+y)$ for $(\zeta,g) \in \Delta^4\times \mathcal{G}$;
  \item $U\bc{\zeta,g}=g(x+y+u)-g(x+y)$ for $(\zeta,g) \in \Delta^4\times \mathcal{G}$;
  \item $V\bc{\zeta,g}=g(x+y+v)-g(x+y)$ for $(\zeta,g) \in \Delta^4\times \mathcal{G}$;
  \item $Z\bc{\zeta,g}=g(y)$ for $(\zeta,g) \in \Delta^4\times \mathcal{G}$.
\end{enumerate}
\end{definition}
Recall that our proof strategy towards \Cref{al11} is via \Cref{prop_allcomestogether}. Given that we have established the hold of assumption \cref{eq_allcomestogether_a1} in \Cref{Sec_stability_type}, it remains to establish assumption \cref{eq_allcomestogether_a2} (or \cref{eq_allcomestogether_a2_2}) for $\vw\in\cbc{\vy,\vu,\vv}$ ($\vw=\vz$):
\begin{align}\label{eq-yuvz-conp-w}
\mathbb{P}\bc{\vnu_{s+1/n}\in\mathcal{W}\bc{\bm{A}_{n,s}[\bm{\theta}]} \mid \vze_{s+1/n}}-W\bc{\vze_{s+1/n},\hat{\psi}_{t_s}}=\oone \qquad (\text{or }\geq \oone,\text{ respectively}).
\end{align}
In order to establish \cref{eq-yuvz-conp-w}, in \Cref{sec_type_events}, we relate the conditional probability that $\vnu_{s+1/n}$ assumes a certain type in $\bm{A}_{n,s}[\bm{\theta}]$ to the types of its neighbors in $\bm{A}_{n,s+1/n}[\bm{\theta}]$ and $\bm{A}_{n,s+1/n}[\bm{\theta}]^T$. We finally prove \Cref{al11} in \Cref{sec_proof_al11}.

\subsection{Basic and type events} \label{sec_type_events}
We first introduce two batches of events, \textit{basic} and \textit{type}. The idea is that the type events are (mostly) derived from different combinations of the basic events, and directly relate to the conditional probability that $\vnu_{s+1/n}$ is of a particular type in $\bm{A}_{n,s}[\bm{\theta}]$.

\begin{definition}[Basic and type events]\label{def_bas_eve}
For $\ranges$, define two basic events as follows:
\begin{align}
\ffB &= \cbc{\supp{\bm{A}_{n,s}(\vnu_{s+1/n},)}\subseteq \mathcal{F}\bc{\bm{A}_{n,s+1/n}[\bm{\theta}]} },\label{def_ffB} \\
\ffBT &= \cbc{ \supp{\bm{A}_{n,s}(\vnu_{s+1/n},)}\subseteq \mathcal{F}\bc{\bm{A}_{n,s+1/n}[\bm{\theta}]^T}} \label{def_ffBT}.
\end{align}
Based on the two basic events, define the following five type events:
\begin{align*}
{\mathfrak Y}_s &= \fB \cap \fBT,  \\
{\mathfrak U}_s &= \fB \cap \ffBT,  \\
{\mathfrak V}_s &= \ffB \cap \fBT, \\
{\mathfrak {XZ}}_s &= \ffB \cap \ffBT
\quad \text{and }\\
{\mathfrak Z}^{\circ}_s &= \{ \supp{\bm{A}_{n,s}(\vnu_{s+1/n},)}\subseteq \mathcal{Y}\bc{\bm{A}_{n,s+1/n}[\bm{\theta}]}\}.
\end{align*}
\end{definition}
The first two intermediate results of this subsection establish that on each type event ${\mathfrak W}_s$, $\vnu_{s+1/n}$ essentially is an element of $\mathcal{W}\bc{\bm{A}_{n,s}[\bm{\theta}]}$.
\begin{lemma}\label{lem_typ_equ}
For any $W \in \{Y,U,V\}$ and $\ranges$,
\begin{align} \label{bound_QY}
\mathbb{P}\left(\vnu_{s+1/n}\notin\mathcal{W}\bc{\bm{A}_{n,s}[\bm{\theta}]}, {\mathfrak W}_s\right) =\bar{o}_{n,P}(1),
\end{align}
as well as
\begin{align} \label{bound_QXZ}
\mathbb{P}\left(\vnu_{s+1/n}\not\in\mathcal{X}\bc{\bm{A}_{n,s}[\bm{\theta}]}\cup \mathcal{Z}\bc{\bm{A}_{n,s}[\bm{\theta}]}, {\mathfrak{XZ}}_s\right) = \bar{o}_{n,P}(1).
\end{align}
\end{lemma}
The proof of \Cref{lem_typ_equ} is given in Appendix \ref{sec_proof_typ_equ}. 
\begin{lemma}\label{lem_typ_zin}
For any $\ranges$,
\[ \mathbb{P}\left(\vnu_{s+1/n}\not\in\mathcal{Z}\bc{\bm{A}_{n,s}[\bm{\theta}]}, {\mathfrak Z}_s^{\circ} \right) = \bar{o}_{n,P}(1).
\]
\end{lemma}
The proof of \Cref{lem_typ_zin} is given in Appendix \ref{sec_proof_typ_zin}. Having identified the close connection between the type events and the type of $\vnu_{s+1/n}$, it remains to estimate their (conditional) probabilities, which we do in the next lemma.

\begin{lemma}\label{lem_eq_tye_tyf}
Fix $\rangee$. For any $W\in\cbc{Y,U,V}$ and $\ranges$, 
\begin{equation}\label{eq_css1_7}
        \mathbb{P}\bc{{\mathfrak W}_s \mid \vze_{s+1/n}}-W(\vze_{s+1/n},\hat{\psi}_{t_s})=\oone,
\end{equation}
and
\begin{equation}\label{eq_css1_7_new}
        \mathbb{P}\bc{{\mathfrak Z}_s^{\circ} \mid \vze_{s+1/n}}-\hat{\psi}_{t_s}\bc{\bm{y}} =\oone.
\end{equation}
\end{lemma}
\begin{proof}
We carry out the proof of \cref{eq_css1_7} for $W=Y$ in detail; the proofs for the other choices of $W$ and \cref{eq_css1_7_new} proceed along the same lines. 
As ${\mathfrak Y}_s = \fB \cap \fBT$, by the inclusion-exclusion principle,
\begin{align}\label{eq_al1_5}
    \PP\bc{{\mathfrak Y}_s}= \PP\bc{\fB} + \PP\bc{\fBT} - \PP\bc{\fB \cup \fBT} = 1 - \PP\bc{\ffB} - \PP\bc{\ffBT} + \PP\bc{\ffB \cap \ffBT}.
\end{align}
Moreover, by Definitions \ref{dxyzuv} and \ref{def_bas_eve},
\begin{enumerate}
  \item[(a)] $\ffB=\{\supp{\bm{A}_n(\vnu_{s+1/n},)}\subseteq \mathcal{X}\bc{\bm{A}_{n,s+1/n}[\bm{\theta}]}\cup \mathcal{Y}\bc{\bm{A}_{n,s+1/n}[\bm{\theta}]}\cup \mathcal{V}\bc{\bm{A}_{n,s+1/n}[\bm{\theta}]}\}$;
  \item[(b)] $\ffBT=\{\supp{\bm{A}_n(\vnu_{s+1/n},)}\subseteq \mathcal{X}\bc{\bm{A}_{n,s+1/n}[\bm{\theta}]}\cup \mathcal{Y}\bc{\bm{A}_{n,s+1/n}[\bm{\theta}]}\cup \mathcal{U}\bc{\bm{A}_{n,s+1/n}[\bm{\theta}]}\}$;
  \item[(c)] $\ffB \cap \ffBT=\{\supp{\bm{A}_n(\vnu_{s+1/n},)}\subseteq \mathcal{X}\bc{\bm{A}_{n,s+1/n}[\bm{\theta}]}\cup \mathcal{Y}\bc{\bm{A}_{n,s+1/n}[\bm{\theta}]}\}$.
\end{enumerate}
Let $f_d:[0,1] \to [0,1]$ be defined by $f_d(x)=x^d$. Instead of starting at \eqref{eq_css1_7} directly, we first show that  
\begin{align}\label{eq_bac_y}
\Erw\abs{\PP\bc{\mathfrak{Y}_s \mid \fscm{n}{s+1/n},\THETA,\bar{\vd}_{\vnu_{s+1/n},s}}-Y\bc{\vze_{s+1/n},f_{\bar{\vd}_{\vnu_{s+1/n},s}}}}= \bar{o}_n(1),
\end{align}
where as always, $\THETA$ is the pair of perturbation matrices. 

Next, recall from \cref{event_gs} that $\mathfrak{G}_{\varepsilon,s+1/n}$ denotes the event that the number of $j\in \verst{s+1/n}$ such that $\bar{\vd}_{j,s+1/n}<d_j$ is greater than $cn$ and from \cref{ineq_dcs} that $\PP(\mathfrak{G}_{\varepsilon,s+1/n})=1-\bar{o}_n(1)$. In what follows, we will thus assume that  \scm{n}{s+1/n}  has the corresponding degree properties.  
Conditionally on such \scm{n}{s+1/n}, $\THETA$ and $\bar{\vd}_{\vnu_{s+1/n},s}$, and additionally on the event that $\vnu_{s+1/n}$ does not have any self-loops, by \Cref{claim_extra}, the neighboring half-edges of $\vnu_{s+1/n}$ in \scm{n}{s} are chosen uniformly at random without replacement from all the half-edges adjacent to vertices in $\verst{s+1/n}$ which have not been paired in \scm{n}{s+1/n}. Denote the number of those possible half-edges by $\bm m_{s+1/n}=\sum_{i\in \verst{s+1/n}}(d_i-\bar{\vd}_{i,s+1/n})$. Out of these, $\sum_{i\in \verst{s+1/n}}(d_i-\bar{\vd}_{i,s+1/n})\teo{i\in \mathcal{Y}\bc{\bm{A}_{n,s+1/n}[\bm{\theta}]}}=\bm m_{s+1/n} \vy_{s+1/n}$ belong to a vertex with type $Y$ in $\bm{A}_{n,s+1/n}[\bm{\theta}]$. Therefore, by \Cref{def_typ_fun} and \Cref{eq_al1_5},
\begin{align*}
    &\abs{\PP\bc{\mathfrak{Y}_s \mid \fscm{n}{s+1/n},\THETA,\bar{\vd}_{\vnu_{s+1/n},s}}-Y\bc{\vze_{s+1/n},f_{\bar{\vd}_{\vnu_{s+1/n},s}}}}\\
&\qquad\leq  \Bigg|  \frac{\binom{(\vx_{s+1/n}+\vy_{s+1/n}+\vv_{s+1/n})\bm m_{s+1/n}}{\bar{\vd}_{\vnu_{s+1/n},s}}}{\binom{\bm m_{s+1/n}}{\bar{\vd}_{\vnu_{s+1/n},s}}} - (\vx_{s+1/n}+\vy_{s+1/n}+\vv_{s+1/n})^{\bar{\vd}_{\vnu_{s+1/n},s}} + \nonumber\\&\qquad\quad\frac{\binom{(\vx_{s+1/n}+ \vy_{s+1/n}+\vu_{s+1/n})\bm m_{s+1/n}}{\bar{\vd}_{\vnu_{s+1/n},s}}}{\binom{\bm m_{s+1/n}}{{\bar{\vd}_{\vnu_{s+1/n},s}}}}  -(\vx_{s+1/n}+\vy_{s+1/n}+\vu_{s+1/n})^{\bar{\vd}_{\vnu_{s+1/n},s}} \nonumber\\&\qquad\quad- \frac{\binom{(\vx_{s+1/n}+\vy_{s+1/n})\bm m_{s+1/n}}{\bar{\vd}_{\vnu_{s+1/n},s}}}{\binom{\bm m_{s+1/n}}{{\bar{\vd}_{\vnu_{s+1/n},s}}}} 
+ (\vx_{s+1/n}+\vy_{s+1/n})^{\bar{\vd}_{\vnu_{s+1/n},s}} \Bigg| \nonumber\\
&\qquad\quad+\PP\bc{\mbox{$\vnu_{s+1/n}$ has self-loop(s)} \mid \fscm{n}{s+1/n},\bar{\vd}_{\vnu_{s+1/n},s}}.
\end{align*}
We next bound the terms above separately.
\begin{enumerate}
    \item For each integer $0\leq k\leq N$,
\begin{align*}
    0\leq 1-\frac{\prod_{i=0}^{k-1}(N-i)}{N^k}\leq\sum_{i=0}^{k-1}\frac{i}{N}= \frac{k(k-1)}{2N}.
\end{align*}
As a consequence, for each integer $0\leq k\leq N_1\leq N_2$,
\begin{align*}
    \abs{\frac{\binom{N_1}{k}}{\binom{N_2}{k}}-\bc{\frac{N_1}{N_2}}^k}\leq \frac{N_1^k}{\prod_{i=0}^{k-1}(N_2-i)}\frac{k(k-1)}{2N_1}.
\end{align*}
Hence, for $a\in[0,1]$ and $k\leq K$, as $\bm m_{s+1/n}\geq cn$ on $\mathfrak{G}_{\eps,s+1/n}$, 
\begin{align*}
    \abs{\frac{\binom{a \bm m_{s+1/n}}{k}}{\binom{\bm m_{s+1/n}}{k}}-a^k}\leq \frac{a^{k-1}}{\prod_{i=0}^{k-1}(1-i/\bm m_{s+1/n})}\frac{k(k-1)}{2\bm m_{s+1/n}}=\bar{o}_n(1).
\end{align*}

    \item On $\mathfrak{G}_{\eps,s+1/n}$, there are at least $cn$ living half-edges when $\vnu_{s+1/n}$ is awakened. Consequently, 
by \Cref{claim_extra}, 
\begin{align*}
    \PP\bc{\mbox{$\vnu_{s+1/n}$ has at least a self-loop},\mathfrak{G}_{\eps,s+1/n}}
    \leq K\frac{K}{cn}=\bar{o}_n(1).
\end{align*}
\end{enumerate}
 Consequently, restricted to $\mathfrak{G}_{\eps,s+1/n}$, \Cref{eq_bac_y} holds, from which the unrestricted version follows using \Cref{ineq_dcs}. 
Hence,
\begin{align}\label{eq_bac_y2}
&\Erw\abs{\PP\bc{\mathfrak{Y}_s \mid \vze_{s+1/n}}-\Erw\brk{Y\bc{\vze_{s+1/n},f_{\bar{\vd}_{\vnu_{s+1/n},s}}}\mid \vze_{s+1/n}}}\\
&\qquad=\Erw\abs{\Erw\brk{\PP\bc{\mathfrak{Y}_s \mid \fscm{n}{s+1/n},\THETA,\bar{\vd}_{\vnu_{s+1/n},s}}-Y\bc{\vze_{s+1/n},f_{\bar{\vd}_{\vnu_{s+1/n},s}}}\mid \vze_{s+1/n}}}\nonumber\\
 &\qquad\leq\Erw\abs{\PP\bc{\mathfrak{Y}_s \mid \fscm{n}{s+1/n},\THETA,\bar{\vd}_{\vnu_{s+1/n},s}}-Y\bc{\vze_{s+1/n},f_{\bar{\vd}_{\vnu_{s+1/n},s}}}}=\bar{o}_n(1).\nonumber
\end{align}
On the other hand, recall $q_k(s)$ from \eqref{def_q_k}, which is the coefficient of $\alpha^m$ in $\hat{\psi}_{t_s}(\alpha)$. 
By the fact that $$\PP\bc{\bar{\vd}_{\vnu_{s+1/n},s}=k \mid \fscm{n}{s+1/n}}=\PP\bc{\bar{\vd}_{\vnu_{s+1/n},s}=k \mid \fscm{n}{s+1/n},\THETA},$$
\Cref{lem_vscm} and the tower property, we have
\begin{align*}
\mathbb{E}\abs{\PP\bc{\bar{\vd}_{\vnu_{s+1/n},s}=k \mid \vze_{s+1/n}}-q_k(s)}\leq \mathbb{E}\abs{\PP\bc{\bar{\vd}_{\vnu_{s+1/n},s}=k \mid \fscm{n}{s+1/n},\THETA}-q_k(s)}=\bar{o}_n(1).
\end{align*}
Consequently,
\begin{align}\label{eq_bac_y_2}
&\Erw\abs{\Erw\brk{Y(\vze_{s+1/n},f_{\bar{\vd}_{\vnu_{s+1/n},s}}) \mid \vze_{s+1/n}}-Y(\vze_{s+1/n},\hat{\psi}_{t_s})}\\
    =&\Erw\abs{\sum_{k=0}^K Y(\vze_{s+1/n},f_k)\PP\bc{\bar{\vd}_{\vnu_{s+1/n},s}=k \mid \vze_{s+1/n}} -Y(\vze_{s+1/n},\hat{\psi}_{t_s})}\nonumber\\
    =&\Erw\abs{ Y\bc{\vze_{s+1/n},\sum_{k=0}^K q_k(s) f_k}-Y(\vze_{s+1/n},\hat{\psi}_{t_s})}+\bar{o}_n(1)=\bar{o}_n(1).\nonumber
\end{align}
Combining \Cref{eq_bac_y2} and \Cref{eq_bac_y_2} gives \Cref{eq_css1_7} for $W=Y$.
\end{proof}

\subsection{Proof of \Cref{al11}} \label{sec_proof_al11}
With \Cref{lem_typ_equ,lem_typ_zin,lem_eq_tye_tyf} in hand, we are now able to prove \Cref{al11}.
It follows directly 
from the combination of \Cref{lem_typ_equ,lem_typ_zin,lem_eq_tye_tyf} that for  any $W \in \cbc{Y,U,V}$ and $\ranges$,
\begin{equation}\label{eq_fixyuv}
    \mathbb{P}\bc{\vnu_{s+1/n}\in\mathcal{W}\bc{\bm{A}_{n,s}[\bm{\theta}]} \mid \vze_{s+1/n}}-W\bc{\vze_{s+1/n},\hat{\psi}_{t_s}}=\oone,
\end{equation}
while
\begin{equation}\label{eq_ineqz}
\mathbb{P}\bc{\vnu_{s+1/n}\in\mathcal{Z}\bc{\bm{A}_{n,s}[\bm{\theta}]} \mid \vze_{s+1/n}}-\hat{\psi}_{t_s}\bc{\vy_{s+1/n}} \geq \oone.
\end{equation}
Combining \eqref{eq_fixyuv} and \eqref{eq_ineqz} with \Cref{prop_allcomestogether} yields the conclusion of \Cref{al11}.

\section{Proof of \Cref{t_main}: Lower bound on the rank}\label{sec_rank_est}
Recall that we fix $\rangee$ and let $\iota=1-\sigma(-\ln \xi)-\varepsilon$. As outlined, our strategy to lower-bound the asymptotic rank is a telescoping decomposition according to the graph exploration. In particular, recall from \Cref{sec_rank_types} that we build on the lower bound
\begin{align}
\Erw\brk{\rk{\bm{A}_n\pth}}\geq&\sum_{ns=\lfloor n\varepsilon \rfloor}^{\lceil n\iota\rceil-1}\bc{\Erw\brk{\rk{\bm{A}_{n,s}\pth}-\rk{\bm{A}_{n,s+1/n}\pth}}}+\Erw\brk{\rk{\bm{A}_{n,\iota}\pth}}.
\end{align}
In this section, we first derive a lower bound on the individual summands $\Erw\brk{\rk{\bm{A}_{n,s}[\vth]}-\rk{\bm{A}_{n,s+1/n}[\vth]}}$ in \Cref{sec_newapp} and then prove \Cref{t_main} in \Cref{sec_proof_t_main}.

\subsection{Rank-difference: From fixed-point equations to lower bound}\label{sec_newapp}
Already in \cref{eq_re_rd}, we derived the expression
\begin{align}\label{eq_re_rd2_0}
\Erw\brk{\rk{\bm{A}_{n,s}[\vth]}-\rk{\bm{A}_{n,s+1/n}[\vth]}}=&\Erw\brk{\vx_s+2\vy_s+\vu_s+\vv_s}+\bar{o}_{n,P}(1)
\end{align}
of the expected rank decrease in terms of the different variable types. Combining \eqref{eq_re_rd2_0} with the fixed-point equations \cref{eq_fpy} and \cref{eq_fpu} and the deterministic relation
 \begin{align}\label{eq_val_xyv}
 \val_s=\vx_s+\vy_s+\vv_s,  
 \end{align}
we obtain 
\begin{align}\label{eq_re_rd2}
\Erw\brk{\rk{\bm{A}_{n,s}[\vth]}-\rk{\bm{A}_{n,s+1/n}[\vth]}}
=\Erw\brk{\val_s+1-\hat{\psi}_{t_s}(\val_s)}+\bar{o}_{n,P}(1) 
=:\Erw\brk{h_{t_s}(\val_s)}+\bar{o}_{n,P}(1).
\end{align}
In the last step, we have defined the function $ h_t:[0,1] \to \RR$ given by
\begin{align}\label{def_h}
h_t\bc{\alpha}:=\alpha+1-\hat{\psi}_t\bc{\alpha}.   
\end{align}
Unfortunately, $\val_s$ in \cref{eq_re_rd2} is a rather complicated random variable and might not even converge in an appropriate sense.  Nevertheless, the type fixed-point equations in \Cref{al11} also open the door towards the derivation of a lower bound on \( \Erw[h_{t_s}(\val_s)] \), 
up to an \( \oone\) error. 
To formulate the lower bound, we define the continuous function $ G_t:[0,1]\to \mathbb R$ given by
\begin{align}\label{def_G}
G_t(\alpha) := \alpha + \hat{\psi}_t(1 - \hat{\psi}_t(\alpha)) - 1 
\end{align}
and let $\alpha^\star(t)$ be the largest zero of $G_t$ in $[0,1]$. Then we have the following lower bound on the rank difference:

\begin{lemma}\label{imp_pro}
Fix $\rangee$.
For any $\ranges$ and $h_t$ as in \cref{def_h},
\begin{align}
\Erw\brk{\rk{\bm{A}_{n,s}\pth}-\rk{\bm{A}_{n,s+1/n}\pth}}
\geq h_{t_s}(\alpha^\star(t_s))+\bar{o}_{n,P}(1).
\end{align}
\end{lemma}
The proofs of \Cref{imp_pro} and later results in the current section heavily depend on the properties of the functions $G_t$ and their zeroes. We therefore take a closer look at these now. Let
\begin{enumerate}
\item $\alpha_0(t)$ be the unique zero of the increasing function $\Xi_t:[0,1]\to \RR$ given by  $\Xi_t(\alpha) = \alpha+\hat{\psi}_{t}(\alpha)-1$. It is straightforward to check that $\alpha_0(t)$ is always a zero of $G_{t}$.
    \item $\alpha_\star(t)$ and $\alpha^\star(t)$  be the smallest and  largest (non-necessarily distinct) zeroes of $G_{t}(\alpha)$ in $[0,1]$, respectively.
\end{enumerate}
Finally, set
\begin{align}\label{def_kappa}
     \kappa:=-\ln{((\hat{\psi}')^{-1}(1)+\hat{\psi}((\hat{\psi}')^{-1}(1)))}.
\end{align}
Then the following results hold:
\begin{restatable}[Properties of $G_{t}$ and its zeroes]{lemma}{lemproal}\label{lem_proal}
Assume that the probability distribution $(p_k)_{k\geq 0}$  satisfies Assumptions \ref{assumption_weaker}, \ref{tech_assumption} and \ref{assumption_stronger}. Then for $t \in [0,1]$, the function $G_t$ defined in \cref{def_G} has the following properties:
    \begin{enumerate}
        \item 
        $G_{t}$ has at most $3$ zeroes, and $G'_{t}$ has at most $2$ zeroes. \label{ap_it1}
        \item \label{ap_it2} If $G_{t}(\alpha)=0$, then also $G(1-\hat{\psi}_{t}(\alpha))=0$. Specifically, $\alpha_\star(t)=1-\hat{\psi}_{t}\bc{\alpha^\star(t)}$ and $\alpha^\star(t)=1-\hat{\psi}_{t}\bc{\alpha_\star(t)}$. 
        \item $G_{t}$ has $1$ or $3$ zeroes. If $G_t$ has $3$ zeroes, then $\alpha_\star(t)<\alpha_0(t)<\alpha^\star(t)$. \label{ap_it3}
 \item For any $\alpha> \alpha^\star(t)$, $G_{t}(\alpha)> 0$. \label{ap_it4}
 
        \item For $t< \kappa$, $G_t'(\alpha_0(t))<0$, $G_t$ has precisely $3$  zeroes and both $G_t'(\alpha_\star(t))$ and $G_t'(\alpha^\star(t))$ are not equal to $0$; for $t> \kappa$, $G_t'(\alpha_0(t))>0$ and $G_t$ has only $1$ zero; for $t=\kappa$, $G_t'(\alpha_0(t))=0$ and $G_t$ has only $1$ zero. As a consequence, for $t\neq \kappa$, $G_t$ and $G'_t$ have no common zero while, for $t=\kappa$, the only common zero is $\alpha_0(\kappa)$.\label{ap_it5}
       
        \item The functions $t \mapsto \alpha_\star(t)$, $t \mapsto \alpha_0(t)$ and $t \mapsto \alpha^\star(t)$ are continuous on $[0,\infty)$ and continuously differentiable on $[0,\infty)\backslash\cbc{\kappa}$.\label{ap_it6}
        \item For a random variable $\bm{b}$, if $ G_t(\bm{b})=\bar{o}_{\mathbb{P}}(1)$, then also
\begin{align*}
  \min\cbc{\abs{\bm{b}-\alpha_\star(t)},\abs{\bm{b}-\alpha_0(t)},\abs{\bm{b}-\alpha^\star(t)}}=\bar{o}_{\mathbb{P}}(1).
    \end{align*}\label{ap_it7}

    \item The function $R_{\psi_t}$ obtains its minimum on $[0,1]$ for $\alpha\in \{\alpha_\star(t),\alpha^\star(t)\}$.\label{ap_it8}
    \end{enumerate}
\end{restatable}
The proof of \Cref{lem_proal} heavily depends on \Cref{tech_assumption} and is given in Appendix \ref{sec_progt}. Given \Cref{al11} and \Cref{lem_proal}, the proof of \Cref{imp_pro} is then almost identical to the proof of \cite[Proposition 2.14]{HofMul25}, which is why we defer it to Appendix \ref{sec_proof_imp_pro}.

\subsection{Overview over the proof of \Cref{t_main}} \label{sec_proof_t_main}
An application of the lower bound of \Cref{imp_pro} on the telescoping decomposition \Cref{eq_re_rd2_0} yields
\begin{align}\label{ineq_rk_decompo}
\Erw\brk{\rk{\bm{A}_n\pth}} 
\geq& \sum_{s=n\varepsilon}^{n\iota-1}h_{t_s}(\alpha^\star(t_s))+\Erw\brk{\rk{\bm{A}_{n,\iota}\pth}}+\bar{o}_{n,P}(n)\nonumber\\
=&n\int_{\varepsilon}^{\iota}h_{t_s}(\alpha^\star(t_s))\dif s+\Erw\brk{\rk{\bm{A}_{n,\iota}\pth}}+\bar{o}_{n,P}(n).
\end{align}
Recall that $\sigma(t)=\sum_{k\geq 0}p_k \eul^{-kt}$.
Equation \cref{ineq_rk_decompo} splits the proof of \Cref{t_main} naturally into the following two lemmas:

\begin{lemma}\label{eq_int_alpha}
For any $S\in [0,\iota]$, with $h_t$ as defined in \eqref{def_h}, 
\begin{align*}
    \int_{0}^{S} h_{t_s}(\alpha^\star(t_s)) \dif s  =\sigma\bc{t_0} R_{\psi_{t_0}}(\alpha^\star(t_0)) - \sigma\bc{t_S}R_{\psi_{t_S}}(\alpha^\star(t_S)).
\end{align*}
\end{lemma}

\begin{lemma}\label{lem_dec_2}
Let $R_{\phi}$ be as defined in \cref{def_rd}. Then, for any $\varepsilon'>0$,
\begin{align}\label{eq_rk_aniota}
    \lim_{n\to\infty}\sup_{J_n\in \syn}\PP\bc{\abs{\frac{1}{n}\rank_{\FF}\bc{\bm{A}_{n,\iota}}-\sigma(t_\iota)R_{\psi_{t_\iota}}(\alpha^\star(t_\iota))}\geq \varepsilon'}=0.
\end{align}
\end{lemma}

We next prove \Cref{eq_int_alpha} and \Cref{lem_dec_2} in their order of appearance, and then conclude with the proof of \Cref{t_main}.

\subsubsection{Proof of \Cref{eq_int_alpha}}
Fix $S \in [0,\iota]$. For any $s \in [0,S]$, the definition of $t_s$ via \cref{def_H_s} ensures that $\sigma(t_s)=1-s$ and $t_0=0$. In particular, 
\begin{align*}
    \dif s=\dif (1-\sigma(t_s))=\lambda(t_s)\dif t_s.
\end{align*}
Moreover,
\begin{flalign*}
    \sigma\bc{t_0} R_{\psi_{t_0}}(\alpha^\star(t_0)) - \sigma\bc{t_S}R_{\psi_{t_S}}(\alpha^\star(t_S))
    =  2-\psi(1-\hat{\psi}(\alpha^\star(0)))-\psi(\alpha^\star(0))-\psi'(\alpha^\star(0))(1-\alpha^\star(0))\\
     -\sigma\bc{t_S} \bc{2-\psi_{t_S}(1-\hat{\psi}_{t_S}(\alpha^\star(t_S)))-\psi_{t_S}(\alpha^\star(t_S))-\psi'_{t_S}(\alpha^\star(t_S))(1-\alpha^\star(t_S))}.
\end{flalign*}
Now let $T=T(S)=t_S$ such that $T \in [0, \sigma^{-1}(1-\iota)]$. We hence aim to prove that
\begin{align}\label{equation_int_t}
    &\int_{0}^{T} \bc{\alpha^\star(t)+1-\hat{\psi}_{t}\bc{\alpha^\star(t)} }\lambda(t)\dif t\\
    &\qquad=2-\psi(1-\hat{\psi}(\alpha^\star(0)))-\psi(\alpha^\star(0))-\psi'(\alpha^\star(0))(1-\alpha^\star(0))\nonumber\\
    &\qquad\quad-\sigma\bc{T}\bc{2-\psi_{T}(1-\hat{\psi}_{T}(\alpha^\star(T)))-\psi_{T}(\alpha^\star(T))-\psi'_{T}(\alpha^\star(T))(1-\alpha^\star(T))}.\nonumber
\end{align}
With the abbreviations 
\begin{align}
    q(T):=\int_{0}^T h(t)\lambda(t)\dif t=\int_{0}^{T} \bc{\alpha^\star(t)+1-\hat{\psi}_{t}\bc{\alpha^\star(t)} }\lambda(t)\dif t
\end{align}
and
\begin{align}
    r(T):=\sigma\bc{T}\bc{2-\psi_{T}(1-\hat{\psi}_{T}(\alpha^\star(T)))-\psi_{T}(\alpha^\star(T))-\psi'_{T}(\alpha^\star(T))(1-\alpha^\star(T))}.
\end{align}
this reduces to showing that $q(T) = r(0)-r(T)$. 

    As $T\downarrow 0$, both sides of \cref{equation_int_t} tend to $0$. Hence, it is sufficient to prove that $-q'(T)=r'(T)$. Looking at the left-hand side, it is immediate that
    \begin{align}
       - q'(T)=-\bc{\alpha^\star(T)+1-\hat{\psi}_{T}\bc{\alpha^\star(T)} }\lambda(T).
    \end{align}
For the right-hand side, let $\hat{\psi}(t,\alpha) := \hat{\psi}_t(\alpha)$,  $f(t,\alpha):=\sigma(t)\psi_t(\alpha)$ and $g(t,\alpha):=\lambda(t)\hat{\psi}_t(\alpha)$. Then 
\begin{align}
 \frac{\partial f}{\partial t}(t,\alpha)=-\bc{1+(1-\alpha)\frac{(\lambda'(t)+\lambda(t))\eul^{2t}}{\lambda(0)}}g(t,\alpha),
\end{align}
\begin{align}
    \frac{\partial f}{\partial \alpha}(t,\alpha)=\frac{\lambda(t)\eul^{2t}}{\lambda(0)}g(t,\alpha),
\end{align}
and
\begin{align}
r(T)=2\sigma\bc{T}-f(T,1-\hat{\psi}_T(\alpha^\star(T)))-f(T,\alpha^\star(T))-\frac{\lambda(T)\eul^{2T}}{\lambda(0)} g(T,\alpha^\star(T))(1-\alpha^\star(T)).\end{align}
Hence,
\begin{align}
r'(T)=&-2\lambda(T)+\bc{1+\hat{\psi}_T(\alpha^\star(T))\frac{(\lambda'(T)+\lambda(T))\eul^{2T}}{\lambda(0)}}g(T,1-\hat{\psi}_T(\alpha^\star(T)))\\
&+\frac{\lambda(T)\eul^{2T}}{\lambda(0)}g(T,1-\hat{\psi}_T(\alpha^\star(T)))\bc{\frac{\partial \hat{\psi}}{\partial t}(T,\alpha^\star(T))+\frac{\partial \hat{\psi}}{\partial \alpha}(T,\alpha^\star(T)){\alpha^\star}'(T)}\nonumber\\
&+\bc{1+(1-\alpha^\star(T))\frac{(\lambda'(T)+\lambda(T))\eul^{2T}}{\lambda(0)}}g(T,\alpha^\star(T))-\frac{\lambda(T)\eul^{2T}}{\lambda(0)}g(T,\alpha^\star(T)){\alpha^\star}'(T)\nonumber\\
&-\frac{(\lambda'(T)+2\lambda(T))\eul^{2T}}{\lambda(0)} g(T,\alpha^\star(T))(1-\alpha^\star(T))-\frac{\lambda(T)\eul^{2T}}{\lambda(0)} \frac{\partial g}{\partial t}(T,\alpha^\star(T))(1-\alpha^\star(T))\nonumber\\
&-\frac{\lambda(T)\eul^{2T}}{\lambda(0)} \frac{\partial g}{\partial \alpha}(T,\alpha^\star(T))(1-\alpha^\star(T)){\alpha^\star}'(T)+\frac{\lambda(T)\eul^{2T}}{\lambda(0)} g(T,\alpha^\star(T)){\alpha^\star}'(T).\nonumber
\end{align}
Since $g(t,\alpha)=\lambda(t)\hat{\psi}_t(\alpha)$ and $\hat{\psi}_t(1-\hat{\psi}_t(\alpha^\star(t_s)))=1-\alpha^\star(t_s)$, we have
\begin{align}
\frac{\partial \hat{\psi}}{\partial t}(t,\alpha)=\frac{\partial (g/\lambda)}{\partial t}(t,\alpha)=\frac{1}{\lambda(t)}\frac{\partial g}{\partial t}(t,\alpha)-\frac{\lambda'(t)}{\lambda^2(t)}g(t,\alpha)\quad\mbox{and}\quad    \frac{\partial \hat{\psi}}{\partial \alpha}(t,\alpha) =\frac{1}{\lambda(t)}\frac{\partial g}{\partial \alpha}(t,\alpha),
\end{align}
while $g(T,1-\hat{\psi}_T(\alpha^\star(T)))=\lambda(T)(1-\alpha^\star(T))$. Hence,
\begin{align}
r'(T)=&-2\lambda(T)+\bc{\lambda(T)+g(T,\alpha^\star(T))\frac{(\lambda'(T)+\lambda(T))\eul^{2T}}{\lambda(0)}}(1-\alpha^\star(T))\\
&+\frac{\lambda(T)\eul^{2T}}{\lambda(0)}(1-\alpha^\star(T))\bc{\frac{\partial g}{\partial t}(T,\alpha^\star(T))-\frac{\lambda'(T)}{\lambda(T)}g(T,\alpha^\star(T))+\frac{\partial g}{\partial \alpha}(T,\alpha^\star(T)){\alpha^\star}'(T)}\nonumber\\
&+\bc{1+(1-\alpha^\star(T))\frac{(\lambda'(T)+\lambda(T))\eul^{2T}}{\lambda(0)}}g(T,\alpha^\star(T))-\frac{\lambda(T)\eul^{2T}}{\lambda(0)}g(T,\alpha^\star(T)){\alpha^\star}'(T)\nonumber\\
&-\frac{(\lambda'(T)+2\lambda(T))\eul^{2T}}{\lambda(0)} g(T,\alpha^\star(T))(1-\alpha^\star(T))-\frac{\lambda(T)\eul^{2T}}{\lambda(0)} \frac{\partial g}{\partial t}(T,\alpha^\star(T))(1-\alpha^\star(T))\nonumber\\
&-\frac{\lambda(T)\eul^{2T}}{\lambda(0)} \frac{\partial g}{\partial \alpha}(T,\alpha^\star(T))(1-\alpha^\star(T)){\alpha^\star}'(T)+\frac{\lambda(T)\eul^{2T}}{\lambda(0)} g(T,\alpha^\star(T)){\alpha^\star}'(T)\nonumber\\
=&-\bc{\alpha^\star(T)+1}\lambda(T)+g\bc{T,\alpha^\star(T)}= -q'(T),\nonumber
\end{align}
as desired.

\subsubsection{Proof of \Cref{lem_dec_2} subject to  \Cref{coro_main_prop_1}}
Observe that $\bm{A}_{n,\iota}$ can naturally be regarded as the adjacency matrix of the graph union of \scm{n}{\iota} and the empty graph on vertices $\{\vnu_{1/n}, \ldots, \vnu_{\iota}\}$. Call the resulting graph $\bm H'$. Then $\bm H'$ is a configuration model, and we will show that $\bm H'$ is subcritical for small enough $\eps$ in the following.
 
By 
\Cref{lem_bvks} and the definition following its proof, as $n \to \infty$, for any $k = 0, \ldots, K$,
the proportion of vertices with \textit{current} degree $k$ in \scm{n}{\iota} converges in probability to $q_k(\iota)$, the coefficient of $\alpha^k$ in the generating function
\begin{align}
    \psi_{t_\iota}(\alpha) = \sigma(t_\iota)^{-1}\sum_{m=0}^K p_m \eul^{-mt_\iota}\bc{1+\frac{\lambda(t_\iota)\eul^{2t_\iota}}{\lambda(0)}(\alpha-1)}^m.
\end{align}
Now, as $\varepsilon\downarrow 0$, $\lambda(t_\iota)\eul^{2t_\iota}\downarrow \lambda(0) = \sum_{k=0}^K k p_k\xi^{k-2}$ and $t_\iota\uparrow -\ln(\xi)$. 
Thus, for $\eps \downarrow 0$, $\psi_{t_\iota}(\alpha)$ converges to 
\begin{align*}
    \frac{\sum_{m=0}^K p_m \xi^m\bc{1+\frac{\lambda(-\ln \xi)}{\xi^{2}\lambda(0)}(\alpha-1)}^m}{\sigma(-\ln \xi)} = 
    \frac{\sum_{m=0}^K p_m\bc{\xi+\frac{\sum_{j=0}^K jp_k\xi^{j-1}}{\sum_{j=0}^Kjp_j}(\alpha-1)}^m}{\sum_{m=0}^{K}p_m\xi^m}=\frac{\psi\bc{\xi+\hat{\psi}(\xi)(\alpha-1)}}{\psi(\xi)}. 
\end{align*}
We next observe that $\alpha \mapsto \hat{\psi}(\alpha)-\alpha = \psi'(\alpha)/\psi'(1) - \alpha$ is a strictly convex function with zeros $\xi$ and $1$. Since $\xi<1$, its derivative in $\xi$ is negative, such that also $\psi'(1)(\hat{\psi}'(\xi)-1) = \psi''(\xi)-\psi'(1) <0$. 
Thus, 
\begin{align*}
 \sum_{k=0}^K k(k-2) q_k(\iota) = \psi_{t_\iota}''(1) - \psi'_{t_\iota}(1) \stackrel{\eps \downarrow 0}{\longrightarrow}  \frac{\psi''(\xi)}{\psi(\xi)}\hat{\psi}(\xi)^2-\frac{\psi'(\xi)}{\psi(\xi)}\hat{\psi}(\xi)=(\psi''(\xi)-\psi'(1))\frac{\hat{\psi}(\xi)^2}{\psi(\xi)}<0.
\end{align*}
Therefore, there exists an $\varepsilon_0$ that only depends on $(p_k)_{k=0}^K$ such that $\sum_{k=0}^K k(k-2) q_k(\iota)<0$ as long as $\varepsilon<\varepsilon_0$. We conclude that, for $\varepsilon$ small enough, the degree distribution of \scm{n}{\iota} is  subcritical.

As a consequence, the proportion of vertices of degree $k$ in $\bm H'$ converges to 
\begin{align*}    \tilde{p}_k(\iota)=\begin{cases}
        q_k(\iota) (1-\iota),\quad&k\geq 1;\\
        q_0(\iota) (1-\iota) + \iota,\quad&k=0.
    \end{cases}
\end{align*}
By the previous considerations and \Cref{lem_bvks}, for $\eps$ small enough, $(\tilde p_k)_{k\geq 0}$ and $\bm H'$ satisfy \Cref{assumption_weaker} with  $\sum_{k\geq 0}k(k-2)\tilde{p}_k(\iota)< 0$.
Moreover, let $\tilde{\psi}$ denote the p.g.f. of  $(\tilde p_k)_{k\geq 0}$. Then for all $\alpha \in [0,1]$, 
\begin{align*}
   \tilde{\psi}(\alpha)=\sigma(t_\iota)\psi_{t_\iota}(\alpha)+1- \sigma(t_\iota), \qquad \text{so that} \qquad R_{\tilde{\psi}} =\sigma(t_\iota)R_{\psi_{t_\iota}}. 
\end{align*}
By \cref{ap_it8} of \Cref{lem_proal}  and \Cref{coro_main_prop_1}, which is proven in \Cref{subsec_subcritical}, \eqref{eq_rk_aniota} holds.

\subsection{Proof of \Cref{t_main}}
By the combination of \cref{ineq_rk_decompo} and \Cref{eq_int_alpha,lem_dec_2}
\begin{align*}
\Erw\brk{\frac{1}{n}\rk{\bm{A}_n\pth}} &\geq 
\sigma\bc{t_\varepsilon} R_{\psi_{t_\varepsilon}}(\alpha^\star(t_\varepsilon)) - \sigma\bc{t_\iota}R_{\psi_{t_\iota}}(\alpha^\star(t_\iota))+\sigma\bc{t_\iota}R_{\psi_{t_\iota}}(\alpha^\star(t_\iota))-\varepsilon'+\bar{o}_{n,P}(1) 
\\ & =  \sigma\bc{t_\varepsilon} R_{\psi_{t_\varepsilon}}(\alpha^\star(t_\varepsilon)) - \varepsilon'+\bar{o}_{n,P}(1). 
\end{align*}
As $\lim_{\varepsilon\downarrow 0}t_\varepsilon=t_0=0$ and $\sigma(0)=1$, taking $\liminf_{\varepsilon'\downarrow 0}\liminf_{\varepsilon\downarrow 0}\liminf_{n\to\infty}\inf_{J_n\in \syn}$ on both sides yields
\begin{align*}
\liminf_{P\to\infty}\liminf_{n\to\infty}\inf_{J_n\in \syn}\Erw\brk{\frac{1}{n}\rk{\bm{A}_n\pth}}
    \geq& R_{\psi_{t_0}}(\alpha^\star(0))  
    = \min_{\alpha\in [0,1]}R_\psi(\alpha).
\end{align*}
The claim follows from the fact that $0\leq\rk{\bm{A}_n\pth}-\rk{\bm{A}_n}\leq 2P$.

\section{{Proof of \Cref{coro_main}}}\label{sec_proof_main}
In this section, we prove \Cref{coro_main_prop_1,coro_main_prop_2,coro_main_prop_3} in their order of appearance. Taken together, they yield \Cref{coro_main}. While the subcritical and critical cases \Cref{coro_main_prop_1} and \Cref{coro_main_prop_2} are relatively independent of the previous sections and given for completeness, removing the finite-degree condition to go from \Cref{t_main} to \Cref{coro_main_prop_3} requires some additional thought.

\subsection{Proof of \Cref{coro_main_prop_1}}\label{sec_rank_subc_c}
We use the Karp-Sipser leaf-removal algorithm, whose relation to the rank of the adjacency matrix is well-known in the literature \cite{bauer2001exactly,bordenave2011rank,karp1981maximum}, to estimate the rank in the subcritical and critical cases.

Given \cm{n} such that $\sum_{k\geq 0}k(k-2)p_k\leq 0$ and $p_2\neq 1$, we aim to compare the rank of the weighted adjacency matrix $\bm{A}_n$ over $\FF$ to that of the unweighted one over $\RR$, for which there is a rank formula available. To this end, we perform the following peeling process, which is known as the Karp-Sipser leaf-removal algorithm: 
\textit{Choose an arbitrary vertex of degree at most $1$ and remove it, along with its unique neighbor (if existent). Repeat the process until all remaining vertices have degree at least $2$. }

The process stops after at most $n$ rounds. In the end, all trees of the graph will have been removed. We call the graph that remains after the peeling process the \textit{pruned} Karp-Sipser core, which is the Karp-Sipser core without isolated vertices. 

Let $\bm{A}^{\text{pKS}}_n$ be the adjacency matrix of the pruned Karp-Sipser core and let $\bm{A}^{0/1}_n$ and $\bm{A}^{\text{pKS},0/1}_{n}$ be the matrices obtained from $\bm{A}_n$ and $\bm{A}^{\text{pKS}}_n$ by replacing all nonzero entries by $1$. Then $\bm{A}^{0/1}_n$ and $\bm{A}^{\text{pKS},0/1}_{n}$ are the adjacency matrices of the unweighted version of \cm{n} and its pruned Karp-Sipser core, respectively.

We next claim that
\begin{align}\label{eqw:difference-matrix-peeling}
    \rank_{\FF}\bc{\bm{A}_{n}}-\rank_{\FF}\bc{\bm{A}^{\text{pKS}}_{n}}=\rank_{\RR}\bc{\bm{A}^{0/1}_{n}}-\rank_{\RR}\bc{\bm{A}^{\text{pKS},0/1}_{n}}.
\end{align}
Indeed, \cref{eqw:difference-matrix-peeling} holds true because in each step of the peeling process: 
\begin{itemize}
    \item If a vertex of degree $0$ is removed, then both the ranks of the adjacency matrix of the weighted and unweighted graph will stay unchanged, since we delete a zero row and a zero column;

    \item If a vertex of degree $1$ and its unique neighbor are removed, then the rank of the adjacency matrix of the weighted and unweighted graph decrease by $2$ (see \cite{bauer2001exactly} and \Cref{f1}).
\end{itemize}

We next show that the size of the pruned Karp-Sipser core is $\ooone{n}$. The main idea here is that the local structure around a uniformly chosen vertex in \cm{n} is tree-like: 
Let $\bm{u}_n$ be a uniform random variable on $[n]$. By \cite[Theorem 4.1]{van2022random} and \Cref{assumption_weaker}, $(\fcm{n},\bm{u}_n)$ converges locally weakly to the unimodular branching process $(\bm{G},\bm{o})$ with root offspring distribution $(p_k)_{0\leq k\leq K}$. Let $\mathfrak{T}_m$ be the subset of all rooted trees with depth at most $m$. Denote by $\bm{C}_n(\bm{u}_n)$ the connected component of $\bm{u}_n$ in \cm{n}. If the $(m+1)$- neighborhood of $\bm{u}_n$, i.e., the subgraph induced in \cm{n} by vertices at distance at most $m+1$ from $\bm{u}_n$, is in $\mathfrak{T}_m$, then $\bm{C}_n(\bm{u}_n)$ is contained in the $m$-neighborhood of $\bm{u}_n$ and thus $\bm{C}_n(\bm{u}_n)$ is in $\mathfrak{T}_m$ as well. In other words, the event that $\bm{C}_n(\bm{u}_n)$ is in $\mathfrak{T}_m$ is determined by the local structure of the graph. Then by \cite[Theorem 2.15]{van2022random} and the local convergence,
\begin{align*}
    \PP\bc{(\bm{C}_n(\bm{u}_n),\bm{u}_n)\in\mathfrak{T}_m}=&\PP\bc{\mbox{the $(m+1)$-neighborhood of $\bm{u}_n$ in \cm{n} rooted at $\bm{u}_n$ is in $\mathfrak{T}_m$}}\\
    =&\PP\bc{\mbox{the $(m+1)$-neighborhood of $\bm{o}$ in $\bm{G}$ rooted at $\bm{o}$ is in $\mathfrak{T}_m$}} +o_n(1)\\
    =&\PP\bc{(\bm{G},\bm{o})\in \mathfrak{T}_m}+o_n(1).
\end{align*}
Moreover, by \cite[Theorem 3.1]{van2017random}, in the subcritical and critical case when $\sum_{k\geq 0}k(k-2)p_k\leq 0$ and $p_2\neq 1$, the unimodular branching process $(\bm{G},\bm{o})$ with root offspring distribution $(p_k)_{0\leq k\leq K}$ dies out almost surely.
Therefore,
\begin{align*}
    \lim_{m\to \infty}\PP\bc{(\bm{G},\bm{o})\in \mathfrak{T}_m}=1.
\end{align*}
As a consequence,
\begin{align*}
    \liminf_{n\to\infty}\PP\bc{\mbox{$\bm{C}_n(\bm{u}_n)$ is a finite} tree}
    \geq \limsup_{m\to \infty}\liminf_{n\to\infty}\PP\bc{(\bm{C}_n(\bm{u}_n),\bm{u}_n)\in\mathfrak{T}_m}=1,
\end{align*}
i.e., $\PP\bc{\mbox{$\bm{C}_n(\bm{u}_n)$ is a finite tree}}=1+o_n(1)$. On the other hand, for each vertex $v\in[n]$, if the connected component in \cm{n} containing $v$ is a tree, then $v$ will be removed in the peeling process. Therefore,
\begin{align*}
    \sum_{v\in[n]}\PP\bc{\mbox{$v$ is in the pruned Karp-Sipser core of \cm{n}}}=&n\PP\bc{\mbox{$\bm{u}_n$ is in the pruned Karp-Sipser core of \cm{n}}}\\
    =&n-n\times\PP\bc{\mbox{$\bm{C}_n(\bm{u}_n)$ is a tree}}=o_n(n).
\end{align*}
As a consequence, the expected number of vertices in the pruned Karp-Sipser core, which is equal to  the number of rows of $\bm{A}^{\text{pKS},0/1}_n$ (or $\bm{A}^{\text{pKS}}_n$), is upper bounded by $o_n(n)$ as well. Hence, \eqref{eqw:difference-matrix-peeling} yields that
\begin{align}\label{eq_uwa_1}
    \rank_{\FF}\bc{\bm{A}_{n}}=\rank_{\RR}\bc{\bm{A}^{0/1}_{n}}-\rank_{\RR}\bc{\bm{A}^{\text{pKS},0/1}_{n}}+\rank_{\FF}\bc{\bm{A}^{\text{pKS}}_{n}}=\rank_{\RR}\bc{\bm{A}^{0/1}_{n}}+\ooone{n}.
\end{align}

Finally, by \cite[Theorem 2]{bordenave2011rank}, \cite[Theorem 4]{abert2013benjamini} and the local convergence of $(\fcm{n},\bm{u}_n)$,
\begin{align}\label{eq_uwa_2}
\lim_{n\to\infty}\PP\bc{\abs{\frac{1}{n}\rank_{\RR}\bc{\bm{A}^{0/1}_{n}}-\min_{\alpha\in [0,1]}R_\psi(\alpha)}\geq \varepsilon/2}=0.
\end{align}
The desired result then follows directly from the combination of \eqref{eq_uwa_1}, \eqref{eq_uwa_2} and Markov's inequality.

\subsection{Proof of \Cref{coro_main_prop_2}}
In the special case $p_2=1$, we remove a small proportion of edges uniformly at random so that the degree distribution becomes subcritical and \Cref{coro_main_prop_1} applies. More specifically, if we remove a small proportion $\eta$ of the edges, the resulting degree distribution satisfies the conditions of \Cref{coro_main_prop_1}, while the normalized rank of the adjacency matrix will only change marginally.

So let $p_2=1$ and fix $\eta>0$. We consider the following graph $\fcm{n}^\eta$ constructed from \cm{n}: Given \cm{n}, remove each edge independently with probability $\eta$. Denote by $\bm{d}^\eta(n,i)$ the degree of vertex $i$ in $\fcm{n}^\eta$. \footnote{Indeed, given the number of edges to be removed from the graph $q$, $\fcm{n}^\eta$ can be constructed as follows: we first pair $2q$ half-edges uniformly at random and remove them. Each such removal yields a realization of the degree sequence $\bm{d}^\eta(n, i)$. Then, pair the remaining half-edges uniformly at random. As a consequence, $\fcm{n}^\eta$ is a configuration model.
}  Then $\fcm{n}^\eta$ is a still a configuration model with a random degree sequence $(\bm{d}^\eta(n,i))_{n\geq 1,i\in[n]}$ satisfying \Cref{assumption_weaker} with  probability distribution $(p_k^\eta)_{k\geq 0}$ defined as 
\begin{align*}
    p_k^\eta=\begin{cases}
       \eta^2,\quad& k=0;\\
       2\eta(1-\eta),\quad& k=1;\\
       (1-\eta)^2,\quad& k=2;\\
        0,\quad& k\geq 3.
    \end{cases}
\end{align*}
Moreover, 
\begin{align*}
    \sum_{k\geq 0}k(k-2)p_k^\eta< 0 \quad\mbox{and}\quad p_2^\eta\neq 1.
\end{align*}
Let $\bm{A}_n^\eta$ denote the adjacency matrix of $\fcm{n}^\eta$ and $\psi^\eta(\alpha)=\sum_{k\geq 0}p_k^\eta \alpha^k$ the p.g.f. of $(p_k^\eta)_{k\geq 0}$. By \Cref{coro_main_prop_1}, 
\begin{equation}\label{r_p2_1}
\lim_{n\to\infty}\sup_{J_n\in \syn}\PP\bc{\abs{\frac{1}{n}\rank_{\FF}\bc{\bm{A}^\eta_{n}}-\min_{\alpha\in [0,1]}R_{\psi^\eta}(\alpha)}\geq \varepsilon/2}=0.
\end{equation}
Moreover, since $R_{\psi^\eta}(\alpha)$ is a polynomial of $(\alpha,\eta)$, $R_{\psi^\eta}(\alpha)$ is uniformly continuous with respect to $\eta$. It follows that 
\begin{align}\label{r_p2_2}
    \lim_{\eta\to 0}\abs{\min_{\alpha\in [0,1]}R_{\psi^\eta}(\alpha)-\min_{\alpha\in [0,1]}R_\psi(\alpha)}=0.
\end{align}
Finally, by \Cref{remark_basic_remove_add}, the rank difference between $\bm{A}_n^\eta$ and $\bm{A}_n$ is upper bounded by $2$ times the number of removed edges, whose expectation divided by $n$ converges to $2\eta$ as $n$ tends to infinity. By Markov's inequality, for any $\varepsilon>0$,
\begin{align}\label{r_p2_3}
    \lim_{\eta\to 0}\lim_{n\to\infty}\sup_{J_n\in \syn}\PP\bc{\abs{\frac{1}{n}\rank_{\FF}\bc{\bm{A}^\eta_{n}}-\frac{1}{n}\rank_{\FF}\bc{\bm{A}_{n}}}\geq \varepsilon/4}=0.
\end{align}
Combining \eqref{r_p2_1}, \eqref{r_p2_2} and \eqref{r_p2_3} gives the desired result.

\subsection{Proof of \Cref{coro_main_prop_3}}\label{subsec_subcritical}
As mentioned earlier, our proof of \Cref{coro_main_prop_3} splits into an upper bound that is independent of \Cref{tech_assumption}, and a lower bound that heavily depends on said assumption.  
Indeed, while there is a way to prove the upper bound using \Cref{tech_assumption} (see \Cref{sec-tech_assumption}), by adopting the approach of~\cite{Bypassing2013}, we can in fact remove the reliance on this assumption as follows:
\begin{proposition}
    \label{coro_upper}
Assume that the degree sequence $\vd$ satisfies \Cref{assumption_weaker}. For any field $\FF$, 
\begin{align*}
\lim_{n\to\infty}\sup_{J_n\in \syn}\PP\bc{\frac{1}{n}\rank_{\FF}\bc{\bm{A}_{n}}\leq \min_{\alpha\in [0,1]}R_\psi(\alpha)- \varepsilon}=0.
\end{align*}
\end{proposition}
The proof of \Cref{coro_upper} is deferred to \Cref{ap_1}. The main work of this section will go into the proof of the following lower bound:
\begin{proposition}\label{coro_lower}
    Assume that the degree sequence $\vd$ satisfies \Cref{assumption_weaker} with a probability  distribution $(p_k)_{k\geq 0}$  satisfying \Cref{tech_assumption} and $\sum_{k\geq 0}k(k-2)p_k>0$. For any field $\FF$,
\begin{align*}
\liminf_{n\to\infty}\inf_{J_n\in \syn}\Erw\brk{\frac{1}{n}\rank_{\FF}\bc{\bm{A}_{n}}}\geq \min_{\alpha\in [0,1]}R_\psi(\alpha).
\end{align*}
\end{proposition}
\begin{proof}[Proof of \Cref{coro_main_prop_3}  subject to \Cref{coro_upper,coro_lower}]
    The proof is identical to the proof of \cite[Theorem 1.2]{HofMul25} subject to \cite[Theorems 2.1 and 2.2]{HofMul25}.
\end{proof}
Hence, we are left with the proof of \Cref{coro_lower}.

\subsubsection{{Proof of \Cref{coro_lower} via coupling}}\label{sec_proof9.1}
With \Cref{t_main} in hand, 
it might appear that \Cref{coro_lower} follows directly 
through a truncation argument. Unfortunately, truncating the tail of the probability distribution and re-introducing a finite upper bound on the degrees might have the effect that the second derivative of the p.g.f. of the resulting probability distribution is no longer log-concave -- even if the original distribution satisfied \Cref{tech_assumption}. 
Thus, in the following proof, to preserve log-concavity the truncation, we will follow a slightly different route and modify the probability distribution in advance.

We also utilize the following result that derives uniform convergence from pointwise convergence:
\begin{restatable}{lemma}{lembasic}
\label{lem_basic}
Let $(\bm{a}_{n,k})_{n,k\geq 0}$ be an array of non-negative, integrable random variables and $(b_k)_{k\geq 0}$ be a sequence of non-negative numbers with $\sum_{k\geq 0}b_k<\infty$. Assume that
    \begin{enumerate}
        \item  For all $k\geq 0$: $\lim_{n\to\infty}\Erw\brk{\abs{\bm{a}_{n,k} - b_k}}$; 
        \item $\lim_{n\to\infty}\Erw\brk{\sum_{k\geq 0} \bm{a}_{n,k}}=\sum_{k\geq 0}b_k$.
    \end{enumerate}
Then $\lim_{n\to\infty}\Erw\brk{\sum_{k\geq 0} \abs{\bm{a}_{n,k}-b_k}}=0$.
\end{restatable}
The proof of \Cref{lem_basic} is given in Appendix \ref{sec-lem_basic}. 
\begin{proof}[Proof of \Cref{coro_lower}]
Let $(p_k)_{k\geq 0}$ be a probability distribution, potentially with infinite support, that satisfies $\sum_{k\geq 0}k(k-2)p_k>0$,  \Cref{assumption_weaker} and \Cref{tech_assumption}. Let $\psi$ denote its p.g.f..

\paragraph{Construction of limiting probability distribution $(p_k(r,\ell))_{k\geq 0}$.} 
For each $r\in [0,1)$, 
we first construct a probability distribution $(p_k(r))_{k\geq 1}$ from $(p_k)_{k\geq 0}$ 
that has the property that the second derivative of its probability generating function is \textit{strictly} log-concave for $r>0$, but that is not yet necessarily finitely supported or satisfies \Cref{assumption_stronger}.

Observe that because $(p_k)_{k \geq 0}$ is supercritical, there exists $k\geq 3$ such that $p_k>0$. For $k\geq 0$, set
\begin{align}\label{eq_def_pkr}
    p_k(r):=b(r)^{-1}\sum_{i\geq k}\binom{i}{k}\eul^{-ir}\bc{p_i+r\teo{i\geq 2}\frac{i-2}{i}p_{i-1}}r^{i-k}(1-r)^k,
\end{align}
 where
\begin{align*}
    b(r):=\sum_{k\geq 0}\eul^{-kr}\bc{p_k+r\teo{k\geq 2}\frac{k-2}{k}p_{k-1}}.
\end{align*}
In particular, $b(r)\in [\sigma(r),1+r]$, where $\sigma(r)=\sum_{k\geq 0}\eul^{-kr}p_k$ was defined in \cref{def_sigma_lambda}. It is straightforward to check that $\lim_{r\downarrow 0}b(r)=1$ and $\lim_{r\downarrow 0}p_k(r)=p_k$ for every $k$. Next, let $\varphi_r$ denote the p.g.f. of $(p_k(r))_{k \geq 0}$, such that  
\begin{align*}
\varphi_r(\alpha):=\sum_{k\geq 0}p_{k}(r)\alpha^k=b(r)^{-1}\sum_{k\geq 0}\eul^{-kr}\bc{p_k+r\teo{k\geq 2}\frac{k-2}{k}p_{k-1}}(r+(1-r)\alpha)^k, \qquad \alpha \in [0,1].
\end{align*}
In particular, $\varphi_0(\alpha)=\psi(\alpha)$, 
and by \Cref{lem_basic},
\begin{align}\label{eq_dif_phir}
    \lim_{r\downarrow 0}\sup_{\alpha\in[0,1]}\abs{\varphi_r(\alpha)-\psi(\alpha)}\leq \lim_{r\downarrow 0}\sum_{k\geq 0}\abs{p_k(r)-p_k}=0.
\end{align}
On the other hand,
\begin{align*}
    \sum_{k\geq 0}kp_k(r)=&b(r)^{-1}\sum_{k\geq 0}\sum_{i\geq k}k\binom{i}{k}\eul^{-ir}\bc{p_i+r\teo{i\geq 2}\frac{i-2}{i}p_{i-1}}r^{i-k}(1-r)^k\\
    =&b(r)^{-1}\sum_{i\geq 0}\sum_{k=1}^i i\binom{i-1}{k-1}\eul^{-ir}\bc{p_i+r\teo{i\geq 2}\frac{i-2}{i}p_{i-1}}r^{i-k}(1-r)^k\\
    =&b(r)^{-1}\sum_{i\geq 0} i\eul^{-ir}\bc{p_i+r\teo{i\geq 2}\frac{i-2}{i}p_{i-1}}(1-r),
\end{align*}
from which it becomes apparent that $\lim_{r\downarrow 0}\sum_{k\geq 0}kp_k(r)=\sum_{k\geq 0} kp_k$. Then, once more by \Cref{lem_basic}, 
\begin{align}\label{eq_kpkpkr}
\lim_{r\downarrow 0}\sum_{k\geq 0}\abs{kp_k(r)- kp_k}=0.    
\end{align}
As a consequence, 
\begin{align}\label{eq_dif_phird}
    \lim_{r\downarrow 0}\sup_{\alpha\in[0,1]}\abs{\varphi_r'(\alpha)-\psi'(\alpha)}\leq \lim_{r\downarrow 0}\sum_{k\geq 0}k\abs{p_k(r)-p_k}=0.
\end{align}
Next, observe that the second derivative of $\varphi_r$ is strictly log-concave  precisely when $(\ln(\varphi''_r(\alpha)))''< 0$ for all $\alpha\in[0,1]$.  For $r\in(0,1)$, 
\begin{align*}
\varphi''_r(\alpha)=&b(r)^{-1}\eul^{-2r}(1-r)^2\Bigg(\sum_{k\geq 0}k(k-1) p_k \bc{\eul^{-r}(r+(1-r)\alpha)}^{k-2}\\
&+r\sum_{k\geq 2}(k-1)(k-2)p_{k-1}\bc{\eul^{-r}(r+(1-r)\alpha)}^{k-2}\Bigg)\\
=&b(r)^{-1}\eul^{-2r}(1-r)^2(1+r\eul^{-r}(r+(1-r)\alpha))\psi''(\eul^{-r}(r+(1-r)\alpha)).
\end{align*}
Finally, for all $r \in (0,1)$ and $\alpha \in [0,1]$, $\varphi''_r(\alpha) \geq \varphi''_r(0) = 2p_2(r) >0$. Hence, taking the second derivative of $\ln(\varphi''_r(\alpha))$ with respect to $\alpha$ gives that
\begin{align}\label{eq_upper_lnphidd}
  \frac{\dif^2 \ln(\varphi''_r(\alpha))}{\dif \alpha^2}=&\frac{\dif^2\ln(\psi''(\eul^{-r}(r+(1-r)\alpha)))}{\dif \alpha^2}+\frac{\dif^2\ln(1+r\eul^{-r}(r+(1-r)\alpha))}{\dif \alpha^2}\\
  =&\eul^{-2r}(1-r)^2\frac{\dif^2\ln(\psi''(y))}{\dif y^2}\bigg\vert_{y=\eul^{-r}(r+(1-r)\alpha)}-\frac{r^2(1-r)^2\eul^{-2r}}{(1+r\eul^{-r}(r+(1-r)\alpha))^2}\leq -\frac{r^2(1-r)^2\eul^{-2r}}{(1+r\eul^{-r})^2}.\nonumber
\end{align}
Hence, $\varphi''_r$ is strictly log-concave, which finishes our considerations on $(p_k(r))_{k \geq 0}$.

For each $\ell\geq 3$ and $r \in (0,1)$, using  $(p_k(r))_{k \geq 0}$, we next define a finitely supported probability distribution $(p_k(r,\ell))_{k\geq 0}$ that satisfies the third condition of \Cref{assumption_stronger} and still has the property that the second derivative of its p.g.f. is \textit{strictly} log-concave. Its definition depends on the support of $(p_k)_{k\geq 0}$: If $(p_k)_{k\geq 0}$ has finite support, we can take $(p_k(r,\ell))_{k\geq 0} = (p_k(r))_{k \geq 0}$, while if $(p_k)_{k\geq 0}$ has infinite support, additional truncation at $\ell$ becomes necessary.
\begin{enumerate}
    \item \textbf{$(p_k)_{k\geq 0}$ is finitely supported:} In this case, 
    we set $p_k(r,\ell):=p_k(r)$ for $k\geq 0$, and check that this choice of $(p_k(r,\ell))_{k \geq 0}$ is finitely supported and satisfies condition 3. of  \Cref{assumption_stronger}. Let  $ m:=\sup\cbc{k\geq 0:p_k>0}\geq 3$. Then: 
    \begin{enumerate}
        \item For $k\geq m+2$, $p_k(r)=0$, so $(p_k(r,\ell))_{k \geq 0}$ is finitely supported.

        \item For $0\leq k\leq m$, 
    \begin{align} \label{pkr_pos}
        p_k(r)\geq b(r)^{-1}\binom{m}{k}\eul^{-mr}p_m r^{m-k}(1-r)^k>0,
    \end{align}
    while
 \begin{align*}
        p_{m+1}(r)\geq \eul^{-(m+1)r}r^2(1-r)^{m}\frac{m-1}{m+1}p_m>0.
    \end{align*}
    \end{enumerate}

    \item \textbf{$(p_k)_{k\geq 0}$ is infinitely supported:} 
    In this case, we set
    \begin{align*}
    p_k(r,\ell):=\begin{cases}
    p_k(r)+\frac{1}{1+\ell}\sum_{j\geq \ell+1} p_j(r),\quad&0\leq k\leq \ell;\\
    0,\quad&k\geq \ell+1.
\end{cases}
\end{align*}
Hence, $(p_k(r,\ell))_{k\geq 0}$ is supported on the finite set $\cbc{0}\cup[\ell]$. As $\sum_{j\geq \ell+1} p_j(r) >0$ for each $\ell \geq 3$, also this choice of $(p_k(r,\ell))_{k \geq 0}$ satisfies condition 3. of  \Cref{assumption_stronger}.
\end{enumerate}

Next, let $\varphi_{r,\ell}: [0,1] \to [0,1]$,
$$\varphi_{r,\ell}(\alpha):=\sum_{k\geq 0}p_k(r,\ell)\alpha^k,$$
denote the p.g.f. of $(p_k(r,\ell))_{k \geq 0}$.  Our final task is to prove that $\varphi_{r,\ell}''$ is strictly log-concave for $r\in (0,1)$ and $\ell$ large enough. Observe that if $(p_k)_{k\geq 0}$ is finitely supported, then $\varphi_{r,\ell} = \varphi_r$ for all $r \in (0,1)$, $\ell \geq 3$, and log-concavity was shown in \cref{eq_upper_lnphidd}. We thus assume that $(p_k)_{k\geq 0}$ is infinitely supported in the following argument. 

For $u \in \NN_0$ and a $u$-times differentiable function $f:[0,1] \to \RR$, denote by  $f^{(u)}$ the $u$th  derivative of $f$. We then have the following estimate.

\begin{claim}\label{claim_ell_r}
For all $r \in (0,1)$ and $u \in \NN_0$,
\begin{align}\label{eq_lim_phiurl}
    \lim_{\ell\to\infty}\sup_{\alpha\in[0,1]}\abs{\varphi^{(u)}_{r,\ell}(\alpha)-\varphi_r^{(u)}(\alpha)}=0.
\end{align}
\end{claim}

\begin{proof}[Proof of \Cref{claim_ell_r}]
First note that
\begin{align}\label{eq_rev_1}
    \abs{\varphi^{(u)}_{r,\ell}(\alpha)-\varphi_r^{(u)}(\alpha)}=&\abs{\sum_{k\geq u}\bc{\prod_{i=0}^{u-1}(k-i)}(p_k(r,\ell)-p_k(r))\alpha^{k-u}}
    \leq \sum_{k\geq 0} k^u \abs{p_k(r,\ell)-p_k(r)}.
\end{align}
 We proceed to further bound the right hand side of \cref{eq_rev_1}:
    \begin{align*}
        \sum_{k\geq 0} k^u \abs{p_k(r,\ell)-p_k(r)}=\sum_{k=0}^\ell \frac{k^u}{1+\ell}\sum_{j\geq \ell+1} p_j(r)+\sum_{k\geq \ell+1}k^u p_k(r)\leq 2\sum_{k\geq \ell+1}k^u p_k(r).
    \end{align*}
    On the other hand, 
    \begin{align*}
       \sum_{k\geq 0}k^u p_k(r)=&\sum_{k\geq 0}\sum_{i\geq k} k^u b(r)^{-1}\binom{i}{k}\eul^{-ir}\bc{p_i+r\teo{i\geq 2}\frac{i-2}{i}p_{i-1}}r^{i-k}(1-r)^k\\
       \leq&\sum_{i\geq 0}\sum_{k=0}^i i^u b(r)^{-1}\binom{i}{k}\eul^{-ir}\bc{p_i+r\teo{i\geq 2}\frac{i-2}{i}p_{i-1}}r^{i-k}(1-r)^k\\
       =&\sum_{i\geq 0} i^u b(r)^{-1}\eul^{-ir}\bc{p_i+r\teo{i\geq 2}\frac{i-2}{i}p_{i-1}}\leq b(r)^{-1}(1+r)\sum_{i\geq 0} i^u \eul^{-ir}.
    \end{align*}
    Since $\sum_{i\geq 0} i^u \eul^{-ir}<\infty$ for $r>0$, we conclude that $\sum_{k\geq 0}k^u p_k(r) < \infty$ for $r\in(0,1)$, so \Cref{eq_lim_phiurl} holds.
\end{proof}

Using \Cref{claim_ell_r}, we next prove that for any $r \in (0,1)$, there exists $L_r\geq 3$ such that $\varphi''_{r,\ell}$ is log-concave for $\ell>L_r$. Since $(\ln(f^{(2)}))^{(2)}=\frac{f^{(4)}f^{(2)}-(f^{(3)})^2}{(f^{(2)})^2}$, we first argue that $\varphi_{r,\ell}^{(2)}(\alpha)$ is uniformly positive for $r \in (0,1)$. 
Note that $p_k(r,\ell)\geq p_k(r)$ for $0 \leq k \leq \ell$. Let $m'\geq 3$ be an integer such that $p_{m'}(r)>0$. Then as in  \cref{pkr_pos}, 
\begin{align*}
  p_2(r,\ell)\geq p_2(r)\geq b(r)^{-1}\binom{m'}{r}\eul^{-m' r}p_{m'} r^{m'-2}(1-r)^{2}>0.  
\end{align*}
Hence, $\varphi''_{r,\ell}(\alpha)\geq \varphi''_{r,\ell}(0)=2 p_{2}(r,\ell)$ is uniformly positive, and 
by \Cref{claim_ell_r},
\begin{align*}
    \lim_{\ell\to \infty}\sup_{\alpha\in[0,1]}\abs{(\ln(\varphi''_{r,\ell}(\alpha)))''-(\ln(\varphi''_r(\alpha)))''}=0.
\end{align*}
Hence, by \Cref{eq_upper_lnphidd}, for each $r>0$, there exists $L_r\geq 3$ such that for $\ell > L_r$,
\begin{align*}
    \sup_{\alpha\in[0,1]}(\ln(\varphi''_{r,\ell}(\alpha)))''\leq -\frac{r^2(1-r)^2\eul^{-2r}}{2(1+r\eul^{-r})^2} <0.
\end{align*}

\paragraph{Uniform approximation of the rank function.}
For  $\phi:[0,1]\mapsto[0,1]$ differentiable, recall the rank function $R_{\phi}$ defined in \cref{def_rd}. 
We next show that also the rank function of $\psi$ is uniformly approximated by the rank function of $\varphi_{r,\ell}$:

\begin{claim}\label{claim_rdrl_rd}
\begin{align}\label{eq_rdrl_rd}
  \lim_{r\downarrow 0}  \limsup_{\ell\to\infty} \sup_{\alpha\in[0,1]}\abs{R_{\varphi_{r,\ell}}(\alpha)-R_{\psi}(\alpha)}=0.
\end{align}
\end{claim}

\begin{proof}[Proof of \Cref{claim_rdrl_rd}]
By (\ref{eq_dif_phir}), (\ref{eq_dif_phird}) and \Cref{claim_ell_r},
\begin{align}\label{eq_phirl_psi_de}
     \lim_{r\downarrow 0} \limsup_{\ell\to\infty}\sup_{\alpha\in[0,1]}\abs{\varphi_{r,\ell}(\alpha)-\psi(\alpha)}=0\quad\mbox{and}\quad   \lim_{r\downarrow 0} \limsup_{\ell\to\infty}\sup_{\alpha\in[0,1]}\abs{\varphi'_{r,\ell}(\alpha)-\psi'(\alpha)}=0.
\end{align}
Note that $1-\frac{\varphi_{r,\ell}'(\alpha)}{\varphi_{r,\ell}'(1)},1-\frac{\psi'(\alpha)}{\psi'(1)}\in[0,1]$. Again by (\ref{eq_dif_phir}) and \Cref{claim_ell_r},
\begin{align}\label{eq_phi_psi1}
     \lim_{r\downarrow 0} \limsup_{\ell\to\infty}\sup_{\alpha\in[0,1]}\abs{\varphi_{r,\ell}\bc{1-\frac{\varphi_{r,\ell}'(\alpha)}{\varphi_{r,\ell}'(1)}}-\psi\bc{1-\frac{\varphi_{r,\ell}'(\alpha)}{\varphi_{r,\ell}'(1)}}}=0,
\end{align}
and because $\psi$ is  uniformly continuous on $[0,1]$,
\begin{align}\label{eq_phi_psi2}
     \lim_{r\downarrow 0} \limsup_{\ell\to\infty}\sup_{\alpha\in[0,1]}\abs{\psi\bc{1-\frac{\varphi_{r,\ell}'(\alpha)}{\varphi_{r,\ell}'(1)}}-\psi\bc{1-\frac{\psi'(\alpha)}{\psi'(1)}}}=0.
\end{align}
The combination of \cref{eq_phirl_psi_de}-\cref{eq_phi_psi2} gives \cref{eq_rdrl_rd}.
    
\end{proof}

\paragraph{Construction of a deterministic degree sequence with limiting distribution $(p_k(r,\ell))_{k\geq 0}$.} Given $r \in (0,1)$ and $\ell \geq 3$, let $m_{r,\ell}:=\sup\cbc{k\geq 0:p_k(r,\ell)>0}\in \NN_{\geq 3}$.  
For $k\geq 0$, let 
\begin{align*}
    n_{r,\ell,k}'=\lfloor (n-1)\sum_{j=0}^k p_j(r,\ell)\rfloor-\lfloor (n-1)\sum_{j=0}^{k-1} p_j(r,\ell)\rfloor.
\end{align*}
Then $\sum_{k\geq 0} n_{r,\ell,k}'=n-1$. To make the degree-sum even, we set
\begin{align}\label{def_nrkell}
    n_{r,\ell,k}=n_{r,\ell,k}'+\teo{k=0,\sum_{j=0}^{m_{r,\ell}}j n_{r,\ell,j}'\mbox{ is even}}+\teo{k=1,\sum_{j=0}^{m_{r,\ell}}j n_{r,\ell,j}'\mbox{ is odd}}.
\end{align}
With this choice, $\sum_{k=0}^{m_{r,\ell}} n_{r,\ell,k}=n$, $\sum_{k=0}^{m_{r,\ell}}k n_{r,\ell,k}$ is even and, for each
$k\geq 0$,
\begin{align}\label{ineq_dkrl_pkrl}
   \abs{\frac{ n_{r,\ell,k}}{n}-p_k(r,\ell)}\leq \frac{3}{n}.
\end{align}
Finally, for $i \in [n]$ and $k \in \{0, \ldots, m_{r,\ell}\}$, let $d_{r,\ell}(n,i)=k$ if and only if $1+\sum_{j=0}^{k-1}n_{r,\ell,j} \leq i \leq \sum_{j=0}^{k}n_{r,\ell,j}$. 
Then for all $\ell>L_r$, the deterministic degree sequence $\cbc{d_{r,\ell}(n,i)}_{i\in[n]}$ satisfies Assumptions \ref{assumption_weaker} and \ref{assumption_stronger} with a probability  distribution $(p_{k}(r,\ell))_{k\geq 0}$  satisfying \Cref{tech_assumption}. 

Finally, let \cm{n,r,\ell} denote the configuration model constructed from the degree sequence $\cbc{d_{r,\ell}(n,i)}_{i\in[n]}$.  
For any symmetric matrix $J_{n}\in \syn$,
we may then define the weighted adjacency matrix 
\begin{equation}\label{ec1rl}
  \bm{A}_{n,r,\ell}(i,j):=\begin{cases}
      \teo{{\mbox{there is at least an edge between $i$ and $j$ in \cm{n,r,\ell}}} }J_{n}(i,j),\quad &i\neq j; \\
      0,&i=j.
  \end{cases}
\end{equation}
Depending on the sign of $\sum_{k\geq 0}k(k-2)p_k(r,\ell)$, as $p_2(r,\ell) \not=0$, by \Cref{coro_main_prop_1} or \Cref{t_main}, 
\begin{equation}\label{e01arl}
\liminf_{n\to\infty}\inf_{J_n\in \syn}\Erw\brk{\frac{1}{n}\rank_{\FF}\bc{\bm{A}_{n,r,\ell}}}\geq\min_{\alpha\in [0,1]}R_{\varphi_{r,\ell}}(\alpha).
\end{equation} 
We remark that even though \eqref{e01arl} was derived for models with fixed vertex degrees in the case $\sum_{k\geq 0}k(k-2)p_k(r,\ell)>0$, $\inf_{J_n\in \syn}\Erw\brk{\rank_{\FF}\bc{\bm{A}_{n,r,\ell}}}$ is invariant under vertex relabeling. In particular, \cref{e01arl} continues to hold true. In the following coupling argument, we may thus randomly chose $n_{r,\ell,k}$-subsets of the vertex set $[n]$ to have degree $k$.

\paragraph{Coupling the original and the truncated graph.} 
We next construct a coupling between the original configuration model \cm{n} with degree sequence $(\vd_i)_{i \geq 0}$ and limiting degree distribution $(p_k)_{k \geq 0}$ and a configuration model with $n_{r,\ell,k}$ vertices of degree $k \geq 0$, where $n_{r,\ell,k}$ was defined in \cref{def_nrkell}.  
Recall that we use $\mathcal N_k:=\cbc{i\in [n]:\vd_i=k}$ to denote the (random) set of vertices of degree $k$ in \cm{n}, and $n_k=|\mathcal{N}_k|$ to denote the cardinality of $\mathcal N_k$.

Let $n_k^{\wedge}:=\min\cbc{ n_k ,n_{r,\ell,k}}$ for $k\geq 0$. Given $(\vd_i)_{i \geq 0}$, for each $k\geq 0$, we then choose an arbitrary subset of $\mathcal N_k$ of size $n_k^{\wedge}$. This subset will be denoted by $\mathcal{N}_k^{\wedge}$.
In the first step, we construct a graph $\fcm{n}^{\wedge}$ with vertex set $\cup_{k\geq 0}\mathcal{N}_k^{\wedge}$. To do so, for each $k \geq 0$, assign $k$ half-edges to each vertex in $\mathcal{N}_k^{\wedge}$. We then pair the half-edges in a uniformly random order, each time choosing a uniformly random neighbor among the unpaired half-edges, 
until no unpaired half-edge is left if $\sum_{k\geq  0} kn_k^{\wedge}$ is even or there is only one unpaired half-edge left if $\sum_{k\geq  0} kn_k^{\wedge}$ is odd. 
For a given matrix $J_n \in \FF^{n \times n}$, let $\bm{A}^\wedge_n\in\FF^{n\times n}$ denote the $J_n$-weighted adjacency matrix of $\fcm{n}^{\wedge}$:
\begin{equation}\label{ec1r2}
  \bm{A}^\wedge_{n}(i,j):=\begin{cases}
      \teo{{\mbox{there is at least an edge between $i$ and $j$ in $\fcm{n}^{\wedge}$}} }J_n(i,j),\quad &i\neq j; \\
      0,&i=j.
  \end{cases}
\end{equation}

We next construct \cm{n} from $\fcm{n}^{\wedge}$ by adding the `missing' vertices
$[n]\backslash\cup_{k\geq 0}\mathcal{N}_k^{\wedge}$ one by one according the the following procedure, e
until all $\sum_{k\geq 0}\bc{ n_k -n_k^{\wedge}}$ vertices have been added.  
For each $k\geq 0$, suppose that $v\in \mathcal N_k\backslash \mathcal{N}_k^{\wedge}$ is about to be added to the graph. We then, one after the other, choose the neighbors of its $k$ half-edges in the following way:

\begin{enumerate}
    \item\label{it1_pairing} 
    Choose an unpaired half-edge $\bm h$ of $v$ uniformly at random. If there is no unpaired half-edge adjacent to $v$, 
    the process ends.

    \item Choose a half-edge $\bm h'$ uniformly at random among all half-edges in the graph, excluding the half-edges of $v$ that were paired in the previous steps and their half-edge neighbors.
 
    \begin{itemize}
        \item If $\bm h'$ is an unpaired half-edge of $v$, connect $\bm h$ and $\bm h'$.
        \item If $\bm h'$ is not an unpaired half-edge of $v$,
        \begin{itemize}
            \item and there is an unpaired half-edge $h''$ of a vertex other than $v$ in the graph, connect $\bm h$ and $h''$;
            \item and all the unpaired half-edges in the graph belong to $v$, break the connection between $\bm h'$ and its current neighbor and connect $\bm h$ and $\bm h'$. This renders the previous neighbor of $\bm h'$ unpaired.
        \end{itemize}
    \end{itemize} 
    Repeat from \ref{it1_pairing}.
\end{enumerate}
The above construction ensures that \cm{n} is a configuration model, since half-edges are paired uniformly at random. Indeed, our construction ensures that before a new vertex is added, if there is an unpaired half-edge in the graph, it is chosen uniformly at random. Hence, when a new vertex $v$ is added, each half-edge, whether already in the graph or attached to $v$, has an equal probability of being paired with the currently unpaired half-edge of $v$. Finally, completely analogously, we construct \cm{n,r,\ell} from $\fcm{n}^{\wedge}$ by successively adding the missing vertices and half-edges. Moreover, we use the same edge weight matrix $J_n$ in all three graphs \cm{n}, \cm{n,r,\ell} and $\fcm{n}^{\wedge}$.

We now upper bound the number of edges that were in changed in the processes of constructing \cm{n} from $\fcm{n}^{\wedge}$ and \cm{n,r,\ell} from $\fcm{n}^{\wedge}$, respectively. Observe that each time a half-edge of an additional vertex is paired, at most one previous edge is broken and a new one created, so that in the end, the number of differing edges in $\fcm{n}^{\wedge}$ and \cm{n} is bounded from above by $2\sum_{k\geq 0}k\bc{ n_k -n_k^{\wedge}}\leq 2\sum_{k\geq 0}k\abs{ n_k -n_{r,\ell,k}}$. Via \Cref{remark_basic_remove_add}, this implies that the expected rank difference of $\bm{A}_{n}$ and $\bm{A}_n^\wedge$ is upper bounded by $\Erw[2\sum_{k\geq 0}k\abs{ n_k -n_{r,\ell,k}}]$. 
Similarly, the expected rank difference between $\bm{A}_{n,r,\ell}$ and $\bm{A}_n^\wedge$ is upper bounded by $\Erw[2\sum_{k\geq 0}k\abs{ n_k -n_{r,\ell,k}}]$. We conclude that 
\begin{align*}
    \Erw\brk{\abs{\rk{\bm A_n} - \rk{\bm A_{n,r,\ell}}} } \leq \Erw\brk{4\sum_{k\geq 0}k\abs{ n_k -n_{r,\ell,k}}}.
\end{align*}
On the other hand, by \cref{ineq_dkrl_pkrl}, \Cref{lem_basic} and \Cref{assumption_weaker},
\begin{align*}
    \Erw\brk{\sum_{k\geq 0}k\abs{\frac{n_{r,\ell,k}}{n}-\frac{ n_k }{n}}} & \leq\sum_{k=0 }^{m_{r,\ell}}k\abs{\frac{n_{r,\ell,k}}{n}-p_k(r,\ell)}+\Erw\brk{\sum_{k\geq 0}k\abs{\frac{ n_k }{n}-p_k}}+\sum_{k\geq 0}k\abs{p_k-p_k(r,\ell)}\\
    &\leq\frac{3m_{r,\ell}(m_{r,\ell}+1)}{n}+\sum_{k\geq  0}k\abs{p_k-p_k(r,\ell)}+o_n(1).
\end{align*}
Hence, the combination of \Cref{eq_kpkpkr}, the proof of \Cref{claim_ell_r} and \Cref{eq_rev_1} yields that
\begin{align*}
    \lim_{r\downarrow 0}\limsup_{\ell\to\infty}\lim_{n\to\infty}\Erw\brk{\sum_{k\geq 0}k\abs{\frac{n_{r,\ell,k}}{n}-\frac{ n_k }{n}}}=0.
\end{align*}
Therefore, 
\begin{align}\label{e01arl2}
&\lim_{r\downarrow 0}\limsup_{\ell\to\infty}\lim_{n\to\infty}\sup_{J_n\in \syn}\Erw\brk{\abs{\frac{\rank_{\FF}\bc{\bm{A}_{n,r,\ell}}}{n}-\frac{\rank_{\FF}\bc{\bm{A}_{n}} }{n}}}\\
&\qquad\leq  \lim_{r\downarrow 0}\limsup_{\ell\to\infty}\lim_{n\to\infty}\Erw\brk{4\sum_{k\geq 0}k\abs{\frac{n_{r,\ell,k}}{n}-\frac{ n_k }{n}}}=0.\nonumber
\end{align}
On the other hand, for any $r\in (0,1)$ and $\ell \geq 3$,
\begin{align}\label{eq_the_end}
&\liminf_{n\to\infty}\bc{\inf_{J_{n}\in \syn}\Erw\brk{\frac{\rank_{\FF}\bc{\bm{A}_{n}}}{n}} - \min_{\alpha\in [0,1]}R_{\psi}(\alpha)} \nonumber\\
\geq &  \liminf_{n\to\infty}\bc{\inf_{J_{n}\in \syn}\Erw\brk{\frac{1}{n}\rank_{\FF}\bc{\bm{A}_{n,r,\ell}}}-\min_{\alpha\in [0,1]}R_{\varphi_{r,\ell}}(\alpha)} \nonumber \\
& - \limsup_{n\to\infty}\bc{\sup_{J_{n}\in \syn}\Erw\brk{\abs{\frac{\rank_{\FF}\bc{\bm{A}_{n,r,\ell}}}{n}-\frac{\rank_{\FF}\bc{\bm{A}_{n}} }{n}}} } - \abs{\min_{\alpha\in [0,1]}R_{\psi}(\alpha) - \min_{\alpha\in [0,1]}R_{\varphi_{r,\ell}}(\alpha)}. 
\end{align} 
Combining  \cref{eq_rdrl_rd}, \cref{e01arl} and \cref{e01arl2} gives the desired lower bound.
\end{proof}

\paragraph{\bf Acknowledgement.}The authors are supported by Netherlands Organisation for Scientific Research
(NWO) through the Gravitation NETWORKS grant 024.002.003.  The work of Haodong Zhu is further supported by the European
Union’s Horizon 2020 research and innovation programme under the Marie
Skłodowska-Curie grant agreement no. 945045.

\bibliographystyle{abbrv}
\DeclareRobustCommand{\VAN}[3]{#3}
\bibliography{references}

\begin{appendices}

\section{From pointwise to uniform convergence}\label{sec-lem_basic}
In this section, we prove \Cref{lem_basic}:
\lembasic*
\begin{proof}
For any positive integer $K$,
    \begin{align*}
        \Erw\brk{\sum_{k\geq 0}\abs{\bm{a}_{n,k}-b_k}}\leq &\Erw\brk{\sum_{k= 0}^K\abs{\bm{a}_{n,k}-b_k}}+\Erw\brk{\sum_{k\geq 0}\bm{a}_{n,k}}-\sum_{k=0}^K \Erw\brk{\bm{a}_{n,k}}+\sum_{k\geq K+1}b_k\\
        \leq &2\sum_{k\geq K+1}b_k+ o_n(1).
    \end{align*}
    We conclude that for any $K$, $\limsup_{n \to \infty}\Erw[\sum_{k\geq 0}\abs{\bm{a}_{n,k}-b_k}]\leq 2\sum_{k\geq K+1}b_k$. As the left-hand side does not depend on $K$, the claim follows.
\end{proof}

\begin{remark}\label{result_p1}
Let $\vd=(\vd_i)_{i\in[n]}$ be a degree sequence and $(p_k)_{k\geq 0}$ be a sequence of deterministic nonnegative numbers that together satisfy \Cref{assumption_weaker}. Then $(p_k)_{k\geq 0}$ is a probability distribution: With the choice $\bm{a}_{n,k} = k N_K/n$, $b_k=kp_k$, all assumptions of \Cref{lem_basic} are satisfied and thus
 \begin{align*}
     \Erw\brk{\sum_{k\geq 0}k \abs{\frac{ n_k }{n}-p_k}}=o_n(1).
 \end{align*}
On the other hand for any $n$, as $\sum_{k\geq 0} n_k =n$, 
 \begin{align*}
     \abs{1-\sum_{k\geq 0}p_k}\leq \Erw\brk{\sum_{k\geq 0} \abs{\frac{ n_k }{n}-p_k}} \leq \Erw\brk{\abs{ \frac{n_0}{n}-p_0}} +\Erw\brk{\sum_{k\geq 0}k \abs{\frac{ n_k }{n}-p_k}}.
 \end{align*}
As the left-hand side does not depend on $n$, and $\Erw\brk{\abs{ \frac{n_0}{n}-p_0}}=o_n(1)$, we conclude that $\sum_{k\geq 0}p_k=1$.
\end{remark}

\section{Continuous-time graph exploration}\label{sec_app_graphexp}
In this section, we summarize how \Cref{lem_janson_s,lem_janson_l_2} follow from the results in \cite{janson2009new}.

\begin{proof}[Proof of \Cref{lem_janson_l_2}]
First observe that for $k=0$, the claim reduces to \Cref{assumption_weaker} (a). For $k \geq 1$, let $\tilde{\bm{V}}_k(t)$ be the number of vertices of degree $k$, all of whose half-edges have lifetimes greater than $t$. By \cite[Lemma 5.2]{janson2009new}, for all $k \geq 0$,
\begin{align}\label{eq_janson_vtilde}
    \sup_{t \geq 0}\abs{n^{-1}\tilde{\bm{V}}_k(t) - p_k\eul^{-kt}} \stackrel{\mathbb{P}}{\longrightarrow} 0.
\end{align}
Let $\tilde{\bm{S}}(t):=\sum_{k=1}^K k\tilde{\bm{V}}_k(t)$. Then by \cite[(5.17)]{janson2009new}, for all $t \geq 0$,
\begin{align}\label{eq_janson_vstilde}
   \abs{\tilde{\bm{V}}_k(t) - \bm{V}_{k}(t)} \leq  k^{-1}\abs{\tilde{\bm{S}}(t) - \bm{S}(t)} \leq \abs{\tilde{\bm{S}}(t) - \bm{S}(t)}.
\end{align}
Moreover, by \cite[(5.9)]{janson2009new},
\begin{align}\label{eq_janson_stilde}
    \sup_{ \in [0, -\ln \xi]}n^{-1}\abs{\tilde{\bm{S}}(t) - \bm{S}(t)} \stackrel{\mathbb{P}}{\longrightarrow} 0.
\end{align}
Combining \eqref{eq_janson_vtilde} - \eqref{eq_janson_stilde}, we obtain that
\begin{align*}
    \sup_{t \in [0, -\ln \xi]}\abs{n^{-1}\bm{V}_{k}(t) - p_k\eul^{-kt}} \leq \sup_{t \geq 0}\abs{n^{-1}\tilde{\bm{V}}_k(t) - p_k\eul^{-kt}} +  \sup_{ \in [0, -\ln \xi]}n^{-1}\abs{\tilde{\bm{S}}(t) - \bm{S}(t)} \stackrel{\mathbb{P}}{\longrightarrow} 0.
\end{align*}
\end{proof}

\begin{proof}[Proof of \Cref{lem_janson_s}]
Let $\eps \in (0, -(\ln \xi)/2)$. For $t \geq 0$, let $\bm{A}(t) = \bm{L}(t) - \bm{S}(t)$ denote the number of active half-edges at time $t$. By \cite[(5.10)]{janson2009new},
\begin{align}
    \sup_{t \in [0, -\ln\xi]} \abs{n^{-1} \bm{A}(t) - H(\eul^{-t})} = 0,
\end{align}
where $H(x) = \lambda(0)x(x - \hat\psi(x))$ in our notation. Moreover, \cite[Lemma 5.5 (i)]{janson2009new} states that under Condition \cref{con_super_c}, $H(x) < 0$ for $x \in (0,\xi)$, from which we infer the existence of a $\delta>0$ such that w.h.p., $\bm{A}(t)/n \geq \delta$ on $[\eps, -\ln \xi - \eps]$. As a new component is started and \textbf{Step 1} is performed only when the number of active  half-edges drops to zero, the claim follows.
\end{proof}

\section{A versatile approximation}\label{app_needproof}
In this section, we prove \Cref{lem_app_divide}:
\appdivide*
\begin{proof}
    By Markov's inequality and \cref{ass_lem_app_divide},
    \begin{align*}
        \mathbb{P}\bc{\ensuremath{\mathds{1}}_{\mathfrak{H}_n}\abs{\bm{X}_n-x n^{\beta}}\geq \sqrt{a_n} n^{\beta}}\leq \sqrt{a_n},
    \end{align*}
and
    \begin{align*}
        \mathbb{P}\bc{\ensuremath{\mathds{1}}_{\mathfrak{H}_n}\abs{\bm{Y}_n-y n^{\beta}}\geq \sqrt{a_n} n^{\beta}}\leq \sqrt{a_n}.
    \end{align*}
Since $\bm{X}_n/\bm{Y}_n\leq C$, 
\begin{align*}
    \mathbb{E}\abs{\frac{\bm{X}_n}{\bm{Y}_n}-\frac{x}{y}}\leq &\mathbb{E}\abs{\ensuremath{\mathds{1}}_{\mathfrak{H}_n}\teo{\abs{\bm{X}_n-x n^{\beta}}\leq \sqrt{a_n} n^{\beta},\abs{\bm{Y}_n-y n^{\beta}}\leq \sqrt{a_n} n^{\beta}}\bc{\frac{\bm{X}_n}{\bm{Y}_n}-\frac{x}{y}}}\\
    &+\bc{C+\frac{x}{y}}\bc{\mathbb{P}\bc{\ensuremath{\mathds{1}}_{\mathfrak{H}_n}\abs{\bm{X}_n-x n^{\beta}}\geq \sqrt{a_n} n^{\beta}}+\mathbb{P}\bc{\ensuremath{\mathds{1}}_{\mathfrak{H}_n}\abs{\bm{Y}_n-y n^{\beta}}\geq \sqrt{a_n} n^{\beta}}+\PP\bc{\mathfrak{H}_n^c}}\\
    \leq &\mathbb{E}\brk{\teo{\abs{\bm{X}_n-x n^{\beta}}\leq \sqrt{a_n} n^{\beta},\abs{\bm{Y}_n-y n^{\beta}}\leq \sqrt{a_n} n^{\beta}}\frac{y\ensuremath{\mathds{1}}_{\mathfrak{H}_n}\abs{\bm{X}_n-x n^{\beta}}+x\ensuremath{\mathds{1}}_{\mathfrak{H}_n}\abs{\bm{Y}_n-y n^{\beta}}}{y\bm{Y}_n}}\\
    &+\bc{C+\frac{x}{y}}(2\sqrt{a_n}+a_n)\\
    \leq & \frac{(x+y) \sqrt{a_n}}{y(y-\sqrt{a_n})} +\bc{C+\frac{x}{y}}(2\sqrt{a_n}+ a_n).
\end{align*}
\end{proof}

\section{Crossing half-edges: The event $\mathfrak G_{\varepsilon,n}$}\label{ap_eventg}
Recall $c=c(\eps)$ from \eqref{def_c} and that $\mathfrak{G}_{\varepsilon,s}$ is the event that the number of $j\in \verst{s}$ such that $\bar{\vd}_{j,s}<d_j$ is larger than $cn$. By using a second-moment method, we can prove that $\mathfrak{G}_{\varepsilon,s}$ occurs with high probability:
\begin{lemma}\label{l_dcs}
Fix $\rangee$. Then, uniformly for $\ranges$, 
\begin{align*}  \PP\bc{\mathfrak{G}_{\varepsilon,s}^c}=\bar{o}_n(1).
\end{align*}
\end{lemma}
\begin{proof}
    Let 
    \begin{align*}
        \bcabc{s}=\sum_{k=0}^K k \verslk{k}{s}-\sum_{k=0}^K k \bverslk{k}{s}
    \end{align*}
   denote the number of half-edges connecting vertices in $\verst{s}$ to vertices in $[n]\setminus \verst{s}$. By the combination of \Cref{l_concentration,lem_bvks} and the dominated convergence theorem,
    \begin{align}\label{eq_est_c}
        \mathbb{E}\abs{n^{-1}\bcabc{s}-\lambda(t_s)+\frac{\lambda^2(t_s)\eul^{2t_s}}{\lambda(0)}}=\bar{o}_n(1).
    \end{align}
Since there are $\bcabc{s}$ edges between $\verst{s}$ and $[n]\backslash \verst{s}$, the number of vertices in $\verst{s}$ which are connected to at least one vertex in $[n]\backslash \verst{s}$ is lower bounded by $\bcabc{s}/K$. As a consequence, the definition of $c$ ensures that
\begin{align*}
\PP\bc{\mathfrak{G}_{\varepsilon,s}^c}\leq \PP\bc{\bcabc{s}\leq Kc n}=\bar{o}_n(1).
\end{align*}
\end{proof}

\section{Proof of \Cref{lem_proal}: Properties of $G_{t}$ and its zeroes}\label{sec_progt}

\lemproal*

\begin{proof}
    \begin{enumerate}
        \item First, observe that by \Cref{assumption_stronger} and \eqref{def_psi_hat}, $\hat{\psi}_t$ is a polynomial of degree $K-1 \geq 2$. Hence, $G_t$ is a polynomial of degree $(K-1)^2\geq 4$. As a consequence, none of the functions $G_t, G'_t, G''_t$ is constant.
        Next, it is straightforward to check that
         \begin{align}\label{eq_psi0t}
             \hat{\psi}_t(\alpha)=\frac{\lambda(0)}{\lambda(t)\eul^t}\hat{\psi}\bc{\eul^{-t}+\frac{\lambda(t)\eul^t}{\lambda(0)}(\alpha-1)}.
             \end{align}
        Let $\beta=\beta_t(\alpha)=\eul^{-t}+\frac{\lambda(t)\eul^{t}}{\lambda(0)}(\alpha-1) \in [0,\eul^{-t}]$. Substituting \cref{eq_psi0t} into the definition of $G_t$ gives that
        \begin{align*}
            G_t(\alpha)=\frac{\lambda(0)}{\lambda(t)\eul^{t}}\bc{\beta-\eul^{-t} + \hat{\psi}\bc{\eul^{-t}-\hat{\psi}(\beta)}}.
        \end{align*}
        Hence,
        \begin{align}\label{eq_dgt}
             G_t'(\alpha)=1-\hat{\psi}'\bc{\eul^{-t}-\hat{\psi}(\beta)}\hat{\psi}'(\beta)
        \end{align}
        and
        \begin{align}\label{eq_ddgt}
             G_t''(\alpha)=&\frac{\lambda(t)\eul^t}{\lambda(0)}\bc{\hat{\psi}''\bc{\eul^{-t}-\hat{\psi}(\beta)}\bc{\hat{\psi}'(\beta)}^2-\hat{\psi}'\bc{\eul^{-t}-\hat{\psi}(\beta)}\hat{\psi}''(\beta)}\\
             =&\frac{\lambda(t)\eul^t}{\lambda(0)}\hat{\psi}'\bc{\eul^{-t}-\hat{\psi}(\beta)}\hat{\psi}'(\beta)\bc{\frac{\hat{\psi}''\bc{\eul^{-t}-\hat{\psi}(\beta)}}{\hat{\psi}'\bc{\eul^{-t}-\hat{\psi}(\beta)}}\hat{\psi}'(\beta)-\frac{\hat{\psi}''(\beta)}{\hat{\psi}'(\beta)}}\nonumber.
        \end{align}
        Thus, by our assumptions on $(p_k)_{k\geq 0}$, $G''_t(\alpha)=0$ only if $$\frac{\hat{\psi}''\bc{\eul^{-t}-\hat{\psi}(\beta)}}{\hat{\psi}'\bc{\eul^{-t}-\hat{\psi}(\beta)}}\hat{\psi}'(\beta)-\frac{\hat{\psi}''(\beta)}{\hat{\psi}'(\beta)}=0.$$ Now, since $\hat{\psi}'$ is a log-concave function, $\hat{\psi}''/\hat{\psi}'$ is a decreasing function. Therefore, $$\alpha\mapsto\frac{\hat{\psi}''\bc{\eul^{-t}-\hat{\psi}(\beta)}}{\hat{\psi}'\bc{\eul^{-t}-\hat{\psi}(\beta)}}\hat{\psi}'(\beta)-\frac{\hat{\psi}''(\beta)}{\hat{\psi}'(\beta)}$$ is an increasing function. Since $G''_t$ is a polynomial and is not a constant, 
we conclude that $G''_t$ has at most one zero. As a consequence, $G'_t$ has at most two zeroes and $G_t$ has at most three zeroes.

        \item It is straightforward to check that
        \begin{align*}
            G_t(1-\hat{\psi}_t(\alpha))=-\hat{\psi}_t(\alpha)+\hat{\psi}_t\bc{\alpha-G_t(\alpha)}.
        \end{align*}
        Hence, if $\alpha$ zeros $G_{t}$, so does $1-\hat{\psi}_{t}(\alpha)$. Specifically, $1-\hat{\psi}_{t}(\alpha^\star(t))$ and $1-\hat{\psi}_{t}(\alpha_\star(t))$ zero $G_{t}$. Since $\alpha_\star(t)$ is the smallest zero of $G_t$, we have $1-\hat{\psi}_{t}(\alpha^\star(t))\geq \alpha_\star(t)$. If $1-\hat{\psi}_{t}(\alpha^\star(t))> \alpha_\star(t)$, then 
        $\alpha^\star(t)=1-\hat{\psi}_{t}\bc{1-\hat{\psi}_{t}(\alpha^\star(t))}< 1-\hat{\psi}_{t}(\alpha_\star(t))$, which contradicts the fact that $\alpha_\star(t)$ is the largest zero of $G_t$. As a result, $\alpha_\star(t)=1-\hat{\psi}_{t}(\alpha^\star(t))$. Analogously, it can be argued that $\alpha^\star(t)=1-\hat{\psi}_{t}(\alpha_\star(t))$.

        \item As remarked earlier, $G_t$ always has the zero $\alpha_0(t)$ (which satisfies $\alpha_0(t)=1-\hat{\psi}_{t}(\alpha_0(t))$ by its definition). If both $\alpha_\star(t)$ and $\alpha^\star(t)$ are equal to $\alpha_0(t)$, then $G_t$ has $1$ zero.
        If one of $\alpha_\star(t)$ and $\alpha^\star(t)$ is different from $\alpha_0(t)$, item \ref{ap_it2} shows that the other is also different from $\alpha_0(t)$ and hence $\alpha_\star(t)<\alpha_0(t)<\alpha^\star(t)$. In this case, using item \ref{ap_it1}, we conclude that $G_t$ has $3$ zeros. 

\item This follows directly from the fact that $G_t(1)=p_0>0=G_t(\alpha^\star(t))$ and there is no zero of $G_t$ in $(\alpha^\star(t),1)$. 

        \item  Recall that $\kappa=-\ln{((\hat{\psi}')^{-1}(1)+\hat{\psi}((\hat{\psi}')^{-1}(1)))}$ and let $\beta_0(t)=\eul^{-t}+\frac{\lambda(t)\eul^{t}}{\lambda(0)}(\alpha_0(t)-1)$. Then the combination of \cref{eq_psi0t} and $\alpha_0(t)=1-\hat{\psi}_{t}(\alpha_0(t))$ gives that
        \begin{align}\label{eq_beta0t}
            \beta_0(t)= \eul^{-t} - \hat{\psi}_{t}(\beta_0(t)).
        \end{align}
        Plugging \cref{eq_beta0t} into \cref{eq_dgt}, we obtain
        \begin{align}\label{eq_dgtab}
            G'_t(\alpha_0(t))=1-\hat{\psi}'\bc{\eul^{-t}-\hat{\psi}(\beta_0(t))}\hat{\psi}'(\beta_0(t))=1-\bc{\hat{\psi}'(\beta_0(t))}^2.
        \end{align}
       We now come to the case distinction in terms of the position of $t$ relative to $\kappa$. If $t<\kappa$, then \cref{eq_beta0t} yields that
        \begin{align*}
    \beta_0(t)+\hat{\psi}(\beta_0(t))>\eul^{-\kappa}=(\hat{\psi}')^{-1}(1)+\hat{\psi}((\hat{\psi}')^{-1}(1))).
        \end{align*}
        By monotonicity of $x \mapsto x+\hat\psi(x)$, we conclude that $\beta_0(t)>(\hat{\psi}')^{-1}(1)$, and by \cref{eq_dgtab} that $G'_t(\alpha_0(t))<0$. In this case, as $G_t(1)=\hat{\psi}_{t}(1-\hat{\psi}_{t}(1))>0$, there is at least one zero in $(\alpha_0(t),1)$. By item \ref{ap_it3}, $G_t$ has $3$ zeroes. Moreover, by Rolle's theorem, there exist $\alpha_1\in (\alpha_0(t),\alpha^\star(t))$ and $\alpha_2\in (\alpha_\star(t),\alpha_0(t))$ such that $G'_t(\alpha_1)=G_t'(\alpha_2)=0$. As $G_t'$ has at most 2 zeroes by item \ref{ap_it1}, neither $\alpha_\star(t)$ nor $\alpha^\star(t)$ are zeroes of $G_t'$.

        For $t>\kappa$, using \cref{eq_beta0t} and \cref{eq_dgtab}, we can analogously deduce that $G'_t(\alpha_0(t))>0$. In this case, we claim that $G_t$ has only $1$ zero. In fact, if $G_t$ has more than $1$ zero, then item \ref{ap_it3} shows that there is a zero $\alpha^\star(t)$ in $(\alpha_0(t),1)$. 
        Since $G'_t(\alpha_0(t))>0$ and $G_t(1)>0$, there exist $\alpha_0(t)<\alpha_3<\alpha^\star(t)<\alpha_4$ such that $G_t'(\alpha_3)<0$ and $G_t'(\alpha_4)>0$. As a consequence, there are at least $2$ zeroes of $G_t'$ in $(\alpha_0(t),1)$. Analogously, there are at least $2$ zeroes of $G_t'$ in $(0,\alpha_0(t))$. Hence, $G_t'$ has at least $4$ zeroes, which contradicts item \ref{ap_it1}. Thus, we conclude that $G_t$ has $1$ zero.

        Finally, for $t=\kappa$, using \cref{eq_beta0t} and \cref{eq_dgtab}, it follows that $\beta_0(\kappa)=(\hat{\psi}')^{-1}(1)$ and $G'_\kappa(\alpha_0(\kappa))=0$. Moreover, plugging \cref{eq_beta0t} into \cref{eq_ddgt} gives that $G''_\kappa(\alpha_0(\kappa))=0$ as well. Since $\alpha \mapsto \frac{\hat{\psi}''\bc{\eul^{-t}-\hat{\psi}(\beta)}}{\hat{\psi}'\bc{\eul^{-t}-\hat{\psi}(\beta)}}\hat{\psi}'(\beta)-\frac{\hat{\psi}''(\beta)}{\hat{\psi}'(\beta)}$ is an increasing function as shown in the proof of item \ref{ap_it1}, there exists $\varepsilon>0$ such that $G''_\kappa(\alpha)<0$ for $\alpha \in (\alpha_0(\kappa)-\varepsilon,\alpha_0(\kappa))$ and
        $G''_\kappa(\alpha)>0$ for $\alpha \in (\alpha_0(\kappa),\alpha_0(\kappa)+\varepsilon)$. Hence, $G'_\kappa(\alpha)>0$ on $(\alpha_0(\kappa)-\varepsilon,\alpha_0(\kappa))\cup(\alpha_0(\kappa),\alpha_0(\kappa)+\varepsilon)$. We can now deduce that $G_t$ has only one zero by the same line of argument as in the case $t>\kappa$.

        \item 
        Provided item \ref{ap_it5} and using the implicit function theorem, the proof is a slight variation of the proof of item 7 in \cite[Lemma 6.1]{HofMul25}.

        \item The proof is identical to the proof of item 8 of \cite[Lemma 6.1]{HofMul25}.

        \item Recall that $R_{\psi_t}(\alpha)=2-\psi_t(1-\hat{\psi}_t(\alpha))-\psi_t(\alpha)-\psi'_t(\alpha)(1-\alpha)$. Since $\hat{\psi}_t(x)=\psi_t'(x)/\psi'(1)$,
        \begin{align*}
            R_{\psi_t}'(\alpha)=\psi_t'(1-\hat{\psi}_t(\alpha))\hat{\psi}_t'(\alpha)-\psi_t'(\alpha)-\psi''_t(\alpha)(1-\alpha)+\psi'_t(\alpha)=\psi_t''(\alpha)G_t(\alpha).
        \end{align*}
        As $\psi_t''(\alpha)>0$, the sign of $G_t$ determines whether $R_{\psi_t}$ is increasing or decreasing.  As $G_t(0)<0$ and $G_t(1)>0$, $R_{\psi_t}$ obtains its minimum at the zeroes of $G_t$. In the case where $G_t$ has three zeroes, tt is straightforward to check that $R_{\psi_t}(\alpha^\star(t))=R_{\psi_t}(\alpha_\star(t))$, and by the proof of \cref{ap_it5}, $G'_t(\alpha_0(t))<0$. In this case, $G_t$ is strictly decreasing around  $\alpha_0(t)$, so $\alpha_0(t)$ is a local maximum of $R_{\psi_t}$. As a consequence, $R_{\psi_t}(\alpha)$ obtains its minimum in $[0,1]$ for $\alpha=\alpha_\star(t)$ or $\alpha^\star(t)$.
    \end{enumerate}
\end{proof}

\section{Proof of \Cref{coro_upper}: Upper bound on the rank}\label{ap_1}
\Cref{coro_upper} is a special case of the following more general result. For more background on local weak convergence, see \cite{van2022random}.

\begin{theorem}\label{at1}
Let $\bm{u}_n$ be chosen uniformly at random from $[n]$. Let $\bc{G_n}_{n\geq 1} = \bc{([n],E_n)}_{n \geq 1}$ be a sequence of finite (multi-)graphs such that $(G_n,\bm{u}_n)$ converges locally weakly to a unimodular Galton-Watson tree with root degree distribution $(p_k)_{k\geq 0}$. 
Then,
\begin{equation}\label{e1}
\lim_{n\to\infty}\limsup_{n\to \infty} \sup_{\substack{A_n\in\text{Sym}_n(\FF)\\A_n(i,j)\not=0 \Leftrightarrow i\neq j,\cbc{i,j}\in E_n}}\PP\bc{\frac{1}{n}\rank_{\FF}\bc{\bm{A}_{n}}\leq \min_{\alpha\in [0,1]}R_\psi(\alpha)- \varepsilon}=0,
\end{equation}
where $\psi$ is the probability generating function of $(p_k)_{k\geq 0}$ and $R_{\psi}$ is as in \Cref{def_rd}.
\end{theorem}

\begin{proof}[Proof of \Cref{coro_upper} subject to \Cref{at1}]
Since \cm{n} is a configuration model with degree sequence $\vd$ satisfying \Cref{assumption_weaker}, by \cite[Theorem 4.1]{van2022random}, for $\bm{u}_n$ chosen uniformly at random from $[n]$, $(\fcm{n},\bm{u}_n)$ converges locally weakly to a unimodular Galton-Watson tree with root degree generating function $\psi$. Recall the definition of  the random variable $\bm{A}_n(i,j)$ in \Cref{ec1}. Then by \Cref{at1},
\begin{align*}
\lim_{n\to\infty}\sup_{J_n\in \syn}\PP\bc{\frac{1}{n}\rank_{\FF}\bc{\bm{A}_{n}}\leq \min_{\alpha\in [0,1]}R_\psi(\alpha)- \varepsilon}=0.
\end{align*}
as desired.
\end{proof}
Thus, it only remains to prove \Cref{at1}.
\subsection{Proof of \Cref{at1}}
The proof of \Cref{at1} relies on a result by Bordenave, Lelarge and Salez \cite{bordenave2013matchings} on the normalized matching numbers of locally weakly convergent sequences of graphs. For a graph $G$, let $\nu(G)$ denote its matching number, i.e. the maximal size of a collection of pairwise non-adjacent edges in $G$.

\begin{theorem}[{\cite[Theorem 2]{bordenave2013matchings}}]\label{at2}
Let $\bm{u}_n$ be chosen uniformly at random from $[n]$. Let $\bc{G_n}_{n\geq 1} = \bc{([n],E_n)}_{n \geq 1}$ be a sequence of finite (multi-)graphs such that $(G_n,\bm{u}_n)$ converges locally weakly to a unimodular Galton-Watson tree with root degree distribution $(p_k)_{k\geq 0}$. Then 
\begin{equation}\label{e3}
  \lim_{n\to \infty} \frac{\nu(G_n)}{n}\stackrel{\mathbb{P}}{\longrightarrow} \frac{1}{2}\min_{\alpha\in [0,1]}R_{\psi}(\alpha),
\end{equation}
where $\psi$ is the probability generating function of $(p_k)_{k\geq 0}$ and $R_{\psi}$ is as in \Cref{def_rd}.
\end{theorem}

To build the connection between the rank of the symmetric matrices $A_{n}$ from \Cref{at1} and the matching number $\nu(G_{n})$ from \Cref{at2}, we adapt Lelarge's idea from \cite{Bypassing2013} to our symmetric setting. We start from the observation that for any field $\FF$, any symmetric matrix $A \in \FF^{n \times n}$ contains a full-rank principal  submatrix of dimension $\rk{A}$.

\begin{lemma}[Any symmetric matrix $A$ has a full rank principal submatrix of dimension $\rk{A}$]\label{l1}
Let $A=(A\bc{i,j})_{i,j\in[n]}\in \FF^{n\times n}$ be a symmetric matrix with $\rk{A}=m$. Then there exist $1\leq \ell_1<\ldots <\ell_m\leq n$ such that $B=(A\bc{\ell_i,\ell_j})_{i,j\in[m]}\in \FF^{m\times m}$ is a symmetric matrix of rank $m$.
\end{lemma}
\begin{proof}
For $i \in [n]$, denote by $r_i=A(i,)$ the $i$th row of $A$. Since the rank of $A$ over $\FF$ is $m$, there exist $1\leq \ell_1<\ldots< \ell_m\leq n$ such that the row vectors $(r_{\ell_i})_{i\in [m]}$ are linearly independent. The rank of the submatrix $C=(A\bc{\ell_i,j})_{i\in[m],j\in[n]}$ is $m$ and in particular, $C$ has the same rank as $A$. Hence, for any $i\in[n]$, there exist $k_{i, 1},\ldots, k_{i, m}\in \FF$ such that
\begin{equation}\label{e4}
r_i=\sum_{j=1}^{m}k_{i, j}r_{\ell_j}.
\end{equation}
For $i \in [n]$, denote by $c_i=C(,i)$ the $i$th column of $C$. For any $i\in[n]$, by (\ref{e4}) and symmetry of $A$,
\[c_i=\sum_{j=1}^{m}k_{i, j}c_{\ell_j}.\]
Therefore, the set $\{c_{\ell_j}:j\in [m]\}$ spans the column space of $C$. Since the rank of $C$ is $m$, the symmetric matrix $B=(A\bc{\ell_i,\ell_j})_{i,j\in[m]}\in\FF^{m\times m}$ has rank $m = \rk{A}$, which is the claim.
\end{proof}

\begin{proof}[Proof of \Cref{at1} subject to \Cref{at2}]
Let now $A_n$ be as in the statement of \Cref{at1}, and denote $\rk{A_n}$ by $m$. By Lemma \ref{l1}, there exist indices $1\leq \ell_1<\ldots <\ell_m\leq n$ such that the symmetric submatrix $B=(B\bc{i,j})_{i,j\in[m]}=(A_n\bc{\ell_i,\ell_j})_{i,j\in[m]}\in \FF^{m\times m}$ has full rank. 
In particular, its determinant is non-zero, so by the Leibniz formula for determinants
\[{\rm det}(B)=\sum_{\pi\in\mathcal{S}_m} {\rm sign}\bc{\pi}\prod_{i\in[m]}B\bc{i,\pi(i)}\neq 0,\]
where $\mathcal{S}_m$ is the set of all permutations of $[m]$ and ${\rm sign}(\pi)$ is the sign of $\pi \in \mathcal{S}_m$. We conclude that there exists a permutation $\tau\in\mathcal{S}_m$ such that $\prod_{i\in[m]}B\bc{i,\tau(i)}\neq 0$. Hence, by definition of $A_n$, for all $i \in [m]$, there is an edge between $\ell_i$ and $\ell_{\tau(i)}$ in the graph $G_n$.

Now consider the directed graph $L_m = (V(L_m), E(L_m)) = (\cbc{\ell_i: i \in [m]}, \{(\ell_i, \ell_{\tau(i)}): i \in [m]\})$, where $(a,b)$ denotes an edge that is directed from $a$ to $b$. 
Since $\tau$ is a permutation, each vertex of $L_m$ has out-degree $1$ and in-degree $1$, even though $L_m$ may have self-loops (corresponding to the fixed-points of $\tau$). 
The graph $L_m$ has the following structure: We call vertices $\ell_i$ with $(\ell_i,\ell_i) \in E(L_m)$ $1$-polygons. For $i\neq j$, if $(\ell_i,\ell_j) \in E(L_m)$ and $(\ell_j,\ell_i) \in E(L_m)$ (corresponding to the transpositions of $\tau$), we call $\{\ell_i,\ell_j\}$ a $2$-polygon. For $k\geq 3$, if there exist $k$ different vertices $\ell_{i_1},\ell_{i_2},\ldots, \ell_{i_k}$ such that $(\ell_{i_s},\ell_{i_{s+1}}) \in E(L_m)$ for all $s\in[k-1]$ and $(\ell_{i_k}, \ell_{i_1}) \in E(L_m)$, we call $\{\ell_{i_1},\ell_{i_2},\ldots, \ell_{i_k}\}$ a $k$-polygon. Since each vertex has out-degree $1$ and in-degree $1$, the vertices of $L_m$ can be partitioned into polygons as described.

The notion of polygons helps to lower bound the matching number $\nu(L_m)$, since each $k$-polygon contributes $\lfloor k/2 \rfloor$ edges to any maximum matching on $L_m$. For $k \in \NN$, denote by $s_k$ the number of $k$-polygons in $L_m$. Moreover, for $f \in \NN$ fixed, denote by $t_f$ the number of vertices belonging to $k$-polygons for $k \in 2\NN \cup \NN_{\geq f}$. Then by definition, 
\begin{align*}
\sum_{2 \mid  k \text{ or } k\geq f}k s_k=t_f \qquad \text{and} \qquad \sum_{2\nmid k \text{ and } k<f}k s_k=m-t_f. 
\end{align*}
Since $\lfloor k/2 \rfloor\geq \frac{f-1}{2f}k$ if $2 \mid k$ or $k\geq f$, writing $L_m^{\bullet}$ for the undirected version of $L_m$, we have
\[\begin{aligned}
\nu(L_m^{\bullet})\geq & \sum_{k}s_k\left\lfloor \frac{k}{2} \right\rfloor
=\sum_{2\nmid k \text{ and } k<f}s_k\left\lfloor \frac{k}{2}\right\rfloor+\sum_{2 \mid  k \text{ or } k\geq f}s_k\left\lfloor \frac{k}{2}\right\rfloor
\geq \sum_{2\nmid k \text{ and } k<f}\frac{s_k(k-1) }{2}+\sum_{2 \mid  k \text{ or } k\geq f} \frac{k s_k (f-1)}{2f}\\
= & \frac{m-t_f}{2}-\sum_{2\nmid k \text{ and } k<f}\frac{s_k}{2}+\frac{t_f(f-1)}{2f}
\geq  \frac{m(f-1)}{2f}-\sum_{2\nmid k \text{ and } k<f}\frac{s_k}{2}.
\end{aligned}\]

We now go back to the original graph $G_n$. For $k \in \NN$, denote by $g_{k,n}$ the number of cycles of length $k$ in $G_n$. Since there is an edge between $\ell_i$ and $\ell_j$ in $L_m$ only if there is an edge between $\ell_i$ and $\ell_j$ in $G_n$, for $k\neq 2$, $g_{k,n}\geq s_k$. Therefore, 
\begin{equation}\label{aeq_1}
    \begin{aligned}
\nu(G_n)\geq & \nu(L_m^{\bullet})
\geq  \frac{m(f-1)}{2f}-\sum_{2\nmid k\text{ and }k<f}\frac{s_k}{2}
\geq   \frac{m(f-1)}{2f}-\sum_{2\nmid k\text{ and }k<f}\frac{g_{k,n}}{2}.
\end{aligned}
\end{equation}

On the other hand,  for any $i\in [n]$, if $i$ is in a cycle of length $k$, its $k$-neighborhood \footnote{I.e., the induced subgraph of $G_n$ composed by the vertices of graph distance at most $k$ from $i$.} is not a tree. By \cite[Theorem 2.7]{van2022random}, since $(G_n,\bm{u}_n)$ converges locally weakly to a unimodular Galton-Watson tree\footnote{We note that convergence in distribution to a constant $c$ is equivalent to convergence in probability to $c$.
},
\[\lim_{n\to\infty}\frac{|\cbc{i\in[n]:\mbox{the $k$-neighborhood of $i$ in $G_n$ is not a tree}}|}{n}\stackrel{\mathbb{P}}{\longrightarrow}0.\]
As a consequence, $\lim_{n\to\infty}g_{k,n}/n=0$ for all $k\geq 1$. Hence, \cref{aeq_1} gives that
\[\begin{aligned}
\limsup_{n\to \infty}\frac{\nu(G_n)}{n} \geq& \limsup_{n\to \infty} \sup_{\substack{A_n\in\text{Sym}_n(\FF)\\A_n(i,j)\not=0 \Leftrightarrow \cbc{i,j}\in E_n}} \bc{\frac{ \bc{f-1}\rk{A_n}}{2nf}-\sum_{2\nmid k\text{ and }k<f}\frac{g_{k,n}}{2n}} \\
\stackrel{\mathbb{P}}{\longrightarrow} &\frac{f-1}{2f}\limsup_{n\to \infty} \sup_{\substack{A_n\in\text{Sym}_n(\FF)\\A_n(i,j)\not=0 \Leftrightarrow \cbc{i,j}\in E_n}} \frac{ \rk{A_n}}{n}.
\end{aligned}\]
Taking $f \to \infty$ on both sides of the equation above and combining \Cref{at2} yields \Cref{at1}.
\end{proof}

\section{Minor adaptions from \cite{HofMul25}}\label{app_pre_results}

\subsection{Proof of \Cref{prop_allcomestogether}}\label{sec_proof_allct}

\begin{proposition}\cite[Proposition C.1]{HofMul25}\label{pexp}
Fix a dimension $k \in \NN$ and $K>1$. Let $\bm{Z}_1,\bm{Z}_2$ and $\bm{X}$ be defined on the same probability space with convex codomains
$\mathcal{R}_{\bm{Z}_1},\mathcal{R}_{\bm{Z}_2}\subset \RR^k$ and $\mathcal{R}_{\bm{X}} \subset \RR$ respectively, such that $\mathcal{R}_{\bm{X}}$ is bounded. Then for any differentiable functions $f,g:\mathbb{R}^k\to\mathbb{R}$,
{\begin{equation*}
\begin{aligned}
&\mathbb{E}\abs{f(\bm{Z}_2)-g(\bm{Z}_2)}\\
\leq& \bc{\sup_{\zeta\in\mathcal{R}_{\bm{Z}_1}}\abs{f(\zeta)}+\sup_{x\in \mathcal{R}_{\bm{X}}}\abs{x}}\bc{4K^2\mathbb{E}\|\bm{Z}_1-\bm{Z}_2\|_\infty+2-2(1-1/K)^k}
+k\sup_{\zeta\in\mathcal{R}_{\bm{Z}_2}}\dabs{\nabla f(\zeta)}_\infty\mathbb{E}\|\bm{Z}_1-\bm{Z}_2\|_\infty\\
&+\mathbb{E}\abs{f(\bm{Z}_1)-\mathbb{E}\brk{\bm{X}|\bm{Z}_1}}
+\mathbb{E}\abs{\mathbb{E}\brk{\bm{X}|\bm{Z}_2}-g(\bm{Z}_2)}+\frac{2k}{K}\bc{\sup_{\zeta\in\mathcal{R}_{\bm{Z}_1}}\dabs{\nabla f(\zeta)}_\infty+\sup_{\zeta\in\mathcal{R}_{\bm{Z}_2}}\dabs{\nabla g(\zeta)}_\infty},
\end{aligned}
\end{equation*}}
and
{\begin{equation*}
\begin{aligned}
&\mathbb{E}\brk{\bc{f(\bm{Z}_2)-g(\bm{Z}_2)}^{-}}\\
\leq& \bc{\sup_{\zeta\in\mathcal{R}_{\bm{Z}_1}}\abs{f(\zeta)}+\sup_{x\in \mathcal{R}_{\bm{X}}}\abs{x}}\bc{4K^2\mathbb{E}\|\bm{Z}_1-\bm{Z}_2\|_\infty+2-2(1-1/K)^k}
+k\sup_{\zeta\in\mathcal{R}_{\bm{Z}_2}}\dabs{\nabla f(\zeta)}_\infty\mathbb{E}\|\bm{Z}_1-\bm{Z}_2\|_\infty\\
&+\mathbb{E}\brk{\bc{f(\bm{Z}_1)-\mathbb{E}\brk{\bm{X}|\bm{Z}_1}}^-}
+\mathbb{E}\brk{\bc{\mathbb{E}\brk{\bm{X}|\bm{Z}_2}-g(\bm{Z}_2)}^-}+\frac{2k}{K}\bc{\sup_{\zeta\in\mathcal{R}_{\bm{Z}_1}}\dabs{\nabla f(\zeta)}_\infty+\sup_{\zeta\in\mathcal{R}_{\bm{Z}_2}}\dabs{\nabla g(\zeta)}_\infty},
\end{aligned}
\end{equation*}}
where $a^-=0\vee (-a)$.
\end{proposition}
\begin{proof}[Proof of \Cref{prop_allcomestogether}]
    We first note from \Cref{lem_con_typp} that 
    \begin{align}\label{eq-allc-con-typp}
\PP\bc{\vnu_{s+1/n}\in\mathcal{W}\bc{\bm{A}_{n,s}[\bm{\theta}]} \mid \vze_s}=\vw_s,
    \end{align}
and from \eqref{eq_allcomestogether_a1} in combination with the assumption that the gradient of $W$ w.r.t. $\zeta$ is uniformly bounded that 
 \begin{align}\label{eq-allc-con-typp-2}
 W\bc{\vze_{s},\hat{\psi}_{t_s}}=W\bc{\vze_{s+1/n},\hat{\psi}_{t_s}}+\oone.
     \end{align}
If further
\begin{align}\label{eq_allcomestogether_e1}
\mathbb{P}\bc{\vnu_{s+1/n}\in\mathcal{W}\bc{\bm{A}_{n,s}[\bm{\theta}]}\mid \vze_s}-\mathbb{P}\bc{\vnu_{s+1/n}\in\mathcal{W}\bc{\bm{A}_{n,s}[\bm{\theta}]}\mid \vze_{s+1/n}}=\oone,
    \end{align}
then \eqref{eq_allcomestogether_a2}, \eqref{eq-allc-con-typp} and \eqref{eq-allc-con-typp-2} together yield
\eqref{eq_allcomestogether_re}. 

Let $M=\sup_{W\in\cbc{X,Y,Z,U,V}}\sup_{\ranges}\dabs{\nabla W(\cdot,\hat{\psi}_{t_s})}_\infty<\infty$.
To certify the missing piece \eqref{eq_allcomestogether_e1} under the hold of \eqref{eq_allcomestogether_a2}, we apply \Cref{pexp}. Taking 
\begin{itemize}
    \item $\bm{Z}_1=\vze_{s}$, $\bm{Z}_2=\vze_{s+1/n}$, $\bm{X}=\teo{\vnu_{s+1/n}\in\mathcal{W}\bc{\bm{A}_{n,s}[\bm{\theta}]}}$;
    \item $f(\zeta)=w$, $g(\zeta)=W(\zeta,\hat{\psi}_{t_s})$ for $\zeta=(x,y,z,u,v)\in [0,1]^5$ (note from \eqref{eq-allc-con-typp} that $f(\bm{Z}_1)=\mathbb{E}\brk{\bm{X}|\bm{Z}_1}$);
\end{itemize}
then \Cref{pexp} yields that
\begin{align*}
&\Erw\abs{\mathbb{E}\brk{\teo{\vnu_{s+1/n}\in\mathcal{W}\bc{\bm{A}_{n,s}[\bm{\theta}]}}\mid \vze_s}-\mathbb{E}\brk{\teo{\vnu_{s+1/n}\in\mathcal{W}\bc{\bm{A}_{n,s}[\bm{\theta}]}}\mid \vze_{s+1/n}}}\\
&\qquad \leq 4\bc{1-(1-1/K)^k}
+\frac{10}{K}\bc{1+M}+\oone.
\end{align*}
Letting $K\to\infty$, we conclude that
\begin{align*}
&\Erw\abs{\mathbb{P}\bc{\vnu_{s+1/n}\in\mathcal{W}\bc{\bm{A}_{n,s}[\bm{\theta}]}\mid \vze_s}-\mathbb{P}\bc{\vnu_{s+1/n}\in\mathcal{W}\bc{\bm{A}_{n,s}[\bm{\theta}]}\mid \vze_{s+1/n}}}\\
&\qquad=\Erw\abs{\mathbb{E}\brk{\teo{\vnu_{s+1/n}\in\mathcal{W}\bc{\bm{A}_{n,s}[\bm{\theta}]}}\mid \vze_s}-\mathbb{E}\brk{\teo{\vnu_{s+1/n}\in\mathcal{W}\bc{\bm{A}_{n,s}[\bm{\theta}]}}\mid \vze_{s+1/n}}}=\oone.
\end{align*}
Then \eqref{eq_allcomestogether_e1} —-- and thus, given \eqref{eq_allcomestogether_a2}, also \eqref{eq_allcomestogether_re} —-- hold. Analogously, given \eqref{eq_allcomestogether_a2_2}, \eqref{eq_allcomestogether_re2} holds.
\end{proof}

\subsection{Implications in the proof of \Cref{l7}}

In the subsequent proofs, we build on the following collection of results from \cite{HofMul25}:
\begin{lemma}[\cite{HofMul25}]\label{lem_frozen_zeroing}
Let $A\in \mathbb{F}^{m\times n}$, $i,j\in [n]$ with $i\neq j$ and $k\in [m]$. Then
\begin{enumerate}[label=(\roman*)]
    \item \makebox[6em][l]{$i\in \mathcal{F}(A)$}$ \quad \Longleftrightarrow  \quad \rank\bc{A}-\rank\bc{A\abc{;i}}=1$ {\rm ({\cite[Lemma 4.7]{DemGla24}}, {\cite[
Lemma 4.1]{HofMul25}})}\label{lem52_1}

\makebox[6em][l]{}$ \quad \Longleftrightarrow  \quad \mbox{the $i$th column of $A$ is not in the column space of $A\abc{;i}$}$;
  \item \makebox[6em][l]{$i\in \mathcal{F}(A)$}$ \quad \Longrightarrow \quad i\in \mathcal{F}(A\abc{;j})$  {\rm{(\cite[
Lemma 4.2]{HofMul25})}};\label{lem52_2}
  \item \makebox[6em][l]{$i\in \mathcal{F}(A\abc{k;})$}$ \quad \Longrightarrow \quad  i\in \mathcal{F}(A)$ {\rm ({\cite[
Lemma 4.2]{HofMul25}})}. 
\end{enumerate}
If further $i\in [m\land n]$, then
\begin{enumerate}
    \item $i \text{ is  frailly frozen in }A \qquad  \Longleftrightarrow \qquad i \text{ is  frailly frozen in } A^T$ {\rm ({\cite[
Proposition 4.5]{HofMul25}})};\label{eq_frailly}
\item for a uniformly chosen $k$-subset $\bm{\mathcal J} \subseteq [m]\backslash\cbc{i}$,
\begin{align*}
    \mathbb{P}\left(i\in\mathcal{F}\bc{A[\bm{\theta}]}\Delta \mathcal{F}\bc{A[\bm{\theta}]\abc{\bm{\mathcal J};}}\right)\leq \frac{k}{P}+\bc{1-\frac{1}{m}}^k\frac{k(k-1)}{2m} + \frac{k}{m}\ \mbox{{\rm ({\cite[
Corollary 4.20]{HofMul25}})}.}
\end{align*}
\end{enumerate}
\end{lemma}

\subsubsection{Proof of \Cref{implication_unfreeze}}
We restrict to the event $\mathfrak P_{n,s}$. On this event, $\bm{A}_{n,s+1/n}[\bm{\theta}]$ has a zero column in position $\vnu_{s+1/n}$, so that $\vnu_{s+1/n}$ is not frozen in $\bm{A}_{n,s+1/n}[\bm{\theta}]$. 
Assume now that $i \in \mathcal{F}\bc{\bm{A}_{n,s+1/n}[\bm{\theta}]}\backslash \mathcal{F}(\bm{A}_{n,s}[\bm{\theta}])$. As $i$ is frozen in $\bm{A}_{n,s+1/n}[\bm{\theta}]$, we first conclude that $i\neq \vnu_{s+1/n}$. 
Next, $\bm{A}_{n,s+1/n}[\bm{\theta}]$ arises from $\bm{A}_{n,s}[\bm{\theta}]$ by zeroing row and column $\vnu_{s+1/n}$, which we break into two steps. By \Cref{lem_frozen_zeroing}, on $\mathfrak P_{n,s}$, replacing a zero row by an arbitrary row cannot unfreeze $i$, so
  \[i\in \mathcal{F}\bc{\bm{A}_{n,s+1/n}[\bm{\theta}]} 
  \qquad  \Longrightarrow \qquad i\in \mathcal{F}\bc{\bm{A}_{n,s}[\bm{\theta}]\abc{;\vnu_{s+1/n}}}.\]
  In particular, there also exists a representation $y$ of $\{i\}$ in $\bm{A}_{n,s}[\bm{\theta}]\abc{;\vnu_{s+1/n}}$. Applying this representation $y$ to $\bm{A}_{n,s}[\bm{\theta}]$ yields
  \[\cbc{i}\subseteq \supp{y\bm{A}_{n,s}[\bm{\theta}]} \subseteq \{i,\vnu_{s+1/n}\}.\]
  Since $i\not\in \mathcal{F}\bc{\bm{A}_{n,s}[\bm{\theta}]}$ by assumption, we conclude that
  \[\supp{y\bm{A}_{n,s}[\bm{\theta}]} = \{i,\vnu_{s+1/n}\}.\]
The just discovered relation $\{i,\vnu_{s+1/n}\}$ is proper in $\bm{A}_{n,s}[\bm{\theta}]$, since any representation of $\{i\}$ or $\{\vnu_{s+1/n}\}$ could be combined with $y$ to yield a representation of $\{\vnu_{s+1/n}\}$ or $\{i\}$, in contradiction to the assumption that  
$i$ is not frozen in $\bm{A}_{n,s}[\bm{\theta}]$.
\qed
\subsubsection{Proof of \cref{implication_freeze}}
We restrict to the event $\mathfrak P_{n,s}$. Assume that 
$i \in \mathcal{F}(\bm{A}_{n,s}[\bm{\theta}])\backslash \mathcal{F}\bc{\bm{A}_{n,s+1/n}[\bm{\theta}]}$ and $i\notin \{\vnu_{s+1/n}\}\cup\supp{\bm{A}_{n,s}(\vnu_{s+1/n},)}$. The matrix $\bm{A}_{n,s+1/n}[\bm{\theta}]$ arises from $\bm{A}_{n,s}[\bm{\theta}]$ by replacing row and column $\vnu_{s+1/n}$ by a zero row and column. By \Cref{lem_frozen_zeroing}, replacing column $\vnu_{s+1/n}$ by a zero row cannot unfreeze $i$, so that
  \[i\in \mathcal{F}\bc{\bm{A}_{n,s}[\bm{\theta}]} \qquad  \Longrightarrow \qquad i\in \mathcal{F}\bc{\bm{A}_{n,s}[\bm{\theta}]\abc{;\vnu_{s+1/n}}}.\]
As we assume that $i$ is frozen in $\bm{A}_{n,s}[\bm{\theta}]$, this implies the existence of a representation $y$ of $\{i\}$ in $\bm{A}_{n,s}[\bm{\theta}]\abc{;\vnu_{s+1/n}}$. Moreover, any such representation $y$ has to use row $\vnu_{s+1/n}$
(i.e., $y_{\vnu_{s+1/n}} \not=0$), as otherwise $y$ would also be a representation of $\{i\}$ in $\bm{A}_{n,s+1/n}[\bm{\theta}]$, in contrast to our assumption that $i$ is not frozen in $\bm{A}_{n,s+1/n}[\bm{\theta}]$.
Hence, 
multiplying $y$ from the left with $\bm{A}_{n,s+1/n}[\bm{\theta}]$ yields non-zero coordinates in places
  \[\supp{y\bm{A}_{n,s+1/n}[\bm{\theta}]} = \cbc{i}\cup\supp{\bm{A}_{n,s}[\bm{\theta}](\vnu_{s+1/n},)}.\]
  Thus, $\cbc{i}\cup\supp{\bm{A}_{n,s}[\bm{\theta}](\vnu_{s+1/n},)}$ is a relation in $\bm{A}_{n,s+1/n}[\bm{\theta}]$. It is proper, since it contains $i$, which is not frozen in $\bm{A}_{n,s+1/n}[\bm{\theta}]$. 
  \qed
\subsection{End of proof of \Cref{pro_sta_vze}}\label{sec-proof-dis-w-s-s+1/n}
Recall the bound \cref{eq_l91_x1}. Starting from there, we now follow a case distinction according to the different types.

\textbf{Case 1}: \textit{Frailly frozen variables}. 
\Cref{dxyzuv} yields the identity
\begin{align}\label{eq_l91_x2}
    \ind\cbc{i \in \mathcal X\bc{\bm{A}_{n,s}[\bm{\theta}]}} =  \ind\cbc{i \in \mathcal F\bc{\bm{A}_{n,s}[\bm{\theta}]} \Delta \mathcal F\bc{\bm{A}_{n,s}[\bm{\theta}] \abc{i;}} }.
\end{align}
Using $(\mathfrak B_1\Delta \mathfrak B_2)\Delta(\mathfrak B_3\Delta \mathfrak B_4)\subseteq (\mathfrak B_1\Delta \mathfrak B_3)\cup (\mathfrak B_2\Delta \mathfrak B_4)$ for sets $\mathfrak B_1, \ldots, \mathfrak B_4$ and plugging (\ref{eq_l91_x2}) into (\ref{eq_l91_x1}) yields 
\begin{equation}\label{eq_l91_x3}
\begin{aligned}
 \mathbb{E}|\vx_s-\vx_{s+1/n}| \leq &\frac{K}{cn}\sum_{i\in [n]}\Big(\mathbb{P}\bc{i\in\mathcal{F}\bc{\bm{A}_{n,s}[\bm{\theta}]}\Delta \mathcal{F}\bc{\bm{A}_{n,s+1/n}[\bm{\theta}]}} \\
& \quad +\mathbb{P}\bc{i\in\mathcal{F}\bc{\bm{A}_{n,s}[\bm{\theta}]\abc{i;}}\Delta \mathcal{F}\bc{\bm{A}_{n,s+1/n}[\bm{\theta}]\abc{i;}}}\Big)+\bar{o}_n(1).
\end{aligned}
\end{equation}
\Cref{l7,l71} now imply that
\begin{align}\label{eq_xn0}
\mathbb{E}|\vx_s-\vx_{s+1/n}|\leq  8\delta (K+1)K^{2K+3}c^{-K-2}+\bar{o}_{n,P}(1).
\end{align}
In particular,
\begin{align}\label{eq_xn_lim}
\limsup_{P\to\infty}\limsup_{n\to\infty}\sup_{J_n\in\syn,\ranges} \mathbb{E}|\vx_s-\vx_{s+1/n}|\leq  8\delta (K+1)K^{2K+3}c^{-K-2}.
\end{align}
Since the {left} hand side of \cref{eq_xn_lim} does not depend on $\delta$, we can send $\delta\downarrow 0$ to conclude that
\begin{align}\label{eq_xn01}
\limsup_{P\to\infty}\limsup_{n\to\infty}\sup_{J_n\in\syn,\ranges} \mathbb{E}|\vx_s-\vx_{s+1/n}|= 0
\end{align}
or equivalently, that $\mathbb{E}|\vx_s-\vx_{s+1/n}|=\bar{o}_{n,P}(1)$.

\textbf{Case 2}: \textit{Variables that are firmly frozen in original and transposed matrix.} 
\Cref{dsc} of firmly frozen variables and \Cref{dxyzuv} of the set $\mathcal Y$ yield the identity
\begin{align}\label{eq_l91_y2}
    \ind\cbc{i \in \mathcal Y\bc{\bm{A}_{n,s}[\bm{\theta}]}} =  \ind\cbc{i \in \mathcal F\bc{\bm{A}_{n,s}[\bm{\theta}]\abc{i;}} \cap \mathcal F\bc{\bm{A}_{n,s}[\bm{\theta}]^T\abc{i;}} }.
\end{align}
Using $(\mathfrak B_1\cap \mathfrak B_2)\Delta(\mathfrak B_3\cap \mathfrak B_4)\subseteq (\mathfrak B_1\Delta \mathfrak B_3)\cup (\mathfrak B_2\Delta \mathfrak B_4)$ for sets $\mathfrak B_1, \ldots, \mathfrak B_4$ and plugging (\ref{eq_l91_y2}) into (\ref{eq_l91_x1}) yields 
\begin{equation}\label{eq_l91_y3}
\begin{aligned}
\mathbb{E}|\vy_s-\vy_{n,s+1/n}|&\leq \frac{K}{cn}\sum_{i\in [n]}\Big(\mathbb{P}\bc{i\in\mathcal{F}\bc{\bm{A}_{n,s}[\bm{\theta}]\abc{i;}}\Delta \mathcal{F}\bc{\bm{A}_{n,s+1/n}[\bm{\theta}]\abc{i;}}} \\
&+\mathbb{P}\bc{i\in\mathcal{F}(\bm{A}_{n,s}[\bm{\theta}]^T\abc{i;})\Delta \mathcal{F}(\bm{A}_{n,s+1/n}[\bm{\theta}]^T\abc{i;})}\Big)+\bar{o}_n(1).
\end{aligned}
\end{equation}
\Cref{l71} then implies that 
\begin{align*}
\mathbb{E}|\vy_s-\vy_{n,s+1/n}|\leq  8\delta (K+1)K^{2K+3}c^{-K-2} +\bar{o}_{n,P}(1).
\end{align*}
Now the same limiting argument as in \textbf{Case 1} yields that
$\mathbb{E}|\vy_s-\vy_{n,s+1/n}|=\bar{o}_{n,P}(1).$

\textbf{Case 3}: \textit{Variables that are neither frozen in original nor transposed matrix.}  
\Cref{dxyzuv} of the set $\mathcal Z$ yields the identity
\begin{align}\label{eq_l91_z2}
    \ind\cbc{i \in \mathcal Z\bc{\bm{A}_{n,s}[\bm{\theta}]}} =  \ind\cbc{i \notin \mathcal F\bc{\bm{A}_{n,s}[\bm{\theta}]} \cup \mathcal F\bc{\bm{A}_{n,s}[\bm{\theta}]^T} }.
\end{align}
Using $(\mathfrak B_1\cup \mathfrak B_2)\Delta(\mathfrak B_3\cup \mathfrak B_4)\subseteq (\mathfrak B_1\Delta \mathfrak B_3)\cup (\mathfrak B_2\Delta \mathfrak B_4)$ for sets $\mathfrak B_1, \ldots, \mathfrak B_4$ and plugging (\ref{eq_l91_z2}) into (\ref{eq_l91_x1}) yields the upper bound
\begin{equation}\label{eq_l91_z3}
\begin{aligned}
 \frac{K}{cn}\sum_{i\in [n]}\bc{\mathbb{P}\bc{i\in\mathcal{F}\bc{\bm{A}_{n,s}[\bm{\theta}]}\Delta \mathcal{F}\bc{\bm{A}_{n,s+1/n}[\bm{\theta}]}} +\mathbb{P}\bc{i\in\mathcal{F}\bc{\bm{A}_{n,s}[\bm{\theta}]^T}\Delta \mathcal{F}\bc{\bm{A}_{n,s+1/n}[\bm{\theta}]^T}}}+\bar{o}_n(1)
\end{aligned}
\end{equation}
for $\mathbb{E}|\vz_s-\vz_{n,s+1/n}|$.
\Cref{l7} then implies that
\begin{align*}
\mathbb{E}|\vz_s-\vz_{n,s+1/n}|\leq  8\delta (K+1)K^{2K+3}c^{-K-2}+\bar{o}_{n,P}(1).
\end{align*}
Again, the same limiting argument as in \textbf{Case 1} yields $\mathbb{E}|\vz_s-\vz_{n,s+1/n}|=\bar{o}_{n,P}(1).$

\textbf{Case 4}: \textit{Variables that are not frozen in original but firmly frozen in transposed matrix.}  
\Cref{dxyzuv} of the set $\mathcal U$ yields the identity
\begin{align}\label{eq_l91_u2}
    \ind\cbc{i \in \mathcal U\bc{\bm{A}_{n,s}[\bm{\theta}]}} =  \ind\cbc{i \in \mathcal F\bc{\bm{A}_{n,s}[\bm{\theta}]^T\abc{i;}} \setminus \mathcal F\bc{\bm{A}_{n,s}[\bm{\theta}]} }.
\end{align}
Using $(\mathfrak B_1\setminus \mathfrak B_2)\Delta(\mathfrak B_3\setminus \mathfrak B_4)\subseteq (\mathfrak B_1\Delta \mathfrak B_3)\cup (\mathfrak B_2\Delta \mathfrak B_4)$ for sets $\mathfrak B_1, \ldots, \mathfrak B_4$ and plugging (\ref{eq_l91_u2}) into (\ref{eq_l91_x1}) yields 
\begin{equation}\label{eq_l91_u3}
\begin{aligned}
\mathbb{E}\abs{\vu_s-\vu_{n,s+1/n}} & \leq \frac{K}{cn}\sum_{i\in [n]}\Big(\mathbb{P}\bc{i\in\mathcal{F}\bc{\bm{A}_{n,s}[\bm{\theta}]}\Delta \mathcal{F}\bc{\bm{A}_{n,s+1/n}[\bm{\theta}]}}\\
&+\mathbb{P}\bc{i\in\mathcal{F}(\bm{A}_{n,s}[\bm{\theta}]^T\abc{i;})\Delta \mathcal{F}(\bm{A}_{n,s+1/n}[\bm{\theta}]^T\abc{i;})}\Big)+\bar{o}_n(1).
\end{aligned}
\end{equation}
\Cref{l7,l71} then give
\begin{align*}
\mathbb{E}|\vu_s-\vu_{n,s+1/n}|\leq
8\delta (K+1)K^{2K+3}c^{-K-2}+\bar{o}_{n,P}(1).
\end{align*}
Again, the same limiting argument as in \textbf{Case 1} yields $\mathbb{E}|\vu_s-\vu_{n,s+1/n}|=o_{n,P}(1).$

\textbf{Case 5}: \textit{Variables that are firmly frozen in original but not frozen in transposed matrix.}  
Analogous to \textbf{Case 4}.

\textbf{Case 6}: \textit{Vector of type proportions.} 
This directly follows from \textbf{Cases 1-5}.
\qed

\subsection{Proof of \Cref{lem_typ_equ}} \label{sec_proof_typ_equ}
We first prove an intermediate result, which establishes the desired overall equivalence of the first and last statements of \cref{im_fnf} for $\vnu_{s+1/n}$:  
\begin{lemma}\label{l10}
For any $\delta>0$ and $\ranges$,
\begin{align} \label{eq_lem_l10_3}
\mathbb{P}\left(\text{$\vnu_{s+1/n}$ is firmly frozen in $\bm{A}_{n,s}[\bm{\theta}]$}, \ffBT \right) = \bar{o}_{n,P}(1),
\end{align}
and
\begin{align} \label{eq_lem_l10_1}
 \mathbb{P}\left(\text{$\vnu_{s+1/n}$ is not firmly frozen in $\bm{A}_{n,s}[\bm{\theta}]$}, \fBT\right)=\bar{o}_{n,P}(1).
\end{align}
\end{lemma}
\begin{proof}
\begin{enumerate}[label=(\roman*)]
\item We first show (\ref{eq_lem_l10_3}). 
By definition, $\vnu_{s+1/n}$ is firmly frozen in $\bm{A}_{n,s}[\bm{\theta}]$ if and only if it is frozen in $\bm{A}_{n,s}[\bm{\theta}]\abc{\vnu_{s+1/n};}$. On the good event $\mathfrak{P}_{n,s}$ from \eqref{event_pns}, zeroing out row $\vnu_{s+1/n}$ leaves us with the matrix $\bm{A}_{n,s+1/n}[\bm \theta]$ except the  $\vnu_{s+1/n}$th  column being $\bm{A}_{n,s}[\bm \theta](,\vnu_{s+1/n})$. By \Cref{im_fnf}, 
\begin{align*}
    \text{$\vnu_{s+1/n}$ is firmly frozen in $\bm{A}_{n,s}[\bm{\theta}]$},\mathfrak{P}_{n,s} \quad  \Longrightarrow\quad \text{$\supp{\bm{A}_{n,s}(,\vnu_{s+1/n})}\not\subseteq \mathcal{F}\bc{\bm{A}_{n,s+1/n}[\bm{\theta}]^T}$.}
\end{align*}
On the other hand, on $\ffBT$, $\supp{\bm{A}_{n,s}(,\vnu_{s+1/n})}=\supp{\bm{A}_{n,s}(\vnu_{s+1/n},)}\subseteq\mathcal{F}\bc{\bm{A}_{n,s+1/n}[\bm{\theta}]^T}$, 
so that
\begin{align*}
    &\mathbb{P}\left(\text{$\vnu_{s+1/n}$ is firmly frozen in $\bm{A}_{n,s}[\bm{\theta}]$}, \ffBT \right)\\
&= \mathbb{P}\left(\text{$\vnu_{s+1/n}$ is firmly frozen in $\bm{A}_{n,s}[\bm{\theta}]$}, \ffBT, \mathfrak P_{n,s} \right) + \bar{o}_{n,P}(1) = \bar{o}_{n,P}(1).
\end{align*}
 \item We next prove (\ref{eq_lem_l10_1}).
    By definition, if $\vnu_{s+1/n}$ is not firmly frozen in $\bm{A}_{n,s}[\bm{\theta}]$, it is not frozen in $\bm{A}_{n,s}[\bm{\theta}]\abc{\vnu_{s+1/n};}$. 
    On the good event $\mathfrak{P}_{n,s}$, zeroing out row $\vnu_{s+1/n}$ leaves us with the matrix $\bm{A}_{n,s+1/n}[\bm \theta]$ except the  $\vnu_{s+1/n}$th  column being $\bm{A}_{n,s}[\bm \theta](,\vnu_{s+1/n})$.  Since $\vnu_{s+1/n}$ is not frozen in this matrix, all diagonal entries are $0$ and $\bm{A}_{n,s+1/n}[\vth] = \bm{A}_{n,s}[\vth]\abc{\vnu_{s+1/n},\vnu_{s+1/n}}$ on $\mathfrak{P}_{n,s}$, by \Cref{im_fnf},
\begin{align*}
    &\text{$\vnu_{s+1/n}$ is not firmly frozen in $\bm{A}_{n,s}[\bm{\theta}]$},\mathfrak{P}_{n,s} \\
    \Longrightarrow\ &  e_{n+\vth_c}(\vnu_{s+1/n}) \mbox{ is not in the row space of $\bm{A}_{n,s}[\bm{\theta}]\abc{\vnu_{s+1/n};}$}, \mathfrak{P}_{n,s}\\
     \Longrightarrow\ & 
\text{ $\bm{A}_{n,s}[\vth](,\vnu_{s+1/n})$ is in the column space of $\bm{A}_{n,s+1/n}[\bm{\theta}]$.}
\end{align*}

On the other hand, on $\fBT$,  $\supp{\bm{A}_{n,s}(,\vnu_{s+1/n})}\not\subseteq\mathcal{F}\bc{\bm{A}_{n,s+1/n}[\bm{\theta}]^T}$. This implies that both $\supp{\bm{A}_{n,s}(,\vnu_{s+1/n})}$ and $\supp{\bm{A}_{n,s}(,\vnu_{s+1/n})}\backslash \mathcal{F}\bc{\bm{A}_{n,s+1/n}[\bm{\theta}]^T}$ are non-empty. If additionally, $\bm{A}_{n,s}(,\vnu_{s+1/n})$ is in the column space of $\bm{A}_{n,s+1/n}[\bm{\theta}]$, then $\supp{\bm{A}_{n,s}(,\vnu_{s+1/n})}\backslash \mathcal{F}\bc{\bm{A}_{n,s+1/n}[\bm{\theta}]^T}$ is a relation of $\bm{A}_{n,s+1/n}[\bm{\theta}]^T$. Hence, by \Cref{d2} (iii),
\begin{align*}
    &\text{$\vnu_{s+1/n}$ is not firmly frozen in $\bm{A}_{n,s}[\bm{\theta}]$}, \fBT, \mathfrak{P}_{n,s}\\
    &\Longrightarrow\quad 
\text{$\supp{\bm{A}_{n,s}(,\vnu_{s+1/n})}$ is a proper relation of $\bm{A}_{n,s+1/n}[\bm{\theta}]^T$.}
\end{align*}
By \Cref{eq_con_halfedge} and the same argument deriving \Cref{eq_useinapp},
\begin{align*}
    \PP\left(\supp{\bm{A}_{n,s}(,\vnu_{s+1/n})} \in \PR\bc{\bm{A}_{n,s+1/n}[\bm{\theta}]^T}\right)\leq \delta K^{2K+1} c^{-K}+\bar{o}_{n,P}(1).
\end{align*}
Hence,
\[\begin{aligned}
&\mathbb{P}\left(\text{$\vnu_{s+1/n}$ is not firmly frozen in $\bm{A}_{n,s}[\bm{\theta}]$}, \fBT, \mathfrak{P}_{n,s} \right)\\
&\leq \PP\left(\supp{\bm{A}_{n,s}(,\vnu_{s+1/n})} \in \PR\bc{\bm{A}_{n,s+1/n}[\bm{\theta}]^T}\right)\leq \delta K^{2K+1} c^{-K}+\bar{o}_{n,P}(1).
\end{aligned}\]
Consequently, 
\begin{align*}
    &\limsup_{P\to\infty}\limsup_{n\to\infty}\sup_{J_n\in\syn,\ranges}\mathbb{P}\left(\text{$\vnu_{s+1/n}$ is not firmly frozen in $\bm{A}_{n,s}[\bm{\theta}]$}, \fBT  \right)\\
    &\qquad\leq \delta K^{2K+1} c^{-K}.
\end{align*}
Since the left-hand side is independent of $\delta$, letting $\delta \downarrow 0$, we conclude that
\begin{align*}
    \mathbb{P}\left(\text{$\vnu_{s+1/n}$ is not firmly frozen in $\bm{A}_{n,s}[\bm{\theta}]$}, \fBT  \right)=o_{n,P}(1).
\end{align*}
\end{enumerate}
\end{proof}

Using the definition of the variable types in terms of frozen variables, we can now easily derive \Cref{lem_typ_equ} from \Cref{l10}:

\begin{proof}[Proof of \Cref{lem_typ_equ}]
We show the claim for each of the possible variable types separately.

    \textit{Two-sided firmly frozen variables - (\ref{bound_QY}) for $W=Y$}:
    By definition, if $\vnu_{s+1/n}\not\in\mathcal{Y}\bc{\bm{A}_{n,s}[\bm{\theta}]}$, then it is not firmly frozen in $\bm{A}_{n,s}[\bm{\theta}]$ or not firmly frozen in $\bm{A}_{n,s}[\bm{\theta}]^T$. Since \Cref{l10} also applies to $\bm{A}_{n,s}[\bm \theta]^T$, a union bound gives
\[\begin{aligned}
&\mathbb{P}\left(\vnu_{s+1/n}\not\in\mathcal{Y}\bc{\bm{A}_{n,s}[\bm{\theta}]}, {\mathfrak Y}_s\right)\\
 \leq& \mathbb{P}\left(\text{$\vnu_{s+1/n}$ not firmly frozen in $\bm{A}_{n,s}[\bm{\theta}]$}, \fBT \right)+\mathbb{P}\left(\text{$\vnu_{s+1/n}$ not firmly frozen in $\bm{A}_{n,s}[\bm{\theta}]^T$}, \fB \right)\\
 =&  \bar{o}_{n,P}(1).
\end{aligned}\]

\textit{One-sided firmly frozen variables  - (\ref{bound_QY}) for $W \in \{U,V\}$}:
    If $\vnu_{s+1/n}\not\in\mathcal{U}\bc{\bm{A}_{n,s}[\bm{\theta}]}$, then, by definition, either $\vnu_{s+1/n}$ is not firmly frozen in $\bm{A}_{n,s}[\bm{\theta}]^T$, or, if this is not the case, it is  frozen in $\bm{A}_{n,s}[\bm{\theta}]$ \textit{and} firmly frozen in $\bm{A}_{n,s}[\bm{\theta}]^T$. In the latter case, the symmetry of frailly frozen variables under transposition (see \Cref{lem_frozen_zeroing}) implies that $\vnu_{s+1/n}$ is also firmly frozen in $\bm{A}_{n,s}[\bm{\theta}]$. We conclude that if $\vnu_{s+1/n}\not\in\mathcal{U}\bc{\bm{A}_{n,s}[\bm{\theta}]}$, then either $\vnu_{s+1/n}$ is not firmly frozen in $\bm{A}_{n,s}[\bm{\theta}]^T$ or $\vnu_{s+1/n}$ is firmly frozen in $\bm{A}_{n,s}[\bm{\theta}]$. Again, by a union bound and \Cref{l10},
\[\begin{aligned}
&\mathbb{P}\left(\vnu_{s+1/n}\not\in\mathcal{U}\bc{\bm{A}_{n,s}[\bm{\theta}]}, {\mathfrak U}_s\right)\\
\leq& \mathbb{P}\left(\text{$\vnu_{s+1/n}$ not firmly frozen in $\bm{A}_{n,s}[\bm{\theta}]^T$}, \fB \right)+\mathbb{P}\left(\text{$\vnu_{s+1/n}$ firmly frozen in $\bm{A}_{n,s}[\bm{\theta}]$}, \ffBT \right)\\
=&\bar{o}_{n,P}(1).
\end{aligned}\]

The claim for $W = V$ follows analogously.

\textit{Frailly frozen or two-sided not frozen variables - (\ref{bound_QXZ})}:
If $\vnu_{s+1/n}\not\in \mathcal{X}\bc{\bm{A}_{n,s}[\bm{\theta}]}\cup \mathcal{Z}\bc{\bm{A}_{n,s}[\bm{\theta}]}$, then by definition, $\vnu_{s+1/n}$ is firmly frozen in $\bm{A}_{n,s}[\bm{\theta}]$ or $\bm{A}_{n,s}[\bm{\theta}]^T$. By a union bound and \Cref{l10},
\[\begin{aligned}
&\mathbb{P}\left(\vnu_{s+1/n}\not\in\mathcal{X}\bc{\bm{A}_{n,s}[\bm{\theta}]}\cup \mathcal{Z}\bc{\bm{A}_{n,s}[\bm{\theta}]}, {\mathfrak{XZ}}_s\right) \\
\leq & \hspace{0.1 cm} \mathbb{P}\left(\text{$\vnu_{s+1/n}$ firmly frozen in $\bm{A}_{n,s}[\bm{\theta}]$}, \ffBT \right)+\mathbb{P}\left(\text{$\vnu_{s+1/n}$ firmly frozen in $\bm{A}_{n,s}[\bm{\theta}]^T$}, \ffB \right)=\bar{o}_{n,P}(1).
\end{aligned}\]
\end{proof}

\subsection{Proof of \Cref{lem_typ_zin}} \label{sec_proof_typ_zin}
The first (and main) step is to prove that on the intersection of $\mathfrak Z^\circ_s$ with a sufficiently likely event,
the $\vnu_{s+1/n}$th row in $\bm{A}_{n,s}[\bm{\theta}]$ can be linearly combined by the other rows of $\bm{A}_{n,s}[\bm{\theta}]$, from which it follows through \Cref{lem_frozen_zeroing} (i) that $\vnu_{s+1/n}$ is not frozen in $\bm{A}_{n,s}[\bm{\theta}]^T$.

On ${\mathfrak Z}_s^{\circ} \cap \mathfrak P_{n,s}$, $\supp{\bm{A}_{n,s}[\bm \theta](\vnu_{s+1/n},)}= \supp{\bm{A}_{n,s}(\vnu_{s+1/n},)}$ and all variables in $\supp{\bm{A}_{n,s}[\bm \theta](\vnu_{s+1/n},)}$ are firmly frozen in $\bm{A}_{n,s+1/n}[\bm{\theta}]$. Ideally, to derive the desired linear combination of $\bm{A}_{n,s}[\bm \theta](\vnu_{s+1/n},)$  by the other rows of $\bm{A}_{n,s}[\bm{\theta}]$, we would like to take one representation for each $i \in \supp{\bm{A}_{n,s}[\bm \theta](\vnu_{s+1/n},)}$ in $\bm{A}_{n,s+1/n}[\bm{\theta}]$, and then simply take the $\bm{A}_{n,s}(\vnu_{s+1/n},i)$-weighted sum over these representations. Alas, the $\vnu_{s+1/n}$th column of $\bm{A}_{n,s}[\bm{\theta}]$ might contain nonzero entries, and it is not clear that for the existing representations, also the entries of column $\vnu_{s+1/n}$ sum to zero. Therefore, we are looking for representations of $i \in \supp{\bm{A}_{n,s}[\bm \theta](\vnu_{s+1/n},)}$ that expressly do not use one of the rows in $\supp{\bm{A}_{n,s}(,\vnu_{s+1/n})}=\supp{\bm{A}_{n,s}(\vnu_{s+1/n},)}$, if such representations exist.

In fact, on ${\mathfrak Z}_s^{\circ} \cap \mathfrak P_{n,s}$, since any $i\in \supp{\bm{A}_{n,s}[\bm \theta](\vnu_{s+1/n},)}$ is firmly frozen in $\bm{A}_{n,s+1/n}[\bm{\theta}]$, it is frozen in $\bm{A}_{n,s+1/n}[\bm{\theta}]\abc{i;}$, so there exists a representation of $i$ that does not use row $i$. To take care of the other rows corresponding to elements of $\supp{\bm{A}_{n,s}[\bm \theta](\vnu_{s+1/n},)}$, we define the event 
\begin{align*}
\mathfrak C=&\cbc{\text{$\forall i\in \supp{\bm{A}_{n,s}[\bm \theta](\vnu_{s+1/n},)}$: $i \notin \mathcal{F}\bc{\bm{A}_{n,s+1/n}[\bm{\theta}]\abc{\supp{\bm{A}_{n,s}[\bm \theta](\vnu_{s+1/n},)};}}\Delta \mathcal{F}\bc{\bm{A}_{n,s+1/n}[\bm{\theta}]\abc{i;}}$}}.
\end{align*}
The event $\mathfrak C$ is sufficiently likely for our purposes, as the combination of \cref{eq_con_halfedge} and \Cref{lem_frozen_zeroing} gives that
\begin{align*}
    \mathbb{P}\bc{\mathfrak C^c}
\leq&\sum_{k=1}^K\sum_{1\leq j_1<\ldots<j_k\leq  n}\sum_{\ell=1}^k\PP\bc{\supp{\bm{A}_{n,s}(\vnu_{s+1/n},)}=\cbc{j_1,\ldots,j_k}}\\
& \hspace{-0.1 cm} \cdot\PP(j_\ell\in \mathcal{F}\bc{\bm{A}_{n,s+1/n}[\bm{\theta}]\abc{\cbc{j_1,\ldots,j_k};}}\Delta \mathcal{F}\bc{\bm{A}_{n,s+1/n}[\bm{\theta}]\abc{j_\ell;}}\mid \supp{\bm{A}_{n,s}(\vnu_{s+1/n},)}=\cbc{j_1,\ldots,j_k})\\
\leq & K^{2K+3}c^{-K}P^{-1}+\bar{o}_{n,P}(1).
\end{align*}
By design, on the event
$\mathfrak C \cap {\mathfrak Z}_s^{\circ}$, any $i \in \supp{\bm{A}_{n,s}(\vnu_{s+1/n},)}$ is frozen in $\bm{A}_{n,s+1/n}[\bm{\theta}]\abc{\supp{\bm{A}_{n,s}[\bm \theta](\vnu_{s+1/n},)};}$. In particular, there exists a representation $b$ of $\{i\}$ in $\bm{A}_{n,s+1/n}[\bm{\theta}]\abc{\supp{\bm{A}_{n,s}[\bm \theta](\vnu_{s+1/n},)};}$ with $b_k=0$ for all $k\in \supp{\bm{A}_{n,s}[\bm \theta](\vnu_{s+1/n},)}$. On the good event $\mathfrak P_{n,s}$, 
\begin{align}\label{eq_l12_1}
b \bm{A}_{n,s}[\bm\theta]\abc{\vnu_{s+1/n};} =e_{n+\bm{\theta}_c}(i).
\end{align}
Thus, on the event ${\mathfrak Z}_s^{\circ}\cap \mathfrak{C}\cap \mathfrak P_{n,s}$, any $i \in \supp{\bm{A}_{n,s}[\bm \theta](\vnu_{s+1/n},)}$ is frozen in $\bm{A}_{n,s}[\bm{\theta}]\abc{\vnu_{s+1/n};}$.
We conclude that the $\vnu_{s+1/n}$th row in $\bm{A}_{n,s}[\bm{\theta}]$ can be linearly combined by the other rows of $\bm{A}_{n,s}[\bm{\theta}]$ (this is also true if $\supp{\bm{A}_{n,s}[\bm \theta](\vnu_{s+1/n},)} = \emptyset$). Therefore, by \Cref{lem_frozen_zeroing}, $\vnu_{s+1/n}$ is not frozen in $\bm{A}_{n,s}[\bm{\theta}]^T$, which only leaves the possibility $\vnu_{s+1/n}\in \mathcal{V}\bc{\bm{A}_{n,s}[\bm{\theta}]}\cup \mathcal{Z}\bc{\bm{A}_{n,s}[\bm{\theta}]}$ on  ${\mathfrak Z}_s^{\circ}\cap \mathfrak{C}\cap \mathfrak P_{n,s}$.

On the other hand, since $\mathfrak Z^{\circ}_s \subseteq \ffBT$, by \Cref{im_fnf}, $\vnu_{s+1/n}$ cannot be firmly frozen in $\bm{A}_{n,s}[\bm{\theta}]$ on the event ${\mathfrak Z}_s^{\circ}\cap \mathfrak{C}\cap \mathfrak P_{n,s}$.
Therefore, $\vnu_{s+1/n}\in  \mathcal{Z}\bc{\bm{A}_{n,s}[\bm{\theta}]}$  and we arrive at
\[
    \mathbb{P}\left(\vnu_{s+1/n}\not\in\mathcal{Z}\bc{\bm{A}_{n,s}[\bm{\theta}]}, {\mathfrak Z}_s^{\circ}, \mathfrak{C} \right)=\bar{o}_{n,P}(1),
\]
i.e.,
\[\begin{aligned}
\mathbb{P}\left(\vnu_{s+1/n}\not\in\mathcal{Z}\bc{\bm{A}_{n,s}[\bm{\theta}]}, {\mathfrak Z}_s^{\circ}\right)
\leq  K^{2K+3}c^{-K}P^{-1}+\bar{o}_{n,P}(1) =  \bar{o}_{n,P}(1).
\end{aligned}\]
\qed
\subsection{Proof of \Cref{imp_pro}} \label{sec_proof_imp_pro}
Recall $h_t:[0,1] \to \RR, h_t\bc{\alpha}=\alpha+1-\hat{\psi}_t\bc{\alpha}$ and that $\alpha_\star(t)\leq \alpha_0(t)\leq \alpha^\star(t)$ denote the zeros of $\alpha \mapsto G_t(\alpha)$. 
\begin{lemma}\label{imp_lem}
For any $t\geq 0$, 
\begin{align}\label{eq_h_al}
h_t\bc{\alpha_\star(t)}=h_t\bc{\alpha^\star(t)}\leq h_t\bc{\alpha_0(t)}.
\end{align}
\end{lemma}

\begin{proof}
By \cref{ap_it2,ap_it3} in \Cref{lem_proal}, $\alpha_\star(t)\leq \alpha_0(t)\leq \alpha^\star(t)$ and $h_t\bc{\alpha_\star(t)}=h_t\bc{\alpha^\star(t)}$. The last inequality in \cref{eq_h_al} follows from the concavity of $h_t$.
\end{proof}

\begin{lemma}\label{imp_lem0}
The following inequality holds for $\ranges$:
\begin{equation}\label{eq_eqineq}
\begin{aligned}
    h_{t_s}\bc{\vx_s+\vy_s+\vu_s}=& \hspace{0.1 cm}  h_{t_s}\bc{\vx_s+\vy_s+\vv_s}+\oone
    = h_{t_s}\bc{\vx_s+\vy_s}+\oone\\
    \leq& \hspace{0.1 cm}  h_{t_s}\bc{\vy_s}+\oone.
    \end{aligned}
\end{equation}
\end{lemma}
\begin{proof}
The first and second equalities in \cref{eq_eqineq} follow directly from \cref{eq_fpu} and \cref{eq_fpv}.

Note that $\vx_s+\vy_s+\vz_s+\vu_s+\vv_s=1$. Then summing up \cref{eq_fpy}, \cref{eq_fpu}, and \cref{eq_fpv} gives that
\begin{align}\label{eq_fpxz}
    \vx_s+\vz_s=1-(\vy_s+\vu_s+\vv_s)=\hat{\psi}_{t_s}\bc{\vx_s+\vy_s}+\oone.
\end{align}
Then the combination of \cref{eq_fpxz} and \cref{eq_fpz} gives
    \begin{align*}
            h_{t_s}\bc{\vx_s+\vy_s}=& \hspace{0.1 cm} \vx_s+\vy_s+1-\hat{\psi}_{t_s}\bc{\vx_s+\vy_s}=\vy_s+1-\vz_s+\oone\\
            \leq&  \hspace{0.1 cm} \vy_s+1-\hat{\psi}_{t_s}\bc{\vy_s}+\oone=h_{t_s}\bc{\vy_s}+\oone,
    \end{align*}
and thus the last inequality in \cref{eq_eqineq} follows.
\end{proof}

\begin{proof}[Proof of \Cref{imp_pro}]
We aim to show that
\begin{align}\label{eq_eqeqh}
    h_{t_s}(\val_s)\geq h_{t_s}(\alpha^\star(t_s))+\oone.
\end{align}
Then \Cref{imp_pro} follows directly from the combination of \Cref{eq_eqeqh}, \cref{eq_fpy}, \cref{eq_fpu} and \Cref{eq_re_rd}.

Define
    \begin{align}\label{def_vtavet}
    \bvta_s=\teo{h_{t_s}'(\vx_s+\vy_s)\geq 0}\quad \mbox{and} \quad\bvet_s=1-\bvta_s=\teo{h_{t_s}'(\vx_s+\vy_s)<0}.
    \end{align}
Since $\bvta_s+\bvet_s=1$, we divide equation \cref{eq_eqeqh} into two parts as follows:
\begin{align}
   &\bvta_s\bc{ h_{t_s}\bc{\alpha^\star(t)} - h_{t_s}\bc{\val_s}}\leq \oone,\label{eq_taur}\\
\mbox{and}\quad&\bvet_s  \bc{h_{t_s}\bc{\alpha^\star(t)}-h_{t_s}\bc{\val_s}}\leq \oone.\label{eq_etar}
\end{align}
In the absence of the error terms $\oone$ and under the assumption that $\bvta_s\equiv 1$ (or $\bvet_s\equiv 1$), the proof of \Cref{imp_pro} would amount to an analytic treatment of the properties of $h_{t_s}$. Unfortunately, we have to deal with the error terms and both cases. In the ensuing argument, we therefore fall back upon Taylor's Theorem with Lagrange Remainder, \Cref{lem_proal} and the following two facts:
\begin{itemize}\label{eq_01eq}
    \item[(i)] For any function $g$ and $\upsilon\in\cbc{0,1}$,  if $\upsilon a=\upsilon b$, then $\upsilon g(a)=\upsilon g(b)$;
    \item[(ii)] For a family of differentiable functions $\bc{g_s}_{\ranges}$, if there exists a uniform bound $b$ such that $\sup_{\ranges,r\in[0,1]}\abs{g'_s(r)}\leq b$, then $g_s(\bm{a}_s')=g_s(\bm{a}_s)+\oone$ for any random variables $\bm{a}_s,\bm{a}'_s\in [0,1]$, $\bm{a}_s'=\bm{a}_s+\oone\in [0,1]$ and $\ranges$,  since $\abs{g_s(\bm{a}_s')-g_s(\bm{a}_s)}\leq b \abs{\bm{a}_s'-\bm{a}_s}$.
\end{itemize}

\begin{enumerate}
    \item \begin{proof}[Proof of \Cref{eq_taur}]
  By definition of $\bvta_s$, $\bvta_s h_{t_s}'\bc{\vx_s+\vy_s}\geq 0$. 
  For a fixed $\varepsilon>0$, Let $b=\inf_{\ranges}\hat{\psi}''_{t_s}(0)>0$, such that $\sup_{\alpha \in [0,1]} h_{t_s}''\bc{\alpha}\leq-b$ for $\ranges$. Then by Taylor's Theorem with Lagrange remainder, for $\ranges$,
    \begin{align*}
            \bvta_s h_{t_s}(\vy_s) 
            \leq& \hspace{0.1 cm} \bvta_s \bc{ h_{t_s}\bc{\vx_s+\vy_s}- h_{t_s}'\bc{\vx_s+\vy_s} \vx_s-\frac{b}{2} \vx_s^2}\\
            \leq& \hspace{0.1 cm} \bvta_s h_{t_s}\bc{\vx_s+\vy_s}
    -\bvta_s\frac{b}{2}\vx_s^2.
    \end{align*}
On the other hand, \cref{eq_eqineq}  shows that for all $\ranges$,
\begin{align*}
    \bvta_s h_{t_s}(\vy_s)\geq \bvta_sh_{t_s}\bc{\vx_s+\vy_s}+\oone.
\end{align*}
Hence
\begin{align}\label{eq_vtavx}
    \bvta_s\vx_s= \oone.
\end{align}
Recall from \Cref{eq_val_xyv} that $\val_s = \vx_s + \vy_s + \vv_s$. We define $\hval_s = \vx_s + \vy_s + \vu_s$. Then \cref{eq_vtavx} in combination with \cref{eq_fpy} and \cref{eq_fpv} implies that 
\begin{align*}
    \bvta_s\val_s=& \hspace{0.1 cm}\bvta_s\bc{\vy_s+\vv_s+\oone}=\bvta_s\bc{1-\hat{\psi}_{t_s}\bc{\vx_s+\vy_s+\vu_s}+\oone}\\
    =& \hspace{0.1 cm} \bvta_s(1-\hat{\psi}_{t_s}(\hval_s)+\oone).
\end{align*}
Analogously, \cref{eq_vtavx} in combination with \cref{eq_fpy} and \cref{eq_fpu} implies that $\bvta_s\hval_s=\bvta_s\bc{1-\hat{\psi}_{t_s}\bc{\val_s}+\oone}$.
Hence, 
    \begin{align*}
            \bvta_s\val_s=\bvta_s(1-\hat{\psi}_{t_s}(\hval_s)+\oone)=\bvta_s(1-\hat{\psi}_{t_s}\bc{1-\hat{\psi}_{t_s}\bc{\val_s}}+\oone),
    \end{align*}
i.e., 
\begin{align}\label{eq_Gat_al}
\bvta_s G_{ t_s }\bc{\val_s}=\oone.
\end{align}
Let $\vbe_s=\bvta_s\val_s+\bvet_s\alpha^\star(t)$. Since $G_t(\alpha^\star(t))=0$, \cref{eq_Gat_al} implies that
    \begin{align}\label{eq_vbe}
   G_{ t_s }(\vbe_s)= \bvta_s G_{ t_s }\bc{\val_s}+\bvet_s G_{ t_s }\bc{\alpha^\star(t)}=\oone.
    \end{align}
Hence, \cref{ap_it7} in \Cref{lem_proal} implies that
    \[
    \min\cbc{\abs{\vbe_s-\alpha_\star( t_s )},\abs{\vbe_s-{\alpha_0( t_s )}},\abs{\vbe_s-\alpha^\star( t_s )}}=\oone.
    \]
By \Cref{imp_lem}, $h_{t_s}\bc{\alpha_\star( t_s )}=h_{t_s}\bc{\alpha^\star( t_s )}\leq h_{t_s}\bc{\alpha_0( t_s )}$, so 
    \[
    h_{t_s}\bc{\alpha^\star( t_s )}\leq h_{t_s}\bc{\vbe_s}+\oone=\bvta_s h_{t_s}\bc{\val_s}+\bvet_s h_{t_s}\bc{\alpha^\star( t_s  )}+\oone,\]
and \cref{eq_taur} follows immediately.
\end{proof}

\item \begin{proof}[Proof of \Cref{eq_etar}]
By definition of $\bvet_s$, $\bvet_s h_{t_s}'\bc{\vx_s+\vy_s}\leq 0$. Another application of Taylor's Theorem with Lagrange remainder to \cref{eq_fpu} and \cref{eq_fpv} as in the argument leading to \cref{eq_vtavx} 
yields that $\bvet_s\vu_s=\oone$ and $\bvet_s\vv_s=\oone$. 
Hence, by \cref{eq_fpy}, 
    \begin{align}
    \label{eq_fixvet}
    \bvet_s\vy_s=\bvet_s(1-\hat{\psi}_{t_s}(\vx_s+\vy_s)+\oone)\ \mbox{and }\ \bvet_s\val_s=\bvet_s(\vx_s+\vy_s+\oone).
    \end{align}

Let
    \[
    \vbe_s'=\bvet_s\teo{\val_s> {\alpha^\star( t_s )}}\val_s+\bc{1-\bvet_s\teo{\val_s> {\alpha^\star( t_s )}}}\alpha^\star( t_s ).
    \]
Then $\vbe_s'\geq \alpha^\star( t_s )$. 
Since $G_{t_s}(\alpha^\star( t_s ))=0$, by \cref{eq_fixvet}, \cref{eq_fpz} and \cref{eq_fpxz},
    \[
    \begin{aligned}
    G_{ t_s }(\vbe_s')=&\bvet_s\teo{\val_s> {\alpha^\star( t_s )}}\bc{\val_s+\hat{\psi}_{t_s}\bc{1-\hat{\psi}_{t_s}\bc{\val_s}}-1}\\
    \leq& \bvet_s\teo{\val_s> {\alpha^\star( t_s )}}\bc{\vx_s+\vy_s+\vz_s-1+\oone}\leq  \oone.
    \end{aligned}
    \]

On the other hand, by \cref{ap_it4} in \Cref{lem_proal},   $ G_{ t_s }(\vbe_s')\geq G_{ t_s }\bc{\alpha^\star( t_s )}=0$. Hence,   $G_{ t_s }(\vbe_s')=\oone$.

Then the combination of \cref{ap_it7} in \Cref{lem_proal} and $\vbe'_s\geq \alpha^\star( t_s )$ yields that $\vbe'_s=\alpha^\star( t_s )+\oone$, which leads to
    \begin{align}
    \label{eq_valup}
    \bvet_s\val_s\leq \bvet_s\vbe'_s=\bvet_s\bc{\alpha^\star( t_s )+\oone}.
    \end{align}
By the concavity of $h_{t_s}$ and $\bvet_s h'_{t_s}\bc{\vx_s+\vy_s}\leq 0$, $\alpha \mapsto \bvet_s h_{t_s}\bc{\alpha}$ is non-increasing on $[\vx_s+\vy_s,1]$. Hence, by 
\cref{eq_valup}, 
    \[
    \bvet_s h_{t_s}\bc{\val_s}\geq \bvet_s h_{t_s}\bc{\vbe_s'}=\bvet_s\bc{h_{t_s}\bc{\alpha^\star( t_s )}+\oone},
    \]
and thus \cref{eq_etar} holds.
\end{proof}
\end{enumerate}
The combination of \cref{eq_taur,eq_etar} gives \cref{eq_eqeqh}.

\end{proof}
\end{appendices}

\end{document}